\documentclass[a4paper,11pt]{article}
\usepackage{tgpagella,eulervm}
\usepackage{multicol} %
\usepackage{tablefootnote}
\usepackage[english]{babel}
\usepackage[normalem]{ulem}
\usepackage[utf8x]{inputenc}
\usepackage[T1]{fontenc}
\usepackage{amsmath}
\usepackage{amssymb}
\usepackage{mathrsfs}
\usepackage{float}
\usepackage{wrapfig}
\usepackage{amsthm}
\usepackage{tikz-cd}
\usepackage{tikz}
\usepackage{tikz-3dplot} 
\usetikzlibrary{shapes.geometric, graphs, arrows, arrows.meta, bending, positioning}
\usepackage{subcaption}
\usepackage{pgfplots}
\pgfplotsset{compat=1.17}
\usepackage[nottoc]{tocbibind}%
\usepackage{tikz}
\usetikzlibrary{calc}

\definecolor{moduleRed}{RGB}{210,45,45}
\definecolor{moduleBlue}{RGB}{35,70,200}
\definecolor{moduleGreen}{RGB}{72,145,55}
\definecolor{moduleOrange}{RGB}{240,125,20}
\definecolor{regionBlue}{RGB}{210,225,250}

\usepackage{savesym}
\savesymbol{faCompress}
\usepackage{fontawesome}

\usepackage[colorlinks=true, allcolors=darkblue]{hyperref}
\usepackage[nameinlink]{cleveref}
\usepackage{nomencl}
\usepackage{multicol}
\usepackage{makecell}
\usepackage{multirow}
\usepackage{caption}

\usepackage[many]{tcolorbox}
\usetikzlibrary{decorations.pathreplacing}

\newtcolorbox{leftbrace}{%
    enhanced jigsaw, 
    breakable, %
    frame hidden, %
    overlay={%
        \draw [
            decoration={brace,amplitude=0.5em},
            decorate,
            ultra thick,
        ]
        (frame.north west)--(frame.south west);
    },
    parbox=false,
}

\renewcommand{\nomgroup}[1]{%
\ifthenelse{\equal{#1}{B}}{\item[\textbf{Invariants}]}{%
\ifthenelse{\equal{#1}{C}}{\item[\textbf{Functions}]}{%
\ifthenelse{\equal{#1}{D}}{\item[\textbf{Other}]}{%
\ifthenelse{\equal{#1}{A}}{\item[\textbf{Collections of subposets}]}{%
}}}}}

\usepackage{thm-restate}
\usepackage{enumerate}
\usepackage[inline]{enumitem}
\usepackage{xcolor}
\usepackage{colortbl} %
\usepackage{bbold}
\usepackage{color}
\usepackage{soul}
\usepackage{pdflscape} %
\usepackage[margin=0.95in]{geometry}
\usepackage{authblk}
\usepackage{framed}

\numberwithin{equation}{section}

\newcommand\freefootnote[1]{%
  \let\thefootnote\relax%
  \footnotetext{#1}%
  \let\thefootnote\svthefootnote%
}

\usepackage{amsmath}
\usepackage{graphicx}
\usepackage{url}

\definecolor{darkblue}{rgb}{0.0, 0.0, 0.8}
\definecolor{darkred}{rgb}{0.8, 0.0, 0.0}
\definecolor{darkgreen}{rgb}{0.0, 0.65, 0.0}

\DeclareMathAlphabet{\mathcal}{OMS}{cmsy}{m}{n}

\newcommand{\white}[1] {{\textcolor{white}{#1}}}

\newcommand{\Lan}{\mathrm{Lan}}

\newcommand{\pfin}{pointwisely finite}
\newcommand{\invertible}{M\"obius invertible{}}
\newcommand{\convolvable}{convolvable}

\newcommand{\abs}[1]{\lvert{#1}\rvert}

\newcommand{\catQ}{\mathrm{Q}}
\newcommand{\R}{\mathbb{R}}
\newcommand{\N}{\mathbb{N}}
\newcommand{\Z}{\mathbb{Z}}
\newcommand{\kf}{k}
\newcommand{\im}{\mathrm{im}}

\newcommand{\K}{\mathbb{K}}
\newcommand{\Cvec}{\mathbf{Vec}_k}
\newcommand{\cvec}{\mathbf{vec}_k}

\newcommand{\inv}{^{-1}}
\newcommand{\id}{\mathrm{id}}
\newcommand{\dgm}{\mathrm{\mathrm{dgm}}}
\newcommand{\rk}{\mathrm{rk}}
\newcommand{\eps}{\epsilon}
\newcommand{\catP}{\mathrm{P}}
\newcommand{\catS}{\mathrm{S}}
\newcommand{\rank}{\mathrm{rk}}

\newcommand{\Int}{\mathrm{Int}}

\newcommand{\Con}{\mathrm{Con}}
\newcommand{\Jcal}{\mathcal{J}}
\newcommand{\Ical}{\mathcal{I}}
\newcommand{\Scal}{\mathcal{S}}

\newcommand{\Rcal}{\mathcal{R}}

\newcommand{\Bcal}{\mathcal{B}}

\newcommand{\barc}{\mathrm{barc}}

\newcommand{\bv}{\mathbf{v}}
\newcommand{\bu}{\mathbf{u}}
\newcommand{\bw}{\mathbf{w}}

\newcommand{\seg}{\mathrm{Seg}}

\newcommand{\supp}{\mathrm{supp}}
\newcommand{\mult}{\mathrm{mult}}

\newcommand{\drm}{\mathrm{d}}

\newcommand{\CC}{\mathbf{C}}
\newcommand{\Iast}{L}
\newenvironment{claimproof}
  {\par\noindent\textit{Proof. }}
  {\hfill$\mbox{\emph{(End of claim)}}\ \blacksquare$\par}

\tikzcdset{scale cd/.style={every label/.append style={scale=#1},
    cells={nodes={scale=#1}}}}

\newtheorem{theorem}{Theorem}%

\newtheorem{proposition}{Proposition}[section]%
\newtheorem{oldtheorem}[proposition]{Theorem}
\newtheorem{lemma}[proposition]{Lemma}
\newtheorem{corollary}[proposition]{Corollary}

\theoremstyle{definition}
\newtheorem{definition}[proposition]{Definition}
\newtheorem{example}[proposition]{Example}

\theoremstyle{remark}
\newtheorem{remark}[proposition]{Remark}

\newcounter{claimgroup}
\newtheorem{claim}{Claim}

\newcommand{\resetclaimcounter}{\stepcounter{claimgroup}\setcounter{claim}{0}}

\newtheorem*{lemma*}{Lemma}

\title{M\"obius Invertibility of Rank Invariants \freefootnote{\emph{2020 Mathematics Subject Classification}:\ \ 55N31, 55U99} \freefootnote{\emph{Keywords}: M\"obius inversion, persistence diagrams, multi-parameter persistent homology, topological data analysis}}

\date{}

\author[1]{Nathaniel Clause} 

\author[2]{Woojin Kim} 

\author[3]{Facundo M\'{e}moli}

\affil[1]{The Ohio State University, Columbus OH \thanks{\texttt{nate.clause@gmail.com}}}

\affil[2]{KAIST, Daejeon, South Korea 	\thanks{\texttt{woojin.kim.math@gmail.com}}}

\affil[3]{Rutgers University, Piscataway NJ \thanks{\texttt{facundo.memoli@gmail.com}}}

\immediate\pdfobj stream attr{
 /Type /XObject /Subtype /Form /FormType 1
/BBox [ 0 0 791 227 ]
/Group <<
/Type /Group
/S /Transparency
/I true
/CS /DeviceRGB
>>
/Resources <<
/ExtGState <<
/a0 <<
/CA 1
/ca 1
>>
>>
>>
}{
q
0 0 0 RG /a0 gs
2.928274 w
0 J
0 j
[] 0.0 d
4 M q 1 0 0 1 0 0 cm
237.668 26.363 145.02 171.398 re S Q
Q
}
\edef\mobiusintervalobjecta{\the\pdflastobj}

\immediate\pdfobj stream attr{
 /Type /XObject /Subtype /Form /FormType 1
/BBox [ 0 0 791 227 ]
/Group <<
/Type /Group
/S /Transparency
/I true
/CS /DeviceRGB
>>
/Resources <<
/ExtGState <<
/a0 <<
/CA 1
/ca 1
>>
>>
>>
}{
q
0 0 0 RG /a0 gs
2.834646 w
0 J
0 j
[] 0.0 d
4 M q 1 0 0 1 0 0 cm
229.172 18.102 341.133 191.117 re S Q
Q
}
\edef\mobiusintervalobjectb{\the\pdflastobj}

\immediate\pdfobj stream attr{
 /Type /XObject /Subtype /Form /FormType 1
/BBox [ 0 0 791 227 ]
/Group <<
/Type /Group
/S /Transparency
/I true
/CS /DeviceRGB
>>
/Resources <<
/ExtGState <<
/a0 <<
/CA 1
/ca 1
>>
>>
>>
}{
q
0 0 0 RG /a0 gs
2.684614 w
0 J
0 j
[] 0.0 d
4 M q 1 0 0 1 0 0 cm
1.344 9.332 579.523 215.934 re S Q
Q
}
\edef\mobiusintervalobjectc{\the\pdflastobj}

\immediate\pdfobj stream attr{
 /Type /XObject /Subtype /Form /FormType 1
/BBox [ 0 0 791 227 ]
/Group <<
/Type /Group
/S /Transparency
/I true
/CS /DeviceRGB
>>
/Resources <<
/ExtGState <<
/a0 <<
/CA 1
/ca 1
>>
>>
>>
}{
q
0 0 0 RG /a0 gs
2.802008 w
0 J
0 j
[] 0.0 d
4 M q 1 0 0 1 0 0 cm
221.406 1.402 568 215.527 re S Q
Q
}
\edef\mobiusintervalobjectd{\the\pdflastobj}
\newsavebox{\mobiusintervalbox}
\sbox{\mobiusintervalbox}{\pdfliteral{
1 0 0 -1 0 226.606277 cm
q
0.709804 0.709804 0.709804 rg /a0 gs
264.746 84.027 m 264.746 186.664 l 364.719 186.664 l 364.719 84.848 l h
264.746 84.027 m f
417.266 73.609 m 477.195 73.609 l 477.195 138.152 l 555.32 138.152 l 555.32
 187.73 l 412.988 187.73 l 412.852 73.609 l h
417.266 73.609 m f
16.23 77.563 m 16.23 136.887 l 78.41 136.887 l 78.41 195.07 l 205.617 195.07
 l 205.617 131.184 l 148.574 131.184 l 148.574 76.992 l h
16.23 77.563 m f
596.051 64.965 m 596.051 192.465 l 770.695 192.465 l 770.695 139.965 l 
710.695 139.965 l 710.695 107.285 l 654.98 107.285 l 654.98 60.676 l 596.371
 60.676 l h
596.051 64.965 m f
q
1 0 0 1 0 0 cm
/a1 gs /x6 Do
Q
q
1 0 0 1 0 0 cm
/a1 gs /x7 Do
Q
q
1 0 0 1 0 0 cm
/a1 gs /x8 Do
Q
q
1 0 0 1 0 0 cm
/a1 gs /x9 Do
Q
0 0 0 rg 254.996 48.074 m 254.996 46.961 255.047 46.621 257.52 46.621 c 258.293 
46.621 l 258.293 45.648 l 257.227 45.746 254.801 45.746 253.59 45.746 c 
252.426 45.746 249.953 45.746 248.887 45.648 c 248.887 46.621 l 249.664 
46.621 l 252.137 46.621 252.234 46.961 252.234 48.074 c 252.234 64.367 l
 252.234 65.48 252.137 65.82 249.664 65.82 c 248.887 65.82 l 248.887 66.793
 l 249.953 66.695 252.426 66.695 253.59 66.695 c 254.801 66.695 257.227 
66.695 258.293 66.793 c 258.293 65.82 l 257.52 65.82 l 255.047 65.82 254.996
 65.48 254.996 64.367 c h
254.996 48.074 m f
262.574 56.172 m 262.574 64.465 l 262.574 65.82 262.234 65.82 260.148 65.82
 c 260.148 66.793 l 261.262 66.742 262.816 66.695 263.688 66.695 c 264.465
 66.695 266.063 66.742 267.129 66.793 c 267.129 65.82 l 265.047 65.82 264.707
 65.82 264.707 64.465 c 264.707 58.742 l 264.707 55.543 266.938 53.797 268.875
 53.797 c 270.863 53.797 271.203 55.445 271.203 57.238 c 271.203 64.465 
l 271.203 65.82 270.863 65.82 268.781 65.82 c 268.781 66.793 l 269.848 66.742
 271.445 66.695 272.27 66.695 c 273.047 66.695 274.695 66.742 275.715 66.793
 c 275.715 65.82 l 274.113 65.82 273.336 65.82 273.336 64.902 c 273.336 
58.984 l 273.336 56.316 273.336 55.395 272.367 54.281 c 271.93 53.746 270.914
 53.117 269.117 53.117 c 266.84 53.117 265.387 54.426 264.512 56.367 c 264.512
 53.117 l 260.148 53.457 l 260.148 54.426 l 262.332 54.426 262.574 54.621
 262.574 56.172 c h
262.574 56.172 m f
280.824 54.426 m 285.285 54.426 l 285.285 53.457 l 280.824 53.457 l 280.824
 47.785 l 280.047 47.785 l 280.047 50.305 279.129 53.602 276.07 53.746 c
 276.07 54.426 l 278.691 54.426 l 278.691 62.961 l 278.691 66.742 281.602
 67.133 282.715 67.133 c 284.898 67.133 285.77 64.949 285.77 62.961 c 285.77
 61.215 l 284.996 61.215 l 284.996 62.914 l 284.996 65.191 284.074 66.355
 282.91 66.355 c 280.824 66.355 280.824 63.543 280.824 63.008 c h
280.824 54.426 m f
294.746 57.656 m 294.746 57.074 294.746 57.074 294.113 57.074 c 292.707
 58.434 290.77 58.434 289.895 58.434 c 289.895 59.207 l 290.379 59.207 291.836
 59.207 293 58.625 c 293 69.633 l 293 70.359 293 70.652 290.867 70.652 c
 290.043 70.652 l 290.043 71.43 l 290.43 71.43 293.098 71.332 293.871 71.332
 c 294.551 71.332 297.266 71.43 297.75 71.43 c 297.75 70.652 l 296.926 70.652
 l 294.746 70.652 294.746 70.359 294.746 69.633 c h
294.746 57.656 m f
304.402 71.09 m 304.402 72.254 304.211 73.66 302.707 75.016 c 302.66 75.113
 302.563 75.211 302.563 75.309 c 302.563 75.453 302.754 75.598 302.852 75.598
 c 303.191 75.598 305.035 73.852 305.035 71.281 c 305.035 69.973 304.5 68.953
 303.484 68.953 c 302.805 68.953 302.27 69.488 302.27 70.168 c 302.27 70.848
 302.754 71.43 303.531 71.43 c 304.066 71.43 304.402 71.09 304.402 71.09
 c h
304.402 71.09 m f
314.414 57.656 m 314.414 57.074 314.414 57.074 313.785 57.074 c 312.379
 58.434 310.438 58.434 309.566 58.434 c 309.566 59.207 l 310.051 59.207 
311.504 59.207 312.668 58.625 c 312.668 69.633 l 312.668 70.359 312.668 
70.652 310.535 70.652 c 309.711 70.652 l 309.711 71.43 l 310.098 71.43 312.766
 71.332 313.539 71.332 c 314.219 71.332 316.934 71.43 317.422 71.43 c 317.422
 70.652 l 316.594 70.652 l 314.414 70.652 314.414 70.359 314.414 69.633 
c h
314.414 57.656 m f
331.289 74.211 m 331.289 74.113 331.289 74.066 330.754 73.531 c 326.926
 69.652 325.906 63.785 325.906 59.035 c 325.906 53.652 327.07 48.27 330.902
 44.438 c 331.289 44.051 331.289 44 331.289 43.906 c 331.289 43.66 331.191
 43.566 330.996 43.566 c 330.66 43.566 327.895 45.699 326.051 49.625 c 324.5
 53.02 324.113 56.465 324.113 59.035 c 324.113 61.457 324.453 65.191 326.148
 68.73 c 328.039 72.512 330.66 74.551 330.996 74.551 c 331.191 74.551 331.289
 74.453 331.289 74.211 c h
331.289 74.211 m f
376.488 59.035 m 376.488 56.656 376.148 52.922 374.449 49.383 c 372.605
 45.602 369.941 43.566 369.648 43.566 c 369.457 43.566 369.309 43.711 369.309
 43.906 c 369.309 44 369.309 44.051 369.891 44.633 c 372.945 47.688 374.691
 52.586 374.691 59.035 c 374.691 64.367 373.578 69.797 369.746 73.676 c 
369.309 74.066 369.309 74.113 369.309 74.211 c 369.309 74.406 369.457 74.551
 369.648 74.551 c 369.941 74.551 372.754 72.418 374.547 68.488 c 376.148
 65.094 376.488 61.652 376.488 59.035 c h
376.488 59.035 m f
630.457 23.289 m 630.457 22.176 630.504 21.836 632.98 21.836 c 633.754 
21.836 l 633.754 20.867 l 632.688 20.961 630.262 20.961 629.051 20.961 c
 627.887 20.961 625.414 20.961 624.348 20.867 c 624.348 21.836 l 625.125
 21.836 l 627.598 21.836 627.691 22.176 627.691 23.289 c 627.691 39.582 
l 627.691 40.699 627.598 41.035 625.125 41.035 c 624.348 41.035 l 624.348
 42.008 l 625.414 41.91 627.887 41.91 629.051 41.91 c 630.262 41.91 632.688
 41.91 633.754 42.008 c 633.754 41.035 l 632.98 41.035 l 630.504 41.035 
630.457 40.699 630.457 39.582 c h
630.457 23.289 m f
638.035 31.387 m 638.035 39.68 l 638.035 41.035 637.695 41.035 635.609 
41.035 c 635.609 42.008 l 636.723 41.957 638.277 41.91 639.148 41.91 c 639.926
 41.91 641.523 41.957 642.59 42.008 c 642.59 41.035 l 640.508 41.035 640.168
 41.035 640.168 39.68 c 640.168 33.957 l 640.168 30.758 642.398 29.012 644.336
 29.012 c 646.324 29.012 646.664 30.66 646.664 32.453 c 646.664 39.68 l 
646.664 41.035 646.324 41.035 644.238 41.035 c 644.238 42.008 l 645.309 
41.957 646.906 41.91 647.73 41.91 c 648.508 41.91 650.156 41.957 651.176
 42.008 c 651.176 41.035 l 649.574 41.035 648.797 41.035 648.797 40.117 
c 648.797 34.199 l 648.797 31.531 648.797 30.613 647.828 29.496 c 647.391
 28.965 646.375 28.332 644.578 28.332 c 642.301 28.332 640.844 29.641 639.973
 31.582 c 639.973 28.332 l 635.609 28.672 l 635.609 29.641 l 637.789 29.641
 638.035 29.836 638.035 31.387 c h
638.035 31.387 m f
656.285 29.641 m 660.746 29.641 l 660.746 28.672 l 656.285 28.672 l 656.285
 23 l 655.508 23 l 655.508 25.52 654.586 28.816 651.531 28.965 c 651.531
 29.641 l 654.152 29.641 l 654.152 38.176 l 654.152 41.957 657.063 42.348
 658.176 42.348 c 660.359 42.348 661.23 40.164 661.23 38.176 c 661.23 36.43
 l 660.453 36.43 l 660.453 38.129 l 660.453 40.406 659.535 41.57 658.371
 41.57 c 656.285 41.57 656.285 38.758 656.285 38.227 c h
656.285 29.641 m f
670.207 32.871 m 670.207 32.289 670.207 32.289 669.574 32.289 c 668.168
 33.648 666.23 33.648 665.355 33.648 c 665.355 34.426 l 665.84 34.426 667.297
 34.426 668.461 33.844 c 668.461 44.848 l 668.461 45.578 668.461 45.867 
666.324 45.867 c 665.5 45.867 l 665.5 46.645 l 665.891 46.645 668.555 46.547
 669.332 46.547 c 670.012 46.547 672.727 46.645 673.211 46.645 c 673.211
 45.867 l 672.387 45.867 l 670.207 45.867 670.207 45.578 670.207 44.848 
c h
670.207 32.871 m f
679.863 46.305 m 679.863 47.469 679.672 48.875 678.168 50.23 c 678.121 
50.328 678.023 50.426 678.023 50.523 c 678.023 50.668 678.215 50.813 678.313
 50.813 c 678.652 50.813 680.496 49.066 680.496 46.496 c 680.496 45.188 
679.961 44.172 678.945 44.172 c 678.266 44.172 677.73 44.703 677.73 45.383
 c 677.73 46.063 678.215 46.645 678.992 46.645 c 679.523 46.645 679.863 
46.305 679.863 46.305 c h
679.863 46.305 m f
688.566 39.418 m 690.262 39.418 691.473 40.582 691.473 42.91 c 691.473 
45.578 689.875 46.402 688.66 46.402 c 687.789 46.402 685.852 46.16 684.977
 44.848 c 685.996 44.848 686.238 44.121 686.238 43.637 c 686.238 42.957 
685.703 42.473 685.023 42.473 c 684.445 42.473 683.813 42.859 683.813 43.734
 c 683.813 45.77 686.043 47.078 688.66 47.078 c 691.668 47.078 693.754 45.043
 693.754 42.91 c 693.754 41.211 692.395 39.516 690.02 39.031 c 692.25 38.207
 693.074 36.605 693.074 35.25 c 693.074 33.551 691.137 32.289 688.711 32.289
 c 686.336 32.289 684.492 33.453 684.492 35.199 c 684.492 35.926 684.977
 36.316 685.605 36.316 c 686.285 36.316 686.723 35.828 686.723 35.25 c 686.723
 34.617 686.285 34.133 685.605 34.086 c 686.383 33.164 687.836 32.922 688.66
 32.922 c 689.633 32.922 690.988 33.406 690.988 35.25 c 690.988 36.168 690.699
 37.188 690.117 37.816 c 689.438 38.641 688.809 38.691 687.742 38.789 c 
687.207 38.836 687.16 38.836 687.063 38.836 c 687.012 38.836 686.82 38.887
 686.82 39.129 c 686.82 39.418 687.012 39.418 687.402 39.418 c h
688.566 39.418 m f
706.75 49.426 m 706.75 49.328 706.75 49.281 706.215 48.746 c 702.387 44.867
 701.367 39 701.367 34.25 c 701.367 28.867 702.531 23.484 706.363 19.652
 c 706.75 19.266 706.75 19.215 706.75 19.121 c 706.75 18.879 706.652 18.781
 706.457 18.781 c 706.117 18.781 703.355 20.914 701.512 24.84 c 699.961 
28.234 699.574 31.68 699.574 34.25 c 699.574 36.672 699.914 40.406 701.609
 43.945 c 703.5 47.73 706.117 49.766 706.457 49.766 c 706.652 49.766 706.75
 49.668 706.75 49.426 c h
706.75 49.426 m f
751.945 34.25 m 751.945 31.871 751.609 28.141 749.91 24.598 c 748.066 20.816
 745.402 18.781 745.109 18.781 c 744.914 18.781 744.77 18.926 744.77 19.121
 c 744.77 19.215 744.77 19.266 745.352 19.848 c 748.406 22.902 750.152 27.801
 750.152 34.25 c 750.152 39.582 749.039 45.012 745.207 48.891 c 744.77 49.281
 744.77 49.328 744.77 49.426 c 744.77 49.621 744.914 49.766 745.109 49.766
 c 745.402 49.766 748.215 47.633 750.008 43.703 c 751.609 40.309 751.945
 36.867 751.945 34.25 c h
751.945 34.25 m f
58.902 30.457 m 58.902 29.344 58.949 29.004 61.422 29.004 c 62.199 29.004
 l 62.199 28.031 l 61.133 28.129 58.707 28.129 57.496 28.129 c 56.332 28.129
 53.859 28.129 52.793 28.031 c 52.793 29.004 l 53.566 29.004 l 56.039 29.004
 56.137 29.344 56.137 30.457 c 56.137 46.75 l 56.137 47.867 56.039 48.203
 53.566 48.203 c 52.793 48.203 l 52.793 49.176 l 53.859 49.078 56.332 49.078
 57.496 49.078 c 58.707 49.078 61.133 49.078 62.199 49.176 c 62.199 48.203
 l 61.422 48.203 l 58.949 48.203 58.902 47.867 58.902 46.75 c h
58.902 30.457 m f
66.477 38.555 m 66.477 46.848 l 66.477 48.203 66.137 48.203 64.055 48.203
 c 64.055 49.176 l 65.168 49.125 66.719 49.078 67.594 49.078 c 68.367 49.078
 69.969 49.125 71.035 49.176 c 71.035 48.203 l 68.949 48.203 68.609 48.203
 68.609 46.848 c 68.609 41.125 l 68.609 37.926 70.84 36.18 72.781 36.18 
c 74.77 36.18 75.109 37.828 75.109 39.621 c 75.109 46.848 l 75.109 48.203
 74.77 48.203 72.684 48.203 c 72.684 49.176 l 73.75 49.125 75.352 49.078
 76.176 49.078 c 76.953 49.078 78.602 49.125 79.617 49.176 c 79.617 48.203
 l 78.02 48.203 77.242 48.203 77.242 47.285 c 77.242 41.367 l 77.242 38.699
 77.242 37.781 76.273 36.664 c 75.836 36.129 74.816 35.5 73.023 35.5 c 70.746
 35.5 69.289 36.809 68.418 38.75 c 68.418 35.5 l 64.055 35.84 l 64.055 36.809
 l 66.234 36.809 66.477 37.004 66.477 38.555 c h
66.477 38.555 m f
84.73 36.809 m 89.191 36.809 l 89.191 35.84 l 84.73 35.84 l 84.73 30.168
 l 83.953 30.168 l 83.953 32.688 83.031 35.984 79.977 36.129 c 79.977 36.809
 l 82.598 36.809 l 82.598 45.344 l 82.598 49.125 85.504 49.516 86.621 49.516
 c 88.801 49.516 89.676 47.332 89.676 45.344 c 89.676 43.598 l 88.898 43.598
 l 88.898 45.297 l 88.898 47.574 87.977 48.738 86.813 48.738 c 84.73 48.738
 84.73 45.926 84.73 45.391 c h
84.73 36.809 m f
102.336 49.883 m 101.609 49.883 l 101.559 50.367 101.316 51.629 101.027
 51.824 c 100.879 51.969 99.23 51.969 98.891 51.969 c 94.918 51.969 l 97.195
 49.98 97.973 49.352 99.23 48.332 c 100.832 47.07 102.336 45.715 102.336
 43.676 c 102.336 41.059 100.055 39.457 97.293 39.457 c 94.625 39.457 92.781
 41.348 92.781 43.336 c 92.781 44.402 93.703 44.551 93.945 44.551 c 94.43
 44.551 95.063 44.16 95.063 43.387 c 95.063 42.996 94.918 42.223 93.801 
42.223 c 94.48 40.719 95.934 40.234 96.953 40.234 c 99.133 40.234 100.25
 41.93 100.25 43.676 c 100.25 45.566 98.891 47.023 98.215 47.797 c 93.023
 52.988 l 92.781 53.18 92.781 53.23 92.781 53.813 c 101.703 53.813 l h
102.336 49.883 m f
108.309 53.473 m 108.309 54.637 108.117 56.043 106.613 57.398 c 106.563
 57.496 106.465 57.594 106.465 57.691 c 106.465 57.836 106.66 57.98 106.758
 57.98 c 107.098 57.98 108.941 56.234 108.941 53.664 c 108.941 52.355 108.406
 51.34 107.387 51.34 c 106.711 51.34 106.176 51.871 106.176 52.551 c 106.176
 53.23 106.66 53.813 107.438 53.813 c 107.969 53.813 108.309 53.473 108.309
 53.473 c h
108.309 53.473 m f
122.004 49.883 m 121.277 49.883 l 121.227 50.367 120.984 51.629 120.695
 51.824 c 120.551 51.969 118.902 51.969 118.563 51.969 c 114.586 51.969 
l 116.863 49.98 117.641 49.352 118.902 48.332 c 120.5 47.07 122.004 45.715
 122.004 43.676 c 122.004 41.059 119.727 39.457 116.961 39.457 c 114.293
 39.457 112.453 41.348 112.453 43.336 c 112.453 44.402 113.371 44.551 113.613
 44.551 c 114.102 44.551 114.73 44.16 114.73 43.387 c 114.73 42.996 114.586
 42.223 113.469 42.223 c 114.148 40.719 115.602 40.234 116.621 40.234 c 
118.805 40.234 119.918 41.93 119.918 43.676 c 119.918 45.566 118.563 47.023
 117.883 47.797 c 112.695 52.988 l 112.453 53.18 112.453 53.23 112.453 53.813
 c 121.375 53.813 l h
122.004 49.883 m f
135.195 56.594 m 135.195 56.496 135.195 56.449 134.66 55.914 c 130.828 
52.035 129.813 46.168 129.813 41.418 c 129.813 36.035 130.977 30.652 134.805
 26.82 c 135.195 26.434 135.195 26.383 135.195 26.289 c 135.195 26.043 135.098
 25.949 134.902 25.949 c 134.563 25.949 131.801 28.082 129.957 32.008 c 
128.406 35.402 128.016 38.848 128.016 41.418 c 128.016 43.84 128.355 47.574
 130.055 51.113 c 131.945 54.895 134.563 56.934 134.902 56.934 c 135.098
 56.934 135.195 56.836 135.195 56.594 c h
135.195 56.594 m f
180.391 41.418 m 180.391 39.039 180.051 35.305 178.355 31.766 c 176.512
 27.984 173.844 25.949 173.555 25.949 c 173.359 25.949 173.215 26.094 173.215
 26.289 c 173.215 26.383 173.215 26.434 173.797 27.016 c 176.852 30.07 178.598
 34.969 178.598 41.418 c 178.598 46.75 177.48 52.18 173.652 56.059 c 173.215
 56.449 173.215 56.496 173.215 56.594 c 173.215 56.789 173.359 56.934 173.555
 56.934 c 173.844 56.934 176.656 54.801 178.453 50.871 c 180.051 47.477 
180.391 44.035 180.391 41.418 c h
180.391 41.418 m f
423.727 39.484 m 423.727 38.324 423.777 37.969 426.352 37.969 c 427.16 
37.969 l 427.16 36.961 l 426.051 37.059 423.523 37.059 422.262 37.059 c 
421.051 37.059 418.477 37.059 417.363 36.961 c 417.363 37.969 l 418.172 
37.969 l 420.746 37.969 420.848 38.324 420.848 39.484 c 420.848 56.453 l
 420.848 57.613 420.746 57.965 418.172 57.965 c 417.363 57.965 l 417.363
 58.977 l 418.477 58.875 421.051 58.875 422.262 58.875 c 423.523 58.875 
426.051 58.875 427.16 58.977 c 427.16 57.965 l 426.352 57.965 l 423.777 
57.965 423.727 57.613 423.727 56.453 c h
423.727 39.484 m f
431.617 47.918 m 431.617 56.551 l 431.617 57.965 431.262 57.965 429.09 
57.965 c 429.09 58.977 l 430.254 58.926 431.867 58.875 432.777 58.875 c 
433.586 58.875 435.254 58.926 436.363 58.977 c 436.363 57.965 l 434.191 
57.965 433.84 57.965 433.84 56.551 c 433.84 50.594 l 433.84 47.262 436.16
 45.441 438.18 45.441 c 440.25 45.441 440.605 47.16 440.605 49.027 c 440.605
 56.551 l 440.605 57.965 440.25 57.965 438.078 57.965 c 438.078 58.977 l
 439.191 58.926 440.855 58.875 441.715 58.875 c 442.523 58.875 444.242 58.926
 445.301 58.977 c 445.301 57.965 l 443.633 57.965 442.828 57.965 442.828
 57.008 c 442.828 50.848 l 442.828 48.07 442.828 47.109 441.816 45.949 c
 441.363 45.391 440.301 44.734 438.434 44.734 c 436.059 44.734 434.547 46.098
 433.637 48.121 c 433.637 44.734 l 429.09 45.09 l 429.09 46.098 l 431.363
 46.098 431.617 46.301 431.617 47.918 c h
431.617 47.918 m f
450.625 46.098 m 455.27 46.098 l 455.27 45.09 l 450.625 45.09 l 450.625
 39.18 l 449.816 39.18 l 449.816 41.809 448.855 45.242 445.676 45.391 c 
445.676 46.098 l 448.402 46.098 l 448.402 54.988 l 448.402 58.926 451.43
 59.328 452.594 59.328 c 454.863 59.328 455.773 57.059 455.773 54.988 c 
455.773 53.168 l 454.965 53.168 l 454.965 54.938 l 454.965 57.309 454.008
 58.523 452.793 58.523 c 450.625 58.523 450.625 55.594 450.625 55.039 c 
h
450.625 46.098 m f
465.121 49.465 m 465.121 48.855 465.121 48.855 464.465 48.855 c 463 50.27
 460.98 50.27 460.07 50.27 c 460.07 51.078 l 460.574 51.078 462.09 51.078
 463.301 50.473 c 463.301 61.938 l 463.301 62.695 463.301 62.996 461.082
 62.996 c 460.223 62.996 l 460.223 63.805 l 460.625 63.805 463.402 63.703
 464.211 63.703 c 464.918 63.703 467.746 63.805 468.25 63.805 c 468.25 62.996
 l 467.391 62.996 l 465.121 62.996 465.121 62.695 465.121 61.938 c h
465.121 49.465 m f
475.18 63.449 m 475.18 64.664 474.977 66.129 473.41 67.543 c 473.363 67.641
 473.262 67.742 473.262 67.844 c 473.262 67.996 473.461 68.148 473.563 68.148
 c 473.918 68.148 475.836 66.328 475.836 63.652 c 475.836 62.289 475.281
 61.23 474.219 61.23 c 473.512 61.23 472.957 61.785 472.957 62.492 c 472.957
 63.199 473.461 63.805 474.27 63.805 c 474.824 63.805 475.18 63.449 475.18
 63.449 c h
475.18 63.449 m f
489.441 59.715 m 488.684 59.715 l 488.633 60.219 488.379 61.531 488.078
 61.734 c 487.926 61.887 486.211 61.887 485.855 61.887 c 481.715 61.887 
l 484.09 59.816 484.895 59.16 486.211 58.098 c 487.875 56.785 489.441 55.371
 489.441 53.25 c 489.441 50.523 487.066 48.855 484.188 48.855 c 481.41 48.855
 479.492 50.828 479.492 52.898 c 479.492 54.008 480.453 54.16 480.703 54.16
 c 481.211 54.16 481.867 53.754 481.867 52.949 c 481.867 52.543 481.715 
51.734 480.555 51.734 c 481.262 50.172 482.773 49.664 483.836 49.664 c 486.109
 49.664 487.27 51.434 487.27 53.25 c 487.27 55.219 485.855 56.734 485.148
 57.543 c 479.746 62.945 l 479.492 63.148 479.492 63.199 479.492 63.805 
c 488.785 63.805 l h
489.441 59.715 m f
503.176 66.703 m 503.176 66.602 503.176 66.551 502.621 65.996 c 498.633
 61.957 497.57 55.844 497.57 50.898 c 497.57 45.293 498.785 39.688 502.773
 35.695 c 503.176 35.293 503.176 35.242 503.176 35.141 c 503.176 34.891 
503.074 34.789 502.875 34.789 c 502.52 34.789 499.641 37.012 497.723 41.102
 c 496.105 44.637 495.703 48.219 495.703 50.898 c 495.703 53.422 496.059
 57.309 497.824 60.996 c 499.793 64.934 502.52 67.055 502.875 67.055 c 503.074
 67.055 503.176 66.953 503.176 66.703 c h
503.176 66.703 m f
550.246 50.898 m 550.246 48.422 549.891 44.535 548.125 40.848 c 546.207
 36.91 543.43 34.789 543.125 34.789 c 542.922 34.789 542.773 34.941 542.773
 35.141 c 542.773 35.242 542.773 35.293 543.379 35.898 c 546.559 39.082 
548.375 44.18 548.375 50.898 c 548.375 56.453 547.215 62.105 543.227 66.148
 c 542.773 66.551 542.773 66.602 542.773 66.703 c 542.773 66.902 542.922
 67.055 543.125 67.055 c 543.43 67.055 546.355 64.832 548.227 60.742 c 549.891
 57.207 550.246 53.625 550.246 50.898 c h
550.246 50.898 m f
147.645 41.16 m 148.777 41.16 l 153.5 48.441 l 153.797 48.934 154.535 50.164
 154.879 50.605 c 155.125 51 155.223 51 155.91 51 c 160.242 51 l 160.832
 51 161.324 51 161.324 50.457 c 161.324 50.164 161.125 49.969 160.832 49.918
 c 159.648 49.672 158.125 47.555 157.387 46.57 c 157.141 46.277 155.566 
44.211 153.5 40.816 c 156.254 40.273 158.961 39.094 158.961 35.355 c 158.961
 31.023 154.387 29.5 150.648 29.5 c 140.809 29.5 l 140.215 29.5 139.723 
29.5 139.723 30.039 c 139.723 30.582 140.316 30.582 140.563 30.582 c 142.383
 30.582 142.527 30.828 142.527 32.402 c 142.527 48.098 l 142.527 49.672 
142.383 49.918 140.563 49.918 c 140.316 49.918 139.723 49.918 139.723 50.457
 c 139.723 51 140.215 51 140.809 51 c 149.418 51 l 149.957 51 150.449 51
 150.449 50.457 c 150.449 49.918 149.91 49.918 149.563 49.918 c 147.793 
49.918 147.645 49.672 147.645 48.098 c h
153.695 39.586 m 154.582 38.453 154.68 36.781 154.68 35.402 c 154.68 33.879
 154.484 32.207 153.5 30.926 c 154.781 31.223 157.879 32.207 157.879 35.355
 c 157.879 37.422 156.945 38.996 153.695 39.586 c h
147.645 32.305 m 147.645 31.617 147.645 30.582 149.613 30.582 c 152.27 
30.582 153.598 31.715 153.598 35.402 c 153.598 39.438 152.613 40.027 147.645
 40.027 c h
143.363 49.918 m 143.609 49.328 143.609 48.441 143.609 48.195 c 143.609
 32.305 l 143.609 32.008 143.609 31.172 143.363 30.582 c 146.957 30.582 
l 146.563 31.125 146.563 31.715 146.563 32.207 c 146.563 48.195 l 146.563
 48.492 146.563 49.328 146.809 49.918 c h
150.055 41.16 m 150.301 41.109 150.402 41.063 150.648 41.063 c 151.141 
41.063 151.828 41.012 152.27 40.965 c 152.762 41.703 156.156 47.309 158.766
 49.918 c 155.715 49.918 l h
150.055 41.16 m f
172.961 35.633 m 172.223 35.633 l 172.176 36.125 171.93 37.402 171.633 
37.602 c 171.484 37.746 169.813 37.746 169.469 37.746 c 165.434 37.746 l
 167.746 35.73 168.531 35.09 169.813 34.059 c 171.438 32.777 172.961 31.402
 172.961 29.336 c 172.961 26.676 170.648 25.055 167.844 25.055 c 165.137
 25.055 163.27 26.973 163.27 28.988 c 163.27 30.074 164.203 30.219 164.449
 30.219 c 164.941 30.219 165.582 29.828 165.582 29.039 c 165.582 28.645 
165.434 27.859 164.301 27.859 c 164.992 26.332 166.465 25.84 167.5 25.84
 c 169.715 25.84 170.844 27.563 170.844 29.336 c 170.844 31.254 169.469 
32.73 168.777 33.516 c 163.516 38.781 l 163.27 38.977 163.27 39.027 163.27
 39.617 c 172.32 39.617 l h
172.961 35.633 m f
342.43 59.309 m 343.563 59.309 l 348.285 66.59 l 348.582 67.082 349.316
 68.313 349.664 68.754 c 349.91 69.148 350.008 69.148 350.695 69.148 c 355.023
 69.148 l 355.617 69.148 356.109 69.148 356.109 68.609 c 356.109 68.313 
355.91 68.117 355.617 68.066 c 354.434 67.82 352.91 65.707 352.172 64.723
 c 351.926 64.426 350.352 62.359 348.285 58.965 c 351.039 58.422 353.746
 57.242 353.746 53.504 c 353.746 49.176 349.172 47.648 345.43 47.648 c 335.59
 47.648 l 335 47.648 334.508 47.648 334.508 48.191 c 334.508 48.73 335.098
 48.73 335.344 48.73 c 337.164 48.73 337.313 48.977 337.313 50.551 c 337.313
 66.246 l 337.313 67.82 337.164 68.066 335.344 68.066 c 335.098 68.066 334.508
 68.066 334.508 68.609 c 334.508 69.148 335 69.148 335.59 69.148 c 344.203
 69.148 l 344.742 69.148 345.234 69.148 345.234 68.609 c 345.234 68.066 
344.691 68.066 344.348 68.066 c 342.578 68.066 342.43 67.82 342.43 66.246
 c h
348.48 57.734 m 349.367 56.602 349.465 54.93 349.465 53.555 c 349.465 52.027
 349.27 50.355 348.285 49.074 c 349.563 49.371 352.664 50.355 352.664 53.504
 c 352.664 55.57 351.73 57.145 348.48 57.734 c h
342.43 50.453 m 342.43 49.766 342.43 48.73 344.398 48.73 c 347.055 48.73
 348.383 49.863 348.383 53.555 c 348.383 57.586 347.398 58.176 342.43 58.176
 c h
338.148 68.066 m 338.395 67.477 338.395 66.59 338.395 66.344 c 338.395 
50.453 l 338.395 50.156 338.395 49.32 338.148 48.73 c 341.742 48.73 l 341.348
 49.273 341.348 49.863 341.348 50.355 c 341.348 66.344 l 341.348 66.641 
341.348 67.477 341.594 68.066 c h
344.84 59.309 m 345.086 59.262 345.184 59.211 345.43 59.211 c 345.922 59.211
 346.613 59.16 347.055 59.113 c 347.547 59.852 350.941 65.461 353.551 68.066
 c 350.5 68.066 l h
344.84 59.309 m f
367.746 53.781 m 367.008 53.781 l 366.957 54.273 366.711 55.551 366.418
 55.75 c 366.27 55.898 364.598 55.898 364.254 55.898 c 360.219 55.898 l 
362.531 53.879 363.316 53.238 364.598 52.207 c 366.219 50.926 367.746 49.551
 367.746 47.484 c 367.746 44.828 365.434 43.203 362.629 43.203 c 359.922
 43.203 358.055 45.121 358.055 47.141 c 358.055 48.223 358.988 48.367 359.234
 48.367 c 359.727 48.367 360.367 47.977 360.367 47.188 c 360.367 46.793 
360.219 46.008 359.086 46.008 c 359.773 44.48 361.25 43.992 362.285 43.992
 c 364.5 43.992 365.629 45.711 365.629 47.484 c 365.629 49.402 364.254 50.879
 363.563 51.664 c 358.301 56.93 l 358.055 57.125 358.055 57.176 358.055 
57.766 c 367.105 57.766 l h
367.746 53.781 m f
515.816 50.664 m 516.945 50.664 l 521.672 57.945 l 521.965 58.438 522.703
 59.668 523.047 60.113 c 523.293 60.504 523.391 60.504 524.082 60.504 c 
528.41 60.504 l 529 60.504 529.492 60.504 529.492 59.965 c 529.492 59.668
 529.297 59.473 529 59.422 c 527.82 59.176 526.297 57.063 525.559 56.078
 c 525.313 55.781 523.738 53.715 521.672 50.32 c 524.426 49.781 527.133 
48.598 527.133 44.859 c 527.133 40.531 522.555 39.004 518.816 39.004 c 508.977
 39.004 l 508.387 39.004 507.895 39.004 507.895 39.547 c 507.895 40.086 
508.484 40.086 508.73 40.086 c 510.551 40.086 510.699 40.332 510.699 41.906
 c 510.699 57.602 l 510.699 59.176 510.551 59.422 508.73 59.422 c 508.484
 59.422 507.895 59.422 507.895 59.965 c 507.895 60.504 508.387 60.504 508.977
 60.504 c 517.586 60.504 l 518.129 60.504 518.621 60.504 518.621 59.965 
c 518.621 59.422 518.078 59.422 517.734 59.422 c 515.965 59.422 515.816 
59.176 515.816 57.602 c h
521.867 49.09 m 522.754 47.961 522.852 46.285 522.852 44.91 c 522.852 43.383
 522.656 41.711 521.672 40.434 c 522.949 40.727 526.051 41.711 526.051 44.859
 c 526.051 46.926 525.113 48.5 521.867 49.09 c h
515.816 41.809 m 515.816 41.121 515.816 40.086 517.785 40.086 c 520.441
 40.086 521.77 41.219 521.77 44.91 c 521.77 48.945 520.785 49.535 515.816
 49.535 c h
511.535 59.422 m 511.781 58.832 511.781 57.945 511.781 57.699 c 511.781
 41.809 l 511.781 41.516 511.781 40.68 511.535 40.086 c 515.125 40.086 l
 514.734 40.629 514.734 41.219 514.734 41.711 c 514.734 57.699 l 514.734
 57.996 514.734 58.832 514.98 59.422 c h
518.227 50.664 m 518.473 50.617 518.57 50.566 518.816 50.566 c 519.309 
50.566 519.996 50.52 520.441 50.469 c 520.934 51.207 524.328 56.816 526.934
 59.422 c 523.883 59.422 l h
518.227 50.664 m f
541.133 45.137 m 540.395 45.137 l 540.344 45.629 540.098 46.91 539.805 
47.105 c 539.656 47.254 537.98 47.254 537.637 47.254 c 533.605 47.254 l 
535.914 45.234 536.703 44.598 537.98 43.563 c 539.605 42.285 541.133 40.906
 541.133 38.84 c 541.133 36.184 538.82 34.559 536.016 34.559 c 533.309 34.559
 531.438 36.477 531.438 38.496 c 531.438 39.578 532.375 39.727 532.621 39.727
 c 533.113 39.727 533.75 39.332 533.75 38.543 c 533.75 38.152 533.605 37.363
 532.473 37.363 c 533.16 35.84 534.637 35.348 535.672 35.348 c 537.883 35.348
 539.016 37.07 539.016 38.84 c 539.016 40.758 537.637 42.234 536.949 43.023
 c 531.684 48.285 l 531.438 48.484 531.438 48.531 531.438 49.121 c 540.492
 49.121 l h
541.133 45.137 m f
718.219 34.191 m 719.352 34.191 l 724.074 41.473 l 724.367 41.965 725.105
 43.195 725.449 43.641 c 725.695 44.031 725.797 44.031 726.484 44.031 c 
730.813 44.031 l 731.402 44.031 731.895 44.031 731.895 43.492 c 731.895 
43.195 731.699 43 731.402 42.949 c 730.223 42.703 728.699 40.59 727.961 
39.605 c 727.715 39.309 726.141 37.242 724.074 33.848 c 726.828 33.309 729.535
 32.125 729.535 28.387 c 729.535 24.059 724.961 22.531 721.219 22.531 c 
711.379 22.531 l 710.789 22.531 710.297 22.531 710.297 23.074 c 710.297 
23.613 710.887 23.613 711.133 23.613 c 712.953 23.613 713.102 23.859 713.102
 25.434 c 713.102 41.129 l 713.102 42.703 712.953 42.949 711.133 42.949 
c 710.887 42.949 710.297 42.949 710.297 43.492 c 710.297 44.031 710.789 
44.031 711.379 44.031 c 719.988 44.031 l 720.531 44.031 721.023 44.031 721.023
 43.492 c 721.023 42.949 720.48 42.949 720.137 42.949 c 718.367 42.949 718.219
 42.703 718.219 41.129 c h
724.27 32.617 m 725.156 31.488 725.254 29.813 725.254 28.438 c 725.254 
26.91 725.059 25.238 724.074 23.961 c 725.352 24.254 728.453 25.238 728.453
 28.387 c 728.453 30.453 727.516 32.027 724.27 32.617 c h
718.219 25.336 m 718.219 24.648 718.219 23.613 720.188 23.613 c 722.844
 23.613 724.172 24.746 724.172 28.438 c 724.172 32.473 723.188 33.063 718.219
 33.063 c h
713.938 42.949 m 714.184 42.359 714.184 41.473 714.184 41.227 c 714.184
 25.336 l 714.184 25.043 714.184 24.207 713.938 23.613 c 717.531 23.613 
l 717.137 24.156 717.137 24.746 717.137 25.238 c 717.137 41.227 l 717.137
 41.523 717.137 42.359 717.383 42.949 c h
720.629 34.191 m 720.875 34.145 720.973 34.094 721.219 34.094 c 721.711
 34.094 722.402 34.047 722.844 33.996 c 723.336 34.734 726.73 40.344 729.336
 42.949 c 726.289 42.949 l h
720.629 34.191 m f
743.535 28.664 m 742.797 28.664 l 742.746 29.156 742.5 30.438 742.207 30.633
 c 742.059 30.781 740.387 30.781 740.043 30.781 c 736.008 30.781 l 738.32
 28.762 739.105 28.125 740.387 27.09 c 742.008 25.813 743.535 24.434 743.535
 22.367 c 743.535 19.711 741.223 18.086 738.418 18.086 c 735.711 18.086 
733.844 20.004 733.844 22.023 c 733.844 23.105 734.777 23.254 735.023 23.254
 c 735.516 23.254 736.152 22.859 736.152 22.07 c 736.152 21.68 736.008 20.891
 734.875 20.891 c 735.563 19.367 737.039 18.875 738.074 18.875 c 740.285
 18.875 741.418 20.598 741.418 22.367 c 741.418 24.285 740.043 25.762 739.352
 26.551 c 734.086 31.813 l 733.844 32.012 733.844 32.059 733.844 32.648 
c 742.895 32.648 l h
743.535 28.664 m f
Q
}}
\wd\mobiusintervalbox=790.807312bp
\ht\mobiusintervalbox=226.606277bp
\dp\mobiusintervalbox=0bp
\immediate\pdfxform attr{/Group <<
/Type /Group
/S /Transparency
/I true
/CS /DeviceRGB
>>} resources{
/ExtGState <<
/a0 <<
/CA 1
/ca 1
>>
/a1 <<
/CA 0.507
/ca 0.507
>>
>>
/XObject << /x6 \mobiusintervalobjecta\space 0 R /x7 \mobiusintervalobjectb\space 0 R /x8 \mobiusintervalobjectc\space 0 R /x9 \mobiusintervalobjectd\space 0 R >>
}\mobiusintervalbox
\edef\mobiusintervalform{\the\pdflastxform}
\newcommand{\mobiusintervalfigure}{\pdfrefxform\mobiusintervalform}

\begin{document}
\maketitle

\begin{abstract}

One of the best known invariants of multi-parameter persistent homology is the rank invariant, which is subsumed by the generalized rank invariant (GRI), obtained by naturally extending its domain.
M\"obius inversion of the GRI not only provides a compact encoding of this invariant but is also useful for clarifying its relationship with other invariants of multi-parameter persistent homology.
However, M\"obius inversion of the GRI has been defined only under suitable finiteness conditions on its domain and support.
This work studies the extent to which M\"obius inversion can be extended beyond these finiteness conditions. 
Namely, we formalize the notion of Möbius invertibility, which captures when Möbius inversion exists beyond these finiteness conditions, and investigate when the GRI is Möbius invertible.

A central result we establish is that M\"obius invertibility of the GRI is guaranteed for finitely presentable $d$-parameter persistence modules when $d\leq 2$, but may fail when $d>2$. Establishing this failure also corrects an error in the existing literature.
In response to this higher-dimensional obstruction, we show that for finitely presentable $d$-parameter persistence modules, restricting the GRI to intervals subordinate to the finite grid supporting a minimal free presentation preserves its M\"obius inversion whenever the latter exists over all intervals. Since computing the GRI over all intervals in $\R^d$ is infeasible in practice, this result provides a principled justification for restricting the GRI to those intervals.

\end{abstract}

\setlength{\nomitemsep}{-\parsep}%

\kern0pt
\tableofcontents

\section{Introduction}\label{section:introduction}

One-parameter persistent homology is a central concept in topological data analysis (TDA) with a wide range of applications \cite{bleher2021topological,cang2015topological,carlsson2021topological,chan2013topology,chowdhury2018importance,dabaghian2012topological,PhysRevC.106.064912,lee2018high,stolz2023relational}. Under fairly general assumptions, a \emph{one-parameter persistence module}, by definition a functor $M:\R\to \Cvec$ from the reals to the category of vector spaces and linear maps over a fixed field $k$, is completely characterized by its  
\emph{persistence diagram} (or \emph{persistence barcode}), which provides a compact and interpretable summary of $M$ \cite{abeasis1981geometry,cohen2007stability,zomorodian2005computing}. 

However, many data sets naturally give rise to more than one filtration parameter, for example scale together with density or another auxiliary measurement. This has motivated extensive research on \emph{multi-parameter} persistent homology, i.e., the setting in which we contend with \emph{multi-parameter} persistence modules $M:\R^d\to \Cvec$ (where $\R^d$ is equipped with the product partial order).
Unfortunately, in contrast with the case of one-parameter persistence modules (where persistence diagrams are full invariants thereof),  there is no complete and simultaneously discrete invariant for multi-parameter persistence modules \cite{carlsson2009theory}.
Accordingly, many discrete and necessarily incomplete invariants have been studied, e.g., \cite{asashiba2023approximation,asashiba2019approximation,blanchette2021homological,botnan2024signed,chacholski2022effective,harrington2019stratifying,landi2018rank,lesnick2015interactive,miller2020homological,oudot2024stability}. 
Then, a central question is how to balance the inherent tension between the discriminating power of an invariant and its computational complexity.

Among these invariants, one that is intimately tied to this work is the classical \emph{rank invariant} (RI). Given a $\catP$-module, i.e., a functor $M:\catP\to \Cvec$ for an arbitrary poset $\catP$, its RI records the ranks of all linear maps in $M$; see \cite{carlsson2009theory,puuska2020erosion}.
For one-parameter persistence modules, the RI contains equivalent information to the classical persistence diagram \cite{abeasis1981geometry,carlsson2009theory}, making it a natural notion to import into the multi-parameter setting.
The RI allows us to quantify persistence over \emph{segments} in the indexing poset (a segment is the set of points lying between a given pair of comparable points), and it is a lossless invariant when applied to so-called rectangle-decomposable modules \cite{botnan2022rectangle}.
As a downside, the RI \emph{only} allows us to quantify persistence over segments, and in the multi-parameter setting there {typically are many}  natural non-segment subsets of the indexing poset $\catP$ over which we desire to quantify persistence.

Motivated by this, the RI has recently evolved into the \emph{generalized rank invariant} (GRI)  by extending the domain of the RI from the set $\seg(\catP)$ of segments of the indexing poset $\catP$ to the set $\Int(\catP)$ of \emph{intervals} of $\catP$ or to the even larger set $\Con(\catP)$ of \emph{connected subposets} of $\catP$ \cite{kim2021generalized} (see also \cite{asashiba2019approximation}).
By recording the rank of the canonical limit-to-colimit map of the diagram over each element of $\Int(\catP)$ (or of $\Con(\catP)$), the GRI  quantifies persistence across a broader range of regions in the indexing poset compared to the RI and, therefore, the GRI acquires more discriminating power than the RI. 

However, the size of $\Int(\catP)$, let alone that of the larger $\Con(\catP)$, is generally formidable, causing a bottleneck for computing the GRI over the entire $\Int(\catP)$ or $\Con(\catP)$. 
For example, if $\catP$ is the 2-dimensional 10-by-10 grid (which serves as a much simplified setup for 2-parameter persistent homology), then $\Int(\catP)$ comprises 1,497,925,315 elements \cite[Theorem 31]{asashiba2022interval}.  %
This has been addressed in the literature through the following strategy:

\begin{enumerate}[label=(S\arabic*),leftmargin=*]
    \item \label{S1}
    restrict the GRI to a 
    tractable family of regions. 
\end{enumerate} 
Indeed, this strategy has already been employed in practice; see, e.g., \cite{xin2023gril,carriere2024sparsification,kim2025super}. 
In these works, however, the choice of restricted domain is driven primarily by computational considerations, rather than by a mathematical justification of the restriction
itself. In this paper we consider the following stricter  strategy which can be seen as an extremal form of \ref{S1}:

\begin{enumerate}[label=(S\arabic*),leftmargin=*,resume]
    \item seek a lossless compact encoding of the full GRI;\label{S2}
\end{enumerate} 
When the relevant poset of regions satisfies suitable local-finiteness
conditions,   \ref{S2} is realized by the \emph{generalized persistence diagram}
\cite{kim2021generalized}, obtained through M\"obius inversion of the GRI.
In that setting,
the generalized persistence diagram and the
GRI determine each other, while the support of the former is contained in
that of the latter (cf. Proposition~\ref{prop:properties of GPD}%
~\ref{item:GPD is more parsimonious}). Thus, the generalized persistence
diagram provides a lossless representation that is at least as parsimonious
as the GRI. This alternative representation has also been used to analyze
the discriminating power of the GRI and compare it with other invariants;
see, e.g., \cite{clause2026generalized,kim2021generalized,kim2021bettis}. 

For many naturally occurring indexing posets, however, the relevant region poset is not locally finite. This leads to the following question.

\begin{quote}
\hypertarget{q:2}{\textbf{Motivating Question.}} 
\emph{When can the RI or GRI be compressed, either through a persistence-diagram-type M\"obius inversion or through a canonical finite restriction of its domain that minimizes information-loss?}
\end{quote}

In order to tackle this question, we formalize the concept of \emph{M\"obius invertibility}.
While M\"obius inversion was first introduced to one-parameter persistence by Patel (see \hyperref[sec:related work]{Related Work} below), the concept of M\"obius \emph{invertibility}—a property, not a procedure—is formalized and systematically studied in this work as a key contribution. 

For persistence modules over finite posets, the GRI is always M\"obius invertible, and its M\"obius inverse is the generalized persistence diagram \cite{kim2021generalized} (see also \cite{asashiba2019approximation}), a parsimonious encoding of the GRI. The situation changes sharply beyond the finite setting. When the indexing poset, and hence the relevant collection of intervals or connected subposets, is infinite, even persistence modules with strong finiteness properties can fail to have M\"obius invertible GRIs. In such cases, the GRI remains defined but its generalized persistence diagram need not exist. 

We give a systematic account of this phenomenon and, in trying to answer \hyperlink{q:2}{Motivating Question}, propose a natural module-dependent finite set on which to restrict the GRI without loss, whenever this is possible.  More precisely, we show that for a finitely presentable $d$-parameter persistence module, restricting its GRI to intervals subordinate to a certain finite grid (i.e., the unique finite grid supporting a minimal free presentation), preserves its Möbius inversion whenever the latter exists over all intervals (Theorem~\ref{thm:unique candidate}). This result provides a mathematically justified restriction of the GRI to a tractable family of regions, thus it also provides a structural justification for methods such as \cite{kim2025super} that implement \ref{S1}.

Beyond its role in answering the motivating question, M\"obius invertibility is also useful for contextualizing known results regarding the RI or GRI within the framework of classical M\"obius inversion theory, which also enables us to simplify existing proofs, rectify errors in the literature, or strengthen statements of existing theorems.

\subsection*{Our Contributions}

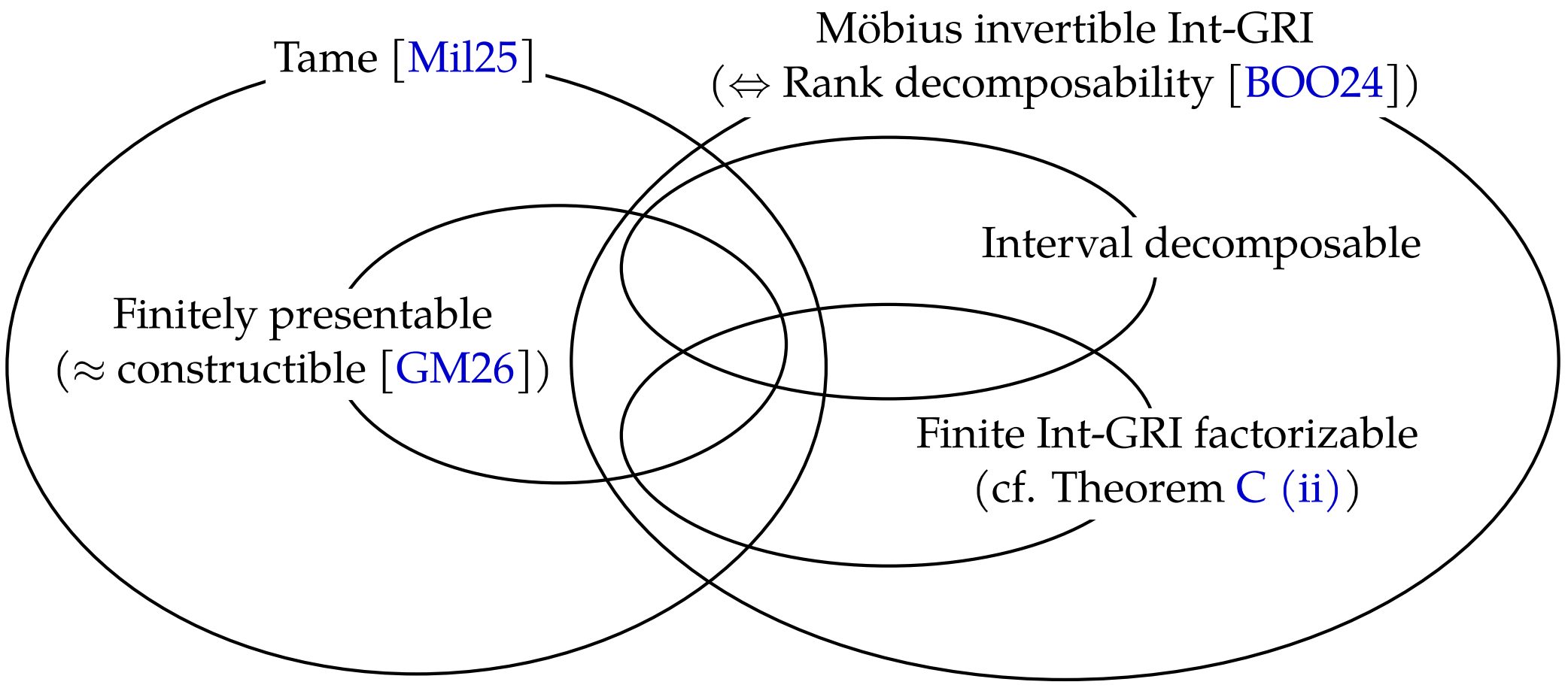
\begin{figure}
\begin{center}
\begin{tikzpicture}[scale=1.0]

\tikzset{
  classlabel/.style={
    align=center,
    fill=white,
    inner sep=2pt
  },
  regnum/.style={
    font=\small\bfseries,
    fill=white,
    inner sep=1pt
  }
}

\draw[thick]%
    (-2.6,0.15) ellipse (3.6 and 2.7);
\node[classlabel]%
    at (-2.7,2.85)
    {Tame \cite{miller2020homological}};

\draw[thick]%
    (3.1,0.20) ellipse (4.34 and 2.8);

\node[classlabel]%
    at (3.1,2.85)
    {M\"obius invertible Int-GRI \\ ($\Leftrightarrow$ Rank decomposability \cite{botnan2024signed})};

\draw[thick]%
    (-1.35,0.35) ellipse (2 and 1.22);

\node[classlabel]      at (-3.6,0.35)
    {Finitely presentable\\
     ($\approx$ constructible \cite{gulen2022galois})};

\draw[thick]%
    (1.55,1.02) ellipse (2.35 and 1.15);
\node[classlabel]
    at (4.3,1.2)
    {Interval decomposable};

\draw[thick]%
    (1.55,-0.45) ellipse (2.35 and 1.15);
\node[classlabel]%
    at (4,-0.70)
    {Finite Int-GRI factorizable\\ (cf. Theorem \ref{thm:sufficient conditions for invertibility}~\ref{item:finitely encodable})};

\end{tikzpicture}
\captionsetup{singlelinecheck=off}
\caption{
Inclusion relations between M\"obius invertibility of the Int-GRI and several notions of structural simplicity for persistence modules over a poset. The indicated implications and non-implications follow from Proposition~\ref{prop:finitely presentable implies tame} (known), Theorems~\ref{thm:rk decomposition is GPD}, \ref{thm:tame does not imply Int-GRI invertible}, and \ref{thm:sufficient conditions for invertibility}, Remark~\ref{rem:implication of the counterexample}~\ref{item:not rank-preserving}, and Corollary~\ref{cor:Mobius invertibility does not imply tameness}. Constructibility and finite presentability are not identical, but they are closely related; see Proposition~\ref{prop:constructible implies finitely presentable} and \cite[Proposition~C.2]{kim2026interleaving}.}
\label{fig:tameness + MI visual}
\end{center}
\end{figure}

    \begin{enumerate}[leftmargin=*,label=(\roman*)] 
    \item 
    We formalize the notion of M\"obius invertibility (\emph{Definition~\ref{def:mobius invertible general}}). Despite its added generality over classical M\"obius inversion, the resulting theory closely parallels the classical theory and also reveals a close connection between \ref{S1} and~\ref{S2}  (see \emph{Proposition~\ref{prop:justifying the term mobius invertibility}}). In particular, by dispensing with the local finiteness assumption on the indexing poset, this framework permits the generalized persistence diagram to be defined for modules indexed by non-discrete posets, such as $\R^d$-modules (\emph{Definition~\ref{def:GGPD}}).%
    \item 
   The extended notion of generalized persistence diagram naturally includes \emph{rank decompositions} \cite{botnan2024signed} as a special case (\emph{Theorem~\ref{thm:rk decomposition is GPD}}). At the same time, we show that M\"obius invertibility of the RI or GRI is distinct from other common finiteness conditions for persistence modules, including Miller's tameness \cite{miller2020homological} and finite presentability (\emph{Proposition~\ref{prop:tame does not imply RI invertible}} and \emph{Theorem \ref{thm:tame does not imply Int-GRI invertible}}). In particular, Theorem~\ref{thm:tame does not imply Int-GRI invertible}~\ref{item:fp not imply MI} corrects \cite[Corollary~3.18]{botnan2024signed}; see Remark~\ref{rem:conflict} below.
\end{enumerate}
We henceforth refer to the GRI over intervals as the \emph{Int-GRI}.
\begin{enumerate}[resume,label=(\roman*)]
      \item  
      We establish several sufficient conditions for M\"obius invertibility of the GRI (\emph{Theorems~\ref{thm:sufficient conditions for invertibility} and \ref{thm:pull-back GPD}} and \emph{Proposition~\ref{prop:finitely presentable invertiver over intmn}}). For instance, the Int-GRI of every finitely presentable $\R^2$-module is M\"obius invertible.
\end{enumerate}    
Some of the results pertaining to Theorems~\ref{thm:rk decomposition is GPD}, \ref{thm:tame does not imply Int-GRI invertible}, and \ref{thm:sufficient conditions for invertibility} are summarized in Figure~\ref{fig:tameness + MI visual}.

\begin{enumerate}[resume,label=(\roman*)]    
\item 
For finitely presentable $\R^d$-modules, we identify an optimal restriction of the domain of the GRI, in a suitable sense. Namely, we show that if the Int-GRI of a finitely presentable $\R^d$-module is M\"obius invertible, then its M\"obius inversion must coincide with the inversion over the intervals subordinate to the finite grid supporting a minimal free presentation (\emph{Definition~\ref{def:GPD for fp modules}} and \emph{Theorem~\ref{thm:unique candidate}}).
    \end{enumerate}

\begin{remark}[\textbf{{A dimension-dependent dichotomy and a correction to the literature}}]
\label{rem:conflict} 
In \cite[Corollary~3.18]{botnan2024signed}, the authors state that, for any finitely presentable $\R^d$-module with $d\geq 1$, the GRI over {half-open} intervals admits a rank decomposition. Theorem~\ref{thm:rk decomposition is GPD} implies that rank decomposability is equivalent to M\"obius invertibility of 
the GRI over half-open intervals.
Thus \cite[Corollary~3.18]{botnan2024signed} would imply that the GRI {over half-open intervals}
of every finitely presentable $\R^d$-module is M\"obius invertible in all dimensions.

Our earlier version
\cite[Theorem~C(iii)]{clause2024generalized_v5} likewise asserted
M\"obius invertibility for finitely presentable $\R^d$-modules in
arbitrary dimension, but over all intervals rather than half-open
intervals. The counterexamples developed here also refute that
assertion.

We show, however, in Appendix~\ref{sec:details on prior errors} that there exists a finitely presentable $\R^3$-module whose GRI over half-open intervals is not M\"obius invertible, using the technique employed in the proof of Theorem~\ref{thm:tame does not imply Int-GRI invertible}~\ref{item:fp not imply MI}. We also trace the discrepancy to an error in the proof of \cite[Corollary~3.18]{botnan2024signed} and identify counterexamples to several statements from \cite{botnan2024signed} on which that corollary depends. On the positive side, Theorem~\ref{thm:sufficient conditions for invertibility}~\ref{item:finitely presentable invertible} establishes M\"obius invertibility for finitely presentable $\R^d$-modules when $d=1,2$. Both the negative statement for $d\geq 3$ and the positive statement for $d=1,2$ require separate geometric and  combinatorial arguments.
\end{remark}

\begin{remark}[Chronology] \label{rem:chrono}
M\"obius invertibility of the GRI, as opposed to the RI, for $\R^d$-modules was first considered in our arXiv preprint \cite{clause2022generalized} and carried forward in \cite{clause2024generalized_v5}, both of which predate \cite{botnan2024signed}. The earlier arXiv version \cite{botnan2021signed} did not contain the developments later appearing in \cite[Section~3]{botnan2024signed}, including \cite[Corollary~3.18]{botnan2024signed}, concerning rank decomposability of the GRI.

We first constructed the counterexamples discussed in
Remark~\ref{rem:conflict}, refuting both
\cite[Corollary~3.18]{botnan2024signed} and our own
\cite[Theorem~C(iii)]{clause2024generalized_v5}, in August 2025.
\end{remark}

\subsubsection*{Related work.}\label{sec:related work} 
Patel first observed that in the one-parameter setting the persistence diagram can be defined as the M\"obius inversion of the RI, and used this perspective to introduce generalized persistence diagrams for constructible $\R$-indexed functors with target categories beyond vector spaces \cite{patel2018generalized}. This perspective motivated subsequent work by Kim and M\'emoli \cite{kim2021generalized} and by McCleary and Patel \cite{mccleary2020bottleneck,mccleary2022edit}.
In particular, in \cite{kim2021generalized}, the generalized persistence diagram of a $\catP$-module is defined as the M\"{o}bius inversion of the GRI over $\Con(\catP)$, which can be viewed as a multiset of signed elements of $\Con(\catP)$ {(assuming that $\Con(\catP)$ is locally finite)}. 

\medskip
Several further lines of work are %
connected to the topics of
the present paper:
\begin{itemize}
\item 
Asashiba et al. also invoke M\"{o}bius inversion to devise methods for approximating a persistence module $M$ over a finite two-dimensional grid by an interval-decomposable module \cite{asashiba2019approximation}.
One of their approximation methods yields an invariant that possesses the same amount of information as 
(the M\"obius inversion of) the GRI of $M$ over the intervals in the grid \cite[Remark 2.19]{kim2021bettis}. 
See also \cite{asashiba2024interval,hiraoka2023refinement} for further developments of these ideas.

\item A characterization of the generalized persistence diagram that does not invoke M\"obius inversion is given in \cite[Theorem~4.11]{escolar2024barcoding}.

\item 
Connections between relative homological algebra and rank-based invariants such as the RI and GRI have been explored in \cite{asashiba2023approximation,blanchette2021homological,botnan2024signed,botnan2022bottleneck}.

\item In recent years, M\"obius inversion has been utilized alongside other invariants of persistence modules, including the graded rank function \cite{betthauser2022graded}, the birth-death function \cite{gulen2022galois,mccleary2022edit,morozov2021output}, the meta-rank \cite{clause2023meta}, the Hilbert function \cite{oudot2024stability},
the multi-rank \cite{thomas2019invariants}, the persistent cup length and persistent cup product \cite{memoli2022persistent}, and the Grassmanian persistence diagram \cite{gulen2023orthogonal}.  
\item  
Dey et al. propose a method based on zigzag persistence for computing the rank of the limit-to-colimit map (cf. Section~\ref{sec:GRI}) for persistence modules over posets more general than $\Z^2$ \cite{dey2024computing}.
See also \cite{dey2026limit}.

\item Finally, we remark that our work lays the foundation for two recent studies \cite{carriere2024sparsification,kim2025super} that delve into the computation and estimation of M\"obius inversion of the GRI for multi-parameter persistence.

\end{itemize}

\subsubsection*{Acknowledgements.} This research was partially supported by the NSF through grants 
DMS-1723003, DMS-2301359, CCF-2310412, IIS-1901360, CCF-1740761, DMS-1547357 and by the National Research Foundation of Korea (NRF) grant funded by the Korea government 
(MSIT and MOE) (RS-2025-00515946 and RS-2025-25397599). 
This research was also partially supported by the \emph{2026 KAIST-U.S. Joint Research Collaboration
Open Track Project for Early-Career Researchers} of the International Office at KAIST.
The authors thank anonymous reviewers for sharing their insight. %
AI-assisted tools were used in the preparation of this manuscript for editorial support, including wording suggestions, proofreading, and checking consistency of exposition. All mathematical content, arguments, and conclusions were reviewed and verified by the authors.

\makeatletter
\nom@tempdim=1in\relax
\makeatother

\begin{thenomenclature} 
\nomgroup{A}
  \item [{\([n]\)}]\begingroup The totally ordered set $\{1<\ldots<n\}$.\nomeqref {0.0}\nompageref{1}
  \item [{\(\catP\)}]\begingroup A poset \(\catP=(\catP,\leq)\), regarded as a category.\nomeqref {0.0}\nompageref{1}
  \item [{\(\Con(\catP)\)}]\begingroup Poset of connected subposets of \(\catP\), ordered by containment (pg. \pageref{nom:Con(P)}).\nomeqref {0.0}\nompageref{1}
  \item [{\(\seg(\catP)\)}]\begingroup Poset of segments of the poset \(\catP\), ordered by containment (Equation \ref{eq:segment}).\nomeqref {0.0}\nompageref{1}
  \item [{\(\Int(\catP)\)}]\begingroup Poset of intervals of the poset \(\catP\), ordered by containment (Definition \ref{def:intervals}).\nomeqref {0.0}\nompageref{1}
  \item [{$\Int_{m,n}(\catP)$}]\begingroup Poset of intervals of $\catP$ with at most $m$ minimal points and $n$ maximal points (pg. \pageref{eq:intmn}).\nomeqref {0.0}\nompageref{1}
\nomgroup{B}
  \item [{RI}]\begingroup Rank invariant (Definition \ref{def:generalized rank invariant}).\nomeqref {0.0}\nompageref{1}
  \item [{GRI}]\begingroup Generalized rank invariant (Definition {\ref{def:generalized rank invariant}}).\nomeqref {0.0}\nompageref{1}
  \item [{$\rk(M)$}]\begingroup The (generalized) rank of $M$ (pg. \pageref{nom:generalized rank}).\nomeqref {0.0}\nompageref{1}
  \item [{$\rk_M^\Ical$}]\begingroup The generalized rank invariant of $M$ over $\Ical\subset\Con(\catP)$ (Definition \ref{def:generalized rank invariant}).\nomeqref {0.0}\nompageref{1}
  \item [{$\rk_M^{\Int}$}]\begingroup The generalized rank invariant of $M$ over $\Int(\catP)$ (Definition \ref{def:generalized rank invariant}). \nomeqref {0.0}\nompageref{1}
  \item [{$\dgm_M^\Ical$}]\begingroup The generalized persistence diagram of $M$ over $\Ical$ (Definition \ref{def:GPD over I}).\nomeqref {0.0}\nompageref{1}
\nomgroup{C}
  \item [{$I(\catQ,\kf)$}]\begingroup The incidence algebra of $\catQ$ over $\kf$ (pg. \pageref{nom:incidence algebra}).\nomeqref {0.0}\nompageref{1}
  \item [{$\delta_\catQ$}]\begingroup The Dirac delta function (Equation (\ref{eq:delta function})).\nomeqref {0.0}\nompageref{1}
  \item [{$\zeta_\catQ$}]\begingroup The zeta function (Equation (\ref{eq:zeta function})).\nomeqref {0.0}\nompageref{1}
  \item [{$\mu_\catQ$}]\begingroup The M\"obius function (Equation (\ref{eq:mobius})).\nomeqref {0.0}\nompageref{1}
  \item [{$\mult_{\Rcal}$}]\begingroup The multiplicity function for the multiset $\Rcal$ (pg. \pageref{nom:mult R}).\nomeqref {0.0}\nompageref{1}
\nomgroup{D}
  \item [{$\barc(M)$}]\begingroup The barcode of $M$ (pg. \pageref{nom:barcode}).\nomeqref {0.0}\nompageref{1}
  \item [{$\kf^\catQ$}]\begingroup The vector space of functions $f:\catQ\to \kf$ (pg. \pageref{nom: kq}).\nomeqref {0.0}\nompageref{1}
  \item [{$\kf_I$}]\begingroup Interval module for the interval $I$ (Equation (\ref{eq:interval module})).\nomeqref {0.0}\nompageref{1}
  \item [{$p^\downarrow_\catP$}]\begingroup The set of elements in $\catP$ less than or equal to $p$ (pg. \pageref{nom: down p}).\nomeqref {0.0}\nompageref{1}
  \item [{$p^\uparrow_\catP$}]\begingroup The set of elements in $\catP$ greater than or equal to $p$ (pg. \pageref{nom: up p}).\nomeqref {0.0}\nompageref{1}

\end{thenomenclature}

\section{Preliminaries}\label{section:prelim}

This section is organized as follows: in Section \ref{sec:persistence modules and their decomposition}, we review the concepts of persistence modules and their structures. In Section \ref{sec:mobius inversion}, we review the notion of the incidence algebra as well as the M\"obius inversion formula. In Section \ref{sec:GRI}, we recall the notions of the rank invariant, the generalized rank invariant, and their properties. In Section \ref{sec:GPD}, we recall the notion of the generalized persistence diagram. In Section \ref{sec:GPD and Rk decomposition}, we recall the notion of the rank decomposition and its connection with the notion of the generalized persistence diagram. 

\subsection{Persistence modules}\label{sec:persistence modules and their decomposition}

Throughout this paper, $\catP=(\catP,\leq)$ is a poset, regarded as the category  with objects the elements $p\in \catP$, and a unique morphism $p\to q$ if and only if $p\leq q\in \catP$.
All vector spaces in this paper are over a fixed field $\kf$. 
Let $\cvec$ denote the category of finite-dimensional vector spaces and linear maps over $\kf$.
\paragraph{Persistence modules and their decompositions.} 
 
A persistence module over $\catP$ is a functor $M:\catP\to \cvec$.\footnote{In the literature, $M$ is often referred to as a pointwisely finite-dimensional persistence module.} We refer to $M$ simply as a \textbf{$\catP$-module}. 
For any $p\in \catP$, we denote the vector space $M_p:=M(p)$, and for any $p\leq q\in \catP$, we denote the linear map $\varphi_M(p,q):=M(p\leq q)$. 
Given any $\catP$-modules $M$ and $N$, the \text{direct sum} $M\oplus N$ is defined pointwisely at each $p\in \catP$.
We say that a nontrivial $\catP$-module $M$ is \textbf{decomposable} if $M$ is isomorphic to $N_1\oplus N_2$ for some non-trivial $\catP$-modules $N_1$ and $N_2$, which we denote by $M\cong N_1\oplus N_2$. 
Otherwise, we say that $M$ is \textbf{indecomposable}. 

By the Azumaya-Krull-Remak-Schmidt Theorem \cite[Theorem 1]{azumaya1950corrections}, every $\catP$-module is isomorphic to a direct sum of indecomposable $\catP$-modules; see also {\cite[Theorem 1.1]{botnan2020decomposition}}.  
This direct sum decomposition is unique up to isomorphism and permutations of summands. The multiset of isomorphism classes of indecomposable summands of $M$ is called the \textbf{barcode} of $M$, denoted by $\barc(M)$.

\begin{definition}[\cite{botnan2018algebraic}] \label{def:intervals} An \textbf{interval} of a poset $\catP$ is any nonempty subset $I\subset \catP$ such that 
\begin{enumerate}[label=(\roman*),leftmargin=*]
    \item  (\textbf{convexity}) If $p,r\in I$ and $q\in \catP$ with $p\leq q\leq r$, then $q\in I$, \label{item:convexity}
    \item (\textbf{connectivity}) For any $p,q\in I$, there is a sequence $p=r_0,r_1,\ldots,r_n=q$ of elements of $I$, where $r_i$ and $r_{i+1}$ are comparable for $0\leq i\leq n-1$. {Such a sequence is called a \textbf{path.}}
    \label{item:connectivity}
\end{enumerate}
\end{definition}

Given an interval $I$ of $\catP$, the \textbf{interval module} $\kf_I$ is the $\catP$-module, with
\begin{equation}\label{eq:interval module} (\kf_I)_p := \begin{cases} \kf & \mathrm{if \ } p\in I\\
0 & \mathrm{otherwise}
\end{cases},
\hspace{10mm}
\varphi_{\kf_I}(p,q) := \begin{cases} \id_\kf & \mathrm{if \ } p\leq q\in I\\
0 & \mathrm{otherwise.}\end{cases}\end{equation}

Every interval module is indecomposable \cite[Proposition 2.2]{botnan2018algebraic}. 
{\label{nom:barcode} 
A $\catP$-module $M$ is \textbf{interval-decomposable} if it is isomorphic to a direct sum of interval modules.}
If {the equivalence class of the interval module $\kf_I$ belongs to} $\barc(M)$, then {we simply say that $I$ belongs to $\barc(M)$ and write $I\in \barc(M)$.}  
A \textbf{zigzag poset} of $n$ points is 
    $\bullet_{1}\leftrightarrow \bullet_{2} \leftrightarrow \ldots \bullet_{n-1} \leftrightarrow \bullet_n$
where $\leftrightarrow$ stands for either $\leq$ or $\geq$. A functor from a zigzag poset (of $n$ points) to $\cvec$ is called a \textbf{zigzag module} \cite{carlsson2010zigzag}. 

\begin{oldtheorem}[\cite{carlsson2010zigzag,gabriel1972unzerlegbare}]\label{thm:zigzag modules are interval-decomposable}
Zigzag modules are interval-decomposable.
\end{oldtheorem}

\paragraph{Tame, Finitely presentable, and Constructible $\catP$-modules.} Aside from \newline interval-decomposability, there are other notions of structural simplicity for persistence modules.  

For any element $p$ of a poset $\catP$, let $p^\uparrow_\catP:=\{q\in \catP:p\leq q\}$\label{nom: up p}, which is an interval of $\catP$. 
If the ambient poset $P$ is clear from context, we simply write $p^\uparrow$. 
More generally, for any subset $Q\subset P$, we define $Q^\uparrow_\catP:=\bigcup_{p\in Q} p^\uparrow$, again dropping the subscript $P$ if the ambient poset is clear. %
\begin{definition}[e.g. {\cite{kerber2021fast,lesnick2015interactive,lesnick2022computing,miller2020homological}}]\label{def:finitely presentable}
 A $\catP$-module $M$ is said to be \textbf{finitely presentable} if $M$ is isomorphic to the cokernel of a morphism $\phi:F_1\rightarrow F_0$, where $F_0$ and $F_1$ are direct sums of the form $F_0=\bigoplus_{a\in A} \kf_{a^\uparrow}$ and $F_1=\bigoplus_{b\in B} \kf_{b^\uparrow}$, with $A$ and $B$ finite multisets of elements of $\catP$. The morphism $\phi$ is called a \textbf{presentation} of $M$. We call $\abs{A}$ (resp. $\abs{B}$) the number of generators (resp. the number of relations). 
 The presentation $\phi$ is \textbf{minimal} if among all possible presentations of the $P$-module $M$, it has the smallest number of generators and relations.
\end{definition}

\begin{remark}\label{rem:minimal grid}
Let $\catP\in \{\R^d,\Z^d\}$.
\begin{enumerate}[label=(\roman*)]
\item \label{item:uniqueness of minimal grid} {
Given any two minimal presentations $\phi:\bigoplus_{a\in A} \kf_{a^\uparrow}\rightarrow \bigoplus_{b\in B} \kf_{b^\uparrow}$ and \newline $\phi':\bigoplus_{a'\in A'} \kf_{a'^\uparrow}\rightarrow \bigoplus_{b'\in B'} \kf_{b'^\uparrow}$ for a finitely presentable $\catP$-module $M$, we have that $A=A'$ and $B=B'$, and these multisets $A$ and $B$ are referred to as the \emph{$0$-th and $1$-st graded Betti numbers} of $M$, respectively; see, for example,  \cite[Sections 2.3-2.4]{lesnick2015interactive}.}

\item Let $M$ and $N$ be any two finitely presentable $\catP$-modules. Then, the direct sum $M\oplus N$ is also finitely presentable, and its minimal presentation is the direct sum of minimal presentations of $M$ and $N$ ; see, for example,  \cite[Section 2.3]{lesnick2022computing}.\label{item:betti numbers of direct sum}

\item A $P$-module $M$ is finitely presentable if and only if there exists a finite full subposet $Q\subset P$ such that $M$ is the left Kan extension of a $Q$-module $N$ along the canonical inclusion $\iota:Q\hookrightarrow P$, i.e., for all $p\in \catP$
\[M_p:=\Lan_{\iota}(N)(p)=\varinjlim N|_{\{q\in \catQ\ :\ \iota(q)\leq p\}}, \]
and the structure maps of $M$ are given by universality of colimits; see, for example, \cite[Proposition 10.3(2)]{chacholski2021realisations}. 
\label{item:characterization of finitely presentable modules} 
\end{enumerate}
\end{remark}

\begin{definition}[\cite{miller2020homological}] \label{def:finite encoding}An \textbf{encoding} of a $\catP$-module $M$ by a poset $\catQ$ is an order-preserving map $\pi : \catP \rightarrow \catQ$ together with a $\catQ$-module $N$ such that $M$ is the \textbf{pullback} of $N$ along $\pi$, i.e., for all $p,p'\in \catP$ with $p\leq p'$, we have that
\[M_p=N_{\pi(p)}\ \mbox{and } \varphi_M(p, p')=\varphi_N(\pi(p), \pi(p')).\]\label{def:up-set}
In this case, we write $M=\pi^\ast N$. The encoding is called \textbf{finite} if the poset $\catQ$ is finite.\footnote{In \cite{miller2020homological}, it is also assumed that for all $q \in \catQ$, the vector space $N_q$ has finite dimension. However, in this paper, that must readily be the case since $M_p$ is assumed to be finite-dimensional for all $p\in \catP$.} A $\catP$-module $M$ is \textbf{tame} if $M$ admits a finite encoding \cite[Theorem 6.12]{miller2020homological}.
\end{definition}

The following proposition is a direct corollary of {\cite[Theorem 6.12, Remark 6.15]{miller2020homological}}.
\begin{proposition}\label{prop:finitely presentable implies tame} Finitely presentable persistence modules are tame.
\end{proposition}

There are other notions of structural simplicity of persistence modules. One of them is \emph{$q$-tameness}. An $\R$-module $M$ is considered $q$-tame if, for any $p < q$ in $\R$, the rank of the map $\varphi_M(p, q)$ is finite \cite{chazal2009proximity, chazal2016structure}. This definition can be directly generalized to $\catP$-modules for arbitrary posets $\catP$. Since this work only deals with pointwise finite-dimensional persistence modules, we readily assumed $q$-tameness throughout. Another is \emph{constructibility} \cite{de2016categorified,gulen2022galois,patel2018generalized}:

\begin{definition}\label{def:co-closure and constructible}
An order-preserving map $c:\catP\rightarrow \catP$ such that $c(p) \leq p$ for all $p \in \catP$, $c \circ c = c$ is called a \textbf{co-closure}. Let $S\subset P$ be any finite subset.
A $\catP$-module $M$ (resp. any function $f$ from $\catP$) is {$S$-\textbf{constructible} or simply \textbf{constructible}} if there is a co-closure $c:\catP\rightarrow \catP$ with image $S$, and  $M\circ c=M$ (resp. $f\circ c= f$). 
\end{definition}

\subsection{The M\"obius inversion formula}\label{sec:mobius inversion}
We review the notions of incidence algebra and M\"obius inversion \cite{rota1964foundations, stanley2011enumerative}.
Throughout this section, let $\catQ$ denote a  \textbf{locally finite} poset, i.e., for all $p,q\in \catQ$ with $p\leq q$, the \textbf{segment} \[[p,q]:=\{x\in \catQ:p\leq x\leq q\}\]
is finite. Let 
\begin{equation}\label{eq:segment}
    \seg(\catQ):=\{[p,q]:p\leq q\ \mbox{in $\catQ$}\}.
\end{equation} Given any function $\alpha:\seg(\catQ)\rightarrow \kf$, we write $\alpha(p,q)$ for $\alpha([p,q])$.
\label{nom:incidence algebra} The \textbf{incidence algebra} $I(\catQ,\kf)$ of $\catQ$ over $\kf$ is the $\kf$-algebra of all functions $\seg(\catQ)\rightarrow \kf$ with the usual structure of a vector space over $\kf$, where multiplication is given by convolution:
\begin{equation}\label{eq:convolution}
(\alpha\, \beta)(p,r):=\sum_{q\in [p,r]}\alpha(p,q)\cdot\beta(q,r).
\end{equation}
Since $\catQ$ is locally finite, the above sum is finite and hence $\alpha\beta$ is well-defined. This multiplication is associative and thus $I(\catQ,k)$ is an associative algebra. The \textbf{Dirac delta function} $\delta_\catQ\in I(\catQ,\kf)$ is given by 
\begin{equation}\label{eq:delta function}
\delta_\catQ(p,q):=\begin{cases}1,&p=q\\0,&\mbox{else,}\end{cases}
\end{equation}
which is the two-sided multiplicative identity. 
\begin{remark}[\cite{stanley2011enumerative}]\label{rem:invertibility} An element $\alpha\in I(\catQ,\kf)$ admits a multiplicative inverse if and only if $\alpha(q,q)\neq 0$ for all $q\in \catQ$. 
\end{remark}

Another important element of $I(\catQ,\kf)$ is the \textbf{zeta function}:
\begin{equation}\label{eq:zeta function}
    \zeta_\catQ(p,q):=\begin{cases}1,&p\leq q\\0,&\mbox{else.}\end{cases}
\end{equation}
By Remark \ref{rem:invertibility}, the zeta function $\zeta_\catQ$ admits a multiplicative inverse, which is called the \textbf{M\"obius function} $\mu_{\catQ}\in I(\catQ,\kf)$. 
The M\"obius function can be computed recursively as 
\begin{equation}\label{eq:mobius}
    \mu_{\catQ}(p,q)=\begin{cases}1,& \mbox{$p=q$,}\\ -\sum\limits_{p\leq r< q}\mu_{\catQ}(p,r),& \mbox{$p<q$,} \\ 0,& \mbox{otherwise.} \end{cases}
\end{equation}

\label{nom: kq} {Let $\kf^\catQ$ denote the
vector space of all functions $\catQ\rightarrow \kf$. 
\label{nom: down p}Also, for $q\in \catQ$, let \[q^\downarrow_\catQ:=\{p\in \catQ: p\leq q\},\] called a \textbf{principal ideal}.}  
As before, if the ambient poset is apparent we drop the poset subscript.
Assuming that $q^{\downarrow}$ 
is finite for each $q\in \catQ$, every element in $I(\catQ,\kf)$ acts on  $\kf^\catQ$ by right multiplication: for any $f\in k^\catQ$ and for any $\alpha\in I(\catQ,k)$, we have

\begin{equation}\label{eq:function times matrix}(f\ast\alpha)(q):=\sum_{p\leq q}f(p)\,\alpha(p,q).
\end{equation}

In fact, even when not every principal ideal in a poset $\catQ$ is finite, Equation (\ref{eq:function times matrix}) still specifies a well-defined multiplication between $f$ and $\alpha$ under the weaker assumption that
\begin{equation}\label{eq:convolvable}
\mbox{for every $q\in \catQ$, 
$f(r)=0$ for all but finitely many $r\in q^\downarrow$. 
} 
\end{equation}
\begin{definition}\label{def:convolvable}
Given a locally finite poset $\catQ$,
we call a function $f:\catQ\rightarrow \kf$ \textbf{\convolvable} (over $\catQ$) if $f$ satisfies Condition (\ref{eq:convolvable}).
\end{definition}

For a function $f\in \kf^\catQ$, its \textbf{support} is the set $\{q\in \catQ \ : \ f(q) \neq 0\}$, {denoted by $\supp(f)$.} 

\begin{proposition}[About convolvability]
\label{prop:convolvability}
Let $\catQ$ be a locally finite poset.
\begin{enumerate}[label=(\roman*),leftmargin=*]
\item  Any function $\catQ\rightarrow \kf$ with finite support is convolvable. {Hence, if $\catQ$ is finite, \emph{any} function $\catQ\rightarrow \kf$ is convolvable.} \label{item:convolvability1}
\item If every principal ideal in $\catQ$ is finite, then every $f\in \kf^\catQ$ is \convolvable.\label{item:convolvability2}
\item The collection of all convolvable functions $\catQ\rightarrow \kf$ is a linear subspace of $\kf^\catQ$.\label{item:convolvability3}
\item If $f:\catQ\rightarrow \kf$ is \convolvable\ over $\catQ$, then for any $\catP\subset \catQ$, the restriction $f|_{\catP}$ is  \convolvable.\label{item:convolvability4}
\item Let $P_1,P_2\subset \catQ$. Assume that, for each $i=1,2$, $f_i:\catP_i\rightarrow k$ is convolvable over $\catP_i$, and that $f_1=f_2$ on $\catP_1\cap \catP_2$. Then $f_1\cup f_2:\catP_1\cup \catP_2\rightarrow k$ is convolvable over $\catP_1\cup \catP_2$.
\label{item:convolvability5}
\item {If $f:Q\rightarrow k$ is convolvable over $Q$ and $\alpha\in I(Q,k)$, then $f\ast \alpha$ is convolvable.} \label{item:convolvable6}
\end{enumerate}
\end{proposition}

\begin{proof} {We only prove Item \ref{item:convolvable6}.} From Equation \eqref{eq:function times matrix}, it is clear that  
$\supp(f\ast\alpha)\subset (\supp(f))^\uparrow$ and thus, for any $q\in Q$, 
\[(\supp(f\ast\alpha))\cap q^\downarrow\subset (\supp(f))^\uparrow \cap q^\downarrow=\left(\bigcup_{x\in \supp(f)}x^\uparrow \right)\cap q^\downarrow = \bigcup_{x\in \supp(f)} x^\uparrow \cap q^\downarrow.\]
Since $f$ is convolvable, there exist only finitely many $x\in \supp(f)$ such that $x^\uparrow\cap q^\downarrow$ is nonempty, in which case $x^\uparrow\cap q^\downarrow= [x,q]$ is finite by local finiteness of $Q$. This proves that $\supp(f\ast\alpha)\cap q^\downarrow$ is finite, as desired.
\end{proof}

\begin{remark}\label{rem:right multiplication} 
Let $\catQ$ be a locally finite poset. Let $\kf_c^\catQ$ be the space of convolvable functions $\catQ\rightarrow \kf$, which is a subspace of $\kf^\catQ$. It follows that:

\begin{enumerate}[label=(\roman*),leftmargin=*]
    \item {By Proposition \ref{prop:convolvability} \ref{item:convolvability1}, if $\catQ$ is finite, then $\kf_c^\catQ=\kf^\catQ$.} \label{item:right multiplication-Q is finite}
    \item Let $\alpha\in I(\catQ,\kf)$. {By Proposition \ref{prop:convolvability} \ref{item:convolvable6},} 
    we have the  \textbf{right multiplication map} $\ast \alpha:\kf_c^\catQ\rightarrow \kf_c^\catQ$ given by $f\mapsto f\ast \alpha$. This map is an automorphism if and only if $\alpha$ is invertible.\label{item:right multiplication}
    \item By  Remark \ref{rem:invertibility} and the previous item, the right multiplication map $\ast\zeta_\catQ$ by the zeta function is an automorphism on $\kf_c^\catQ$ with inverse $\ast\mu_{\catQ}.$\label{item:zeta defines an automorphism}
\end{enumerate}  
\end{remark}

The M\"obius inversion formula is a powerful tool that has found widespread applications in combinatorics despite the fact that it has its origins in number theory. It will be a central tool for establishing our main results.

\begin{oldtheorem}[M\"obius Inversion formula]\label{thm:mobius}
Let $\catQ$ be a locally finite poset.
For any pair of \convolvable\  functions $f,g:\catQ\rightarrow k$,   \begin{equation}\label{eq:zeta}\displaystyle g(q)=\sum_{r\leq q} f(r)\ \mbox{for all $q\in \catQ$}
\end{equation} if and only if 
\begin{equation}\label{eq:zeta inverse}\displaystyle f(q)=\sum_{r\leq q} g(r)\cdot \mu_{\catQ}(r,q) \ \mbox{for all } q\in \catQ.
\end{equation}
\end{oldtheorem}

\begin{proof}
    Equation (\ref{eq:zeta}) can be represented as $g=f\ast\zeta_\catQ$. By multiplying both sides by $\zeta_\catQ^{-1}=\mu_\catQ$ on the right, we have $g\ast\mu_\catQ=f$, which is precisely Equation (\ref{eq:zeta inverse}).
\end{proof}

\begin{definition}\label{def:mobius inversion}The function $f=g\ast\mu_Q$ is referred to as the \textbf{M\"obius inversion} of $g$ (over $\catQ$).
\end{definition}

{M\"obius inversion serves as a discrete analogue of the concept of a derivative in calculus, as illustrated by the following example.}

\begin{example}Let $q\in \catQ$ and define the two functions $\mathbf{1}_q,\ \mathbf{1}_{\geq q}:\catQ\rightarrow \kf$ to be 
\[\mathbf{1}_q(p):=\begin{cases}1,&p=q\\0,&\mbox{otherwise.}\end{cases} \ \ \ \ \ \ \mathbf{1}_{\geq q}(p):=\begin{cases}1,&p \geq q\\0,&\mbox{otherwise.}\end{cases}
\]
Then, both functions are convolvable. Indeed, Proposition \ref{prop:convolvability} \ref{item:convolvability1} implies that $\mathbf{1}_q$ is convolvable and the local finiteness of $\catQ$ guarantees that $\mathbf{1}_{\geq q}$ is convolvable. Notice that the M\"obius inversion of $\mathbf{1}_{\geq q}$ is equal to $\mathbf{1}_q$
and that 
$\mathbf{1}_q$ captures where 
 the function values of $\mathbf{1}_{q\geq}$ change. 
\end{example}

The following proposition will be useful for some proofs in the sequel. 

\begin{proposition}\label{prop:mobius inversion of non-decreasing map is convolvable} 
Let $g:\catQ\rightarrow \R_{\geq 0}$ be non-decreasing and convolvable over $\catQ$, and let $f$ be the M\"obius inversion of $g$. Then, $\supp(f)\subset \supp(g)$.
\end{proposition}

\begin{proof} 
We prove that the complement of $\supp(f)$ contains the complement  of $\supp(g)$. Since $g$ is non-decreasing, if $g(r)=0$, then $g(r')=0$ for all $r'\leq r$  and thus $f(r)=\sum_{r'\leq r}g(r')\cdot\mu(r',r)=0$.  %
\end{proof}

A constructible function (Definition \ref{def:co-closure and constructible}) admits M\"obius inversion after its proper restriction:
\begin{proposition}[{\cite[Proposition 3.4]{gulen2022galois}}]\label{prop:Mobius inversion of a contructible function}\label{prop:inversion after restriction}Let $\catS$ be a nonempty finite subset of a given poset $\catP$. Let $m:\catP\rightarrow \Z$ be an $\catS$-constructible function. 
Then,  for all $q\in \catP$,
\[m(q)=\sum_{\substack{p\leq q \\ p\in S}} (m|_S\ast \mu_S)(p).\]
\end{proposition}

\paragraph{A matrix algebra perspective on the incidence algebra} Let $m\in \N=\{1,2,\ldots\}$ and let $\catQ$ be a poset with $m$ elements. 
By the order-extension principle, we can extend the order on $\catQ$ to a total order.
Thus, we fix $\catQ:= \{q_1,\ldots,q_m\}$, where $q_i<q_j$ implies $i<j$.
Then, each element $\alpha$ in the incidence algebra $I(\catQ,\kf)$ is canonically identified with the $(m\times m)$-upper-triangular matrix $(\alpha_{ij})$ whose $(i,j)$-entry is
\[ \alpha_{ij}:=\begin{cases}\alpha(q_i,q_j)&\mbox{if $q_i<q_j$}
\\0,&\mbox{else.}\end{cases}
\]
Then, for $\beta\in I(\catQ,\kf)$, the product $\alpha\, \beta$ in Equation (\ref{eq:convolution}) can be expressed as the multiplication of the upper-triangular matrices $(\alpha_{ij})$ and $(\beta_{ij})$, where \[(\alpha\, \beta)_{ij}= \mbox{(the $i$-th row of $\alpha$)$\cdot$(the $j$-th column of $\beta$$)^t$}.\]
Now let us identify each $f\in k^\catQ$ with the $m$-dimensional row vector $(f_i)$ where $f_i=f(q_i)$ for $i=1,\ldots,m$. 
Then, the multiplication $f\ast\alpha$ in Equation (\ref{eq:function times matrix}) can be expressed as the multiplication of the $(1\times m)$-matrix $(f_i)$ and the $(m\times m)$-matrix $(\alpha_{ij})$.   

\begin{remark} Recall that an upper-triangular matrix is invertible if and only if all of its diagonal entries are nonzero. Remark \ref{rem:invertibility} can be seen as a straightforward from this fact. 
\end{remark}

\subsection{Generalized rank invariant {(GRI)}}\label{sec:GRI}
In this section, we recall the definitions of the rank invariant \cite{carlsson2009theory,puuska2020erosion} and the generalized rank invariant \cite{kim2021generalized}. 

Let $M$ be any $\catP$-module. 
\label{nom:rank invariant}
Then, $M$ admits both a \emph{limit} and a \emph{colimit} \cite[Chapter V]{mac2013categories}: 
A limit of $M$, denoted by $\varprojlim M$, consists of a vector space $L$ together with a collection of linear maps $\{\pi_p:L\to M_p\}_{p\in \catP}$ such that 
\begin{equation}\label{eq:limit}
    \varphi_M(p, q)\circ \pi_p =\pi_q\  \mbox{for every $p\leq q$ in $\catP$}.
\end{equation}
A colimit of $M$, denoted by $\varinjlim M$, consists of a vector space $C$ together with a collection of linear maps $\{i_p:M_p\to C\}_{p\in \catP}$ such that 
\begin{equation}\label{eq:colimit}
    i_q\circ \varphi_M(p, q)=i_p\  \mbox{for every $p\leq q$ in $\catP$}.
\end{equation}
Both $\varprojlim M$ and $\varinjlim M$ satisfy certain universal properties, making them unique up to isomorphism. The following well-known concrete formulation for the limit and colimit will be useful later :

\begin{proposition}[{see, e.g. \cite[Appendix E]{kim2021generalized}}]\phantomsection\label{prop:limit and colimit} 
\begin{enumerate}[label=(\roman*),topsep=0pt,itemsep=-1ex,partopsep=1ex,parsep=1ex,leftmargin=*]
    \item The limit of $M$ is isomorphic to the pair
    $(W,(\pi_p)_{p\in \catP})$ described as:
    \[W:=\left\{ (\ell_p)_{p\in \catP}\in \prod_{p\in \catP} M_p: \ \forall p\leq q\in \catP, \varphi_M(p,q)(\ell_p)=\ell_q\right\}\]
    where for each $p\in \catP$, the map $\pi_p:W\rightarrow M_p$ is the canonical projection. Elements of $W$ are called  \textbf{sections} of $M$.\label{item:limit}
    \item The colimit of $M$ is isomorphic to a pair $(C,(\iota_p)_{p\in \catP})$ that is described as follows.
    For $p\in \catP$, let $\bar{\iota}_p:M_p\to \bigoplus_{p\in \catP}M_p$ be the canonical injection.
    $C$ is the quotient space $\left(\bigoplus_{p\in\catP} M_p\right)/V$, where $V$ is generated by the vectors $\bar{\iota}_p(v_p)-\bar{\iota}_{q}(v_{q})$ over all $p\leq q$ in $\catP$, with $v_p\in M_p,v_{q}\in M_{q}$ , and $\varphi_M(p,q)(v_p)=v_{q}$. 
    Letting $f$ be the quotient map from $\bigoplus_{p\in \catP}M_p$ to $C$, for $p\in \catP$, $\iota_p:M_p\to C$ is the composition $f\circ \bar{\iota}_p$.\label{item:colimit}
\end{enumerate}
\end{proposition}

Let us assume that $\catP$ is connected (Definition \ref{def:intervals} \ref{item:connectivity}). The connectedness of $\catP$ alongside the equalities given in Equations (\ref{eq:limit}) and (\ref{eq:colimit}) imply that $i_p \circ \pi_p=i_q \circ \pi_q:L\rightarrow C$ for any $p,q\in \catP$. 
This fact ensures that the \textbf{canonical limit-to-colimit map} \[\psi_M:\varprojlim M\longrightarrow \varinjlim M\] given by $i_p\circ \pi_p$ for any $p\in \catP$ is well-defined. 
\label{nom:generalized rank} The \textbf{(generalized) rank} of $M$ is defined to be:\footnote{{This construction was considered in the study of quiver representations \cite{kinser2008rank}.}} 
\begin{equation}\label{eq:generalized rank of M}
    \rank(M):=\rank(\psi_M).
\end{equation}

\begin{remark}\label{rem:generalized rank is finite for pfd modules}
    $\rank(M)$ is finite as $\rank(M)=\rank(i_p\circ \pi_p)\leq \dim(M_p)<\infty$ for any $p\in P$.
\end{remark}

The rank of $M$ is a count of the `persistent features' in $M$ that span the entire indexing poset $\catP$. 
{Such `persistent features' appear as summands of $M$ in the form of $\kf_\catP$:}

\begin{oldtheorem}[{\cite[Lemma 3.1]{chambers2018persistent}, \cite[Lemma 5.41]{asashiba2024interval}}]
\label{thm:rkequalsintervals}
Let $\catP$ be a connected poset. Assume that a $\catP$-module $M$ is isomorphic to a direct sum $\bigoplus_{a\in A} M_a$ for some indexing set $A$ where each $M_a$ is indecomposable. Then, the rank of $M$ is equal to the cardinality of the set $\{a\in A: M_a\cong \kf_P\}.$ 
\end{oldtheorem}

We refine the rank of a $\catP$-module, which is a single integer, into an integer-valued function. Possible domain options for the function include the following sets, in addition to $\seg(\catP)$ as shown in Equation (\ref{eq:segment}):
\label{nom:Con(P)}
\begin{align*}
\Con(\catP)&:=\{I\subset\catP:\mbox{$I$ is nonempty and connected}\}, \mbox{}
\\
\Int(\catP)&:=\{I\subset\catP:\mbox{$I$ is an interval}\}.
\end{align*}
Note the inclusions $\seg(\catP)\subset \Int(\catP)\subset \Con(\catP).$

\begin{definition}[\cite{asashiba2019approximation,carlsson2009theory,dey2022computing,kim2021generalized}]\label{def:generalized rank invariant}   
    The \textbf{generalized rank invariant (GRI) of} a $\catP$-module $M$ is the map 
    \begin{align*}
    \rank_M:\Con(\catP)&\to \Z_{\geq 0}\\
    I&\mapsto \rank(M\vert_I)
    \end{align*}
where $M\vert_I$ is the restriction of $M$ to $I$.  The restriction of $\rk_M$ to $\Ical\subset \Con(\catP)$ is denoted by $\rk_M^\Ical$ and is called the \textbf{GRI of $M$ over} $\Ical$.  
When $\mathcal{I}=\Int(\catP)$, the GRI over $\Ical$ is simply called the \textbf{Int-GRI} and denoted by $\rk_M^{\Int}$. When $\Ical=\seg(\catP)$, the GRI of $M$ over $\Ical$ is the \textbf{rank invariant (RI)} of $M$, introduced in \cite{carlsson2009theory}, which will be denoted by $\rk_M^\seg$. The GRI over $\Con(P)$ will also be referred to as the Con-GRI to emphasize the entirety of the domain.
\end{definition} 
We remark that the Int-GRI  has been shown to be strictly stronger than the RI for several indexing posets \cite{kim2021generalized,dey2022computing}, i.e., there exists a pair of persistence modules that have the same RI but different Int-GRI. In addition, when $\catP=\Z^2$, the Con-GRI  has been shown to be strictly stronger than the Int-GRI \cite{kim2021bettis}.

\begin{remark} The following are useful properties of the GRI.\label{rem:properties of rank invariant}
\begin{enumerate}[label=(\roman*),leftmargin=*]
\item (Monotonicity)\label{item:monotonicity} If $I\subset J$ in $\Con(\catP)$, then $\rank_M(I)\geq \rank_M(J)$. This is because the canonical limit-to-colimit map over $I$ is a factor of the canonical
limit-to-colimit map over $J$ \cite[Proposition 3.8]{kim2021generalized}. 

\item (Additivity) If $M\cong \bigoplus_{i=1}^n M_i$, and $I\in \Con(\catP)$, then $\rank_M(I) = \sum_{i=1}^n \rank_{M_i}(I)$. \label{item:additivity}

\item (The GRI of an interval module) Let $J\in \Int(\catP)$. For the interval module $\kf_J$ and any $I\in \Con(\catP)$, we have
\[\rank_{\kf_J}(I)=\begin{cases}1,& J\supset I\\0,&\mbox{else.}\end{cases}\]
\label{item:the GRI of an interval module}
\end{enumerate}
\end{remark}

{We review a simple but useful proposition concerning the colimit of a $\catP$-module, which will be used later.} 
Let $M$ be any $\catP$-module.
Let $\Gamma:p_0,\ldots,p_k$ be a path in $\catP$. 
A $(k+1)$-tuple $\bv\in \bigoplus_{i=0}^k M_{p_i}$ is called a \textbf{section of $M$ along $\Gamma$} if for each $i$, either $v_{p_i}$ is mapped to $v_{p_{i+1}}$ via $\varphi_M(p_i,p_{i+1})$ or $v_{p_{i+1}}$ is mapped to $v_{p_i}$ via $\varphi_M(p_{i+1},p_i)$.
Note that $\bv$ is not necessarily a section of the restriction  $M\vert_{\{p_0,\ldots,p_k\}}$  of $M$ to the subposet $\{p_0,\ldots,p_k\}\subset I$ \cite[Example 21]{dey2022computing}. 
Furthermore, $\Gamma$ can contain multiple copies of the same point in $\catP$.

\begin{proposition}[{\cite[Proposition 22]{dey2022computing}}]\label{prop:path characterization of zero in colim}
Let $p,p'\in \catP$. 
For any vectors $v_{p}\in M_{p}$ and $v_{p'}\in M_{p'}$, it holds that $[v_{p}]=[v_{p'}]$ as elements of\footnote{For simplicity, we write $[v_{p}]$
 and $[v_{p'}]$ instead of $[j_p(v_p)]$ and $[j_q(v_{p'})]$ respectively where $j_p:M_p\rightarrow \bigoplus_{r\in \catP}M_r$  and  $j_{p'}:M_{p'}\rightarrow \bigoplus_{r\in \catP}M_r$ are the canonical inclusion maps.} the colimit $\varinjlim M$ if and only if there exist a path $
 \Gamma$ in $\catP$ connecting $p$ and $p'$
  and a section $\bv$ of $M$ along $\Gamma$ such that $\bv_{p}=v_{p}$ and $\bv_{p'}=v_{p'}$. 
\end{proposition}

\begin{figure}
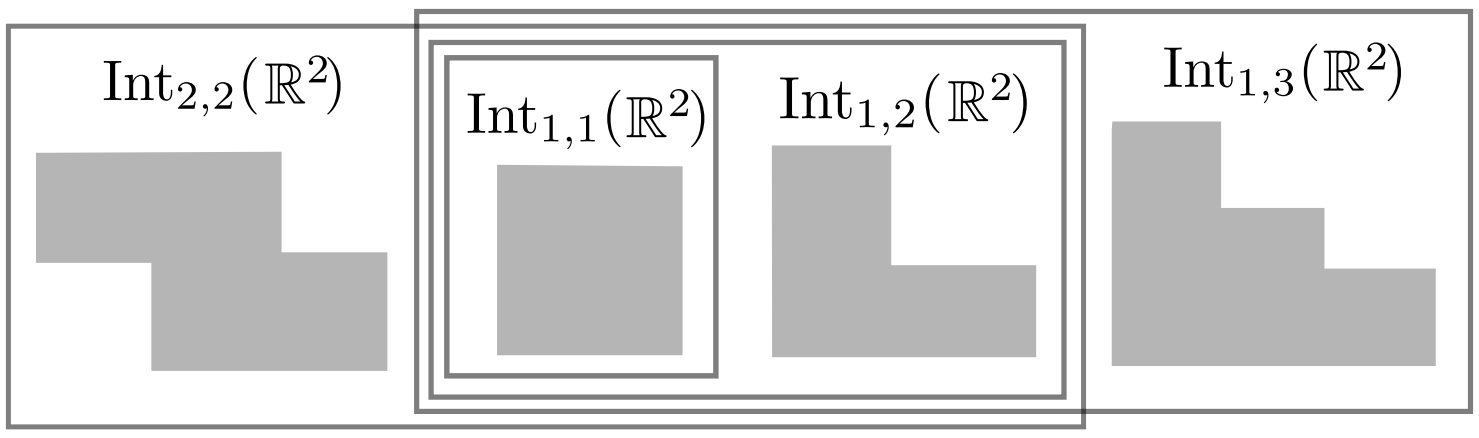

    \centering
    \resizebox{0.6\textwidth}{!}{\mobiusintervalfigure}
    \caption{$\Int_{m,n}(\R^2)$ for some values of $m$ and $n$.  \label{fig:intmnr2}}
\end{figure}

\smallskip 
We close this section by describing a parameterized family of subsets of $\Int(\catP)$ that will be useful later. Given a poset $\catP$, let $a\in \catP$ and let $B\subset \catP$ be an \emph{antichain}, i.e., any two different points in $B$ are \emph{not} comparable. We write $a\leq B$ (resp. $a\geq B$) if there exists $b\in B$ such that $a\leq b$ (resp. $a\geq b$). When two antichains $A,B\subset \catP$ are given, we write $A\leq B$ if $\forall a\in A$, $\forall b\in B$, $a\leq B$ and $A\leq b$. This defines a partial order on the set of antichains in $\catP$ \cite[Section 2.1]{blanchette2021homological}.
When $A\leq B$, define 
 \[[A,B):=\{x\in \catP: \exists a\in A, \exists b\in B, a\leq x <b\},\]
We remark that $[A,B)$ is either empty (when $A=B$), or when $A<B$, $[A,B)$ is either an interval, or a disjoint union of multiple intervals of $\catP$. 
For $m,n\in \N$, let
    \begin{equation}\label{eq:intmn}
    \Int_{m,n}(\catP):=\{[A,B)\in\Int(\catP): A,B\subset \catP\ \mbox{are antichains s.t.  $A<B$, $\abs{A}\leq m$, $\abs{B} \leq n$ \}}
    \end{equation}
See Figure~\ref{fig:intmnr2} for an illustration of $\Int_{m,n}(\catP)$ corresponding to the case $\catP=\R^2$.
We also define, for $m\in \N$,
\[ \Int_{m,0}(\catP):=\{A^\uparrow\in \Int(\catP): A\subset \catP \mbox{ is an antichain s.t. } \abs{A}\leq m \}.
\]

\subsection{Generalized persistence diagrams {(GPD)}}\label{sec:GPD}

In this section we review the notion of the generalized persistence diagram (GPD) introduced in \cite{kim2021generalized}, but we modify the original formulation by adding more flexibility in choosing the domain of the GPD. 
This additional flexibility is useful when clarifying the relationship between the GRI/GPD and other invariants of persistence modules such as graded Betti numbers\footnote{{Graded Betti numbers are classical invariants \cite{eisenbudbetti}, and, to the best of our knowledge, they first appeared in \cite{carlsson2009theory} in the context of multi-parameter persistence.}} \cite{kim2021bettis} and rank decompositions (Section \ref{sec:GPD and Rk decomposition}).

 Let $\Con(\catP)$ and any subcollection of $\Con(\catP)$ be ordered by containment $\supset${, the opposite to the usual inclusion order.} Let $M$ be a $\catP$-module.
If $\Ical\subset \Con(\catP)$ is locally finite and $\rank_M^\Ical$ is convolvable over $\Ical$ (Definition \ref{def:convolvable}), we say that \textbf{the GRI of $M$ is convolvable over $\Ical$}.

\begin{definition}\label{def:GPD over I}
    Let $M$ be a $\catP$-module and let $\Ical\subset \Con(\catP)$.
    Assume that the GRI of $M$ is \convolvable\ over $\Ical$. 
    The \textbf{generalized persistence diagram (GPD) of $M$ over $\Ical$} is the M\"obius inversion of $\rank_M^\Ical$ over the poset $(\Ical,\supset)$, i.e., the function $\dgm_M^{\Ical}:\Ical\rightarrow \Z$ given by:\begin{equation}\label{eq:dgm}\dgm_M^{\Ical}:=\rank_M^{\Ical}\ast\mu_{\Ical}.\end{equation}
   When $\Ical=\Int(\catP)$, the GPD of $M$ over $\Ical$ will be referred to as simply the \textbf{Int-GPD} of $M$, and will be denoted by $\dgm_M^\Int$. 
\end{definition}

\begin{remark}\label{rem:additivity of GPD} When the GRIs of $\catP$-modules $M$ and $N$ are both convolvable over $\Ical\subset \Con(\catP)$, by additivity of the GRI (cf. Remark \ref{rem:properties of rank invariant} \ref{item:monotonicity}), we have $\dgm_{M\oplus N}^\Ical=\dgm_M^\Ical+\dgm_N^\Ical$. 
\end{remark}

We establish a few other desirable properties of the GPD.

\begin{proposition}\label{prop:properties of GPD}  Assume that the GRI of a given $\catP$-module $M$ is convolvable over a locally finite $\Ical\subset \Con(\catP)$. Then the following holds.

\begin{enumerate}[leftmargin=*, label=(\roman*)]
\item The support of $\dgm_M^{\Ical}$ is contained in the support of $\rank_M^{\Ical}$. In other words, for any $I\in \Ical$, if $\rank_M^{\Ical}(I)=0$, then $\dgm_M^{\Ical}(I)=0$. \label{item:GPD is more parsimonious}
\item $\dgm_M^{\Ical}$ is convolvable over $\Ical$.\label{item:convolvability of GPD}
\item The unique function $\drm:\Ical\rightarrow \R$ satisfying 
\[\rank_M(I)=\sum_{\substack{J\supset I\\ J\in \Ical}}\drm(J) \ \ \mbox{for all $I\in \Ical$},\]
where the sum contains only finitely many nonzero summands, is $\dgm_M^{\Ical}$.
In other words, Equation (\ref{eq:dgm}) implies $\rk_M^\Ical = \dgm_M^\Ical \ast \zeta_\Ical$.
\label{item:universal property of GPD}
\item If $\Ical\subset \Int(\catP)$ and $M$ is interval-decomposable with summands supported on intervals in $\Ical$, then $\dgm_M^\Ical(I)$ equals the multiplicity of $I$ in $\barc(M)$.\label{item:GPD captures multiplicity}
\end{enumerate} 
\end{proposition}

We note that \cite[Theorem~4.11]{escolar2024barcoding} gives a  characterization of the GPD different from the one in Item~\ref{item:universal property of GPD}.

\begin{proof}Item \ref{item:GPD is more parsimonious} follows from Remark \ref{rem:properties of rank invariant} \ref{item:monotonicity} and  Proposition \ref{prop:mobius inversion of non-decreasing map is convolvable}.
Item \ref{item:convolvability of GPD} is an immediate consequence of 
Proposition \ref{prop:convolvability} \ref{item:convolvable6}.
Item \ref{item:universal property of GPD} follows from Theorem \ref{thm:mobius} and Item \ref{item:convolvability of GPD}. Item~\ref{item:GPD captures multiplicity} follows from Item~\ref{item:universal property of GPD} and the fact that \[\rank_M(I)=\sum_{\substack{J\supset I\\J\in \Ical}} \mbox{(multiplicity of $J$ in $\barc(M)$)} \qquad \mbox{for all $I\in \Ical$,}\] which itself follows from Remark~\ref{rem:properties of rank invariant}~\ref{item:additivity} and \ref{item:the GRI of an interval module}.
\end{proof}
Item \ref{item:universal property of GPD} implies that $\dgm_M^\Ical$ determines $\rk_M^\Ical$. This fact together with Item \ref{item:GPD is more parsimonious} suggests that $\dgm_M^\Ical$ provides a representation of $\rk_M^\Ical$ with smaller support.

\subsection{Rank decompositions}\label{sec:GPD and Rk decomposition}

In this section, we provide a review of the notion of \emph{rank decomposition} \cite{botnan2024signed}, with a strong focus on its connection with the M\"obius inversion formula. 

For any multiset $\Rcal$ of intervals in a poset $\catP$, let $\mult_{\Rcal}:\Int(\catP)\rightarrow \Z_{\geq 0}$ be defined by $R\mapsto$ (the multiplicity of $R$ in $\Rcal$). The multiset $\Rcal$ is said to be \textbf{\pfin}\ if for all $p\in \catP$, the sum $\sum_{R\ni p}\mult_{\Rcal}(R)$ is finite.\label{nom:mult R} Also, given any $\Ical\subset \Int(\catP)$, the multiset $\Rcal$ is said to be \textbf{$\Ical$-finite}\ if for all $I\in \Ical$, the sum $\sum_{R\supset I}\mult_{\Rcal}(R)$ is finite. We define the direct sum \[\kf_{\Rcal}:=\bigoplus_{R\in \Rcal}k_R.\] Although $(\kf_{\Rcal})_p$ may be infinite-dimensional for some
$p\in\catP$, we use the same definition of the GRI for $\kf_{\Rcal}$ via the
canonical limit-to-colimit map, allowing the value $+\infty$.

\begin{remark}\label{rem:pointwisely finite multiset-new}
Let $\Ical \subset \Int(\catP)$, and let $\Rcal$ be a multiset of elements of $\Ical$. Then $\Rcal$ is $\Ical$-finite if and only if $\rk^{\Ical}_{\kf_{\Rcal}}$
 is finite-valued. If $\Ical$ contains all singleton intervals, then $\Ical$-finiteness is equivalent to pointwise finiteness of $\Rcal$, and hence to pointwise finite-dimensionality of $\kf_\Rcal$.
\end{remark}

\begin{definition}\label{def:rank decomposition}
     Let $M$ be a $\catP$-module and $\Ical\subset \Int(\catP)$. 
          Any pair $(\Rcal,\Scal)$ of $\Ical$-finite multisets of elements in $\Ical$  such that
    \[\rank_M^{\Ical}=\rank_{\kf_{\Rcal}}^{\Ical}-\rank_{\kf_{\Scal}}^{\Ical}\]
    is called a \textbf{rank decomposition} of $\rank_M^\Ical$.

    If $\Rcal$ and $\Scal$ are disjoint, then the rank decomposition is called \textbf{minimal} and we write $\Rcal_M$ and $\Scal_M$. The uniqueness of the pair $(\Rcal_M,\Scal_M)$ follows from the next proposition.
\end{definition}

\begin{proposition}[{\cite[Corollary 2.7]{botnan2024signed}}]\label{prop:if rk decomposes then the minimal one exists}
  Let $M$ be a $\catP$-module such that $\rank_M^\Ical$ admits a rank decomposition $(\Rcal,\Scal)$. Then, the unique minimal rank decomposition of $\rank_M^\Ical$ is given by 
  $(\Rcal_M,\Scal_M)$, where $\Rcal_M:=\Rcal-(\Rcal\cap\Scal)$ and $\Scal_M:=\Scal-(\Rcal\cap\Scal)$.
\end{proposition}

\begin{oldtheorem}[{Restatement of \cite[Corollary 2.13]{botnan2024signed}}]
\label{them:rk decomposition is GPD} 
Assume that $\Ical\subset\Int(\catP)$ is locally finite, and let $M$ be a $\catP$-module.
If $\rk_M^\Ical$ is convolvable over $\Ical$, then the minimal rank decomposition of $\rk^\Ical_M$ is given by
\[\Rcal_M=\{d_I\cdot I:\dgm_M^\Ical(I)>0\}\ \ \ \mbox{and}\ \ \  \Scal_M=\{d_I\cdot I:\dgm_M^\Ical(I)<0\},\] 
where $d_I:=\abs{\dgm_M^\Ical(I)}$ and $d_I\cdot I$ stands for $d_I$ copies of $I$.
\end{oldtheorem}
Since the statement above is a slight reformulation of \cite[Corollary~2.13]{botnan2024signed}, we include a proof adapted to our presentation.

\begin{proof}
Let $\mathbf{1}_{J\supset}:\Int(\catP)\rightarrow \{0,1\}$ be defined by $\mathbf{1}_{J\supset}(I)=1$ if $J\supset I$ and $\mathbf{1}_{J\supset}(I)=0$, otherwise. For every $I\in \Ical$, we have
     \begin{align}        \rk_M^\Ical(I)&=\sum_{\substack{J\supset I \\ J\in\Ical}}\dgm_M^\Ical(J) \label{eq:line1}\\&=\sum_{\substack{J\supset I\\ \dgm_M^\Ical(J)>0}}\dgm_M^\Ical(J)+\sum_{\substack{J\supset I \\ \dgm_M^\Ical(J)<0}}\dgm_M^\Ical(J) \label{eq:line2}\\&=\sum_{\substack{\dgm_M^\Ical(J)>0}}\dgm_M^\Ical(J)\cdot\mathbf{1}_{J\supset}(I)+\sum_{\substack{\dgm_M^\Ical(J)<0}}\dgm_M^\Ical(J)\cdot\mathbf{1}_{J\supset}(I)\label{eq:line3}
     \\&=\sum_{J\supset I} \mult_{\Rcal_M}(J)-\sum_{J\supset I} \mult_{\Scal_M}(J)\\&=\rank_{\kf_{\Rcal_M}}^\Ical(I)
-\rank_{\kf_{\Scal_M}}^\Ical(I),
\qquad\text{by \cite[Proposition~2.1]{botnan2024signed}.}
\label{eq:line5}
    \end{align}
    Since $\Rcal_M\cap\Scal_M=\emptyset$, this rank decomposition is minimal.
\end{proof}

    Theorem \ref{them:rk decomposition is GPD} and its proof implies, under the assumption that $\Ical\subset\Int(\catP)$ is locally finite and $\rank_M^\Ical$ is convolvable over $\Ical$:
    \begin{enumerate}[label=(\roman*)]
        \item  the existence of a rank decomposition is implied by the forward direction of M\"obius Inversion Formula (Theorem \ref{thm:mobius}).
        \item the uniqueness of a minimal rank decomposition is implied by the backward direction of M\"obius Inversion Formula.  
    \end{enumerate}

The following theorem is proved by simply replacing $\Ical$, $\rk_M^\Ical$, and $\dgm_M^\Ical$ in the preceding proof with $\Ical$, $r$, and $r\ast\mu_{\Ical}$, respectively.

 \begin{oldtheorem}[{Restatement of \cite[Theorem 2.12]{botnan2024signed}}]\label{thm:BOO's setup is more general} Let $\Ical\subset \Int(\catP)$ be locally finite and let $r:\Ical\rightarrow \Z$ be convolvable. 
 Then, there exists $a_I\in \Z$ for each $I\in \Ical$ such that for $M:=\bigoplus\limits_{\substack{I\in \Ical\\a_I>0}} (\kf_I)^{a_I}$ and $N:=\bigoplus\limits_{\substack{I\in \Ical\\a_I<0}} (\kf_I)^{-a_I}$, we have  
 $r= \rk_M - \rk_N$, or equivalently,  for every $J\in \Ical$ 
 \begin{equation*}
     r(J)=\sum_{I\in \Ical}a_I\rk_{\kf_I}^{\Ical}(J)
 \end{equation*}
where the RHS contains only finitely many nonzero terms.
  \end{oldtheorem}

\section{M\"obius invertibility of the GRI}
\label{sec:mobius invertible Rk}

In Section~\ref{sec:relative to a dictionary}, we extend the notion of M\"obius inversion beyond locally finite posets by introducing the notion of M\"obius invertibility. 
In particular, M\"obius invertibility of the GRI dissociates the existence of the GPD from that of the GRI as well as extends the notion of GPD, and serves as a novel measure of structural simplicity for persistence modules. 
In Section~\ref{sec:Mobius invertibility and rank decompositions}, 
{we show that the notion of rank decomposition can be viewed as an instance of the extended notion of M\"obius inversion.}
In Section~\ref{sec:On structural simplicity of persistence modules}, we compare M\"obius invertibility of the GRI with other concepts pertaining to the structural simplicity of persistence modules, as depicted in Figure~\ref{fig:tameness + MI visual}. 
In Section \ref{sec:sufficient conditions}, we identify certain sufficient conditions on persistence modules that guarantee the M\"obius invertibility of their GRI. 
In Section~\ref{sec:GPD for finitely presentables}, we introduce a proxy for the GPD that circumvents the fact that finitely presentable $d$-parameter persistence modules with $d\geq 3$ 
may fail to admit a M\"obius invertible GRI.

\subsection{M\"obius invertibility and an extension of the GPD
}\label{sec:relative to a dictionary}

In this section, we generalize the notion of M\"obius inversion (Definition~\ref{def:mobius inversion}) by dropping the local finiteness assumption on the domain of a function.
This enables us to extend the notion of GPD. 

\begin{definition}\label{def:mobius invertible general} Let $\catP$ be a poset (not necessarily locally finite). A function $g:\catP\to \Z$ is said to be \textbf{\invertible} if there exists a function $f:\catP \to \Z$ such that 
\begin{enumerate}[label=(\roman*)]
\item \label{item:mobius invertible general 1} for each $q\in \catP$, the set $q^\downarrow_P \cap \supp(f)$ is finite and
\item  \label{item:mobius invertible general 2}
$g(q)=\sum_{p\leq q}f(p).$
\end{enumerate}
By abuse of terminology, the function $f$ will be called a \textbf{M\"obius inversion} of $g$ (cf. Definition~\ref{def:mobius inversion}). In particular, if $\supp(f)$ is contained in some subset $Q \subset P$, we say that $g$ is \invertible\ \textbf{over} $Q$. 
\end{definition}

The abuse of terminology will be justified by showing that the M\"obius inversion $f$ is unique and is obtained via the original notion of M\"obius inversion applied to a restriction of $g$ to a locally finite subposet of $\catP$.

\begin{proposition}\label{prop:justifying the term mobius invertibility}
Let $g:\catP\to\Z$ be \invertible\ with M\"obius inversion $f$. Then, 
\begin{enumerate}[label=(\roman*)]
    \item $f$ is uniquely determined, \label{item:justifying the term mobius invertibility-unique}
    \item $Q:=\supp (f)$ is a locally finite subposet of $\catP$, \label{item:justifying the term mobius invertibility-locally finite}
    \item  $g|_Q$ is convolvable, and \label{item:justifying the term mobius invertibility-convolvable}
    \item $f|_{Q}=g|_{Q}\ast\mu_{Q}$, 
where $\mu_Q$ denotes the M\"obius function of $Q$. \label{item:justifying the term mobius invertibility-mobius inversion} 
\end{enumerate}
\end{proposition}

We remark that this proposition implies that the restriction of $g$ to $Q$ yields a compact encoding of $g$ via $f|_Q$, without loss of information. %

\begin{proof}
\noindent \ref{item:justifying the term mobius invertibility-locally finite}: Since, for each $q\in\catP$, the set
$q^\downarrow_{\catP}\cap Q=q^\downarrow_{\catP}\cap\supp (f)$ is finite by
Definition~\ref{def:mobius invertible general}~\ref{item:mobius invertible general 1}, the induced subposet $Q$ is
locally finite. Indeed, for $p\le q$ in $Q$, the interval
$[p,q]_Q=[p,q]_{\catP}\cap Q$ is contained in
$q^\downarrow_{\catP}\cap Q$, and hence is finite. 

\noindent \ref{item:justifying the term mobius invertibility-convolvable}: The finiteness of $q^\downarrow_{\catP}\cap \catQ$ for all $q\in \catQ(\subset \catP)$ implies the convolvability of $f|_Q$.
For each $q\in Q$, Definition~\ref{def:mobius invertible general}~\ref{item:mobius invertible general 2} gives
\[
g(q)=\sum_{p\le q}f(p)
      =\sum_{p\in q^\downarrow_{\catP}\cap Q}f(p)
      =\sum_{p\le_Q q}f(p),
\]
where the second equality holds because $f(p)=0$ for $p\notin Q$.
Thus,
\[
g|_Q= f|_Q\ast\zeta_Q,
\]
where $\zeta_Q$ is the zeta function of the locally finite poset $Q$. By Proposition~\ref{prop:convolvability}~\ref{item:convolvable6}, $g|_Q$ is convolvable. 

\noindent \ref{item:justifying the term mobius invertibility-mobius inversion}: From the preceding equation, we obtain
\[
f|_Q=g|_Q\ast\mu_Q.
\]

\noindent \ref{item:justifying the term mobius invertibility-unique}: Suppose that $f'$ is another M\"obius inversion of $g$ in the sense of Definition~\ref{def:mobius invertible general}, and let
$Q'=\supp (f')$. By the argument above, $Q'$ is also locally finite and \[\qquad
f'|_{Q'}=g|_{Q'}\ast\mu_{Q'}.
\]
Let $R:=Q\cup Q'.$
Since, for every $q\in\catP$,
\[
q^\downarrow_{\catP}\cap R
=
(q^\downarrow_{\catP}\cap Q)\cup
(q^\downarrow_{\catP}\cap Q'),
\]
the set $q^\downarrow_{\catP}\cap R$ is finite. Hence $R$ is a locally finite
subposet of $\catP$.

Now define $h:=f-f'$. Since both $f$ and $f'$ vanish outside $R$, we have
$\supp(h)\subseteq R$. Also, the restriction $h|_R$ is convolvable as the set $q^\downarrow_{\catP}\cap R$ is finite for every $q\in R(\subset P)$. Moreover, for every $q\in R$,
\[
\sum_{p\le_R q}h(p)
=
\sum_{p\le q}f(p)-\sum_{p\le q}f'(p)
=
g(q)-g(q)
=
0.
\]
Thus $h|_R$ is a M\"obius inversion of the zero function on the locally finite
poset $R$. By the uniqueness of M\"obius inversion on locally finite posets,
it follows that $h|_R=0$. Hence $f=f'$ on $R$, and therefore on all of
$\catP$.
\end{proof}

\begin{remark}\phantomsection\label{rem:monotonicity of mobius invertibility new}
\begin{enumerate}[label=(\roman*)]
\item (Comparison with Definition~\ref{def:co-closure and constructible})
If $g:\catP\to \Z$ is $\catQ$-constructible for some finite $Q\subset P$, then by Proposition~\ref{prop:Mobius inversion of a contructible function}, $g$ is M\"obius invertible over $\catQ$. However, the converse does not hold; see Example~\ref{ex:mobius invertibility is different from constructibility} below. 

\item (Comparison with Definition~\ref{def:convolvable}) Convolvability of $g:\catP\to\Z$ directly implies M\"obius invertibility of $g$. However, the converse does not hold. Convolvability requires the poset $\catP$ to be locally finite, whereas M\"obius invertibility does not. 
{For example, the floor function $[0,\infty)\to\N$ admits the M\"obius inversion
$f:[0,\infty)\to\N$ defined by $f(n)=1$ for $n\in\N$ and $f(x)=0$ for
$x\in[0,\infty)\setminus\N$. However, the floor function is not convolvable,
since $[0,\infty)$ is not locally finite.}
\label{item:comparison between mobius and convolv}
\item (Inheritance of M\"obius invertibility) 
Let $\catQ \subset \catQ' \subset \catP' \subset \catP$.
If $g:\catP\to \Z$ is M\"obius invertible over $\catQ$, then it directly follows from definition that $g|_{\catP'}$ is M\"obius invertible over $\catQ'$.
\label{item:monotonicity of mobius invertibility}
\end{enumerate}
\end{remark}

Now we further generalize the notion of GPD (Definition~\ref{def:GPD over I}):
\begin{definition}[Generalized persistence diagram]\label{def:GGPD} Let $M$ be a $\catP$-module and $\Jcal\subset \Con(\catP)$. If $\rk_M^\Jcal$ is \invertible{}, then its M\"obius inversion is called the \textbf{generalized persistence diagram (GPD) of $M$ over $\Jcal$}, denoted by $\dgm_M^\Jcal$. 
In particular, when $\Jcal$ is either $\Con(\catP)$ or $\Int(\catP)$, we simply write $\dgm_M^\Con$ or $\dgm_M^{\Int}$.
Also, if $\rk_M^\Jcal$ is M\"obius invertible over $\Ical$ for some $\Ical\subset \Jcal\cap \Int(\catP)$, then $\rk_M^\Jcal$ is said to be M\"obius invertible \textbf{over intervals}. 
 \end{definition}
$\dgm_M^\Jcal$ is a function $\Jcal\to \Z$ such that, for each $I\in \Jcal$, the set $\{J\in \Jcal:\ J\supset I \ \mbox{and } \dgm_M^\Jcal(J)\neq 0\}$ is finite and
\begin{equation}\label{eq:mobius invertible over I_v2}
\rank^\Jcal_M(I)=\sum_{\substack{J\supset I\\ J\in \Jcal}}\dgm_M^\Jcal(J).
\end{equation}
This shows that the preceding notion of GPD is obtained by axiomatizing the \emph{fundamental lemma of persistent homology} \cite[Chapter VII]{edelsbrunner2008computational}, without assuming anything about the indexing poset $\catP$, the persistence module $M$, the domain $\Jcal$ of the GRI, or the structural simplicity of the GRI. Indeed, when $\catP$ is a finite totally ordered set, $\Jcal=\Int(\catP)$, and $M$ arises from applying the homology functor to a filtration over $\catP$ of a simplicial complex, $\dgm_M^\Jcal$ given in Equation (\ref{eq:mobius invertible over I_v2}) coincides with the classical persistence diagram and Equation (\ref{eq:mobius invertible over I_v2}) is known as the fundamental lemma of persistent homology. 

From another viewpoint, M\"obius invertibility of the GRI (i.e., the existence of the GPD) is a relaxation of interval-decomposability of persistence modules, as evidenced by Theorem \ref{thm:sufficient conditions for invertibility} \ref{item:interval-decomposable} below.

The desirable properties of the earlier GPD, presented in Remark~\ref{rem:additivity of GPD} and Proposition~\ref{prop:properties of GPD}, generalize in the following two propositions.

\begin{proposition}\label{prop:properties of GGPD} %
Let $M$ be a $\catP$-module and $\Jcal\subset \Con(\catP)$. If $\rk_M^\Jcal$ is M\"obius invertible, then:
\begin{enumerate}[leftmargin=*, label=(\roman*)]
\item %
$\Ical:=\supp \left( \dgm_M^{\Jcal}\right)\subset  \supp \left(\rank_M^{\Jcal}\right)$. 
\label{item:GGPD is more parsimonious}
\item $\Ical$ is locally finite, and the restriction of $\dgm_M^{\Jcal}$ to $\Ical$ is convolvable over $\Ical$.\label{item:convolvability of GGPD}
\item The unique function $\drm:\Jcal\rightarrow \R$ satisfying 
\[\rank_M^\Jcal(I)=\sum_{\substack{J\supset I\\ J\in \Jcal}}\drm(J) \ \ \mbox{for all $I\in \Jcal$},\]
where the sum contains only finitely many nonzero summands, is $\dgm_M^\Jcal$.
\label{item:universal property of GGPD}
\end{enumerate} 
\end{proposition}
\begin{proof}
\noindent \ref{item:universal property of GGPD}: This directly follows from the definition of M\"obius invertibility and Proposition~\ref{prop:justifying the term mobius invertibility}~\ref{item:justifying the term mobius invertibility-unique}.

\noindent \ref{item:convolvability of GGPD}: M\"obius invertibility of $\rk_M^\Jcal$ implies that 
$\Ical$ is locally finite (by Proposition~\ref{prop:justifying the term mobius invertibility}~\ref{item:justifying the term mobius invertibility-locally finite}), and
that, for all $J\in \Jcal$, there are finitely many $I\in \Ical$ containing $J$. Hence, for all $I\in \Ical (\subset \Jcal)$, there are finitely many $J\in \Ical$ containing $I$, i.e., $I^{\downarrow}_{\Ical}$ is finite. This precisely proves that $\dgm_M^\Jcal|_{\Ical}$ is convolvable over $\Ical$.

\noindent \ref{item:GGPD is more parsimonious}: By the original definition of the GPD (Definition~\ref{def:GPD over I}) and Proposition~\ref{prop:justifying the term mobius invertibility}~\ref{item:justifying the term mobius invertibility-mobius inversion}, \[
\dgm_M^\Ical=\rk_M^\Ical\ast\mu_{\Ical}=\dgm_M^{\Jcal}|_{\Ical}.
\]
Hence, we have $\rk_M^\Ical=\dgm_M^\Ical\ast\zeta_{\Ical}$ by Proposition~\ref{prop:properties of GPD}~\ref{item:universal property of GPD}. Since we saw $\dgm_M^\Jcal|_{\Ical}(=\dgm_M^\Ical)$ is convolvable, Proposition~\ref{prop:convolvability}~\ref{item:convolvable6} implies that $\rk_M^\Ical$ is also convolvable. Then, Proposition~\ref{prop:mobius inversion of non-decreasing map is convolvable} and monotonicity of $\rk_M^\Ical$ imply that
\[\Ical\left(=
\supp \left( \dgm_M^\Ical\right)\right) \subset \supp \left(\rk_M^\Ical\right) \subset \supp \left( \rk_M^\Jcal\right),\]
as desired.
\end{proof}
\begin{proposition}[Additivity]\label{prop:GPD additivity}
{Let $M$ and $N$ be $\catP$-modules and $\Jcal \subset \Con(\catP)$. If $\rk_M^\Jcal$ and $\rk_N^\Jcal$ are M\"obius invertible, then so is $\rk_{M\oplus N}^\Jcal$, and \[\dgm^\Jcal_{M\oplus N}=\dgm^\Jcal_{M }+\dgm^\Jcal_{N}.\]}
\end{proposition}
\begin{proof}
Let $I\in \Jcal$. Then,
\begin{align*}
    \rk^{\Jcal}_{M\oplus N}(I)&=\rk^{\Jcal}_{M}(I)+\rk^{\Jcal}_{N}(I)&\mbox{by Remark~\ref{rem:properties of rank invariant}~\ref{item:additivity}}
    \\ &= \sum_{\substack{J\supset I\\J\in \Jcal}} \dgm_M^\Jcal(J)+\sum_{\substack{J\supset I\\J\in \Jcal}} \dgm_N^\Jcal(J)&\mbox{by Definition~\ref{def:GGPD}}
    \\&=\sum_{\substack{J\supset I\\J\in \Jcal}} \left(\dgm_M^\Jcal(J)+\dgm_N^\Jcal(J)\right).
\end{align*}
Hence, by Definition~\ref{def:GGPD}, $\rk^{\Jcal}_{M\oplus N}$ is M\"obius invertible over $\Jcal$, and the unique function $\dgm^\Jcal_{M\oplus N}$ (cf. Proposition~\ref{prop:properties of GGPD}~\ref{item:universal property of GGPD}) is identical to $\dgm^\Jcal_{M }+\dgm^\Jcal_{N}.$
\end{proof}

\begin{remark}[{Comparing Definition \ref{def:GGPD} with ideas from \cite[Section 8]{botnan2024signed}}]\label{rem:comparison with BOO21}
\white{.} 
\phantomsection

\begin{enumerate}[label=(\roman*)]
\item 
When the set $\Jcal$ in Definition \ref{def:GGPD} consists solely of intervals of $\catP$,
then saying that $\rk_M^\Jcal$ is M\"obius invertible over $\Ical$ is equivalent to saying that the GRI of $M$ with \emph{test set} $\Jcal$ admits a rank decomposition with \emph{dictionary} $\Ical$ \cite[Section 8]{botnan2024signed}, i.e., there exist $\Ical$-decomposable $\catP$-modules $N_1$ and $N_2$ such that
\[\rk_M^{\Jcal}=\rk_{N_1}^{\Jcal}-\rk_{N_2}^\Jcal.\]
This can be proved in a similar way as Theorem \ref{them:rk decomposition is GPD}. 
However, Definition \ref{def:GGPD} slightly generalizes ideas in \cite[Section 8]{botnan2024signed} in that neither the test set nor the dictionary set needs to be a set of intervals.
This level of generality is useful in building a connection between invariants of persistence modules or clarifying the discriminating power of the GRI; {see, for example, \cite{clause2026generalized,kim2021generalized,kim2021bettis}.} 
\label{item:GPD is more general than rank decomposition}

\item 
In contrast to \cite[Section 8]{botnan2021signed}, we consider the case where $\Ical\subsetneq \Jcal$ to be important both in theory and in practice: 
 Since $\dgm_M^\Ical$ and $\rank_M^{\Ical}$ determine each other (cf. Proposition \ref{prop:properties of GPD} \ref{item:universal property of GPD}), the M\"obius invertibility of {$\rk_M^\Jcal$} over $\Ical$ implies that 
{$\rk_M^\Jcal$} can be restricted to the smaller domain $\Ical$ \emph{without losing any information.} In particular, significant discrepancy between $\Ical$ and $\Jcal$ indicates $\rk^\Jcal_M$ is mostly constant.  
Identifying such an $\mathcal{I}$, ideally the minimal one, is crucial to minimize the computational burden as well as facilitating vectorization of the information for usage in machine learning pipeline \cite{loiseaux2024stable,xin2023gril}.  Consider, for instance, a finitely presentable $\R^2$-module $M$. While the poset $\Jcal:=(\Int(\R^2),\supset)$ is not locally finite, $
\rk_M^\Jcal$ is completely encoded into $\dgm_M^{\Ical}$ 
(which follows from Theorem \ref{thm:sufficient conditions for invertibility} \ref{item:finitely presentable invertible} below). 
\end{enumerate}
\end{remark}

The following example demonstrates that M\"obius invertibility is a strictly more general notion than constructibility (Definition \ref{def:co-closure and constructible}).

\begin{example}\label{ex:mobius invertibility is different from constructibility} Consider the $\R$-module $M:=k_{[0,2)}\oplus k_{[1,3)}$. The rank invariant $\rk_M^\Int$ is M\"obius invertible over $\Ical:=\{[0,2),[1,3)\}$ (which follows from Theorem~\ref{thm:sufficient conditions for invertibility}~\ref{item:interval-decomposable} below). Nevertheless, no co-closure $c:\Int(\R)\to\Int(\R)$ with image $\Ical$ can satisfy $\rk_M^\Int\circ c=\rk_M^\Int$, since $\rk_M^\Int\circ c$ has maximum value $1$, while $\rk_M^\Int$ achieves the value $2$.
\end{example}

\subsection{{M\"obius invertibility and rank decomposability}}
\label{sec:Mobius invertibility and rank decompositions}

In this section, we show that rank decomposability of the GRI, introduced in \cite{botnan2024signed}, is subsumed by the framework of M\"obius invertibility of the GRI. One implication is, by Proposition \ref{prop:justifying the term mobius invertibility}~\ref{item:justifying the term mobius invertibility-mobius inversion}, that every rank decomposition of the GRI arises from standard M\"obius inversion (cf. Section \ref{sec:mobius inversion}), \emph{regardless of whether the domain of the GRI is locally finite.}

\begin{theorem}\label{thm:rk decomposition is GPD}  Let $M$ be any $\catP$-module. For any  $\Jcal\subset\Int(\catP)$, not necessarily locally finite, the following are equivalent.
    \begin{enumerate}[label=(\roman*),leftmargin=*]
        \item $\rank_M^\Jcal$ is M\"obius invertible. 
        \label{item:1}
         \item $\rank_M^\Jcal$ admits a rank decomposition.  \label{item:2}     
    \end{enumerate}
\end{theorem}

\begin{proof}
\noindent
\ref{item:1} $\Rightarrow$ \ref{item:2}:
Let $\Ical:=\supp(\dgm_M^\Jcal).$
By Proposition
\ref{prop:justifying the term mobius invertibility}%
~\ref{item:justifying the term mobius invertibility-locally finite},
\ref{item:justifying the term mobius invertibility-convolvable},
and
\ref{item:justifying the term mobius invertibility-mobius inversion},
the poset $\Ical$ is locally finite, $\rk_M^\Ical$ is convolvable,
and
\[
\dgm_M^\Ical
=
\rk_M^\Ical\ast\mu_\Ical
=
\left.\dgm_M^\Jcal\right|_\Ical.
\]
Therefore, for every $I\in\Jcal$,
\[
\rank_M^\Jcal(I)
=
\sum_{\substack{J\supset I\\J\in\Ical}}
\dgm_M^\Jcal(J)
=
\sum_{\substack{J\supset I\\J\in\Ical}}
\dgm_M^\Ical(J).
\]
Applying the argument of
Equations~\eqref{eq:line1}--\eqref{eq:line5} to the preceding
identity yields disjoint multisets $\Rcal_M$ and $\Scal_M$ of elements in $\Jcal$ satisfying
\[
\rank_M^\Jcal
=
\rank_{\kf_{\Rcal_M}}^\Jcal
-
\rank_{\kf_{\Scal_M}}^\Jcal.
\]
Moreover, since $\{J\in \Jcal:\dgm_M^{\Jcal}(J)\neq 0,\ J\supset I\} $ is finite for every $I\in \Jcal$ (cf. Item \ref{item:mobius invertible general 1} of
Definition~\ref{def:mobius invertible general}),
$\Rcal_M$ and $\Scal_M$ are $\Jcal$-finite. Hence,
$(\Rcal_M,\Scal_M)$ is a rank decomposition of $\rank_M^\Jcal$.

\medskip

\noindent
\ref{item:2} $\Rightarrow$ \ref{item:1}:
Let $(\Rcal,\Scal)$ be a rank decomposition of $\rank_M^\Jcal$.
Define
\begin{align*}
\drm_M:\Jcal&\longrightarrow\Z,\\
I&\longmapsto
\mult_\Rcal(I)-\mult_\Scal(I).
\end{align*}
Since $\Rcal$ and $\Scal$ are $\Jcal$-finite, for every
$I\in\Jcal$, the set
\[
\left\{
J\in\Jcal:
J\supset I
\text{ and }
\drm_M(J)\neq 0
\right\}
\]
is finite. Moreover, for every $I\in\Jcal$,
\begin{align*}
\rank_M^\Jcal(I)
&=
\left(
\rank_{\kf_\Rcal}^\Jcal
-
\rank_{\kf_\Scal}^\Jcal
\right)(I)\\
&=
\sum_{\substack{J\supset I\\J\in\Jcal}}
\mult_\Rcal(J)
-
\sum_{\substack{J\supset I\\J\in\Jcal}}
\mult_\Scal(J)\\
&=
\sum_{\substack{J\supset I\\J\in\Jcal}}
\drm_M(J),
\end{align*}
where the second equality follows from
\cite[Proposition~2.1]{botnan2024signed}.
Thus, by Definition~\ref{def:mobius invertible general},
$\rank_M^\Jcal$ is M\"obius invertible, with $\dgm_M^\Jcal=\drm_M.$
\end{proof}

\subsection{M\"obius invertibility and other measures of structural simplicity in persistence modules
}\label{sec:On structural simplicity of persistence modules}

In this section, we elucidate the relationships among various concepts pertaining to the structural simplicity of persistence modules, as depicted in Figure \ref{fig:tameness + MI visual}, extending beyond Proposition~\ref{prop:finitely presentable implies tame}. 
{Also, as a byproduct, we obtain counterexamples to the statement given in \cite[Corollary~3.18]{botnan2024signed}}: see Propositions \ref{prop:constructible implies finitely  presentable} and  \ref{prop:tame does not imply RI invertible} and Theorem \ref{thm:tame does not imply Int-GRI invertible} and Remark~\ref{rem:about tame does not imply Int-GRI mobius invertible}.

The following proposition is  straightforward from Definitions \ref{def:finitely presentable}, \ref{def:finite encoding}, and \ref{def:co-closure and constructible} together with \cite[Proposition 7.5]{gulen2022galois}.
\begin{proposition}\label{prop:constructible implies finitely presentable}
Constructible persistence modules are finitely presentable. 
\end{proposition} 
\begin{proof}
Let $M:\catP\rightarrow \cvec$ be an $\catS$-constructible module with a finite set $\catS \subset \catP$. Then, there exist an $\catS$-constructible $\catP$-module $F_0=\bigoplus_{i=1}^n k_{p_i^\uparrow}$ and a natural transformation $\Psi^0:F_0\Rightarrow M$ such that 
the map $\Psi_a^0:(F_0)_a\rightarrow M_a$ is surjective for each $a\in \catP$ 
\cite[Proposition 7.5]{gulen2022galois}.
Since $F_0$ is $S$-constructible, there exists a co-closure $c:\catP\to\catP$, where $S$ is the image of $c$ with $F_0\circ c=F_0$. Then, we have that $\ker(\Psi^0)\circ c=\ker(\Psi^0)$ since $\ker(\Psi^0)$ is a submodule of $F_0$. This implies that $\ker(\Psi^0)$ is also $\catS$-constructible. Then, again by 
\cite[Proposition 7.5]{gulen2022galois}, there exists another $\catS$-constructible $\catP$-module $F_1=\bigoplus_{i=1}^mk_{q_i^\uparrow}$  with a natural transformation $\Psi^1:F_1\Rightarrow F_0$
such that $\im(\Psi^1)=\ker(\Psi^0)$.
Now we have that $M$ is the cokernel of $\Psi^1$, completing the proof.
\end{proof}

While the converse of the preceding proposition does not hold strictly speaking, constructibility and finite presentability are almost identical concepts in a suitable sense; see \cite[Proposition C.2]{kim2026interleaving}.

In the next proposition, we see that finitely presentability of a $\Z^d$- or $\R^d$-module (Definition~\ref{def:finite encoding}) does not guarantee M\"obius invertibility of its rank invariant.

\begin{proposition}\phantomsection\label{prop:tame does not imply RI invertible}
For each $d\geq 1$, there exist finitely presentable $\Z^d$- or $\R^d$-modules ({and thus tame by Proposition~\ref{prop:finitely presentable implies tame}}) whose rank invariants are not M\"obius invertible. 
\end{proposition}

One way to remedy this failure of M\"obius inversion is, somewhat ironically, to enlarge the domain of the rank invariant.
When $d=1$, extending the domain to $\Int(\Z)$ or $\Int(\R)$ guarantees
M\"obius invertibility of the resulting GRI, as a consequence of
Theorem~\ref{thm:sufficient conditions for invertibility}%
~\ref{item:interval-decomposable}. Note that $\Int(\Z)$ and $\Int(\R)$ are obtained from $\seg(\Z)$ and $\seg(\R)$, respectively, by adding infinite intervals. More generally, for arbitrary $d\geq 1$,
it actually suffices to extend the domain to the collection of ``half-open rectangles'';
see Corollary~\ref{cor:BOO's Cor 3.18} (i.e., \cite[Corollary 2.18]{botnan2024signed}),
Proposition~\ref{prop:invertible over half-open}, and the subsequent remark.

Another remark is that the preceding theorem implies the existence of modules whose rank invariant is not M\"obius invertible over \emph{any} poset containing $\Z^d$ or $\R^d$ as a convex subposet for $d \geq 1$.

\begin{proof}
Let $d\geq 1$, let $\catP\in\{\Z^d,\R^d\}$, and consider the
$\catP$-module $M:=\kf_{0^{\uparrow}}.$
Then $M$ is clearly finitely presentable. We claim, however, its rank invariant $\rk_M^{\seg}$ is not M\"obius invertible.

For every
$u\in\catP$ with $0\leq u$, we have $[0,u]\subseteq0^{\uparrow}$,
and hence $\rk_M^{\seg}([0,u])=1.$
Suppose, by way of contradiction, that
$\rk_M^{\seg}$
is M\"obius invertible, and let $f:\seg(\catP)\to\Z$
be its M\"obius inversion. Set $I_0:=\{0\}$, $\Ical:=\supp(f)$, and
\[
\Ical_{\supset I_0}
:=\left\{
J\in\Ical:
J\supset I_0
\right\} =
\left\{
J\in\seg(\catP):
J\supset I_0
\text{ and }
f(J)\neq 0
\right\}.
\]
By M\"obius invertibility, the set $\Ical_{\supset I_0}$ is finite. It is also
nonempty, since
\[
1
=
g(I_0)
=
\sum_{\substack{J\in\seg(\catP)\\J\supset I_0}}f(J).
\]
For each $J\in\Ical_{\supset I_0}$, write $J=[a_J,b_J]$, where $a_J,b_J\in\catP$. Define $u\in\catP$ coordinatewise by
\[
u_i
:=
1+\max_{J\in\mathcal A}(b_J)_i,
\qquad i=1,\ldots,d,
\]
and let
$K:=[0,u].$
Since $K\subset0^{\uparrow}$, we have $\rk_M^\seg(K)=1.$

On the other hand, suppose that $J\in\seg(\catP)$ satisfies
$J\supset K$ and $f(J)\neq 0$. Since $K\supset I_0$, we must have
$J\in\Ical_{\supset I_0}$. Also, $u\in J=[a_J,b_J]$, and therefore $u\leq b_J.$
This contradicts the choice of $u$, according to which
\[
u_i>(b_J)_i
\qquad
\text{for every }i=1,\ldots,d.
\]
Consequently, we have that
$f(J)=0$
for every $J\in\seg(\catP)\text{ with }J\supset K.
$
It follows that
\[
1
=
g(K)
=
\sum_{\substack{J\in\seg(\catP)\\J\supset K}}f(J)
=
0,
\]
which is a contradiction. Therefore, the rank invariant $\rk_M^{\seg}$ is not
M\"obius invertible.
\end{proof}

The next theorem is the counterpart of the previous proposition for the Int-GRI. 

\begin{theorem}\phantomsection\label{thm:tame does not imply Int-GRI invertible}
\begin{enumerate}[label=(\roman*)]
    \item \label{item:tame not imply MI} For each $d\geq 2$, there exist tame $\Z^d$- and $\R^d$-modules whose Int-GRIs are not M\"obius invertible.
    \item \label{item:fp not imply MI} For each $d\geq 3$, there exist  finitely presentable $\Z^d$- and $\R^d$-modules whose Int-GRIs are not M\"obius invertible. 
\end{enumerate}
\end{theorem}

{The preceding theorem implies the existence of modules whose Int-GRI is not M\"obius invertible over \emph{any} poset containing $\Z^d$ or $\R^d$ as a convex subposet for $d \geq 2$. More remarks follow.} %
\begin{remark}\label{rem:about tame does not imply Int-GRI mobius invertible}
\begin{enumerate}[label=(\roman*)]

\item In contrast to Theorem \ref{thm:tame does not imply Int-GRI invertible}, we will see in Corollary~\ref{cor:Mobius invertibility does not imply tameness}  below that there exists a persistence module that is neither tame nor finitely presentable, yet whose Int-GRI is M\"obius invertible. Also, in contrast to Theorem \ref{thm:tame does not imply Int-GRI invertible}~\ref{item:fp not imply MI}, we will see in Theorem~\ref{thm:sufficient conditions for invertibility}~\ref{item:interval-decomposable} and \ref{item:finitely presentable invertible} that, for $d=1,2$, the Int-GRI of every finitely presentable $\Z^d$- or $\R^d$-module is M\"obius invertible. 

\item 
We refer the reader to Remark~\ref{rem:conflict} for a discussion of how the techniques used in the proof of Theorem~\ref{thm:tame does not imply Int-GRI invertible}~\ref{item:fp not imply MI} rectify 
\cite[Corollary 3.18]{botnan2024signed}.

\end{enumerate}
\end{remark}

For establishing 
the previous two theorems, we introduce a useful proposition. Let $M$ be any given $\catP$-module and $\Ical\subset\Con(\catP)$. We call $I\in \Ical$ \textbf{support-maximal} %
if $I\in \supp(\rk_M^\Ical)$ and any $J \in \Ical$ properly containing $I$ does not belong to $\supp(\rk_M^\Ical)$.
\begin{proposition}\label{prop:property of maximal intervals} Assume that  $\rk_M^{\Ical}$ is M\"obius invertible (i.e., $\dgm_M^\Ical$ exists, as given in Definition \ref{def:GGPD}). If $I\in \Ical$ is support-maximal, then 
$\rk_M^\Ical(I)=\dgm_M^\Ical(I),$ implying that $I\in\supp(\dgm_M^{\Ical})$.
\end{proposition}

\begin{proof}  By Definition~\ref{def:GGPD}, 
\[\rank^\Ical_M(I)=\sum_{\substack{J\supset I\\ J\in \Ical}}\dgm_M^\Ical(J),\]
where the RHS contains only finitely many nonzero summands. In particular,  Proposition~\ref{prop:properties of GGPD}~\ref{item:GGPD is more parsimonious} together with the maximality of $I$ implies that the RHS is precisely $\dgm_M^\Ical(I)$.
\end{proof}

\begin{corollary}\label{cor:intersection}If there exist an infinite collection of support-maximal intervals whose intersection contains an element of $\Ical$, then $\rk_M^\Ical$ is not M\"obius invertible. 
\end{corollary}

\begin{proof} Let $\{I_a\}_a$ be an infinite set of support-maximal intervals in $\Ical$ such that the intersection $\bigcap_a I_a$ contains some $I\in \Ical$. %
Suppose that $\rk_M^\Ical$ is M\"obius invertible.  Then, by the preceding proposition, $\dgm_M^\Ical(I_a)\neq 0$ for all $a$. 
Then $\{I_a\}_a$ is contained in the set $\{J\in \Ical:J\supset I \ \mbox{and } \dgm_M^\Ical(J)\neq 0\}$, which is therefore  infinite. This contradicts to the M\"obius invertibility of $\rk_M^\Ical$ (Definition~\ref{def:mobius invertible general}~\ref{item:mobius invertible general 1}).
\end{proof}

We remark that the converse to Corollary \ref{cor:intersection} does not hold, as shown by the following argument.
For example, consider any $\Z^3$-module $M$ whose Int-GRI is not M\"obius invertible (such an $M$ is constructed in the proof of Theorem~\ref{thm:tame does not imply Int-GRI invertible}~\ref{item:fp not imply MI}).
Let $M':=M\oplus \kf_{\Z^3}$. Then, $\rk_{M'}^\Int$ clearly has only one support-maximal interval $\Z^3$.
We claim, however, that the Int-GRI of $M'$ is not M\"obius invertible.
If it were, then 
it would follow that
\[
\dgm^\Int_{M'}-\dgm^\Int_{\kf_{\Z^3}}
\]
is a valid GPD for $M$, contradicting the fact that the Int-GRI of $M$ is not M\"obius invertible.

\begin{proof}[Proof of Theorem \ref{thm:tame does not imply Int-GRI invertible}  \ref{item:tame not imply MI}]

\begin{figure}[ht]
\centering
\begin{tikzpicture}[
    >=Stealth,
    modulearrow/.style={
        ->,
        line width=0.8pt,
        >={Stealth[length=2.2mm,width=1.5mm]}
    },
    guide/.style={
        densely dotted,
        line width=0.8pt
    },
    obj/.style={
        inner sep=1pt,
        font=\normalsize
    },
    maplabel/.style={
        font=\scriptsize,
        inner sep=1pt,
        fill=white
    }
]

\begin{scope}[shift={(0,7.2)}]
\node[font=\Large] at (3.75,4.70) {$M$};
\node[anchor=west] at (-0.75,4.52) {$y=-x+1$};
\node[anchor=west] at (-0.18,3.63) {$y=-x$};
\draw[guide] (1.30,4.58) -- (1.73,4.20);
\draw[guide] (1.30,3.55) -- (1.73,3.18);
\draw[guide] (4.92,0.12) -- (5.30,-0.25);
\draw[guide] (5.92,0.12) -- (6.30,-0.25);
\node[obj] (A00) at (2.10,0.40) {$k$};
\node[obj,text=moduleOrange] (A10) at (2.10,3.05) {$k^2$};
\node[obj,text=moduleRed] (A20) at (2.10,4.05) {$k$};
\node[obj,text=moduleGreen] (A11) at (2.95,2.15) {$k^2$};
\node[obj,text=moduleBlue] (A21) at (2.95,3.05) {$k$};
\node[obj,text=moduleOrange] (A12) at (3.85,1.30) {$k^2$};
\node[obj,text=moduleRed] (A22) at (3.85,2.15) {$k$};
\node[obj,text=moduleGreen] (A13) at (4.75,0.40) {$k^2$};
\node[obj,text=moduleBlue] (A23) at (4.75,1.30) {$k$};
\node[obj,text=moduleRed] (A14) at (5.75,0.40) {$k$};
\node[obj,text=gray] (Azero) at (5.75,4.05) {$0$};
\draw[modulearrow] (A00) -- (A10);
\draw[modulearrow] (A10) -- (A20);
\draw[modulearrow] (A11) -- (A21);
\draw[modulearrow] (A12) -- (A22);
\draw[modulearrow] (A13) -- (A23);
\draw[modulearrow] (A14) -- (Azero);
\draw[modulearrow] (A10) -- (A21);
\draw[modulearrow] (A11) -- (A22);
\draw[modulearrow] (A12) -- (A23);
\draw[modulearrow] (A13) -- (A14);
\draw[modulearrow]
($(A00.east)+(0.08,0.10)$)
--
($(A11.south west)+(-0.0,-0.15)$);
\draw[modulearrow]
($(A00.east)+(0.12,0.03)$)
--
($(A12.west)+(-0.10,-0.10)$);
\draw[modulearrow]
($(A00.east)+(0.10,-0.03)$)
--
($(A13.west)+(-0.10,0)$);
\draw[modulearrow] (A20) -- (Azero);
\draw[modulearrow]
($(A21.east)+(0.05,0.03)$)
--
($(Azero.west)+(-0.10,-0.20)$);
\draw[modulearrow]
($(A22.east)+(0.05,0.05)$)
--
($(Azero.south west)+(-0.05,-0.05)$);
\draw[modulearrow]
($(A23.east)+(0.04,0.05)$)
--
($(Azero.south)+(-0.12,-0.06)$);
\node[font=\large] at (3.75,-0.60) {(A)};
\end{scope}

\begin{scope}[shift={(8.0,7.2)}]

\node[font=\Large] at (3.10,4.70) {$N$};

\node[obj,text=gray]          (B0)  at (3.10,4.10) {$0$};

\node[obj,text=moduleRed]     (B1L) at (1.30,3.00) {$k$};
\node[obj,text=moduleBlue]    (B1R) at (4.90,3.00) {$k$};

\node[obj,text=moduleOrange]  (B2L) at (1.30,1.40) {$k^2$};
\node[obj,text=moduleGreen]   (B2R) at (4.90,1.40) {$k^2$};

\node[obj]                    (B3)  at (3.10,0.15) {$k$};

\draw[modulearrow] (B1L) -- (B0);
\draw[modulearrow] (B1R) -- (B0);

\draw[modulearrow] (B2L) -- (B1L);
\draw[modulearrow] (B2R) -- (B1R);

\draw[modulearrow]
    (B2L)
    --
    node[pos=.62,above=7pt,xshift=4pt,maplabel]
        {$\begin{bmatrix}0&1\end{bmatrix}$}
    (B1R);

\draw[modulearrow]
    (B2R)
    --
    node[pos=.62,above=7pt,xshift=-4pt,maplabel]
        {$\begin{bmatrix}0&1\end{bmatrix}$}
    (B1L);

\draw[modulearrow]
    (B3)
    --
    node[pos=.48,left=8pt,xshift=-4pt,yshift=-5pt,maplabel]
        {$\begin{bmatrix}1\\[1pt]1\end{bmatrix}$}
    (B2L);

\draw[modulearrow]
    (B3)
    --
    node[pos=.48,right=8pt,xshift=4pt,yshift=-5pt,maplabel]
        {$\begin{bmatrix}1\\[1pt]1\end{bmatrix}$}
    (B2R);

\node[maplabel,anchor=east]
    at (0.95,2.20)
    {$\begin{bmatrix}1&0\end{bmatrix}$};

\node[maplabel,anchor=west]
    at (5.25,2.20)
    {$\begin{bmatrix}1&0\end{bmatrix}$};
\node[font=\large] at (3.10,-0.60) {(B)};

\end{scope}

\begin{scope}[shift={(7.4,1.5)},x=.78cm,y=.78cm]
\begin{scope}[shift={(-0.5,-0.5)}]
\clip (-4.8,-4.25) rectangle (4.9,5.25);

\foreach \x in {-4,...,6} {
    \pgfmathtruncatemacro{\ytop}{-\x+1}
    \fill[regionBlue]
        (\x,-3) rectangle ({\x+1},{\ytop});
}

\foreach \x in {-5,-3,-1} {
    \pgfmathtruncatemacro{\yextra}{-\x+1}
    \fill[regionBlue]
        (\x,\yextra) rectangle ({\x+1},{\yextra+1});
}

\foreach \x in {2,4,6} {
    \pgfmathtruncatemacro{\yextra}{-\x+1}
    \fill[regionBlue]
        (\x,\yextra) rectangle ({\x+1},{\yextra+1});
}
\end{scope}

\draw[
    <->,
    gray!55,
    line width=1pt,
    >={Stealth[length=2.8mm,width=2mm]}
]
(-4.65,0) -- (4.75,0);

\draw[
    <->,
    gray!55,
    line width=1pt,
    >={Stealth[length=2.8mm,width=2mm]}
]
(0,-4.10) -- (0,5.05);

\node[anchor=south west,font=\large]
    at (4.40,0.10) {$x$};

\node[anchor=south west,font=\large]
    at (0.12,4.68) {$y$};

\node[font=\large]
    at (2.45,3.35) {$I_0$};

\node[obj] (L1) at (-4.00,4) {\textcolor{orange}{$k^2$}};
\node[obj] (R1) at (-3.0,4) {\textcolor{blue}{$k$}};

\node[obj] (L2) at (-3.0,3) {\textcolor{moduleGreen}{$k^2$}};
\node[obj] (R2) at (-2.0,3) {\textcolor{red}{$k$}};

\node[obj] (L3) at (-2.0,2) {\textcolor{orange}{$k^2$}};
\node[obj] (R3) at (-1,2) {\textcolor{blue}{$k$}};

\node[obj] (L4) at (-1,1.00) {\textcolor{moduleGreen}{$k^2$}};
\node[obj] (R4) at (-0,1.00) {\textcolor{red}{$k$}};

\node[obj] (L5) at (-0,0.0) {\textcolor{orange}{$k^2$}};
\node[obj] (R5) at (1,0.0) {\textcolor{blue}{$k$}};

\node[obj] (L6) at (1,-1) {\textcolor{moduleGreen}{$k^2$}};
\node[obj] (R6) at (2,-1) {\textcolor{red}{$k$}};

\node[obj] (L7) at (2,-2) {\textcolor{orange}{$k^2$}};
\node[obj] (R7) at (3,-2) {\textcolor{blue}{$k$}};

\node[obj] (L8) at (3,-3) {\textcolor{moduleGreen}{$k^2$}};
\node[obj] (R8) at (4,-3) {\textcolor{red}{$k$}};

\draw[modulearrow] (L1.east) -- (R1.west);
\draw[modulearrow] (L2.east) -- (R2.west);
\draw[modulearrow] (L3.east) -- (R3.west);
\draw[modulearrow] (L4.east) -- (R4.west);
\draw[modulearrow] (L5.east) -- (R5.west);
\draw[modulearrow] (L6.east) -- (R6.west);
\draw[modulearrow] (L7.east) -- (R7.west);
\draw[modulearrow] (L8.east) -- (R8.west);

\draw[modulearrow] (L2.north) -- (R1.south);
\draw[modulearrow] (L3.north) -- (R2.south);
\draw[modulearrow] (L4.north) -- (R3.south);
\draw[modulearrow] (L5.north) -- (R4.south);
\draw[modulearrow] (L6.north) -- (R5.south);
\draw[modulearrow] (L7.north) -- (R6.south);
\draw[modulearrow] (L8.north) -- (R7.south);

\node[obj] at (-4.0,-3) {$k$};

\node[font=\large] at (0.00,-4.75) {(C)};

\end{scope}
\end{tikzpicture}

\caption{(A) and (B) Illustrations of the $\Z^2$-module $M$ and $\catQ$-module $N$ from the proof of Theorem \ref{thm:tame does not imply Int-GRI invertible}. (C) A visualization of $M$ and the interval $I_0$ used in the proof of Theorem \ref{thm:tame does not imply Int-GRI invertible}.}\label{fig:tame does not imply Int-GRI invertible module}
\end{figure}
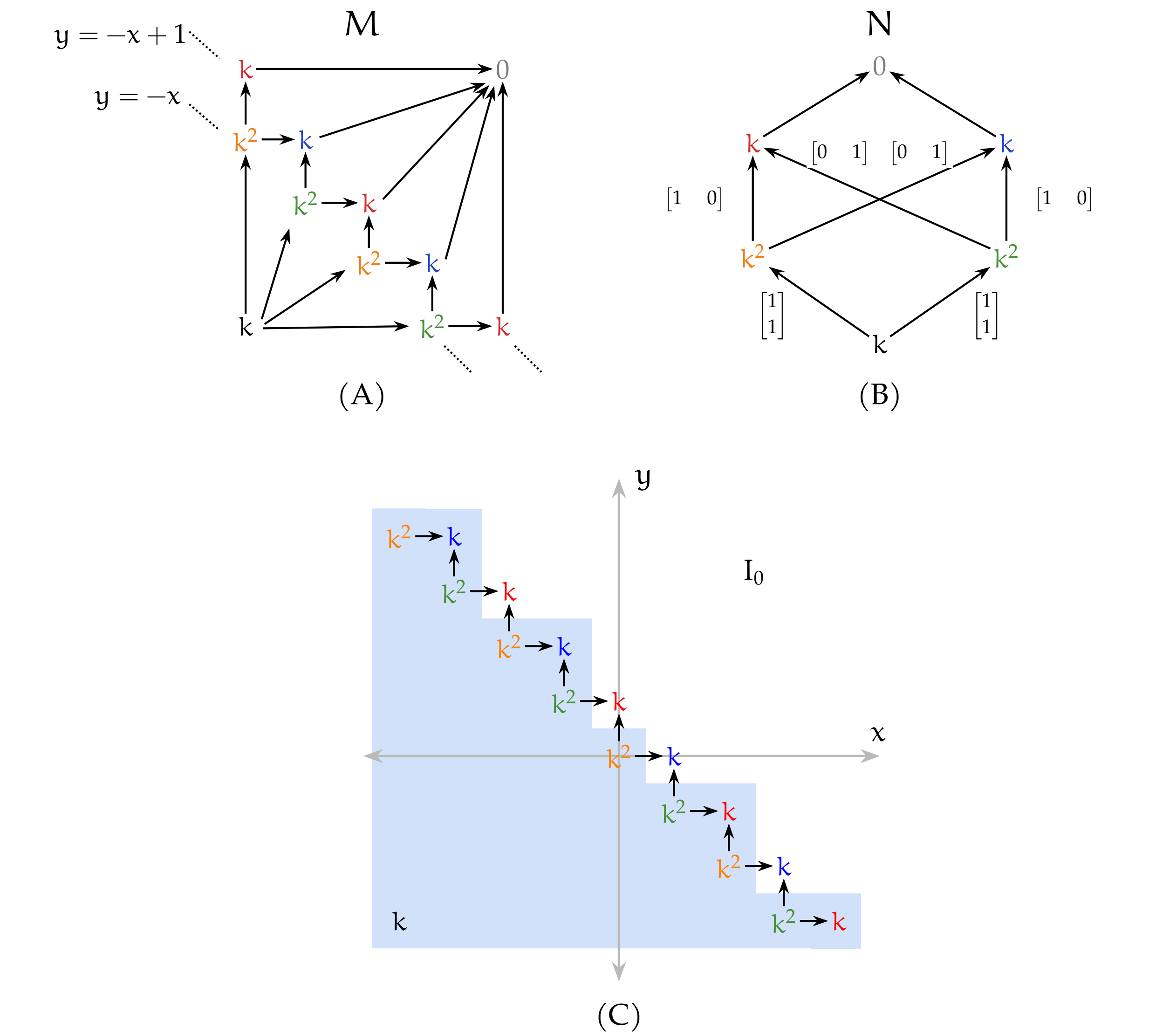

    We first construct a $\Z^2$-module that is tame while its Int-GRI is not M\"obius invertible.
    Let $N$ be the persistence module defined in Figure~\ref{fig:tame does not imply Int-GRI invertible module} (B), over the 6-point poset that we denote by $\catQ$.
    Let $M$ be the $\Z^2$-module that is the pullback $\pi^\ast N$, where $\pi:\Z^2\to \catQ$ is the order-preserving map depicted in Figure \ref{fig:tame does not imply Int-GRI invertible module} (A) and (B); points indexing red vector spaces $k$ in (A) map to points indexing the red vector space $k$ in (B), points indexing blue vector spaces $k$ in (A) map to points indexing the blue vector space $k$ in (B), points below the line $y=-x$ in (A) map to the point indexing the black vector space $k$, and points above the line $y=-x+1$ map to the point indexing the gray vector space $0$.
    Then, by construction, $M$ is tame.

    We show that $\rk^\Int_M$ is \emph{not} M\"obius invertible.    Let $I_0 \subset \Z^2$ consist of the points at or below the line $y=-x$, together with every other point on the line $y=-x+1$, except that the two consecutive points $(0,1)$ and $(1,0)$ are omitted (see Figure~\ref{fig:tame does not imply Int-GRI invertible module} (C)); namely,  
     \[I_0:=\{(x,y)\in \Z^2:y\leq-x\} \cup \{(-2i-1,2i+2)\}_{i=0}^\infty \cup \{(2i+2,-2i-1)\}_{i=0}^\infty.\]
    Then for each $a\in \Z$, let $I_a\subset \Z^2$ be the shift of $I_0$ by $(a,-a)$. %
    Note that $I_a$ is an interval of $\Z^2$. 

In the following two claims, we show that for all $a\in \Z$, $I_a$ is a support-maximal interval in $\Z^2$.
Then, since the singleton interval $\{(-1,-1)\}\in \Int(\Z^2)$ is contained in $I_a$ for all $a\in \Z$, Corollary \ref{cor:intersection} implies that the Int-GRI of $M$ is not M\"obius invertible.

    \begin{claim}\label{claim:I_a nonzero rank}
        For each $a\in \Z$, $\rk^\Int_M(I_a) = 1$.
    \end{claim}
    \begin{claimproof}
        Following Proposition \ref{prop:limit and colimit}, we first describe a nonzero section of  $M\vert_{I_a}$. For $p\in I_a$, let \[\ell_p=\begin{cases}
            1\in \kf,&\mbox{if $M_p=\kf$}\\(1,1)\in \kf^2,&\mbox{if $M_p=\kf^2$.}
        \end{cases}\]
        Then $(\ell_p)_{p \in I_a}\in \prod_{p \in I_a}M_p$ is a section of $M|_{I_a}$, and
        the image of this section in the colimit is nonzero; 
        to see this, pick any $p\in I_a$, and let $[\ell_p]$ denote the image of the corresponding section in the colimit. 
By Proposition~\ref{prop:path characterization of zero in colim}, it suffices to show that
there is no section $\mathbf{v}$ along a path $p$, $\ldots$, $q$ in $I_a$ with $\mathbf{v}_p=\ell_p$ and $\mathbf{v}_q=0\in M_q$. 
Note that every morphism in $M|_{I_a}$ is one of
\begin{equation}\label{eq:possible morphisms}
\kf \xrightarrow{\begin{bmatrix}1\\1\end{bmatrix}} \kf^2,
\qquad
\kf^2\xrightarrow{\begin{bmatrix}1&0\end{bmatrix}}\kf,
\qquad
\kf^2\xrightarrow{\begin{bmatrix}0&1\end{bmatrix}}\kf,
\qquad
\kf\xrightarrow{1}\kf,
\end{equation}
and that a section along any path in $I_a$ starting at $p$ with value $\ell_p$ (i.e., $1$ or $(1,1)$) can take only the nonzero value $1\in\kf$ or one of the nonzero values $(1,0)$, $(0,1)$, and $(1,1)$ in $\kf^2$. Thus, $[\ell_p]$ is nonzero in the colimit.
Hence, $\rk^\Int_M(I_a)\geq 1$.
        Also, as $(-1,-1)\in I_a$ and $M_{(-1,-1)}=\kf$ has dimension 1, we have $\rk^\Int_M(I_a)\leq 1$ (cf. Remark~\ref{rem:generalized rank is finite for pfd modules}), and thus $\rk^\Int_M(I_a)=1$.        
    \end{claimproof}

    \begin{claim}\label{claim:I_a maximal}
        For $a\in \Z$ and $J\in \Con(\catP)$ such that $J\supsetneq I_a$, $\rk_M(J) = 0$.
    \end{claim}
    \begin{claimproof}
        Let $J\in \Con(\catP)$ such that $J\supsetneq I_a$ for some fixed $a\in \Z$.
        If $J$ contains a point $(x,y)\in \Z^2$ with $y>-x+1$, then $M_{(x,y)} = 0$ and thus $\rk_M(J)=0$.
        Hence, suppose $J$ contains  
        no point $(x,y)\in \Z^2$ with $y>-x+1$. Then, since $J$ properly contains $I_a$, there exists $(x,y)\in J\setminus I_a$ with $y=-x+1$.
        Then, by the definition of $I_a$, at least one of $(x+1,y-1)$ or $(x-1,y+1)$ must be in $I_a$ and thus in $J$ as well.
        This ensures that the diagram $M\vert_J$ contains the following diagram as a part of it: \[\kf \xleftarrow{\begin{bmatrix} 0 & 1\end{bmatrix}} \kf^2 \xrightarrow{\begin{bmatrix}1 & 0\end{bmatrix}} \kf.\]
        The rank over this diagram is zero, and thus by monotonicity of the GRI (Remark \ref{rem:properties of rank invariant} \ref{item:monotonicity}), we have $\rk_M(J)=0$.
    \end{claimproof}

Now, let $d\geq 3$, and we construct a $\Z^d$-module that is tame while its Int-GRI is not M\"obius invertible. 
Consider the projection $\pi_{\leq 2}^d:\Z^d\to \Z^2$ given by $(x_1,\ldots,x_d)\mapsto (x_1,x_2)$. Recall the map $\pi:\Z^2\to \catQ$ given in Figure~\ref{fig:tame does not imply Int-GRI invertible module}, which provides a finite encoding of the $\Z^2$-module $M$ (see Definition~\ref{def:finite encoding}). Consider the $\Z^d$-module 
$(\pi\circ \pi_{\leq 2}^d)^\ast N= N \circ \pi \circ \pi_{\leq 2}^d = M\circ \pi_{\leq 2}^d$
, which is tame by construction.
We claim that its Int-GRI is not M\"obius invertible.
{By arguments similar to those above, the intervals $I_a\times\Z^{d-2}$, where $a\in\Z$, form an infinite family of support-maximal intervals in $\Z^d$, all containing the singleton interval $\{(-1,-1,0,\ldots,0)\}$. Hence, by Corollary~\ref{cor:intersection}, the Int-GRI of $(\pi\circ\pi_{\leq2}^d)^\ast N$ is not M\"obius invertible.} %

It remains to show that for each $d\geq 2$, there exists a tame $\R^d$-module whose Int-GRI is not M\"obius invertible.
Fix $d\geq 2$, and let $\lfloor - \rfloor_{\Z^d}:\R^d\to \Z^d$ be the floor function $(x_1,\ldots,x_d)\mapsto (\lfloor x_1\rfloor, \ldots, \lfloor x_d\rfloor)$.
Define the $\R^d$-module {$M':=(\pi\circ\pi_{\leq 2}^d\circ \lfloor -\rfloor_{\Z^d})^\ast N$}, which is tame by construction.
We claim that the Int-GRI of $M'$ is not M\"obius invertible.
Recall the interval $I_a\in \Int(\Z^2)$ defined above for each $a\in \Z$ and let $I'_a$ be the preimage of $I_a$ under the map $\pi_{\leq 2}^d\circ\lfloor-\rfloor_{\Z^d}:\R^d\to\Z^2$.
By arguments similar  to those in the proofs of Claims \ref{claim:I_a nonzero rank} and \ref{claim:I_a maximal}, these $I'_a$ have $\rk_{M'}(I'_a)=1$, and for any $a\in \Z$ and $J\in \Int(\R^d)$ with $J\supsetneq I'_a$, we have $\rk_{M'}(J)=0$.
All $I'_a$ contain the singleton interval $\{(-1,-1,0,\ldots,0)\}$; hence, by Corollary \ref{cor:intersection}, the Int-GRI of $M'$ is not M\"obius invertible. 
\end{proof}

\begin{remark}\label{rem:implication of the counterexample} 
For the $\Z^2$-module $M$ considered in the previous proof, we observe the following.
\begin{enumerate}[label=(\roman*)]
\item {In Claims~\ref{claim:I_a nonzero rank} and \ref{claim:I_a maximal} above, we actually showed that, for any $\Ical\subset \Con(\catP)$ containing $\Int(\catP)$ and $\rk_M^\Ical$, $\Ical$ contains an infinite subcollection of support-maximal intervals that contain a common element of $\Ical$, namely the singleton interval $\{(-1,-1)\}$. Therefore, by Corollary \ref{cor:intersection}, for any such $\Ical$, $\rk_M^\Ical$ is not M\"obius invertible.}
\label{item:implication of the counterexample1} 

\item The finite encoding $\pi:\Z^2\to Q$ in the preceding proof is \emph{not} `rank-preserving', e.g., $\rk_M(I_a)=1\neq 0=\rk_N(\pi(I_a))$ for each $a\in \Z$. {This observation is relevant to Theorem~\ref{thm:sufficient conditions for invertibility} \ref{item:finitely encodable} in the next section.} \label{item:not rank-preserving}
\end{enumerate}
\end{remark}

\begin{proof}[Proof of Theorem \ref{thm:tame does not imply Int-GRI invertible} \ref{item:fp not imply MI}]

We first construct a $\Z^3$-module that is finitely presentable while its Int-GRI is not M\"obius invertible.
Let $N$ be the persistence module defined over the subposet $[0,3] \times [0,3] \times \{0\} \subset \Z^3$ given by:
\[
N:=\begin{tikzcd}
0 \arrow[r] & 0 \arrow[r] & 0 \arrow[r] & 0\\
\kf \arrow[u]\arrow[r,"1"] & \kf \arrow[u]\arrow[r,"1"] & \kf \arrow[u]\arrow[r] & 0\arrow[u] \\
\kf \arrow[u,"1"]\arrow[r,"{\footnotesize \begin{bmatrix}0\\ 1\end{bmatrix}}"] & \kf^2 \arrow[r,"\begin{bmatrix} 0 \ 1\end{bmatrix}"']\arrow[u,"\begin{bmatrix}0 \ 1\end{bmatrix}"'] & \kf\arrow[u,"1"]\arrow[r] & 0 \arrow[u]\\
0\arrow[r]\arrow[u] & \kf\arrow[u,"{\footnotesize\begin{bmatrix}1\\ 1\end{bmatrix}}"]\arrow[r,"1"] & \kf\arrow[u,"1"]\arrow[r] & 0 \arrow[u]
\end{tikzcd}\]

Then let $M$ be the $\Z^3$-module defined as the left Kan extension of $N$ along the canonical inclusion $([0,3]\times [0,3]\times \{0\}) \hookrightarrow \Z^3$.
By Remark \ref{rem:minimal grid} \ref{item:characterization of finitely presentable modules}, $M$ is finitely presentable.
Intuitively, $M$ consists of a ``copy" of $N$ at each $z$-axis value at or above the plane $z=0$, {with $\varphi_M(p,p+(0,0,1))$ being the identity map for $p\geq (0,0,0)$.}
Note the support of $M$ is $([1,2]\times[0,2]\times 0^\uparrow_{\Z}) \cup (\{0\}\times[1,2] \times 0^\uparrow_{\Z})$.

Now, for each $a\in 0^\uparrow_\Z$, define the interval $I_a$ of $\Z^3$ as 
\begin{equation}\label{eq:I_a}
I_a:=\left(\{(0,1)\}\times (a+1)_\Z^\uparrow\right) \cup \left(\{(1,0),(2,0),(0,2)\}\times 0_\Z^\uparrow\right) \cup \left([1,2]\times [1,2] \times [0,a]\right)
\end{equation}

Note that 
$\{(2,2,0)\}\subset \bigcap_{a\in 0_\Z^\uparrow} I_a$. Hence, by Corollary \ref{cor:intersection},
it suffices to show that $\{I_a\}_{a\in 0_\Z^\uparrow}$ is a collection of support-maximal intervals. This follows from the following two claims. 

\begin{claim}\label{claim:z3 ia support}
    For each $a\in 0_\Z^\uparrow$, $\rk_M(I_a)=1$.
\end{claim}
\begin{claimproof}
    It is immediate that $\rk_M(I_a)\leq 1$ as $\rk_M(I_a)\leq \dim(M_{(1,0,0)}) = 1$.
    Let $\ell:=(\ell_p)_{p\in I_a}\in \varprojlim M\vert_{I_a}$ be given by:
    \[ \ell_p := \begin{cases}
    1 & \mbox{if } M_p=\kf\\
    (1,1) & \mbox{if } M_p=\kf^2.
\end{cases}\]
Notably, the morphism:
\begin{equation*}
   \varphi_N((0,1),(1,1)):\ \kf \xrightarrow{\begin{bmatrix} 0 \\ 1 \end{bmatrix}} \kf^2
\end{equation*}
does not appear in $M\vert_{I_a}$ as there is no point in $\{(0,1)\}\times (a+1)_\Z^\uparrow$ below any point in $\{(1,1)\}\times [0,a]$.
Thus, via an argument similar to that in the proof of Claim~\ref{claim:I_a nonzero rank} from Theorem \ref{thm:tame does not imply Int-GRI invertible} \ref{item:tame not imply MI}, one can verify that the image of $\ell$ under the limit-to-colimit map over $I_a$ is nonzero.
Hence, invoking that $0<\rk_M(I_a)\leq 1$, we have that $\rk_M(I_a)=1$.
\end{claimproof}

\begin{claim}\label{claim:z3 ia maximal}
    Fix $a\in 0_\Z^\uparrow$ and let $J\in \Int(\Z^3)$ such that $J\supsetneq I_a$.
    Then $\rk_M(J) = 0$.
\end{claim}
\begin{claimproof}
If $p\in J$ such that $M_p=0$, then the conclusion follows, so assume $M_p\neq 0$. Then, either (i) $p\in [\{(0,1)\}]\times [0,a]$, or (ii) $p\in [1,2]\times [1,2]\times (a+1)^\uparrow$.

In Case (i), by convexity of $J$ and the fact $J\supset I_a$, we  have $p$, $p+(1, 0, 0)$, $p + (1, -1, 0)\in J$.
Thus, 
the diagram 
\begin{equation}\label{eq:sub-diagram}
k\xrightarrow{\begin{bmatrix} 0\\1\end{bmatrix}}\kf^2\xleftarrow{\begin{bmatrix}1\\1\end{bmatrix}}\kf
\end{equation}
appears in $M\vert_J$, and thus $\varprojlim M\vert_J = 0$, implying $\rk_M(J) = 0$.
In Case (ii), $p$ is equal to $(x, y, z)$ with $x,y\in \{1,2\}$ and $z\geq a+1$.
By convexity of $J$ and the fact $J\supset I_a$, this implies that $(0, 1, z), (1, 0, z), (1,1, z)\in J$.
As a result, Diagram \eqref{eq:sub-diagram} appears as a part of the diagram $M\vert_J$, and thus $\varprojlim M\vert_J = 0$, implying $\rk_M(J) = 0$.
\end{claimproof}

Finally, let $d>3$. 
Then, instead of viewing $N$ as a persistence module defined over $[0,3]\times [0,3]\times \{0\}\subset \Z^3$, view $N$ in the obvious way as a persistence module defined over $[0,3]\times [0,3]\times \{0\}^{d-2}\subset \Z^d$.
As before, let $M$ be the finitely presentable $\Z^d$-module defined as the left Kan extension along the canonical inclusion of $[0,3]\times [0,3]\times \{0\}^{d-2}\hookrightarrow \Z^d$.
By replacing the interval $I_a$ in Equation \eqref{eq:I_a} by the interval $I_a\times (0_\Z^\uparrow)^{d-3}$, the direct analogs of Claims \ref{claim:z3 ia support} and \ref{claim:z3 ia maximal} hold, giving us that the Int-GRI of $M$ is not M\"obius invertible.

It remains to show that for each $d\geq 3$, there exists a finitely presentable $\R^d$-modules whose Int-GRI is not M\"obius invertible.
Fix $d\geq 3$, and as before let $\lfloor - \rfloor_{\Z^d}:\R^d\to \Z^d$ be the floor function.
Let $M$ be the $\Z^d$-module with non-M\"obius invertible Int-GRI constructed above, and define the $\R^d$-module $M':=\lfloor - \rfloor_{\Z^d}^\ast M$. %
Equivalently, viewing $N$ as a module supported on $\{0,1,2,3\}\times \{0,1,2,3\} \times \{0\}^{d-2}\subset \R^d$, $M'$ is the left Kan extension of $N$ along the canonical inclusion $\{0,1,2,3\}\times \{0,1,2,3\} \times \{0\}^{d-2}\hookrightarrow \R^d$, so by Remark \ref{rem:minimal grid} \ref{item:characterization of finitely presentable modules}, $M'$ is finitely presentable.

For $a\in 0_\Z^\uparrow$, define $I_a':=\lfloor - \rfloor_{\Z^d}\inv(I_a\times (0_\R^\uparrow)^{d-3})$.
Then, $I_a'\in \Int(\R^d)$ and the statements in Claims \ref{claim:z3 ia support} and \ref{claim:z3 ia maximal} hold when replacing $I_a$ with $I_a'$ and $\Z^3$ with $\R^d$, via directly analogous arguments.
Lastly, note that $\{(2,2,\{0\}^{d-2})\}\subset \bigcap_{a\in 0_\Z^\uparrow} I'_a$, and thus by Corollary \ref{cor:intersection}, the Int-GRI of $M'$ is not M\"obius invertible.
\end{proof}

\begin{example}\label{ex:another R3 module} {We provide another finitely presentable $\R^3$-module whose Int-GRI is not M\"obius invertible; this example is easy to visualize and will be useful later for several purposes.}
For the $\catQ$-module $N$ depicted  in Figure~\ref{fig:M fp not MI}, let $M$ be the left Kan extension of the $Q$-module along the inclusion $\iota:Q\hookrightarrow \R^3$, as depicted  Figure~\ref{fig:ia for fp not imply mi} (left). 
Then, by Remark~\ref{rem:minimal grid}~\ref{item:characterization of finitely presentable modules}, $M$ is finitely presentable.
We also observe that if $p\in \R^3$ does not belong to the following set, then $M_p=0$:  \[[(0,0,0),(3,1,3)):=\{x\in \R^3:(0,0,0) \leq x < (3,1,3)\}\in \Int(\R^3).\]  This set is partitioned into the following nine unit blocks as depicted Figure~\ref{fig:ia for fp not imply mi} (left): 
\begin{equation}\label{eq:unit blocks}
B_{i,j}:=[(i,0,j),(i,0,j)+(1,1,1)), \ \ 0\leq i,j\leq 2.
\end{equation}
Now, for each $a\in (0,1)$, define $I_a\in \Int(\R^3)$ as follows; see the right-hand side of Figure~\ref{fig:ia for fp not imply mi}:
\[I_a:=[(1,0,0), (3,a,3)) \cup [(0,0,2),(1,1,3))\cup [(0,a,1),(1,1,2))\cup [(1,a,0), (3,1,1)).\]
In the following two claims, we show that $\{I_a\}_{a\in (0,1)}$ is a collection of support-maximal intervals with $\rk_M^\Int(I_a)=1$ for all $a\in (0,1)$. Observe that the singleton $I=\{(2,0,2)\}$ is contained in $I_a$ for all $a\in (0,1)$. Thus, by Corollary \ref{cor:intersection}, the Int-GRI of $M$ is not M\"obius invertible.

\setcounter{claim}{0}
\resetclaimcounter
\begin{claim}\label{claim:I_a nonzero rank 2}
For each $a\in (0,1)$, we have that $\rk_M^\Int(I_a)=1$.    
\end{claim}
\begin{claimproof}
       Fix any $a\in (0,1)$. It is immediate that $\rk_M(I_a)\leq 1$ as $\rk_M(I_a)\leq \dim(M_{(0,0,2)})=\dim(\kf)=1$. Observe that the map $1\in \kf\mapsto (0,1)\in \kf^2$ stemming from blocks $B_{0, 1}$ and $B_{1,1}$ does not appear in the diagram $M|_{I_a}$. This observation leads to another observation that $\ell:=(\ell_p)_{p\in I_a}\in \varprojlim M\vert_{I_a}$ (cf. Proposition~\ref{prop:limit and colimit}~\ref{item:limit}), where
\[ \ell_p := \begin{cases}
    1 & \mbox{if } M_p=\kf\\
    (1,1) & \mbox{if } M_p=\kf^2.
\end{cases}\]
Via an argument similar to that in the proof of Claim~\ref{claim:I_a nonzero rank} from Theorem \ref{thm:tame does not imply Int-GRI invertible}  \ref{item:tame not imply MI}, one can see that the image of $\ell$ under the limit-to-colimit map over $I_a$ is nonzero.
Hence, we have $0<\rk_M(I_a)\leq 1$, i.e., $\rk_M(I_a) = 1$. 
\end{claimproof}

\begin{claim}\label{claim:greater interval has zero rank}
    Fix $a\in (0,1)$, and let $J\supsetneq I_a$ be an interval of $\R^3$.
    Then $\rk_M(J) = 0$.
\end{claim}
\begin{claimproof}
    Let $p\in J\setminus I_a$.
    If $M_p = 0$, then the claim is immediate, and thus we assume that $p\in B_{i,j}$ for some $0\leq i,j\leq 2$, with $i,j$ not both $0$. 
    Then, since $B_{0,2}\cup B_{1,0}\cup B_{2,0}\subset I_a$, we have that $p\in B_{i,j}$ for some $(i,j)\in \{(0,1),(1,1),(1,2),(2,1),(2,2)\}$.

    Assume that $p=(x, y, z)\in B_{0,1}\setminus I_a$. Then,
    $(x+1,y,z)\in I_a\cap B_{1,1}$ and $(x+1,y,z-1)\in I_a\cap B_{1,0}$. Hence, over the subposet $\{p<(x+1,y,z)>(x+1,y,z-1)\}\subset J$, Diagram \eqref{eq:sub-diagram}
    appears as a part of the diagram $M\vert_J$. Diagram \eqref{eq:sub-diagram} does not admit any nontrivial section, and thus the limit of $M$ over $J$ is trivial, implying that $\rk_M(J) = 0$.

    Next, assume that $p=(x,y,z)\in B_{i,j}\setminus I_a$ for $1\leq i,j\leq 2$, i.e., $1\leq x < 3$, $a\leq y<1$, and $1\leq z<3$. Observe that  $(0,y,1)\in I_a\cap B_{0,1}\subset J$ and that $(0,y,1)\leq (1,y,1)\leq p\in J$. Hence, by convexity of $J$, we have $(1,y,1)\in J$.
   Therefore, over the subposet \[\{(0,y,1)<(1,y,1)>(1,y,0)\}\subset J,\] we have Diagram \eqref{eq:sub-diagram} again, and hence $\rk_M(J) = 0$, as desired.    
\end{claimproof}
\end{example}

\tdplotsetmaincoords{60}{125}
\tdplotsetrotatedcoords{0}{30}{120} 
\begin{figure}
\begin{center}
\begin{subfigure}[c]{0.4\textwidth}
\centering

\begin{tikzpicture}[scale=1.8,tdplot_rotated_coords,
                    edge/.style={black,thick},
                    edgeDash/.style={black,thick,dashed},
                    centernode/.style={circle,inner sep=1pt,draw=none,fill=none}]
\def\nodegap{0.2}

\tikzset{
  edge/.style={->, >=stealth, line width=0.6pt, shorten >=0.3pt, shorten <=0.3pt, >=stealth, scale=1},
  edgeDash/.style={->, >=stealth, dashed, line width=0.6pt, shorten >=0.3pt, shorten <=0.3pt}
}

\foreach \x in {0,1,2,3}
  \foreach \y in {0,1,2,3}
    \foreach \z in {0,1}{
      \ifthenelse{\x<3}{
        \ifthenelse{\z=1 \OR \y=3}{
          \draw[edgeDash] (\x+\nodegap,\y,\z) -- (\x+1-\nodegap,\y,\z);
        }{
          \ifthenelse{\z=0 \AND \x<2}{
  \draw[edge, red] (\x+\nodegap,\y,\z) -- (\x+1-\nodegap,\y,\z);
}{
  \draw[edge] (\x+\nodegap,\y,\z) -- (\x+1-\nodegap,\y,\z);
}
        }
      }{}
      \ifthenelse{\y<3}{
        \ifthenelse{\z=1 \OR \x=3}{
          \draw[edgeDash] (\x,\y+\nodegap,\z) -- (\x,\y+1-\nodegap,\z);
        }{
          \ifthenelse{\z=0 \AND \y<2}{
  \draw[edge, red] (\x,\y+\nodegap,\z) -- (\x,\y+1-\nodegap,\z);
}{
  \draw[edge] (\x,\y+\nodegap,\z) -- (\x,\y+1-\nodegap,\z);
}
        }
      }{}
      \ifthenelse{\z<1}{
        \ifthenelse{\x=3 \OR \y=3}{
          \draw[edgeDash] (\x,\y,\z+\nodegap) -- (\x,\y,\z+1-\nodegap);
        }{
          \draw[edge] (\x,\y,\z+\nodegap) -- (\x,\y,\z+1-\nodegap);
        }
      }{}
}
\foreach \x in {0,1,2,3}
   \foreach \y in {0,1,2,3}
      \foreach \z in {0,1}{
         \node[centernode] (C-\x-\y-\z) at (\x,\y,\z) {};
      }          

\node[centernode] at (0,0,0) {0};
\node[centernode] at (1,0,0) {$\kf$};
\node[centernode] at (2,0,0) {$\kf$};
\node[centernode] at (0,1,0) {$\kf$};
\node[centernode] at (1,1,0) {$\kf^2$};
\node[centernode] at (2,1,0) {$\kf$};
\node[centernode] at (0,2,0) {$\kf$};
\node[centernode] at (1,2,0) {$\kf$};
\node[centernode] at (2,2,0) {$\kf$};
\node[centernode] at (0,0,1) {0};
\node[centernode] at (1,0,1) {0};
\node[centernode] at (2,0,1) {0};
\node[centernode] at (0,1,1) {0};
\node[centernode] at (1,1,1) {0};
\node[centernode] at (2,1,1) {0};
\node[centernode] at (0,2,1) {0};
\node[centernode] at (1,2,1) {0};
\node[centernode] at (2,2,1) {0};
\node[centernode] at (3,0,0) {0};
\node[centernode] at (3,1,0) {0};
\node[centernode] at (3,2,0) {0};
\node[centernode] at (3,3,0) {0};
\node[centernode] at (3,0,1) {0};
\node[centernode] at (3,1,1) {0};
\node[centernode] at (3,2,1) {0};
\node[centernode] at (3,3,1) {0};
\node[centernode] at (0,3,0) {0};
\node[centernode] at (0,3,1) {0};
\node[centernode] at (1,3,0) {0};
\node[centernode] at (1,3,1) {0};
\node[centernode] at (2,3,0) {0};
\node[centernode] at (2,3,1) {0};

\node[centernode] at (-0.6,-0.1,0.1) {$(0,0,0)$};
\node[centernode] at (3.6,3,1) {$(3,1,3)$};

\begin{scope}[shift={(4,-0.5,0)}]
  \draw[->, thick, black] (0,0,0) -- (0.6,0,0) node[below right=-1pt] {\small $x$};
  \draw[->, thick, black] (0,0,0) -- (0,0.6,0) node[left=-1pt] {\small $z$};
  \draw[->, thick, black] (0,0,0) -- (0,0,0.6) node[above left=-1pt] {\small $y$};
  \end{scope}
\end{tikzpicture}
\end{subfigure}
\hspace{20mm}
\begin{subfigure}[c]{0.38\textwidth}
\centering
\begin{tikzpicture}[scale=1.8,tdplot_rotated_coords,
                    centernode/.style={circle,inner sep=1pt,draw=none,fill=none}]
\def\nodegap{0.2}

\tikzset{
  >={Stealth[length=4pt,width=4pt]},
  edge/.style={->, red, line width=0.6pt, shorten >=0.3pt, shorten <=0.3pt},
  lab/.style={font=\scriptsize, text=black, fill=white, inner sep=0.5pt}
}

\foreach \x in {0,1,2}
  \foreach \y in {0,1,2}{
    \node[centernode] (C-\x-\y-0) at (\x,\y,0) {};
  }

\node[centernode] at (0,0,0) {0};
\node[centernode] at (1,0,0) {$\kf$};
\node[centernode] at (2,0,0) {$\kf$};
\node[centernode] at (0,1,0) {$\kf$};
\node[centernode] at (1,1,0) {$\kf^2$};
\node[centernode] at (2,1,0) {$\kf$};
\node[centernode] at (0,2,0) {$\kf$};
\node[centernode] at (1,2,0) {$\kf$};
\node[centernode] at (2,2,0) {$\kf$};

\draw[edge] (0+\nodegap,2,0) -- node[midway, lab, above=2pt] {$1$} (1-\nodegap,2,0);
\draw[edge] (1+\nodegap,2,0) -- node[midway, lab, above=2pt] {$1$} (2-\nodegap,2,0);

\draw[edge] (0+\nodegap,1,0) -- node[midway, lab, above=2pt] {$\begin{bmatrix}0\\[1pt]1\end{bmatrix}$} (1-\nodegap,1,0);
\draw[edge] (1+\nodegap,1,0) -- node[midway, lab, above=2pt] {$\begin{bmatrix}0\ 1\end{bmatrix}$} (2-\nodegap,1,0);

\draw[edge] (0+\nodegap,0,0) -- (1-\nodegap,0,0);
\draw[edge] (1+\nodegap,0,0) -- node[midway, lab, below=2pt] {$1$} (2-\nodegap,0,0);

\draw[edge] (0,0+\nodegap,0) -- (0,1-\nodegap,0);
\draw[edge] (0,1+\nodegap,0) -- node[midway, lab, left=3pt] {$1$} (0,2-\nodegap,0);

\draw[edge] (1,0+\nodegap,0) -- node[midway, lab, right=3pt] {$\begin{bmatrix}1\\[1pt]1\end{bmatrix}$} (1,1-\nodegap,0);
\draw[edge] (1,1+\nodegap,0) -- node[midway, lab, right=3pt, yshift=1pt] {$\begin{bmatrix}0\ 1\end{bmatrix}$} (1,2-\nodegap,0);

\draw[edge] (2,0+\nodegap,0) -- node[midway, lab, right=3pt] {$1$} (2,1-\nodegap,0);
\draw[edge] (2,1+\nodegap,0) -- node[midway, lab, right=3pt] {$1$} (2,2-\nodegap,0);
\end{tikzpicture}

\end{subfigure}
\end{center}
\vspace{-8mm}
\caption{A \(Q\)-module $N$ for \(Q = \{0,1,2,3\} \times \{0,1\} \times \{0,1,2,3\}\), where the maps corresponding to the red arrows (within the plane $y=0$) are shown on the right.
}
\label{fig:M fp not MI}

\begin{center}
\begin{subfigure}[c]{0.5\textwidth}
 $   
\begin{tikzpicture}[scale=1.8,tdplot_rotated_coords,
                    edge/.style={black,thick},
                    edgeDash/.style={black,thick,dashed},
                    centernode/.style={circle,inner sep=1pt,draw=none,fill=none}]
\tikzset{>={Stealth[length=4pt,width=4pt]}} %
\foreach \x in {0,1,2,3}
   \foreach \y in {0,1,2,3}
      \foreach \z in {0,1}{
           \ifthenelse{\x<3}{
              \ifthenelse{\z=1 \OR \y=3}{ %
                 \draw[edgeDash] (\x,\y,\z) -- (\x+1,\y,\z);
              }{
                 \draw[edge] (\x,\y,\z) -- (\x+1,\y,\z);
              }
           }{}
           \ifthenelse{\y<3}{
              \ifthenelse{\z=1 \OR \x=3}{ %
                 \draw[edgeDash] (\x,\y,\z) -- (\x,\y+1,\z);
              }{
                 \draw[edge] (\x,\y,\z) -- (\x,\y+1,\z);
              }
           }{}
           \ifthenelse{\z<1}{
              \ifthenelse{\x=3 \OR \y=3}{ %
                 \draw[edgeDash] (\x,\y,\z) -- (\x,\y,\z+1);
              }{
                 \draw[edge] (\x,\y,\z) -- (\x,\y,\z+1);
              }
           }{}
      }

\foreach \x in {0,1,2}
   \foreach \y in {0,1,2}
      \foreach \z in {0}{
         \node[centernode] (C-\x-\y-\z) at (\x+0.5,\y+0.5,\z+0.5) {}; %
      }    

\foreach \x in {0,1,2}
   \foreach \y in {0,1,2}
      \foreach \z in {0}{
         \ifthenelse{\x<2}{
            \draw[->,red,thick,shorten >=10pt,shorten <=10pt] 
               (C-\x-\y-\z) -- (C-\the\numexpr\x+1\relax-\y-\z);
         }{}
         \ifthenelse{\y<2}{
            \draw[->,red,thick,shorten >=10pt,shorten <=10pt] 
               (C-\x-\y-\z) -- (C-\x-\the\numexpr\y+1\relax-\z);
         }{}
      }

\node[centernode] at (0.5,0.5,0.5) {0};
\node[centernode] at (1.5,0.5,0.5) {$\kf$};
\node[centernode] at (2.5,0.5,0.5) {$\kf$};
\node[centernode] at (0.5,1.5,0.5) {$\kf$};
\node[centernode] at (1.5,1.5,0.5) {$\kf^2$};
\node[centernode] at (2.5,1.5,0.5) {$\kf$};
\node[centernode] at (0.5,2.5,0.5) {$\kf$};
\node[centernode] at (1.5,2.5,0.5) {$\kf$};
\node[centernode] at (2.5,2.5,0.5) {$\kf$};

\tikzset{
  blockwm/.style={
    text=gray,
    font=\bfseries\footnotesize,
  }
}
\def\wmx{-9pt}  %
\def\wmy{-15pt}  %

\node[blockwm, xshift=\wmx, yshift=\wmy] at (0.5,0.5,0.5) {$B_{0,0}$};
\node[blockwm, xshift=\wmx, yshift=\wmy] at (1.5,0.5,0.5) {$B_{1,0}$};
\node[blockwm, xshift=\wmx, yshift=\wmy] at (2.5,0.5,0.5) {$B_{2,0}$};

\node[blockwm, xshift=\wmx, yshift=\wmy] at (0.5,1.5,0.5) {$B_{0,1}$};
\node[blockwm, xshift=\wmx, yshift=\wmy] at (1.5,1.5,0.5) {$B_{1,1}$};
\node[blockwm, xshift=\wmx, yshift=\wmy] at (2.5,1.5,0.5) {$B_{2,1}$};

\node[blockwm, xshift=\wmx, yshift=\wmy] at (0.5,2.5,0.5) {$B_{0,2}$};
\node[blockwm, xshift=\wmx, yshift=\wmy] at (1.5,2.5,0.5) {$B_{1,2}$};
\node[blockwm, xshift=\wmx, yshift=\wmy] at (2.5,2.5,0.5) {$B_{2,2}$};

\node[centernode] at (-0.4,-0.1,0.3) {$(0,0,0)$};
\node[centernode] at (3.4, 3.1, 0.7) {$(3,1,3)$};
\end{tikzpicture}$
\end{subfigure}
\hfill
    \begin{subfigure}[c]{0.49\textwidth}
\tdplotsetmaincoords{60}{125}
\tdplotsetrotatedcoords{0}{30}{120} 
\begin{tikzpicture}[scale=1.8,tdplot_rotated_coords,
                    edge/.style={black,thick},
                    edgeDash/.style={black,thick,dashed},
                    face/.style={fill=gray!30,opacity=0.6},centernode/.style={circle,inner sep=1pt,draw=none,fill=none}]
\fill[face] (1,0,0) -- (2,0,0) -- (3,0,0) -- (3,0,1) -- (2,0,1) -- (1,0,1) -- cycle;

\fill[face] (1,0,0) -- (1,1,0) -- (1,1,1) -- (1,0,1) -- cycle;

\fill[face] (3,0,0) -- (3,1,0) -- (3,1,1) -- (3,0,1) -- cycle;

\fill[face] (1,0,0) -- (3,0,0) -- (3,1,0)--(1,1,0)--cycle;

\fill[face] (1,1,0) -- (3,1,0) -- (3,1,0.5) -- (1,1,0.5) -- cycle;

\fill[face] (1,1,0.5) -- (3,1,0.5) -- (3,1,1) -- (1,1,1) -- cycle;

\fill[face] (1,1,0)--(3,1,0)--(3,3,0)--(1,3,0)--cycle;
\fill[face] (1,1,0)--(1,2,0)--(1,2,0.5)--(1,1,0.5)--cycle;
\fill[face] (0,1,0.5)--(1,1,0.5)--(1,2,0.5)--(0,2,0.5)--cycle;
\fill[face] (0,2,0)--(1,2,0)--(1,2,0.5)--(0,2,0.5)--cycle;
\fill[face] (0,2,0)--(1,2,0)--(1,3,0)--(0,3,0)--cycle;

\fill[face] (0,3,0)--(1,3,0)--(1,3,1)--(0,3,1)--cycle;
\fill[face] (0,1,0.5)--(1,1,0.5)--(1,1,1)--(0,1,1)--cycle;
\fill[face] (0,1,1)--(0,3,1)--(1,3,1)--(1,1,1)--cycle;
\fill[face] (1,1,0.5)--(1,1,1)--(1,3,1)--(1,3,0.5)--cycle;
\fill[face] (3,1,1)--(1,1,1)--(1,0,1)--(3,0,1)--cycle;
\fill[face] (1,1,0.5)--(3,1,0.5)--(3,3,0.5)--(1,3,0.5)--cycle;
\fill[face] (3,1,0)--(3,3,0)--(3,3,0.5)--(3,1,0.5)--cycle;

\fill[face] (0,3,0) -- (3,3,0) -- (3,3,0.5) -- (1,3,0.5) -- (1,3,1) -- (0,3,1) -- cycle;

\fill[face] (0,1,0.5) -- (0,2,0.5) -- (0,2,0)-- (0,3,0) -- (0,3,1) -- (0,1,1) -- cycle;

\draw[edge] (1, 0, 0) -- (2, 0, 0);
\draw[edge] (1, 0, 0) -- (1, 0, 1);
\draw[edge] (1, 0, 0) -- (1, 1, 0);
\draw[edgeDash] (1, 0, 1) -- (2, 0, 1);
\draw[edgeDash] (1, 0, 1) -- (1, 1, 1);
\draw[edge] (2, 0, 0) -- (3, 0, 0);
\draw[edgeDash] (2, 0, 1) -- (3, 0, 1);
\draw[edgeDash] (3, 0, 0) -- (3, 0, 1);
\draw[edge] (1, 1, 0) -- (1, 1, 0.5);
\draw[edgeDash] (1, 1, 0.5) -- (1, 1, 1);
\draw[edgeDash] (1, 1, 1) -- (3, 1, 1);
\draw[edgeDash] (3, 0, 0) -- (3, 3, 0);
\draw[edgeDash] (3, 0, 1) -- (3, 1, 1);
\draw[edge] (0, 1, 0.5) -- (0, 1, 1);
\draw[edge] (0, 1, 0.5) -- (1, 1, 0.5);
\draw[edgeDash] (1, 1, 0.5) -- (3, 1, 0.5);
\draw[edgeDash] (3, 1, 0) -- (3, 1, 1);
\draw[edgeDash] (0, 1, 1) -- (1, 1, 1);
\draw[edge] (1, 1, 0) -- (3, 1, 0);
\draw[edgeDash] (3, 1, 0.5) -- (3, 3, 0.5);
\draw[edgeDash] (1, 1, 0.5) -- (1, 3, 0.5);
\draw[edgeDash] (1, 3, 0.5) -- (3, 3, 0.5);
\draw[edge] (1, 1, 0) -- (1, 3, 0);
\draw[edgeDash] (1, 3, 0) -- (1, 3, 0.5);
\draw[edgeDash] (3, 3, 0) -- (3, 3, 0.5);
\draw[edge] (0, 1, 0.5) -- (0, 2, 0.5);
\draw[edge] (0, 2, 0.5) -- (1, 2, 0.5);
\draw[edgeDash] (0, 1, 1) -- (0, 3, 1);
\draw[edgeDash] (0, 3, 1) -- (1, 3, 1);
\draw[edgeDash] (1, 3, 1) -- (1, 3, 0.5);
\draw[edge] (0, 2, 0.5) -- (0, 2, 0);
\draw[edge] (0, 2, 0) -- (0, 3, 0);
\draw[edgeDash] (0, 3, 0) -- (0, 3, 1);
\draw[edge] (0, 2, 0) -- (1, 2, 0);
\draw[edgeDash] (0, 3, 0) -- (3, 3, 0);
\draw[edge] (1, 2, 0) -- (1, 2, 0.5);
\draw[edgeDash] (1,1,1) -- (1,3,1);

\node[shape=circle, inner sep=2pt, fill=black] at (0,0,0) {};
\node[centernode] at (0,-0.25,0) {$(0,0,0)$};
\draw[edgeDash] (0,0,0) -- (1,0,0);

\node[shape=circle, inner sep=2pt, fill=black] at (3,3,0.5) {};
\node[centernode] at (3.15,3,1) {$(3,a,3)$};

\node[shape=circle, inner sep=2pt, fill=black] at (1,3,0.5) {};
\node[centernode] at (1.4,3.3,0.5) {$(1,a,3)$};

\node[shape=circle, inner sep=2pt, fill=black] at (0,1,0.5) {};
\node[centernode] at (-0.5,1.1,0) {$(0,a,1)$};

\node[shape=circle, inner sep=2pt, fill=black] at (0,3,1) {};
\node[centernode] at (-0.5,3,1) {$(0,1,3)$};

\node[shape=circle, inner sep=2pt, fill=black] at (1,2,0) {};
\node[centernode] at (1.5,2,0) {$(1,0,2)$};

\node[shape=circle, inner sep=2pt, fill=black] at (1,0,0) {};
\node[centernode] at (1,0,0.-0.5) {$(1,0,0)$};

\end{tikzpicture}
    \end{subfigure}
    \end{center}
    \caption{(Left) The left Kan extension of $N$ depicted in Figure~\ref{fig:M fp not MI} along the inclusion $Q\hookrightarrow \R^3$, which is denoted by $M$. Nine blocks $B_{i,j}$, $0\leq i,j\leq 2$ partition the block $[(0,0,0),(3,1,3))$. For $p\not\in [(0,0,0),(3,1,3))$, we have that $M_p=0$.  (Right) The interval $I_a$ for $a=0.5$.}
    \label{fig:ia for fp not imply mi}
\end{figure}

We close this section with a final remark on an implication of Theorem~\ref{thm:tame does not imply Int-GRI invertible}.
\begin{remark}
Given a $\catP$-module $M$, if $\rk_M^\Int$ is not M\"obius invertible, then there exists \emph{no} finite exact sequence of interval-decomposable $\catP$-modules $\{M_i\}_{i=0}^n$ and morphisms  $\{f_i\}_{i=0}^n$ with
\[0\longrightarrow M_n\stackrel{f_n}{\longrightarrow}\cdots\longrightarrow  M_1\stackrel{f_1}{\longrightarrow} M_0 \stackrel{f_0}{\longrightarrow} M \stackrel{f_{-1}:=0}{\longrightarrow} 0\]
such that $\rk_{M_i}^{\Int}=\rk_{\ker(f_i)}^\Int +\rk_{\ker(f_{i-1})}^\Int$ for each $i=0,
\ldots, n$. Indeed, the existence of such a sequence implies the rank decomposition $\rk_M^\Int=\rk_{\bigoplus_{i}M_{2i}}^{\Int}-\rk_{\bigoplus_{i}M_{2i-1}}^\Int$ and thus the M\"obius invertibility of $\rk_M^{\Int}$. Theorem \ref{thm:tame does not imply Int-GRI invertible}
clarifies that there are $\Z^d$- or $\R^d$-modules that do not admit such finite sequences. This observation might have interesting implications in the perspective of recent studies at the intersection of representation theory and multi-parameter persistence \cite{amiot2024invariants,asashiba2023approximation,blanchette2021homological,botnan2024signed,chacholski2022effective}.
\end{remark}

\subsection{Sufficient conditions for M\"obius invertibility of the GRI}\label{sec:sufficient conditions}

In this section, we explore several conditions on persistence modules $M$ that ensure the M\"obius invertibility of the Con- or Int-GRI of $M$: {Theorem~\ref{thm:sufficient conditions for invertibility}, Proposition \ref{prop:tame+ does imply MI}, and Theorem \ref{thm:pull-back GPD}.} 

Let $\pi:\catP\to \catQ$ be an order-preserving map.
If for all $I\in \Int(\catP)$, the image $\pi(I)$ belongs to $\Int(\catQ)$, then we say $\pi$ is \emph{interval-preserving}.
Let $M$ be a $\catP$-module and $N$ be a $\catQ$-module. We say that $\pi$ is Int\emph{-rank-preserving} (from $M$ to $N$), if for all $I\in \Int(\catP)$,  $\rk_M(I) = \rk_N(\pi(I))$. 

\begin{theorem}\label{thm:sufficient conditions for invertibility}  Each of the following assumptions implies that {the Int-GRI of} a given $\catP$-module $M$ is \invertible{}. 
\begin{enumerate}[label=(\roman*),leftmargin=*]
  \item  $M$ is interval-decomposable. 
    \label{item:interval-decomposable}
    \item (Finite Int-GRI factorization condition) There exist a finite poset $\catQ$, an order-preserving map $\pi:\catP\rightarrow \catQ$, and a $\catQ$-module $N$ such that $\pi$ is interval-preserving and Int-rank-preserving from $M$ to $N$.
   \label{item:finitely encodable}  
   \item $\catP\in \{\Z^2,\R^2\}$ and $M$ is finitely presentable \label{item:finitely presentable invertible}. 
\end{enumerate}
In addition, under the assumption of \ref{item:interval-decomposable}, the GPD sends each interval to its multiplicity in the barcode of $M$.
\end{theorem}
We remark that the fact that the Int-GRI of a finitely presentable $\Z^2$- or $\R^2$-module $M$ is M\"obius invertible will be refined in Proposition~\ref{prop:finitely presentable invertiver over intmn}. In contrast, the Con-GRI of $M$ is not necessarily M\"obius invertible; see Remark~\ref{rem:finitely presented R^d not Con MI}.
We also remark that the finite Int-GRI factorization condition is distinct from the constructibility of 
$\rk_M^\Int$ (Definition~\ref{def:co-closure and constructible}).
For example, let $\catP$ be an infinite antichain, let $\catQ$ be the
one-point poset, and let $
\pi:\catP\to\catQ$ be the unique map. Let $M$ be the constant $\catP$-module with value $\kf$,
and let $N$ be the constant $\catQ$-module with value $\kf$. Then $\pi$ is
interval-preserving and Int-rank-preserving, so the finite Int-GRI
factorization condition holds.
On the other hand, since $\catP$ is an antichain, every interval in
$\catP$ is a singleton. Hence, $
\Int(\catP)=\bigl\{\{p\}:p\in\catP\bigr\}
$
is itself an infinite antichain. Therefore, every co-closure operator on
$\Int(\catP)$ is the identity map, and in particular has infinite image.
Thus, $\rk_M^\Int$ cannot be constructible.

\begin{corollary}\phantomsection\label{cor:Mobius invertibility does not imply tameness} There exists a poset $\catP$ such that:
\begin{enumerate}[label=(\roman*)]
    \item   there exists a non-tame $\catP$-module whose Int-GRI is M\"obius invertible.
    \item there exists a non finitely presentable $\catP$-module whose Int-GRI is M\"obius invertible.
\end{enumerate}
\end{corollary}

\begin{proof}
Let $M$ be a $\Z$-module that has infinitely many summands, e.g. $M=\bigoplus_{n\in \Z}k_{\{n\}}$. 
Then, $M$ cannot be finitely encoded and thus $M$ is not tame. Also, $M$ is not finitely presentable. However, by Theorem \ref{thm:sufficient conditions for invertibility} \ref{item:interval-decomposable}, its Int-GRI is M\"obius invertible.  
\end{proof}

\begin{remark}[Comparison between Definition~\ref{def:finite encoding} and Theorem~\ref{thm:sufficient conditions for invertibility}~\ref{item:finitely encodable}]\phantomsection\label{rem:about thm sufficient conditions for invertibility}
{A map $\pi:\catP \to \catQ$ as in Theorem~\ref{thm:sufficient conditions for invertibility}~\ref{item:finitely encodable} need not be a finite encoding, and vice versa. For the latter statement, see Remark \ref{rem:implication of the counterexample} \ref{item:not rank-preserving}. For the former, let $M$ and $N$ be non-isomorphic persistence modules over the same finite poset $\catP$ with the same Int-GRI. Then the identity map $\id_{\catP}$ is order-preserving, interval-preserving, and Int-rank-preserving, but it is not a finite encoding.}
Nevertheless, natural candidates for 
such $\pi$ {as in Theorem~\ref{thm:sufficient conditions for invertibility}~\ref{item:finitely encodable}}
are finite encodings:
see Example \ref{ex:anti-diagonal}. 
\end{remark}

Now we prove Theorem \ref{thm:sufficient conditions for invertibility}. Item \ref{item:interval-decomposable} is  straightforward: Define $\drm_M:\Int(\catP)\rightarrow \Z$ by sending each interval $I$ to the multiplicity of $I$ in the barcode of $M$. 
Then, {by Remark~\ref{rem:properties of rank invariant}~\ref{item:additivity} and \ref{item:the GRI of an interval module}}, the function $\drm_M$ 
{is the M\"obius inversion of $\rk_M^\Int$ (Definition~\ref{def:GGPD})}, completing the proof. In proving item \ref{item:finitely encodable}, the proposition below is useful. 

\begin{proposition}\label{prop:each component is an interval}
Let $\pi:\catP\rightarrow \catQ$ be an order-preserving map and let $J$ be an interval of $\catQ$. Then $\pi^{-1}(J)$ is convex (i.e., a disjoint union of intervals of $\catP$).
\end{proposition}
\begin{proof}
Let $p,r\in \pi^{-1}(J)$ and $q\in \catP$ with $q\in [p,r]$. Since $\pi$ is order-preserving,  $\pi(q)\in [\pi(p),\pi(r)]$ in $Q$. Then, as $J$ is an interval, we have that $\pi(q)\in  J$, i.e., $q\in \pi^{-1}(J)$, as desired.
\end{proof}

\begin{proof}[Proof of Theorem \ref{thm:sufficient conditions for invertibility}
\ref{item:finitely encodable}] 

Since $\catQ$ is finite, the set $\pi^{-1}\Int(\catQ):=\{\pi\inv(J):J\in \Int(\catQ)\}$ is finite. For $J\in \Int(\catQ)$, let $C(\pi\inv(J))$ be the set of connected components of $\pi\inv(J)$. {This implies that any pair of elements in $C(\pi\inv(J))$ are disjoint. Also,}
by Proposition \ref{prop:each component is an interval}, every element of $C(\pi\inv(J))$ is in $\Int(\catP)$.

 Let \[\Ical := \bigcup_{J\in \Int(\catQ)} C(\pi\inv(J))=\{I\in \Int(\catP):\mbox{$I$ is a connected component of an element in $\pi^{-1}\Int(\catQ)$}\}.\] For each $I\in \Int(\catP)$, the set
\[\Ical_{\supset I}:=\{J\in \Ical:J\supset I\}
\]
is finite, because for each $J\in\Int(\catQ)$, the intersection $\Ical_{\supset I} \cap C(\pi^{-1}(J))$ contains at most one element. This proves that, by {Proposition \ref{prop:convolvability} \ref{item:convolvability2},} $\rk^\Ical_M$ is convolvable. 

Next, for $I\in \Int(\catP)$, let $I^\Ical\in \Ical$ be
the connected component of $\pi\inv(\pi(I))$ containing $I$.
Since $\pi(I^\Ical)=\pi(I)$ {and $\pi$ is $\Int$-rank-preserving,}
we have 
\[
\rk_M(I)
=
\rk_N(\pi(I))
=
\rk_N(\pi(I^\Ical))
=
\rk_M(I^\Ical).
\]
Also, observe that for any $J\in \Ical$, we have $I\subset J$ if and only if $I^\Ical\subset J$. Therefore,

\[\rk_M(I) = \rk_M(I^\Ical) = \sum_{\substack{J\supset I^\Ical\\ J\in \Ical}} \dgm_M^\Ical(J) = \sum_{\substack{J\supset I\\ J\in \Ical}} \dgm_M^\Ical(J), \]
i.e., $\rk^\Int_M$ is M\"obius invertible over $\Ical\subset \Int(\catP)$, as desired.
\end{proof}

Before proving Theorem \ref{thm:sufficient conditions for invertibility} \ref{item:finitely presentable invertible}, we provide 
a proposition and subsequent corollaries about the geometry of intervals in $\R^2$. 
\begin{proposition}\label{prop:path contained in rectangle}
    Let $p,q\in I\in \Int(\R^2)$. 
    Then there exists a path in $I$ connecting $p$ to $q$ that is contained in the segment $B=[p\wedge q, p\vee q]\subset \R^2$.  The same statement holds when $\R^2$ is replaced by $\Z^2$.
\end{proposition}
\begin{proof}
We prove the statement only for $\R^2$, as the proof extends trivially to $\Z^2$.
If we suppose $p$ and $q$ are comparable, the statement trivially holds: both $p,q$ and $q,p$ are such paths. Suppose not.
    Let $p=:(p^x,p^y)$ and $q=:(q^x,q^y)$, and without loss of generality, assume $p^x < q^x$ and $p^y > q^y$.
    Define the map $\varphi:\R^2\to B$ by sending each point in $\R^2$ to its closest point in $B$; see Figure~\ref{fig:projection to B}. Equivalently, we have
    \begin{equation}\label{eq:projection to B}
    \varphi((s,t)) = ((s\vee (p^x\wedge q^x))\wedge (p^x\vee q^x),\qquad (t\vee (p^y\wedge q^y))\wedge(p^y\vee q^y)).
    \end{equation}

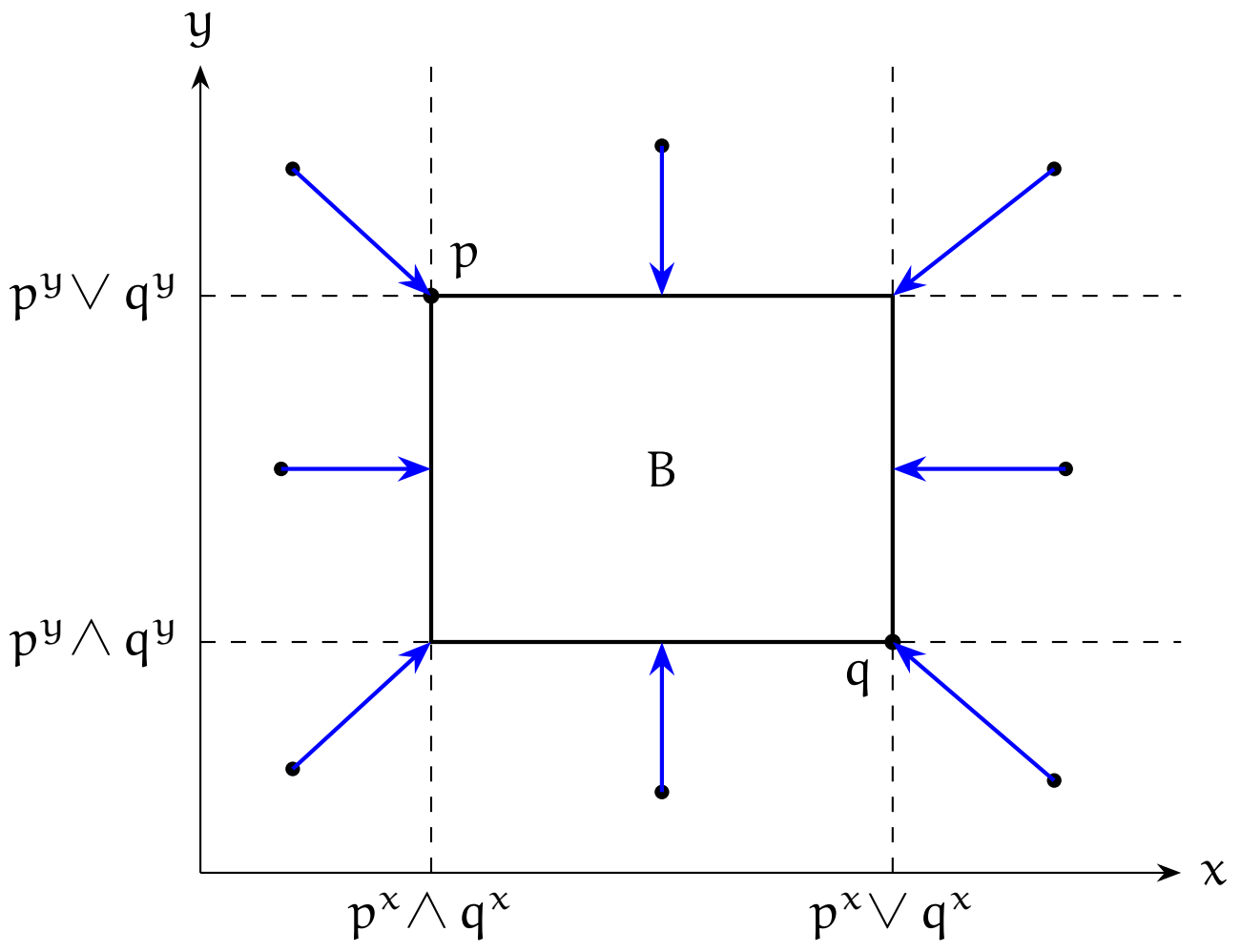
\begin{figure}
    \centering    
\begin{tikzpicture}[
    >=Stealth,
    scale=0.8,
    every node/.style={font=\small}
]

\def\xL{2}
\def\xR{6}
\def\yB{2}
\def\yT{5}

\draw[->] (0,0) -- (8.5,0) node[right] {$x$};
\draw[->] (0,0) -- (0,7) node[above] {$y$};

\draw[thick] (\xL,\yB) rectangle (\xR,\yT);

\draw[dashed] (\xL,0) -- (\xL,7);
\draw[dashed] (\xR,0) -- (\xR,7);
\draw[dashed] (0,\yB) -- (8.5,\yB);
\draw[dashed] (0,\yT) -- (8.5,\yT);

\node at (4,3.5) {$B$};

\fill (\xL,\yT) circle(2pt);
\node[above right] at (\xL,\yT) {$p$};

\fill (\xR,\yB) circle(2pt);
\node[below left] at (\xR,\yB) {$q$};

\node[below] at (\xL,0) {$p^x\!\wedge q^x$};
\node[below] at (\xR,0) {$p^x\!\vee q^x$};

\node[left] at (0,\yB) {$p^y\!\wedge q^y$};
\node[left] at (0,\yT) {$p^y\!\vee q^y$};

\fill (0.8,6.1) circle(1.8pt);
\draw[->,blue,thick] (0.8,6.1)--(\xL,\yT);

\fill (4,6.3) circle(1.8pt);
\draw[->,blue,thick] (4,6.3)--(4,\yT);

\fill (7.4,6.1) circle(1.8pt);
\draw[->,blue,thick] (7.4,6.1)--(\xR,\yT);

\fill (0.7,3.5) circle(1.8pt);
\draw[->,blue,thick] (0.7,3.5)--(\xL,3.5);

\fill (7.5,3.5) circle(1.8pt);
\draw[->,blue,thick] (7.5,3.5)--(\xR,3.5);

\fill (0.8,0.9) circle(1.8pt);
\draw[->,blue,thick] (0.8,0.9)--(\xL,\yB);

\fill (4,0.7) circle(1.8pt);
\draw[->,blue,thick] (4,0.7)--(4,\yB);

\fill (7.4,0.8) circle(1.8pt);
\draw[->,blue,thick] (7.4,0.8)--(\xR,\yB);

\end{tikzpicture}
\caption{The action of the map $\varphi$ given in Equation~\eqref{eq:projection to B} is depicted.}\label{fig:projection to B}    
\end{figure}

    Since equal application of meets and joins preserves order, $\varphi$ is order-preserving. Also observe that $\varphi(p) = p$ and $\varphi(q) = q$. %
    Hence, for any 
    path connecting $p$ to $q$ inside $I$, the image of the path under $\varphi$ is a path connecting $p$ to $q$ in $B$.    
    It suffices to show that 
   the image is in $I$. To this end, it suffices to show that $\varphi(I)\subset I$.
    Let $r=(s,t)\in I$ and $r':=\varphi(r)$.
    If $r\in B$, then $r' = r$, so $r'\in I$.
        Assume $r\notin B$. Then $r$ can be chosen from the eight regions outside $B$, that are depicted in Figure~\ref{fig:projection to B}. In all cases,  $r'$ is either $p$, $q$, a point in a segment whose endpoints are among $p,q$ and $r$. Since $p,q,r\in I$, and $I$ is convex, in all cases, $r'\in I$, as desired.
    In all possible cases, $r'\in I$.
\end{proof}

\begin{corollary}\label{cor:B cap I is an interval}
  Let $I\in \Int(\R^2)$. If the set $B:=[s_1,t_1)\times [s_2,t_2)\subset \R^2$ intersects $I$, then $B\cap I\in \Int(\R^2)$.  
  In addition, an analogous statement holds for $\Z^2$: the intersection of an interval and a segment in $\Z^2$ is either empty or an interval.
\end{corollary}
\begin{proof}
We prove the statement only for $\R^2$, as the proof extends trivially to $\Z^2$.
Since both $B$ and $I$ are convex (Definition~\ref{def:intervals}), $B\cap I$ is convex. It remains to show that $B\cap I$ is connected. Let $p,q\in B\cap I$. By the preceding proposition, there is a path in $[p\wedge q,\ p\vee q]\cap I$ connecting $p$ and $q$. 
It is not difficult to see that $B$ is closed under meet and join, and thus both $p\wedge q$ and $p\vee q$ belong to $B$.
Therefore, $[p\wedge q,\ p\vee q]\subset B$. Hence, the path connecting $p$ and $q$ is in $B\cap I$, as desired.
\end{proof}

\begin{corollary}[{Comparable points across adjacent blocks in an interval of $\R^2$}]\label{cor:adjacent blocks}
{Let \[
B_1:=[a_1,b_1)\times[c_1,d_1)
\qquad\text{and}\qquad
B_2:=[a_2,b_2)\times[c_2,d_2)
\]
be blocks in $\R^2$ such that either
\[
[a_1,b_1)=[a_2,b_2)
\quad\text{and}\quad
c_2=d_1,
\]
or
\[
[c_1,d_1)=[c_2,d_2)
\quad\text{and}\quad
a_2=b_1.
\]}
Let $I\in\Int(\R^2)$ such that $I\cap B_1\neq\emptyset$ and $I\cap B_2\neq\emptyset.$
Then there exist $p\in I\cap B_1$ and $q\in I\cap B_2$ such that
$
p\leq q.
$ 
 In addition, an analogous statement holds for $\Z^2$: if an interval in $\Z^2$ intersects two adjacent segments\footnote{For integers $a_i\leq b_i$ and $c_i\leq d_i$ ($i=1,2$), the segments
$
[a_1,b_1]\times[c_1,d_1]
 \text{ and } 
[a_2,b_2]\times[c_2,d_2]
$
in $\Z^2$ are called \emph{adjacent} if either (1)
$
[a_1,b_1]=[a_2,b_2]
\text{ and }
c_2=d_1+1,
$
or (2) 
$
[c_1,d_1]=[c_2,d_2]
\text{ and }
a_2=b_1+1.
$} in $\Z^2$, then there exist comparable elements, one from each intersection.
\end{corollary}
\begin{proof} 
We prove the statement only for $\R^2$, as the proof extends trivially to $\Z^2$. {We prove the statement only under the former condition.} Let $p=(p^x,p^y)\in I\cap B_1$ and $q=(q^x,q^y)\in I\cap B_2$.
    By Proposition \ref{prop:path contained in rectangle}, there exists a path $p=r_0,r_1,\ldots,r_n=q$ in $I$, which is also contained in %
    $B:=[p\wedge q,p\vee q]$. Note that $B$ is contained in the union $B_1\cup B_2=[a_1,b_1)\times [c_1,d_2)$.
    As $r_0\in B_1$ and $r_n\in B_2$, there exists the smallest %
    $j\in \{1,\ldots,n\}$ such that 
    $r_j\in B_2$.
    By this choice, $r_{j-1}\in B_1$, and as $r_{j-1},r_j$ are comparable we have $r_{j-1}\leq r_j$, as desired.
\end{proof}

\begin{remark}\label{rem:differences in higher dimensions}

    We note that the analogues of Proposition \ref{prop:path contained in rectangle} 
    and Corollaries \ref{cor:B cap I is an interval} and \ref{cor:adjacent blocks}
    in higher dimensions ($\Z^d$ and $\R^d$ for $d\geq 3$) do not hold, which prevents the proof that finitely presentable $\Z^2$- and $\R^2$-modules have M\"obius invertible Int-GRIs (i.e., Theorem \ref{thm:sufficient conditions for invertibility} \ref{item:finitely presentable invertible}) from extending to $\Z^d$- and $\R^d$-modules for $d\geq 3$; see Remark \ref{rem:difference for fp modules in dims 2 and 3+} for further details.

    For an interval in $\R^3$ violating the higher-dimensional analogue of Corollary \ref{cor:B cap I is an interval}, consider the interval $I_a$ shown in Figure \ref{fig:ia for fp not imply mi} (right), and the set $B = [0,2) \times [0,1) \times [1,2)$.
    Then $B\cap I = ([0,1)\times [a,1)\times [1,2)) \cup ([1,2) \times [0,a) \times [1,2))$ is not in $\Int(\R^3)$, as it is disconnected.
    This example extends to dimensions $d>3$ by replacing $I_a$ with $I_a\times [0,1)^{d-3}$ and $B$ with $B\times [0,1)^{d-3}$.%
       
    {For an interval in $\R^3$ violating the higher-dimensional analogue of Corollary \ref{cor:adjacent blocks},
    consider the interval $I_a$ shown Figure \ref{fig:ia for fp not imply mi} (right), and the adjacent unit blocks $B_{0,1}$ and $B_{1,1}$ from Equation~\eqref{eq:unit blocks}.
    Although $I_a$ intersects both blocks,}
    there is no pair of points $p,q\in I_a$ with $p\in B_{0,1}$, $q\in B_{1,1}$, and $p\leq q$.   
\end{remark}

We now prove Theorem \ref{thm:sufficient conditions for invertibility} \ref{item:finitely presentable invertible}, building on Theorem \ref{thm:sufficient conditions for invertibility} \ref{item:finitely encodable}. 
{More specifically, given a finitely presentable $\R^2$- (or $\Z^2$-) module, we identify an order-preserving map from $\R^2$ to a finite poset $Q$ satisfying all the conditions of Theorem~\ref{thm:sufficient conditions for invertibility}~\ref{item:finitely encodable}.
} 

\begin{definition}\label{def:smallest grid}{Let $S\subset \R^d$ (resp. $\Z^d$) be a set or multiset. For $i=1,\ldots,d$, let $S_i\subset \R$ (resp. $\Z$) be the set of the $i$-th coordinates of the points in $S$. Then we call the product $\prod_{i=1}^d S_i \subset \R^d$ (resp. $\Z^d$)  the \textbf{smallest grid containing} $S$.} 
\end{definition}

\begin{proof}[Proof of Theorem \ref{thm:sufficient conditions for invertibility} \ref{item:finitely presentable invertible}] 
{Assume that $P=\R^2$. For the case of $P=\Z^2$, all of the claims below and their proofs follow upon simple replacement of all instances of $\R^2$ with $\Z^2$.} Let $M'$ be a nontrivial finitely presentable $\R^2$-module.
Let $\phi:\bigoplus_{a\in A} \kf_{a^\uparrow}\rightarrow \bigoplus_{b\in B}\kf_{b^\uparrow}$ be a minimal presentation for $M'$ (Definition~\ref{def:finitely presentable}).
Let 
\begin{equation}\label{eq:grid}\catQ:=\{g_1^x<\ldots<g_{n_x}^x\}\times \{g_1^y<\ldots<g_{n_y}^y\}\subset \R^2
\end{equation}
be the smallest grid containing all the elements in $A$ and $B$.
Let $q_\ast$ be the unique minimal element of $\catQ$.
Denote $M:=M'\vert_{q_\ast^\uparrow}$, and define the $\catQ$-module $N$ by $N:= M\vert_{\catQ}$.

Let $\lfloor - \rfloor_{\catQ}:q_\ast^\uparrow\rightarrow \catQ$ be the map sending each $p\in q_\ast^\uparrow$ to the maximal $q\in \catQ$ such that $q\leq p$.
We will show that $\lfloor - \rfloor_{\catQ}$ is interval-preserving, and $\Int$-rank-preserving from $M$ to $N$, which by Theorem \ref{thm:sufficient conditions for invertibility} \ref{item:finitely encodable} implies that the Int-GRI of $M$ is M\"{o}bius invertible. {This will ultimately imply that the Int-GRI of $M'$ is M\"obius invertible, as shown in Claim~\ref{claim:M' is mobius invertible} below.}

\resetclaimcounter
\begin{claim}
$\lfloor - \rfloor_{\catQ}$ is interval-preserving.
\end{claim}
\begin{claimproof}
Fix $I\in \Int(q_\ast^\uparrow)$. Then $I$ is connected. Since $\lfloor - \rfloor_{\catQ}$ is order-preserving, the image $\lfloor I \rfloor_{\catQ}$ is also connected. 
Next, we show that $\lfloor I \rfloor_{\catQ}$ is convex. Let $q_1< q_2<q_3\in \catQ$ such that $q_1,q_3\in \lfloor I\rfloor_{\catQ}$. We wish to show that $q_2\in \lfloor I \rfloor_Q$.
Let $p_1, p_3\in I$ such that $\lfloor p_1\rfloor_{\catQ} = q_1$ and $\lfloor p_3\rfloor_{\catQ} = q_3$.
This implies that 
\begin{equation}\label{eq:inequalities}
\begin{tikzcd}[column sep=-0.5em, row sep=-0.5em]
    p_{1} & {\geq} & q_{1} & {<} & q_{2} & {<} & q_{3} & {\leq} &p_{3}.
\end{tikzcd}
\end{equation}
We now break into cases based on the order relationship between $p_1$ and $q_2$.
If $q_2<p_1$, then $q_2\leq q_1$ by the definition of $\lfloor-\rfloor_Q$, which is a contradiction and so this case is impossible.
If $p_1\leq q_2$, then $p_1\leq q_2\leq p_3$ alongside the convexity of $I$ implies $q_2\in I$.
Then, $\lfloor q_2\rfloor_{\catQ} = q_2$, implying that $q_2\in \lfloor I\rfloor_{\catQ}$.

Lastly, assume that $p_1$ and $q_2$ are incomparable. We write $p_1=(p_1^x,p_1^y)$, and let $(i,j)\in \{1,\ldots,n_x\}\times\{1,\ldots,n_y\}$ such that $q_2=(g_i^x,g_j^y)$ (cf. Equation~\eqref{eq:grid}).
Without loss of generality, assume that 
\begin{equation}\label{eq:p1 and q2}
    p_1^x<g_i^x \mbox{ and } p_1^y>g_j^y\ \mbox{(i.e., $p_1$  is above and to the left of $q_2$ in the plane)}.
\end{equation}
As $\lfloor p_1\rfloor_{\catQ} = q_1 < q_2$, we have that 
\begin{equation}\label{eq:p1y leq gi+1y}
g_j^y \leq p_1^y < g_{j+1}^y,
\end{equation}
{with the convention that $g_{n_y+1}^y:=\infty$.}
Again, we break into cases based on the order relationship between $p_1$ and $p_3$. We note that, $p_3<p_1$  is impossible: assuming $p_3<p_1$ contradicts the earlier assumption that $q_1=\lfloor p_1 \rfloor_Q <q_3= \lfloor p_3 \rfloor_Q$ given in \eqref{eq:inequalities}. Hence, either $p_1<p_3$, or $p_1$ and $p_3$ are incomparable.

Assume that $p_1<p_3$. Then by the inequalities given in \eqref{eq:inequalities}, \[p_1 < p_1 \vee q_2 = (g_i^x,p_1^y) \leq  p_3.\] 
Thus, by convexity of $I$, we have $p_1\vee q_2\in I$. Also, the equality $p_1 \vee q_2 = (g_i^x,p_1^y)$ together with \eqref{eq:p1y leq gi+1y} implies that $\lfloor p_1\vee q_2\rfloor_{\catQ} = (g^x_i, g^y_j) = q_2$, implying that $q_2\in \lfloor I\rfloor_{\catQ}$.

Assume that $p_1=(p_1^x,p_1^y)$ and $p_3:=(p_3^x,p_3^y)$ are incomparable. 
Then, by assumption \eqref{eq:p1 and q2} and 
the inequality $q_2< p_3$ from \eqref{eq:inequalities}, we have that $p_1^x<p_3^x$ and $p_1^y>p_3^y$ (i.e., $p_1$ is above and to the left of $p_3$ in the plane).
By connectivity of $I$, there exists a path $p_1=a_1, a_2, \ldots, a_n=p_3$ in $I$. {There are the following two cases.} 

\noindent \emph{(Case 1) There exists  $1\leq \ell\leq n$ such that  $a_\ell\leq q_2$.} Then by \eqref{eq:inequalities}, we have that $a_\ell\leq q_2 < p_3$, and convexity of $I$ implies that $q_2\in I$. Hence, $q_2=\lfloor q_2 \rfloor_Q$, implying that $q_2 \in \lfloor I \rfloor_Q$.

\noindent \emph{(Case 2) For all $1\leq \ell\leq n$, $a_\ell\not\leq q_2$ .} Let us write $a_\ell=(a_\ell^x,a_\ell^y)$ for each $\ell$. As $I$ is an interval in $\R^2$, we may assume:
\begin{equation}\label{eq:shortest path in R2}
  \mbox{For all $1\leq \ell\leq n$ 
, \ \ \  $p_1^x\leq a_\ell^x\leq p_3^x$ and $p_3^y\leq a_\ell^y\leq p_1^y$.}
\end{equation}
Now let $k:=\min\{\ell\in \{1\ldots,n-1\}: q_2<a_{\ell+1}\}$.  This $k$ is well-defined since $q_2<a_n=p_3$.
By \eqref{eq:inequalities} and \eqref{eq:shortest path in R2},
we have $ g_j^y \leq  p_3^y \leq  a_k^y$.
Since $q_2\not \leq a_k$, we then have $ a_k^x< g_i^x$. 
Further, since $a_k\not\leq q_2=(g_i^x,g_j^y)$, we have $g_j^y< a_k^y$. 
As a result, $a_k\vee q_2 = (g_i^x,a_k^y)$.
{By assumption \eqref{eq:shortest path in R2}}, $a_k^y\leq p_1^y$. 
Thus, invoking \eqref{eq:p1y leq gi+1y}, we have $g_j^y < a_k^y \leq p_1^y < g_{j+1}^y$, and in turn
\[\lfloor a_k \vee q_2 \rfloor_{\catQ} = \lfloor (g_i^x,a_k^y)\rfloor_{\catQ} = (g_i^x, g_j^y) = q_2.\]
By assumption, $a_k, q_2 \leq a_{k+1}$ and $a_k,a_{k+1}\in I$.
Thus we have $a_k\leq a_k\vee q_2\leq a_{k+1}$, and convexity of $I$ implies $a_k\vee q_2\in I$.
Since $q_2 = \lfloor a_k\vee q_2 \rfloor_{\catQ}$, $q_2\in \lfloor I \rfloor_{\catQ}$, and we are done.

\end{claimproof}

\begin{claim} $\lfloor - \rfloor_{\catQ}$ is $\Int$-rank-preserving. 
\end{claim}
\begin{claimproof}
We  wish to prove that
\begin{equation}\label{eq:corollary Rd int}
    \mbox{for all intervals}\ I\subset q_\ast^\uparrow,\ \  \rk_M(I) = \rk_N(\lfloor I \rfloor_{\catQ}).
\end{equation}

Let $I\in \Int(q_\ast^\uparrow)$.
Let $\Psi_M$ be the canonical limit-to-colimit map over $M\vert_I$ and $\Psi_N$ be the canonical limit-to-colimit map over $N\vert_{\lfloor I\rfloor_\catQ}$. By the rank-nullity theorem, it suffices to prove that   $\varprojlim M\vert_I\cong \varprojlim N\vert_{\lfloor I \rfloor_{\catQ}}$ and $\ker(\Psi_M)\cong \ker(\Psi_N)$. 

Consider the map $\Phi:\varprojlim M\vert_I\to \varprojlim N\vert_{\lfloor I \rfloor_{\catQ}}$ given by, for each $\bw:=(\bw_p)_{p\in I}\in \varprojlim M\vert_I$ (cf. Proposition \ref{prop:limit and colimit}),
    \[\Phi(\bw):=\bv:=(\bv_q)_{q\in \lfloor I \rfloor_{\catQ}}\ ,\] where $\bv_q:=\bw_p$ for some $p\in I$ such that $\lfloor p \rfloor_{\catQ} = q$ (see Figure \ref{fig:rd mobius invertible} (A)). %
    It suffices to prove:
    
\emph{Subclaim 1.} $\Phi$ is well-defined. 

\emph{Subclaim 2.} $\Phi$ is a linear isomorphism.

\emph{Subclaim 3.} %
The restriction of $\Phi$ to $\ker(\Psi_M)$ is a bijection between $\ker(\Psi_M)$ and $\ker(\Psi_N)$.

\noindent Before proving these, note:
\begin{enumerate}

 \item   
 $q_\ast^\uparrow$ is partitioned into the preimages $B_q$ of the elements $q\in\catQ$ under $\lfloor-\rfloor_Q$.
 Note that
    \[B_q = [g_i^x, g_{i+1}^x) \times [g_j^y, g_{j+1}^y)\]
    for some $1\leq i \leq n_x$ and $1\leq j \leq n_y$ where $g_{n_x+1}:=\infty$ and $g_{n_{y}+1}:=\infty$. 
    Each $B_q$ will be called a \emph{block}. 
    If $q,q'\in Q$ with $q\leq q'$, then we write $B_q\leq B_{q'}$. In addition, if $[q,q']=\{q,q'\}$, then the blocks $B_q$ and $B_{q'}$ are called \emph{adjacent}.
 \item For any ${p\leq p'\in \R^2}$, if
 $\lfloor p \rfloor_{\catQ} = \lfloor p' \rfloor_{\catQ}$, then $\varphi_M(p,p'):M_p \rightarrow M_{p'}$ is the identity. In other words, $M$ is constant on each block $B_q$. \label{item:floor1}
\end{enumerate}

\begin{proof}[Proof of Subclaim 1]     %
For each $q\in \catQ$, since $M$ is constant on the block $B_q$, the vector $\bv_q$ does not depend on the choice of $p$. 
    It remains to show that 
    $\bv$ belongs to $\varprojlim N\vert_{\lfloor I \rfloor_\catQ}$, i.e., 
    if $q< q'\in \lfloor I \rfloor_{\catQ}$, then $\varphi_N(q, q')(\bv_q) = \bv_{q'}$.
    Let $p,p'\in I$ such that $\lfloor p \rfloor_{\catQ} = q$ and $\lfloor p' \rfloor_{\catQ} = q'$.
    If $p\leq p'$, then the equality $\varphi_N(q, q')(\bv_q) = \bv_{q'}$ directly follows from the facts that $\varphi_M(p, p')(\bw_p) = \bw_{p'}$, that $M$ is constant on each of the blocks $B_q$ and $B_{q'}$, and that $N=M|_{Q}.$

    Thus, 
    suppose $p$ and $p'$ are not comparable (see Figure \ref{fig:rd mobius invertible} (B) for an illustrative example). This, together with the earlier assumption $q<q'$, implies that $q$ and $q'$ have either the same first coordinate or the same second coordinate. Without loss of generality, assume the former. This implies the existence of a finite sequence of adjacent blocks $B_q=B_1\leq B_2\leq \ldots\leq B_n=B_{q'}$ for some $n\in \N$, all of which are of the form $[g_{i}^x,g_{i+1}^x)\times [-)$  for some $1\leq i\leq n_X$ so that the union $B:=\bigcup_{k=1}^n B_k$ is also in the form $[g_{i}^x,g_{i+1}^x)\times [-)$. 
    
    Since $B\cap I$ contains $p$ and $p'$, by Corollary~\ref{cor:B cap I is an interval}, $B\cap I$ is an interval. 
    {Also, convexity of $I$ implies that $I$ intersects each of $B_1,\ldots,B_n$. By Corollary~\ref{cor:adjacent blocks}, for each $k\in \{1,\ldots,n-1\}$, we can pick a pair $(p_k,\bar{p}_{k+1})\in B_k\times B_{k+1}$ such that $p_k\leq \bar{p}_{k+1}$. Let $\bar{p}_1:=p$ and $p_n:=p'$. Let $k\in \{1,\ldots,n\}$. Then,  $B_k\cap I$ contains both $\bar{p}_k$ and $p_k$ and thus, by Corollary~\ref{cor:B cap I is an interval}, $B_k\cap I$ is an interval. Let $\Gamma_k$ be a path in $B_k\cap I$ connecting from $\bar{p}_k$ to $p_k$. Then, the concatenation $\Gamma_1,\Gamma_2,\ldots,\Gamma_n$ forms a path in $B\cap I$ that connects from $p$ to $p'$. We denote this path by $\Gamma$.
    } 
    Note that the restriction of the section $(\bw_p)_{p\in I}$ to $\Gamma$ induces a section of $M$ along $\Gamma$ and in turn, since $M$ is constant on every block, induces a section of $N=M|_{Q}$ along $\lfloor \Gamma \rfloor_Q$. This implies that the composition of the structure maps in $N$ along $\lfloor \Gamma \rfloor_Q$, possibly flipping the direction of identity maps, yields the map $\varphi_N(q, q')$ (since $N$ is a functor $\catQ\to \cvec$) and that $\varphi_N(q, q')(\bv_q)=\bv_{q'}$, as desired.   
\end{proof}

    \begin{figure}
        \begin{center}
            \begin{tikzpicture}[
                scale=0.7,
    x=1.25cm,
    y=1.25cm,
    >=Stealth,
    gridline/.style={
        draw=gray!35,
        line width=3pt
    },
    panelborder/.style={
        draw=gray!70,
        line width=1.2pt,
        rounded corners=4mm
    },
    latticepoint/.style={
        circle,
        fill=black,
        draw=gray!20,
        line width=.8pt,
        minimum size=4mm,
        inner sep=0pt
    },
    region/.style={
        fill=regionBlue,
        fill opacity=.65,
        draw=none
    },
    patharrow/.style={
        ->,
        line width=1pt,
        >={Stealth[length=2.3mm,width=1.6mm]}
    },
    blacklabel/.style={
        font=\Large,
        text=black
    },
    redlabel/.style={
        font=\Large,
        text=red
    }
]

\begin{scope}

\foreach \x in {0,2,4}
    \draw[gridline] (\x,-0.50) -- (\x,6.55);

\foreach \y in {0,2,4,6}
    \draw[gridline] (-0.72,\y) -- (4.70,\y);

\fill[region]
    (0.25,5.40) --
    (0.95,5.40) --
    (0.95,3.75) --
    (2.70,3.75) --
    (2.70,2.65) --
    (3.12,2.65) --
    (3.12,0.82) --
    (1.35,0.82) --
    (1.35,2.65) --
    (0.95,2.65) --
    (0.95,3.15) --
    (0.25,3.15) --
    cycle;

\foreach \x in {0,2,4}
    \foreach \y in {0,2,4,6}
        \node[latticepoint] at (\x,\y) {};

\node[blacklabel] at (0.63,5.02)
    {$\mathbf{w}_{p_1}$};

\node[blacklabel] at (1.48,3.37)
    {$\mathbf{w}_{p_2}$};

\node[blacklabel] at (2.52,3.37)
    {$\mathbf{w}_{p_3}$};

\node[blacklabel] at (1.62,1.35)
    {$\mathbf{w}_{p_4}$};

\node[blacklabel] at (2.62,1.35)
    {$\mathbf{w}_{p_5}$};

\node[redlabel,anchor=east] at (0.15,3.68)
    {$\mathbf{v}_{q_1}$};

\node[redlabel,anchor=east] at (0.15,1.68)
    {$\mathbf{v}_{q_2}$};

\node[redlabel,anchor=east] at (0.15,-0.32)
    {$\mathbf{v}_{q_4}$};

\node[redlabel] at (2.42,2.30)
    {$\mathbf{v}_{q_3}$};

\node[redlabel] at (1.58,-0.33)
    {$\mathbf{v}_{q_5}$};

\draw[panelborder]
    (-0.82,-0.65) rectangle (4.82,6.68);

\node[font=\Large] at (2,-1.12) {(A)};

\end{scope}

\begin{scope}[shift={(5.75,0)}]

\foreach \x in {0,2,4}
    \draw[gridline] (\x,-0.50) -- (\x,6.55);

\foreach \y in {0,2,4,6}
    \draw[gridline] (-0.72,\y) -- (4.70,\y);

\fill[region]
    (0.25,5.40) --
    (0.95,5.40) --
    (0.95,3.75) --
    (2.70,3.75) --
    (2.70,2.65) --
    (3.12,2.65) --
    (3.12,0.82) --
    (1.35,0.82) --
    (1.35,2.65) --
    (0.95,2.65) --
    (0.95,3.15) --
    (0.25,3.15) --
    cycle;

\foreach \x in {0,2,4}
    \foreach \y in {0,2,4,6}
        \node[latticepoint] at (\x,\y) {};

\node[blacklabel] (pprime) at (0.58,5.08)
    {$p'$};

\node[blacklabel] (p2) at (0.58,3.55)
    {$p_1$};

\node[blacklabel] (p1) at (1.7,3.55)
    {$p_2$};

\node[blacklabel] (p) at (1.7,1.35)
    {$p$};

\draw[patharrow]
    (p2.north) -- (pprime.south);

\draw[patharrow]
    (p2.east) -- (p1.west);

\draw[patharrow]
    (p.north) -- (p1.south);

\node[redlabel,anchor=east] at (-0.10,-0.12)
    {$q$};

\node[redlabel,anchor=east] at (-0.10,3.98)
    {$q'$};

\node[blacklabel] at (2.58,1.5)
    {$I$};

\draw[panelborder]
    (-0.82,-0.65) rectangle (4.82,6.68);

\node[font=\Large] at (2,-1.12) {(B)};

\end{scope}

\end{tikzpicture}
            \caption{In each of panels~(A) and~(B), an interval $I\in\R^2$ is depicted as the shaded region. (A) A section $(\bw_p)_{p\in I}$ and its image $(\bv_q)_{q\in \lfloor I \rfloor_\catQ}:= \Phi((\bw_p)_{p\in I})$. 
            (B) Illustration for $q=\lfloor p \rfloor_Q$ , $q'=\lfloor p' \rfloor_Q$. %
            The existence of a path in $Q$ connecting $q$ and $q'$ guarantees the existence of a path in $I$ connecting $p$ and $p'$, which is exemplified by the path $p', p_1, p_2, p$.}
            \label{fig:rd mobius invertible}
        \end{center}
    \end{figure}
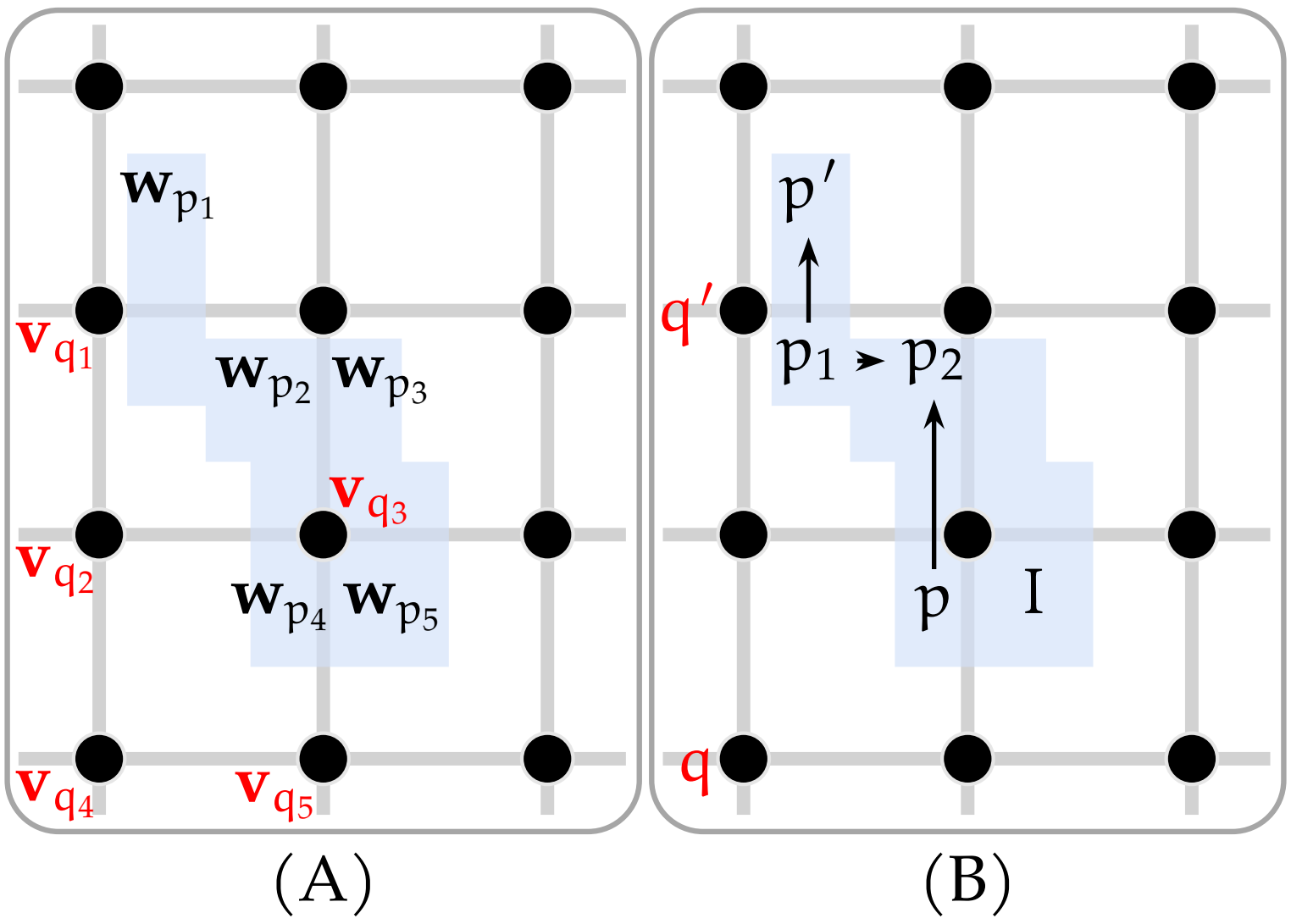

    \begin{proof}[Proof of Subclaim 2] 
    We construct the inverse $\Phi\inv:\varprojlim N\vert_{\lfloor I\rfloor_\catQ}\to \varprojlim M\vert_I$.
    For 
    $\bv\in\varprojlim N\vert_{\lfloor I\rfloor_\catQ}$, define $\Phi\inv(\bv):= \bw$, where $\bw_p :=\bv_q$ for $p\in I$ with $\lfloor p \rfloor_\catQ=q$.   Indeed $\bw$ belongs to $\varprojlim M\vert_I$, as (1)
    if $p\leq p'\in I$, then $\lfloor p \rfloor_\catQ\leq \lfloor p'\rfloor_\catQ$, (2) $M$ is constant on each block $B_q$, and (3) $\bv$ is a section.
    Lastly, it is immediate to see by the definitions that $\Phi\circ \Phi\inv$ is the identity on $\varprojlim N\vert_{\lfloor I\rfloor_\catQ}$ and $\Phi\inv\circ \Phi$ is the identity on $\varprojlim M\vert_I$.
    \end{proof}

    \begin{proof}[Proof of Subclaim 3.]
    
    We have already demonstrated that $\Phi$  injective, so it remains to show that $\Phi(\ker(\Psi_M)) = \ker(\Psi_N)$. 
    To see $\Phi(\ker(\Psi_M))\subset \ker(\Psi_N)$, suppose $(\bw_p)_{p\in I}\in \ker(\Psi_M)$.
    By Proposition \ref{prop:path characterization of zero in colim}, 
    there is a section $\mathbf{v}$ along a path $\Gamma$: $p'$, $\ldots$, $p''$ in $I$ with $\mathbf{v}_{p'}=\bw_{p'}$ and $\mathbf{v}_{p''}=0\in M_{p''}$.
    Define $\Phi_\Gamma:\varprojlim M|_{\Gamma}\to \varprojlim M|_{\lfloor \Gamma \rfloor_Q}$ similarly to $\Phi$, and let
    $\Phi_\Gamma(\mathbf{v})=:\bu$.
    Then, denoting $q':=\lfloor p' \rfloor_\catQ$, $q'':=\lfloor p''\rfloor_\catQ$ and $\Gamma_\catQ:=\lfloor \Gamma \rfloor_\catQ$, we see that $\bu_{q'} = \bw_{p'}$, $\bu_{q''} = 0$, and $\Gamma_\catQ$ is a path 
    in $\lfloor I\rfloor_\catQ$ connecting $q'$ to $q''$. By Proposition \ref{prop:path characterization of zero in colim}, this shows $\Psi_N\left(\Phi((\bw_p)_{p\in I})\right) = 0$, i.e., $\Phi((\bw_p)_{p\in I})\in \ker(\Psi_N)$, as desired.
    
  On the flip side, suppose we have a section $\bv:=(\bv_q)_{q\in \lfloor I\rfloor_\catQ}\in \ker(\Psi_N)$. 
    By Proposition \ref{prop:path characterization of zero in colim}, this means that there is a section $\bu$ along a path $\Delta$: $q_1$, $\ldots$, $q_n$ in $\lfloor I\rfloor_\catQ$ with $\bu_{q_1}=\bv_{q_1}$ and $\bu_{q_n}=0\in N_{q_n}$. 
    We may assume that any two consecutive points in $\Delta$ are distinct and consecutive in $Q$, i.e., for every $1\leq k< n$, the set $\{q_k,q_{k+1}\}$ forms a size-two segment in $Q$. 
    We wish to show that 
    $\Phi\inv(\bv)\in \varprojlim M|_{I}$ belongs to   $\ker(\Psi_M)$. To this end, we again invoke Proposition \ref{prop:path characterization of zero in colim}.  For each $1\leq k\leq n$, pick any $p_k\in I$ with $\lfloor p_k\rfloor_Q=q_k$. To each $p_k$, assign the vector $\mathbf{w}_{p_k}:=\bu_{q_k}$. Since $M$ is constant on the block $B_{q_k}$, we have $\mathbf{w}_{p_k}\in M_{p_k}$. 
    If $p_1,\ldots,p_n$ form a path in $I$, then we are done, as $\mathbf{w}_{p_1}=\bu_{q_1}=\bv_{q_1}=\Phi\inv(\bv)_{p_1}$ and $\mathbf{w}_{p_n}=\bu_{q_n}=0$. 
    Suppose not, i.e., there exists $1\leq k<n$ such that $p_k$ and $p_{k+1}$ are incomparable. 
    Fix such $k$.
    By the construction of the path $\Delta$, $B_{q_k}$ and $B_{q_{k+1}}$ are adjacent blocks.
    Hence, by Corollary~\ref{cor:B cap I is an interval}, $(B_{q_k}\cup B_{q_{k+1}})\cap I$ is an interval containing $p_k$ and $p_{k+1}$. 
    Let $\Gamma$ be a path in $(B_{q_k}\cup B_{q_{k+1}})\cap I$ connecting $p_k$ and $p_{k+1}$. Assigning $\bu_{q_k}$ to points in $\Gamma$ within $B_{q_k}$ and $\bu_{q_{k+1}}$ to points in $\Gamma$ within $B_{q_{k+1}}$ clearly induces a section of $M$ along $\Gamma$. This proves that we can extend $(\mathbf{w}_{p_k})_{1\leq k\leq n}$ to a section of $M$ along a path in $I$, which takes the values $\Phi\inv(\bv)_{p_1}$ and $0$ at its endpoints.
        \end{proof}   
\end{claimproof}

\begin{claim} The Int-GRI of $M'$ is M\"obius invertible. \label{claim:M' is mobius invertible}
\end{claim}

\begin{claimproof}We have shown $\lfloor - \rfloor_{\catQ}$ is interval-preserving and Int-rank-preserving with respect to $M$ and $N$, so by Theorem \ref{thm:sufficient conditions for invertibility} \ref{item:finitely encodable}, $M:= M'\vert_{q_\ast^\uparrow}$ has a M\"obius invertible Int-GRI, i.e.,
$\dgm_{M}^{\Int(q_\ast^\uparrow)}:\Int(q_\ast^\uparrow)\to \Z$ exists.
If $I\in \Int(\R^2)$ with $I\not\subset q_\ast^\uparrow$, then there exists $p\in I$ with $p\not\geq q_\ast$, which implies $M_p=0$.
Then $\rk^\Int_{M'}(\{p\}) = 0$, and so by Remark \ref{rem:properties of rank invariant} \ref{item:monotonicity}, $\rk^\Int_{M'}(I) = 0$.

Define $d_{M'}:\Int(\R^2)\to \Z$ as follows:
\[d_{M'}(I):= \begin{cases}
    0,& \mbox{if } I \not\subset q_\ast^\uparrow\\
    \dgm_M^{\Int(q_\ast^\uparrow)}(I), & \mbox{if } I\subset q_\ast^\uparrow.
\end{cases}\]

It is immediate that for all $I\in \Int(\R^2)$, the set $\{J\supset I: d_{M'}(J)\neq 0 \mbox{ and } J\in \Int(\R^2)\}$ is finite, since it coincides with the set $\{J\supset I: \dgm_M^{\Int(q_\ast^\uparrow)}(J)\neq 0 \mbox{ and } J\in \Int(q_\ast^\uparrow)\}$.
Let $I\in \Int(\R^2)$.
If $I\notin \Int(q_\ast^\uparrow)$, then as $J\supset I$ implies $J\notin \Int(q_\ast^\uparrow)$, we have:
\[0 = \rk_{M'}^{\Int(\R^2)}(I) = \sum_{\substack{J\supset I\\J\in \Int(\R^2)}} d_{M'}(J).\]
If $I\in \Int(q_\ast^\uparrow)$, then as $M = M'\vert_{q_\ast^\uparrow}$ we have:
\[ \rk_{M'}^{\Int(\R^2)}(I) = \rk_M^{\Int(q_\ast^\uparrow)}(I) = \sum\limits_{\substack{J\supset I\\J\in \Int(q_\ast^\uparrow)}} \dgm_M^{\Int(q_\ast^\uparrow)}(J) = \sum\limits_{\substack{J\supset I\\J\in \Int(\R^2)}} d_{M'}(J).\]
Thus, by Definition \ref{def:GGPD}, the GRI of $M'$ is M\"obius invertible, as desired.
\end{claimproof}
\end{proof}

\begin{figure}
\begin{center}
\tdplotsetmaincoords{60}{125}
\tdplotsetrotatedcoords{0}{30}{120} 
\begin{tikzpicture}[scale=1.8,tdplot_rotated_coords,
                    edge/.style={black,thick},
                    edgeDash/.style={black,thick,dashed},
                    face/.style={fill=gray!30,opacity=0.6},centernode/.style={circle,inner sep=1pt,draw=none,fill=none}]

\fill[face] (0, 1, 0.5) -- (1,1,0.5) -- (1,2,0.5) -- (0,2,0.5) -- cycle;
\fill[face] (0,2,0) -- (1,2,0) -- (1,2,0.5) -- (0,2,0.5) -- cycle;
\fill[face] (1,2,0) -- (2,2,0) -- (2,2,0.5) -- (1,2,0.5) -- cycle;
\fill[face] (2,1,0) -- (2,2,0) -- (2,2,0.5) -- (2,1,0.5) -- cycle;
\fill[face] (2,0,0) -- (2,1,0) -- (2,1,0.5) -- (2,0,0.5) -- cycle;
\fill[face] (1,0,0.5) -- (2,0,0.5) -- (2,1,0.5) -- (1,1,0.5) -- cycle;
\fill[face] (1,0,0.5) -- (1,0,1) -- (1,1,1) -- (1,1,0.5) -- cycle;
\fill[face] (0,1,0.5) -- (0,2,0.5) -- (0,2,0) -- (0,3,0) -- (0,3,1) -- (0,1,1) -- cycle;
\fill[face] (1,1,0.5) -- (1,3,0.5) -- (1,3,1) -- (1,1,1) -- cycle;
\fill[face] (1,1,0.5) -- (3,1,0.5) -- (3,1,1) -- (1,1,1) -- cycle;
\fill[face] (0,2,0) -- (2,2,0) -- (2,0,0) -- (3,0,0) -- (3,3,0) -- (0,3,0) -- cycle;
\fill[face] (3,0,0) -- (3,3,0) -- (3,3,0.5) -- (3,1,0.5) -- (3,1,1) -- (3,0,1) -- cycle;
\fill[face] (0,3,0) -- (3,3,0) -- (3,3,0.5) -- (1,3,0.5) -- (1,3,1) -- (0,3,1) -- cycle;

\draw[edgeDash] (2,1,0.5) -- (2,2,0.5);
\draw[edgeDash] (1,2,0.5) -- (2,2,0.5);
\draw[edge] (2,2,0) -- (2,2,0.5);
\draw[edge] (1,0,0.5) -- (1, 1, 0.5);
\draw[edge] (1,0,0.5) -- (2,0,0.5);
\draw[edge] (1,0,0.5) -- (1,0,1);
\draw[edge] (2, 0, 0) -- (2, 0, 0.5);
\draw[edge] (2, 0, 0.5) -- (2, 1, 0.5);
\draw[edgeDash] (1, 0, 1) -- (2, 0, 1);
\draw[edgeDash] (1, 0, 1) -- (1, 1, 1);
\draw[edge] (2, 1, 0) -- (2, 1, 0.5);
\draw[edge] (2, 0, 0) -- (2, 1, 0);
\draw[edge] (2, 0, 0) -- (3, 0, 0);
\draw[edgeDash] (2, 0, 1) -- (3, 0, 1);
\draw[edgeDash] (3, 0, 0) -- (3, 0, 1);
\draw[edgeDash] (1, 1, 0.5) -- (1, 1, 1);
\draw[edgeDash] (1, 1, 1) -- (3, 1, 1);
\draw[edgeDash] (3, 0, 0) -- (3, 3, 0);
\draw[edgeDash] (3, 0, 1) -- (3, 1, 1);
\draw[edge] (0, 1, 0.5) -- (0, 1, 1);
\draw[edge] (0, 1, 0.5) -- (1, 1, 0.5);
\draw[edgeDash] (1, 1, 0.5) -- (3, 1, 0.5);
\draw[edgeDash] (3, 1, 0) -- (3, 1, 1);
\draw[edgeDash] (0, 1, 1) -- (1, 1, 1);
\draw[edgeDash] (3, 1, 0.5) -- (3, 3, 0.5);
\draw[edgeDash] (1, 1, 0.5) -- (1, 3, 0.5);
\draw[edgeDash] (1, 3, 0.5) -- (3, 3, 0.5);
\draw[edge] (1, 2, 0) -- (2, 2, 0);
\draw[edge] (2, 1, 0) -- (2, 2, 0);
\draw[edgeDash] (1, 3, 0) -- (1, 3, 0.5);
\draw[edgeDash] (3, 3, 0) -- (3, 3, 0.5);
\draw[edge] (0, 1, 0.5) -- (0, 2, 0.5);
\draw[edge] (0, 2, 0.5) -- (1, 2, 0.5);
\draw[edgeDash] (0, 1, 1) -- (0, 3, 1);
\draw[edgeDash] (0, 3, 1) -- (1, 3, 1);
\draw[edgeDash] (1, 3, 1) -- (1, 3, 0.5);
\draw[edge] (0, 2, 0.5) -- (0, 2, 0);
\draw[edge] (0, 2, 0) -- (0, 3, 0);
\draw[edgeDash] (0, 3, 0) -- (0, 3, 1);
\draw[edge] (0, 2, 0) -- (1, 2, 0);
\draw[edgeDash] (0, 3, 0) -- (3, 3, 0);
\draw[edge] (1, 2, 0) -- (1, 2, 0.5);
\draw[edgeDash] (1,1,1) -- (1,3,1);

\node[shape=circle, inner sep=2pt, fill=black] at (0,0,0) {};
\node[centernode] at (0,-0.25,0) {$(0,0,0)$};
\draw[edgeDash] (0,0,0) -- (2,0,0);

\node[shape=circle, inner sep=2pt, fill=black] at (3,3,0.5) {};
\node[centernode] at (3.15,3,1) {$(3,0.5,3)$};

\node[shape=circle, inner sep=2pt, fill=black] at (1,3,0.5) {};
\node[centernode] at (1.4,3.3,0.5) {$(1,0.5,3)$};

\node[shape=circle, inner sep=2pt, fill=black] at (1,1,0.5) {};
\node[centernode] at (1.5,0.7,0.5) {$(1,0.5,1)$};

\node[shape=circle, inner sep=2pt, fill=black] at (2,0,0) {};
\node[centernode] at (2,0,-0.35) {$(2,0,0)$};

\node[shape=circle, inner sep=2pt, fill=black] at (0,1,0.5) {};
\node[centernode] at (-0.5,1.1,0) {$(0,0.5,1)$};

\node[shape=circle, inner sep=2pt, fill=black] at (0,3,1) {};
\node[centernode] at (-0.5,3,1) {$(0,1,3)$};

\end{tikzpicture}
\caption{An interval $I$ in $\R^3$ such that $\lfloor I \rfloor_{\catQ}$ is not an interval in $\catQ$. {Corner points at the ends of dashed line segments (e.g. $(1,0.5,1)$) do not belong to the interval.}}
\label{fig:interval not mapped to interval}
\end{center}
\end{figure}
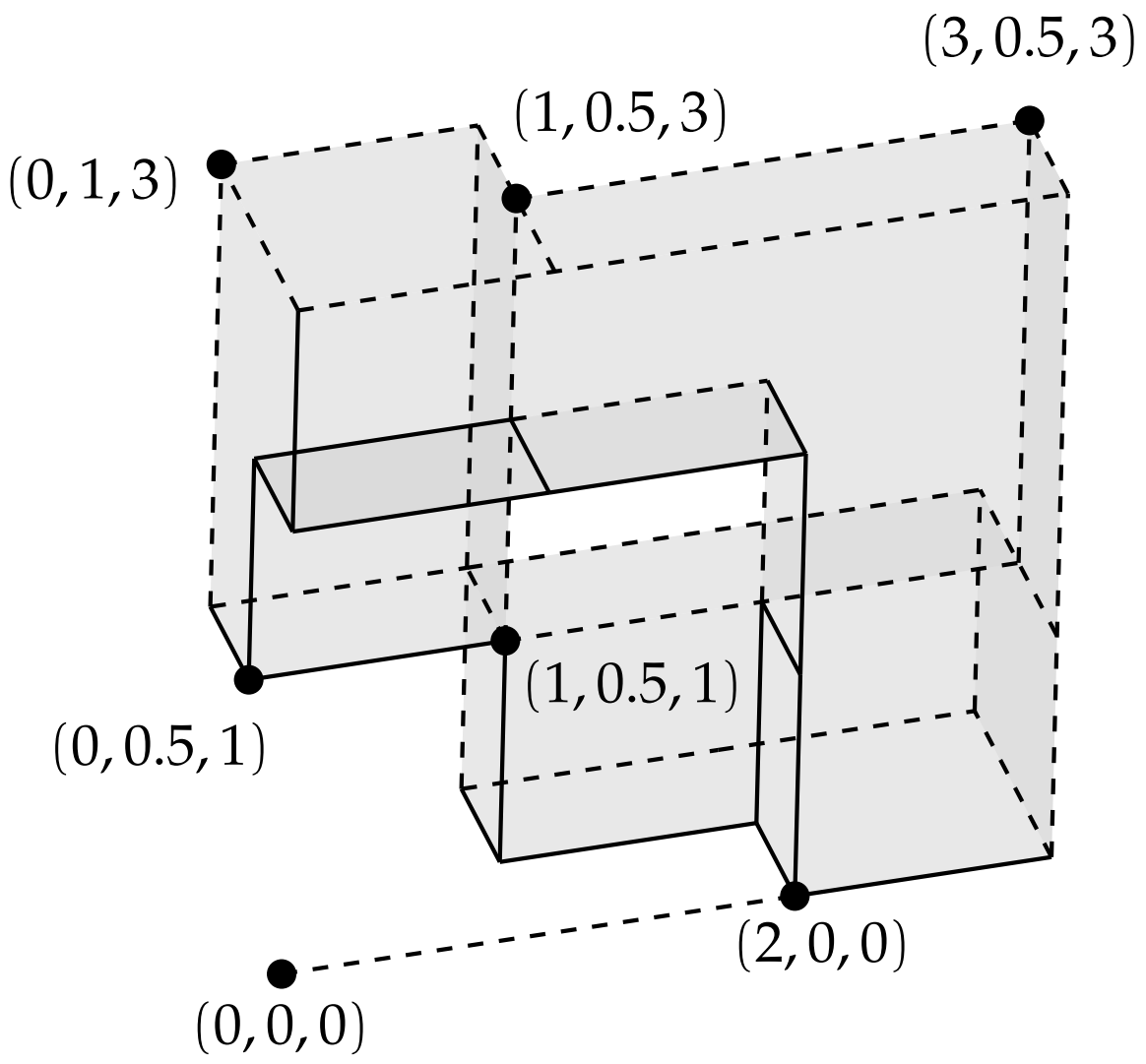

\begin{remark}[Discrepancy between Theorem~\ref{thm:tame does not imply Int-GRI invertible}~\ref{item:fp not imply MI} and Theorem~\ref{thm:sufficient conditions for invertibility}~\ref{item:finitely presentable invertible}]\phantomsection\label{rem:difference for fp modules in dims 2 and 3+}
  Let $d\geq 3$ and consider a finitely presentable $\R^d$-module $M$. From the proof of Theorem \ref{thm:sufficient conditions for invertibility}~\ref{item:finitely presentable invertible}, take the higher-dimensional analog of the grid $\catQ$ and the map $\lfloor - \rfloor_{\catQ}$ for $M$. Unlike the case of $d=2$, the map is not necessarily Int-rank preserving {nor interval-preserving}. Thus, Theorem~\ref{thm:sufficient conditions for invertibility} \ref{item:finitely encodable} does not apply to $\lfloor - \rfloor_{\catQ}$ for establishing the M\"obius invertibility of $\rk_M^\Int$. 
  
  To see the failure of the Int-rank preserving property, from Example~\ref{ex:another R3 module}, 
    recall the $\R^3$-module $M$, and the interval $I_a$ for any fixed $a\in (0,1)$, and the blocks $B_{i,j}$ for $0\leq i,j\leq 2$.
We saw that $\rk_M(I_a)=1$. 
However, it is not difficult to see that $\rk_{M\vert_{\catQ}}(\lfloor I_a\rfloor_{\catQ})=0$. This discrepancy follows from the fact that $\lfloor B_{0,1} \rfloor =(0,0,1)\leq (1,0,1)=\lfloor B_{1,1}\rfloor$ in $\lfloor I_a\rfloor_{\catQ}$, whereas no pair of points $p\in B_{0,1}\cap I_a$ and $q\in B_{1,1}\cap I_a$ satisfies $p\leq q$, as observed in Remark~\ref{rem:differences in higher dimensions}.
    This shows that $\lfloor-\rfloor_Q$ is not a full functor from the category $I_a$ to the category $\lfloor I_a \rfloor$ and thus there exists a structure map in $M\vert_{\lfloor I_a\rfloor_{\catQ}}$ not corresponding to any structure map of $M\vert_{I_a}$. This is the source of the discrepancy.
    
    {To see the failure of the interval-preserving property,}
    consider the following interval of $\R^3$ (see Figure \ref{fig:interval not mapped to interval}):
    \[I:= [(0,0.5,1),(1,1,3)) \cup [(0,0,2),(3, 0.5, 3)) \cup [(1, 0.5, 0), (3, 1, 1)) \cup [(2, 0, 0), (3, 0.5, 2)).\]
    This interval equals $I_{0.5}\setminus [(1,0,0), (2,0.5,2))$, with $I_{0.5}$ being the interval defined in the proof of Example~\ref{ex:another R3 module}
    and depicted in Figure~\ref{fig:ia for fp not imply mi}. 
    Then, its image
{\[\lfloor I\rfloor_{\catQ} = [(0,0,0),(2,0,2)]_{\Z^3}\setminus\{(0,0,0),(1,0,1)\}\]}
    is not an interval, as it is not convex due to missing the point $(1, 0, 1)$.
\end{remark}

For the Con-GRI, the statement analogous to  Theorem \ref{thm:sufficient conditions for invertibility} \ref{item:finitely encodable} also holds:

\begin{proposition}\label{prop:tame+ does imply MI}

    Let $M$ be a $\catP$-module.
    If 
    there exists a finite poset $\catQ$, an order-preserving map $\pi:\catP\to \catQ$ and $\catQ$-module $N$ such that $\pi$ is {Con-}rank-preserving, i.e., for all $I\in \Con(\catP)$, $\rk_M(I) = \rk_N(\pi(I))$, then $\rk_M$ is M\"obius invertible. 
    
\end{proposition}

\begin{proof}
For any $J\subset \catQ$, let $C(\pi\inv(J))$ be the set of connected components of $\pi\inv(J)$.
This implies $C(\pi\inv(J))\subset \Con(\catP)$.
Define $\Ical\subset \Con(\catP)$ as $\Ical:=\bigcup_{J\in \Con(\catQ)} C(\pi\inv(J))$.
From here, the argument that $\rk_M:\Con(\catP)\to \Z$ is M\"obius invertible over $\Ical$ follows the same as the proof of Theorem \ref{thm:sufficient conditions for invertibility}
\ref{item:finitely encodable} upon replacing all instances of $\Int(-)$ with $\Con(-)$.\qedhere
\end{proof}
Interestingly, the statement analogous to Theorem \ref{thm:sufficient conditions for invertibility} \ref{item:finitely presentable invertible} for the Con-GRI is false:

\begin{remark}\label{rem:finitely presented R^d not Con MI}
       There exist finitely presentable $\R^2$-modules whose Con-GRIs are not M\"obius invertible.   The construction of such an example is parallel to the one described in the proof of Theorem \ref{thm:tame does not imply Int-GRI invertible} \ref{item:tame not imply MI}, and thus we omit it. In particular, when constructing such an example, the map $\lfloor - \rfloor_\catQ$ considered in the proof of Theorem \ref{thm:sufficient conditions for invertibility} \ref{item:finitely presentable invertible} is not Con-rank-preserving, and thus we cannot utilize Proposition \ref{prop:tame+ does imply MI} in that case. 
\end{remark}

\paragraph{Pullback of the GPD.} Theorem \ref{thm:sufficient conditions for invertibility} \ref{item:finitely encodable} and Proposition \ref{prop:tame+ does imply MI} show that the presence of a rank-preserving order-preserving map $\pi:\catP\rightarrow\catQ$ with $\abs{\catQ}<\infty$ ensures the M\"obius invertibility of the GRI. Next, we will demonstrate that by imposing an additional assumption on $\pi$, the GPD of $M$ can be fully determined by the GPD of $N$.
 
 An order-preserving map $\pi:\catP\rightarrow \catQ$ is called a \textbf{covering} (resp. \textbf{Int-covering}) if $\pi$ is surjective and for each $J\in \Con(\catQ)$ (resp. $J\in \Int(\catQ)$), 
 each connected component of $\pi^{-1}(J)$ maps onto $J$ under $\pi$. 
If $\pi$ is a covering, then it is an Int-covering by definition. However, the converse does not necessarily hold (see Example~\ref{ex:counterexample} in the appendix). 
 Let $\pi^{-1}\Con(\catQ):=\{\pi^{-1}(J):J\in \Con(\catQ)\}$, and define $\pi^{-1}\Int(\catQ)$ similarly.

\begin{theorem}\label{thm:pull-back GPD} Let $M$ be a $\catP$-module and let $N$ be a $\catQ$-module with $|Q|<\infty$.
\begin{enumerate}[label=(\roman*)]
    \item Assume that a covering $\pi:\catP\rightarrow \catQ$ is Con-rank-preserving. 
Then, the Con-GRI of $M$ is M\"obius invertible, and
 for $I\in \Con(\catP)$, 
\[\dgm_M^\Con(I)=\begin{cases}\dgm_N^\Con(\pi(I)),&\mbox{if $I$ is a connected component of an element in $\pi^{-1}\Con(\catQ)$}\\0,&\mbox{otherwise.}\end{cases}\]
\label{item:pull-back Con-GPD}
\item Assume that an Int-covering $\pi:\catP\rightarrow \catQ$ is Int-rank-preserving and interval-preserving. 
Then, the Int-GRI of $M$ is M\"obius invertible, and
 for $I\in \Int(\catP)$, 
\begin{equation}\label{eq:pullback Int-GPD}
\dgm_M^\Int(I)=\begin{cases}\dgm_N^\Int(\pi(I)),&\mbox{if $I$ is a connected component of an element in $\pi^{-1}\Int(\catQ)$}\\0,&\mbox{otherwise.}\end{cases}
\end{equation}
\label{item:pull-back Int-GPD}
\end{enumerate}
\end{theorem}

This theorem can be useful for significantly reducing the computational complexity of the GPD; such a use of the theorem can be found in \cite[Section 3]{kim2025super}.

\begin{proof} 
We first prove Item~\ref{item:pull-back Con-GPD}. Since $\catQ$ is finite,
$\dgm_N^\Con$ exists. Let
\[
\Ical
:=
\left\{
J\in\Con(\catP):
J\text{ is a connected component of }\pi^{-1}(U)
\text{ for some }U\in\Con(\catQ)
\right\},
\]
and define the map
\begin{align*}
d:\Con(\catP)&\longrightarrow \Z,\\
J&\longmapsto
\begin{cases}
\dgm_N^\Con(\pi(J)), & \text{if } J\in\Ical,\\
0, & \text{otherwise}.
\end{cases}
\end{align*}

We wish to show that $d=\dgm_M^\Con$. Fix $I\in\Con(\catP)$. For every $U\in\Con(\catQ)$ with
$U\supset\pi(I)$, the connected set $I$ is contained in a unique connected
component $J_U$ of $\pi^{-1}(U)$.  Conversely, every $J\in\Ical$ containing
$I$ is of this form, with $U=\pi(J)$. Thus,
we have a surjection
\begin{align*}
\{U\in\Con(\catQ):U\supset \pi(I)\}
&\longrightarrow
\{J\in\Ical:J\supset I\},\\
U&\longmapsto J_U.
\end{align*}
Since $\pi$ is a covering, (1) the above map is also injective (thus bijective), and (2) for every $U\in\Con(\catQ)$ with
$U\supset\pi(I)$, we have $\pi(J_U)=U$ and hence
$d(J_U)=\dgm_N^\Con(U).$
Therefore,
\begin{align*}
\sum_{\substack{J\supset I\\J\in\Con(\catP)}}d(J)
=\sum_{\substack{J\supset I\\J\in\Ical}}d(J)
=
\sum_{\substack{U\supset\pi(I)\\U\in\Con(\catQ)}}
\dgm_N^\Con(U)
=
\rk_N^\Con(\pi(I))
=
\rk_M^\Con(I),
\end{align*}
where the last equality follows from the Con-rank-preserving property of
$\pi$. The sum is finite because $\catQ$ is finite. Hence, $d$ is the
M\"obius inversion of $\rk_M^\Con$, proving
Item~\ref{item:pull-back Con-GPD}.

For Item~\ref{item:pull-back Int-GPD}, the same argument applies with
$\Con$ replaced by $\Int$. The interval-preserving assumption ensures that
$\pi(I)\in\Int(\catQ)$ for every $I\in\Int(\catP)$,
Proposition~\ref{prop:each component is an interval} ensures that each
connected component of $\pi^{-1}(U)$ belongs to $\Int(\catP)$, and the
Int-covering assumption gives $\pi(J_U)=U$.
\end{proof}

\begin{example}\label{ex:anti-diagonal}
Consider the $\Z^2$-module $M$ defined as follows:
\[M_{(x,y)}=\begin{cases} 0,&y<-x \\ k^2,&y=-x \\ k,&y>-x,
\end{cases}
\]
where $M$ is constant on $y>-x$ and each map $k^2\rightarrow k$ in $M$ is the projection $p_1$ %
{$:(x,y)\mapsto x$.}
Although $M$ is not finitely presentable, $M$ is tame: Let $N$ be the $\{1<2<3\}$-module given by
$0\rightarrow k^2\stackrel{p_1}{\rightarrow} k.$
There exists the obvious order-preserving map $\pi:\Z^2\rightarrow \{1,2,3\}$ such that $M=\pi^\ast N$. This $\pi$ is a covering and Int-rank-preserving. Thus, by the previous theorem, not only $\dgm_M^\Int$ exists, but also $\dgm_M^\Int$ is directly obtained from $\dgm_N^\Int$. 

Another way to show that the GRI of $M$ is M\"obius invertible is to prove that $M$ is interval-decomposable, and invoke Theorem \ref{thm:sufficient conditions for invertibility} \ref{item:interval-decomposable}. 
\end{example}

{By exploiting Theorem \ref{thm:pull-back GPD},} we refine Theorem \ref{thm:sufficient conditions for invertibility} \ref{item:finitely presentable invertible} by specifying the precise types of intervals that appear in the support of the Int-GPD of a finitely presentable $\R^2$- or $\Z^2$-module. 

\begin{definition}\label{def:minimal grid}Let $M$ be a finitely presentable $\R^d$- or $\Z^d$-module and $\phi:\bigoplus_{a\in  A}k_{a^\uparrow}\to\bigoplus_{b\in B}k_{b^\uparrow}$ be a 
minimal presentation of $M$. Let $Q_M$ be the smallest grid containing $A\cup B$ (Definition~\ref{def:smallest grid}), which we call the \textbf{minimal grid} for $M$. Let $\Ical_M$
be the set of the preimages of all elements in $\Int(\catQ_M)$ under $\lfloor-\rfloor_{\catQ_M}$. Equivalently, $\Ical_M$ is the subset of $\Int(\R^d)$ whose elements are unions of \textbf{unit blocks} $\lfloor q\rfloor_Q^{-1}$, $q\in Q_{M}$.
Each interval in $\Ical_M$ is said to be \textbf{subordinate} to $\catQ_M$. 
\end{definition}

\begin{remark}\label{rem:about smallest grid} 
\begin{enumerate}[label=(\roman*)]
\item If $M=0$, then $Q_M$ is empty and thus $\Ical_M$ is also empty.\label{item:smallest grid of zero module}
\item Let $M$ and $N$ be any two finitely presentable $\R^d$-modules (or $\Z^d$-modules). Then, by Remark~\ref{rem:minimal grid}~\ref{item:betti numbers of direct sum}, both $\Ical_M$ and $\Ical_N$ are subsets of $\Ical_{M\oplus N}$.\label{item:smallest grid of direct sum} 
\end{enumerate}
\end{remark}

\begin{proposition}\label{prop:preimages under the floor map are connected}
Let $M$ be a finitely presentable $\R^d$-module. Then $\Ical_M$ consists solely of intervals that possess finitely many minimal and maximal elements, i.e., 
$\Ical_M\subset \bigcup_{\substack{m\in\N\\n\in \Z_{\geq 0}}}\Int_{m,n}(\R^d)$ (cf. Figure~\ref{fig:intmnr2}). An analogous statement holds for any finitely presentable $\Z^d$-module.
\end{proposition}
\begin{proof} Let $I\in \Ical_M$. By Proposition~\ref{prop:each component is an interval},  $I$ is convex. We show
that $I$ is connected. Let $a,b\in I$, and assume that $I$ is the preimage of an $I'\in \Int(\catQ_M)$ under $\lfloor-\rfloor_{Q_M}$. 
Let $a':=\lfloor a \rfloor_{Q_M}, b':=\lfloor b \rfloor_{Q_M}\in I'$. Since $I'$ is connected, there exists a path $a',p_1,\ldots, p_n,b'$ in $I'$. By the definition of the map $\lfloor-\rfloor_{Q_M}$, we have $I'\subset I$ and thus this path is contained in $I$. Then, $a,a',p_1,\ldots, p_n,b', b$ is a path in $I$, 
showing that $I$ is connected. Since $I$ is the preimage of an interval of $\catQ_M$ under $\lfloor-\rfloor_{\catQ_M}$, and $\catQ_M$ is finite, we have either $I\in \Int_{m,n}(\R^d)$ with $m=n=\abs{\catQ_M}$ or $I\in \Int_{m,0}(\R^d)$ with $m=\abs{\catQ_M}$, completing the proof.
\end{proof}

\begin{proposition}\label{prop:finitely presentable invertiver over intmn}
The Int-GRI of any finitely presentable $\R^2$-module $M$ is M\"obius invertible over $\Ical_M$. More specifically, for any $I\in \Int(\R^2)$,
\[\dgm_M^\Int(I)=\begin{cases}\dgm_{M|_{\catQ}}^\Int(\lfloor I\rfloor_\catQ),&\mbox{if $I\in \Ical_M$}\\0,&\mbox{otherwise.}\end{cases}\]
 An analogous statement holds for any finitely presentable $\Z^2$-module.
\end{proposition}

\begin{proof}
 Without loss of generality, let $M$ be the $\R^2$-module $M'$ given in the proof of Theorem \ref{thm:sufficient conditions for invertibility}~\ref{item:finitely presentable invertible}. In that proof, we showed that, for the $n_x\times n_y$ grid $Q(=\catQ_M)$ in $\R^2$ with the minimum point $q_\ast$, the map $\lfloor - \rfloor_{\catQ}:q_\ast^\uparrow \rightarrow Q$ is interval-preserving and Int-rank-preserving, and that $\dgm_M^\Int$ coincides with $\dgm_{M|_{q_\ast^\uparrow}}^\Int$ on $\Int(q_\ast^\uparrow)$ and vanishes outside $\Int(q_\ast^\uparrow)$.
 
 Observe that the map $\lfloor - \rfloor_{\catQ}$ is an Int-covering: (1) $\lfloor - \rfloor_{\catQ}$ is surjective by definition; (2) We already saw from the proof of the preceding proposition that the preimage of any interval $J$ in $\catQ$ under $\lfloor - \rfloor_{\catQ}$ is simply an interval $I$ in $\R^2$. 
 {This implies that the unique component $I$ of $\lfloor-\rfloor_{\catQ}^{-1}(J)$ maps onto $J$. Therefore, by applying Theorem~\ref{thm:pull-back GPD}~\ref{item:pull-back Int-GPD} to the module $M|_{q_\ast^\uparrow}$ (corresponding to $M$ in the theorem), the map $\lfloor - \rfloor_{\catQ}:q_\ast^\uparrow \rightarrow Q$ (corresponding to $\pi$), and the module $M|_Q$ (corresponding to $N$), Equation~\eqref{eq:pullback Int-GPD} becomes}
\[\forall I\in \Int(q_\ast^\uparrow),\ \  \dgm_{M|_{q_\ast^\uparrow}}^\Int(I)=\begin{cases}\dgm_{M|_{\catQ}}^\Int(\lfloor I\rfloor_\catQ),&\mbox{if $I\in \Ical_M$}\\0,&\mbox{otherwise.,}\end{cases}\]
This completes the proof of the claimed equation.

One can similarly prove the analogous statement for a finitely presentable $\Z^2$-module.
\end{proof}

\subsection{GPD for finitely presentable $\R^d$-modules}\label{sec:GPD for finitely presentables}
For persistence modules over arbitrary posets, the Int-GRI generalizes the persistence diagram or the barcode in the sense of Theorem~\ref{thm:sufficient conditions for invertibility} \ref{item:interval-decomposable}. However, we saw that for non-interval-decomposable $\R^d$- or $\Z^d$-modules with $d\geq 3$, the Int-GRI might not be M\"obius invertible, even when those modules are finitely presentable (cf. Theorem~\ref{thm:tame does not imply Int-GRI invertible} \ref{item:fp not imply MI} and Theorem~\ref{thm:sufficient conditions for invertibility}~\ref{item:finitely presentable invertible}). 

This fact suggests that the Int-GRI might be an unnecessarily complex invariant for finitely presentable $\R^d$- or $\Z^d$-modules, and motivates us {to ask whether there exists a domain of the GRI such that its M\"obius inversion is always well defined and allows recovery of the GPD over all intervals whenever the latter exists.
In this section, we show that the answer to this question is affirmative (Theorem~\ref{thm:unique candidate}). Interestingly, to this end, we make use of the minimal grid (Definition~\ref{def:minimal grid}).}

From Definition~\ref{def:minimal grid}, recall that, given a finitely presentable $\R^d$-module $M$, $\Ical_M$ denotes the set of intervals of $\R^d$ subordinate to the minimal grid $\catQ_M$. Since $\Ical_M$ is finite, $\rk_M^{\Ical_M}$ is  convolvable, and thus  M\"obius invertible, as noted in  Remark~\ref{rem:monotonicity of mobius invertibility new}~\ref{item:comparison between mobius and convolv}. 

\begin{definition}[GPD for finitely presentable $\R^d$-modules]\label{def:GPD for fp modules}Given any finitely presentable $\R^d$-module $M$, we define the \textbf{GPD of $M$} to be the trivial extension of $\dgm_M^{\Ical_M}$ to $\Int(\R^d)$, i.e., 
\begin{align*}
\dgm_M : \Int(\R^d) &\longrightarrow \Z, \\
I &\longmapsto
\begin{cases}
\dgm_M^{\Ical_M}(I) & \text{if } I \in {\Ical_M}, \\
0               & \text{otherwise.}
\end{cases}
\end{align*}
\end{definition}

We remark that, by Remark~\ref{rem:about smallest grid}~\ref{item:smallest grid of zero module}, if $M=0$, then $\dgm_M$ is constant at zero. The GPD for finitely presentable $\R^d$-modules satisfies the following desirable properties:

\begin{proposition}
\label{prop:pixelized GPD} 
Let $M$ be any finitely presentable $\R^d$-module. Then,
\begin{enumerate}[label=(\roman*)]
\item (computation finiteness) 
For any $I\in \Ical_M$, \[\dgm_M(I)=\dgm_{M|_{Q_M}}^{\Int(Q_M)}(\lfloor I\rfloor_{Q_M}).\]\label{item:computation finiteness}
\item (generalizing barcodes) If $M$ is interval-decomposable, then for each $I\in \Int(\R^d)$, $\dgm_M(I)$ equals the multiplicity of $I$ in $\barc(M)$.  \label{item:generalizing barcode}

\end{enumerate}

\end{proposition}

\begin{proof}
\noindent\ref{item:computation finiteness}: The poset $(\Ical_M,\supset)$ is isomorphic to $(\Int(Q_M),\supset)$ via the poset isomorphism $I\mapsto \lfloor I\rfloor_Q$. Therefore, for any $J,I\in \Ical_M$ with $J\supset I$, we have that \[\mu_{\Ical_M}(J,I)=\mu_{\Int(Q_M)}(\lfloor J\rfloor_{Q},\lfloor I \rfloor_{Q}).\] Also, it is clear that, for any $I\in \Ical_M$, $\rk_M(I)=\rk_{M|_{Q_M}}(\lfloor I\rfloor_{Q_M})$. These observations directly imply the claim.

\noindent \ref{item:generalizing barcode}: Assume that $M$ is interval-decomposable. Then, its restriction to $Q_M$ is also interval-decomposable, i.e., there exist $U_1,\ldots,U_n\in \Int(Q_M)$ such that $M|_{Q_M}\cong \bigoplus_{i=1}^n k_{U_i}$. Let $\iota:Q_M\hookrightarrow \R^d$ be the canonical inclusion. Then, we have
\begin{align*}
M&\cong \Lan_{\iota}\left(M|_{Q_M}\right)&\mbox{by construction} 
\\ &\cong \Lan_{\iota}\left(\bigoplus_{i=1}^n k_{U_i}\right)&\mbox{by assumption}\\
&\stackrel{(\ast)}{\cong} \bigoplus_{i=1}^n \Lan_{\iota}\left(k_{U_i}\right)&\mbox{see below}\\
&\cong \bigoplus_{i=1}^n k_{\lfloor U_i \rfloor_{{Q_M}}^{-1}}.
\end{align*}
Let $\cvec^{\catQ_M}$ and $\cvec^{\R^d}$ be the category of $\catQ_M$- and $\R^d$-modules, respectively. The isomorphism ($\ast$) follows from the two facts: (1) the left Kan extension functor $\Lan_\iota:\cvec^{\catQ_M}\to \cvec^{\R^d}$ is left adjoint to the restriction functor $(-)|_{Q_M}:\cvec^{\R^d}\to \cvec^{\catQ_M}$. (2) Left adjoint functors preserve colimits (especially, coproducts) \cite[Theorem V.5.1]{mac2013categories}.
Since $\lfloor U_i \rfloor_{{Q_M}}^{-1}\in \Ical_M$ for each $i$, the claim directly follows from Item~\ref{item:computation finiteness} and Theorem~\ref{thm:sufficient conditions for invertibility}~\ref{item:interval-decomposable}.
\end{proof}

We show that if $\dgm^\Int_M$ exists, then it must coincide with $\dgm_M$:
\begin{theorem}\label{thm:unique candidate} Let $M$ be any finitely presentable $\R^d$-module. If $\rk_M^\Int$ is M\"obius invertible,   
i.e., $\dgm^{\Int}_M$ exists, 
then $\dgm^\Int_M=\dgm_M$.  %
\end{theorem}
We remark that the analogue of Theorem \ref{thm:unique candidate} for $\Z^d$-modules remains an open question. 
In particular, the proof of Theorem \ref{thm:unique candidate} relies heavily on a construction of an infinite family of intervals using continuous deformations.

{This theorem, together with Proposition~\ref{prop:pixelized GPD}~\ref{item:generalizing barcode}, suggests that the restriction of the GRI to the intervals subordinate to the minimal grid is an optimal choice of domain that both generalizes the concept of the persistence barcode and provides a finite encoding of the entire Int-GRI.} This theorem also implies that $\dgm_M^{\Int}$ can be computed in practice by first computing $\dgm_M$, which only requires a minimal presentation of $M$ and the GPD of the corresponding persistence module over the minimal grid. {The following remark is another consequence of the preceding theorem:}

\begin{remark}[Additivity of $\dgm$ relies on the dimension $d$] Let $M,N$ be any two finitely presentable $\R^d$-modules. Then, $M\oplus N$ is  finitely presentable (cf. Remark~\ref{rem:minimal grid}~\ref{item:betti numbers of direct sum}). If $d\in \{1,2\}$, then 
\begin{align*}
\dgm_{M\oplus N}&=\dgm^{\Int}_{M\oplus N}&&\mbox{by Theorems \ref{thm:sufficient conditions for invertibility}~\ref{item:finitely presentable invertible} and  \ref{thm:unique candidate}}
\\&=\dgm^{\Int}_{M}+\dgm^{\Int}_{N}&&\mbox{by Proposition~\ref{prop:GPD additivity}}
\\&=\dgm_{M}+\dgm_{N}&&\mbox{by Theorems \ref{thm:sufficient conditions for invertibility}~\ref{item:finitely presentable invertible} and \ref{thm:unique candidate}}.
\end{align*}
However, when $d\geq 3$, the equality $\dgm_{M\oplus N}=\dgm_M+\dgm_N$ does not necessarily hold. 
For example, let $M$ be the $\R^3$-module defined in Theorem~\ref{thm:tame does not imply Int-GRI invertible}~\ref{item:fp not imply MI}, which is depicted in Figure~\ref{fig:ia for fp not imply mi}. 
Also, for $I_{0.5}\in \Int(\R^3)$ depicted in Figure~\ref{fig:ia for fp not imply mi}, let $N$ be the interval module $k_{I_{0.5}}$. 
Then, $\dgm_M(I_{0.5})=0$ since $I_{0.5}\not\in \Ical_M$, and  $\dgm_N(I_{0.5})=1$ by Proposition~\ref{prop:pixelized GPD}~\ref{item:generalizing barcode}. This shows that $\dgm_M(I_{0.5})+\dgm_N(I_{0.5})=1$. On the other hand,
\begin{align*}
\dgm_{M\oplus N}(I_{0.5})&=\dgm_{M\oplus N}^{\Ical_{M\oplus N}}(I_{0.5})&&\mbox{since $I\in \Ical_N \subset \Ical_{M\oplus N}$ (cf. Remark~\ref{rem:about smallest grid}~\ref{item:smallest grid of direct sum})}
\\&\stackrel{(\ast)}{=}\rk_{M\oplus N}^{\Ical_{M\oplus N}}(I_{0.5})&&\mbox{see below}
\\&=\rk_{M}^{\Ical_{M\oplus N}}(I_{0.5})+\rk_{N}^{\Ical_{M\oplus N}}(I_{0.5})&&\mbox{by Remark~\ref{rem:properties of rank invariant}~\ref{item:additivity}}
\\&=\rk(M|_{I_{0.5}})+\rk(N|_{I_{0.5}})&&\mbox{by definition}
\\&=2.
\end{align*}
To see the equality ($\ast$), recall from the proof of Theorem~\ref{thm:tame does not imply Int-GRI invertible}~\ref{item:fp not imply MI} that $I_{0.5}$ is a maximal interval for $\rk_M^\Int$. It is also clear that $I_{0.5}$ is a maximal interval for $\rk_N^{\Int}$. These imply that, for all $J\in \Ical_{M\oplus N}$ with $J\supsetneq I_{0.5}$, $\rk(M|_{J})+\rk(N|_J)=0$, and thus $\rk_{M\oplus N}^{\Ical_{M\oplus N}}(J)=0$. Hence, $I_{0.5}$ is a maximal interval for $\rk_{M\oplus N}^{\Ical_{M\oplus N}}$, and the equality ($\ast$) follows from Proposition~\ref{prop:property of maximal intervals}. 
\end{remark}
{Although $\dgm$ generalizes the barcode  for finitely presentable $\R^d$-modules with any $d\geq 1$ (Proposition~\ref{prop:pixelized GPD}~\ref{item:generalizing barcode}), its failure to be additive when $d\geq 3$ precludes it from being a \emph{barcoding invariant} \cite{escolar2024barcoding} in dimensions $d\geq 3$.}

To prove Theorem~\ref{thm:unique candidate}, we make several remarks and establish a number of technical propositions.
In the sequel, let $M$ be a finitely presentable $\R^d$-module, let $Q_M$ be the minimal grid for $M$ with minimal point $q_\ast$, $\lfloor - \rfloor:q_\ast^\uparrow \to Q_M$ be the floor map, and $\Ical_M$ the collection of intervals subordinate to the grid (cf Definition \ref{def:minimal grid}).
We call $B\in \Ical_M$ a \emph{(unit) block} if $B = \lfloor - \rfloor\inv(q)$ for some $q\in Q_M$.
Define $\Bcal$ to be the set of all unit blocks in $\Ical_M$.

\begin{remark}\label{rem:properties of finitely generated modules}
\begin{enumerate}[label=(\roman*)]

\item 
    If $I_1,I_2\in \Con(\R^d)$ with $I_1\cap I_2\neq \emptyset$, then $I_1\cup I_2\in \Con(\R^d)$. \label{item:union of connected sets with overlap is connected}

\item     Let $\catP$ be a poset and $I\subset \catP$. Then, the (order-)convex hull of $I$ is defined as
    \[ \overline{I} := \bigcup_{p\leq q\in I}\, [p,q]. \]
    If $I\in \Con(\catP)$, then $\overline{I}$ is the smallest interval of $P$ containing $I$. \label{item:convex hull is smallest interval}
 \item For all $p,q\in q_\ast^\uparrow$ with $p\leq q$, the map $\varphi_M(p, q)$ depends only on the pair of blocks of $\Bcal$ containing $p$ and $q$.
 \label{item:properties of finitely generated modules-maps between blocks}

 \item Let $I\in \Int(\R^d)$ with $I\subset q_\ast^\uparrow$. 
 By the definition of $\Ical_M$, $I\notin \Ical_M$ if and only if there exists a unit block $B\in \Bcal$ with $\emptyset \subsetneq B\cap I\subsetneq B$.%
\label{item:partially contained block}
\end{enumerate} 
\end{remark}

\begin{proposition}\label{prop:on convex hulls and IcalM}
    Let $I\in \Int(\R^d)$ with $I\subset q_\ast^\uparrow$.
    Then the smallest interval in $\Ical_M$ that contains $I$ is 
    $\lfloor - \rfloor^{-1}\left(\overline{\lfloor I\rfloor} \right)$.
\end{proposition}

\begin{proof} 
First, we claim that $\mathbb{I}:=\lfloor - \rfloor\inv\left(\overline{\lfloor I \rfloor}\right)$ is an element of $\Ical_M$.
Since $\lfloor - \rfloor$ is order-preserving and $I$ is connected, it is clear that $\lfloor I \rfloor\in \Con(Q_M)$. Then, by Remark~\ref{rem:properties of finitely generated modules}~\ref{item:convex hull is smallest interval}, $\overline{\lfloor I \rfloor}\in \Int(Q_M)$. Then by the definition of $\Ical_M$ (Definition~\ref{def:minimal grid}), we have $\mathbb{I}\in \Ical_M$.

    We show that $\mathbb{I}$ is the smallest element of $\Ical_M$ containing $I$.
    Since $I\subset q_\ast^\uparrow$ and $\lfloor-\rfloor$ is a map $q_\ast^\uparrow \to Q_M$,  $\mathbb{I}$ clearly contains $I$.
        Suppose some $\mathbb{J}\in \Ical_M$ contains $I$. Then, by definition, $\mathbb{J} = \lfloor - \rfloor\inv(J)$ for some $J\in \Int(Q_M)$.
    The assumption $\mathbb{J}\supset I$ implies that %
    $J =\lfloor \mathbb{J} \rfloor
    \supset \lfloor I \rfloor$.
    As $J\in \Int(Q_M)$, by Remark~\ref{rem:properties of finitely generated modules}~\ref{item:convex hull is smallest interval}, we have that $J\supset \overline{\lfloor I \rfloor}$, which implies $
    \mathbb{J} \supset \mathbb{I}$, as desired.    
\end{proof}

For $1\leq i\leq d$, let $\pi_i:\R^d\to \R$ be the standard projection onto the $i$-th coordinate, and let $\widehat{\pi}_i:\R^d\to \R^{d-1}$ be the projection onto all \emph{except} the $i$-th coordinate.
Also, for $t\in \R$ and $I\subset \R^d$
let \[I\vert_{\pi_i=t} :=\pi_i^{-1}(t)\cap I=\{x\in I \, : \, \mbox{the } i\mbox{-th coordinate of } x \mbox{ equals } t\}.\]

The following remark follows immediately from the preceding definitions.
\begin{remark}\label{rem:slice projection vs. slice}
For any $t\in \R$, any $I,J\subset \R^d$, and  any $1\leq i\leq d$, \[\widehat{\pi}_i(I\vert_{\pi_i=t}) = \widehat{\pi}_i(J\vert_{\pi_i=t})\iff I\vert_{\pi_i=t} = J\vert_{\pi_i=t}.\]
\end{remark}

\begin{proposition}\label{prop:two unequal slices}
Let $B:=\Pi_{1\leq i\leq d}[a_i,b_i)\subset \R^d$. For a nonempty subset $ U\subset B$, $U=B$ if and only if for all $1\leq i\leq d$ and for all $t_1<t_2\in \pi_i(B)$ \[\widehat{\pi}_i(U\vert_{\pi_i=t_1})= \widehat{\pi}_i(U\vert_{\pi_i= t_2}).\]
\end{proposition}
\begin{proof} %
    The forward direction is direct since, if $U=B$, then for all $1\leq i\leq d$ and for all $t\in \pi_i(B)$, we have \[\widehat{\pi}_i(U\vert_{\pi_i=t})=\prod_{\substack{1\leq j\leq d \\j\neq i}.}[a_j,b_j).\] Let us prove the backward implication.
First, we claim that for all $1\leq i\leq d$ and for all $t\in \pi_i(B)$, $\widehat{\pi}_i(U\vert_{\pi_i=t})$ is nonempty. Indeed, if there exist $1\leq i\leq d$ and $t\in \pi_i(B)$ such that %
$\widehat{\pi}_i(U\vert_{\pi_i=t})$ is empty, 
then by assumption, for \emph{all} $t\in \pi_i(B)$,
$\widehat{\pi}_i(U\vert_{\pi_i=t})$ is empty.
This implies that 
$U\vert_{\pi_i=t}=\emptyset$ for \emph{all} $t\in \pi_i(B)$, which is impossible since $U$ is a nonempty subset of $B$. 

It remains to show that $B\subset U$. Let $(x_1,\ldots,x_d)\in B$.
Since $\widehat{\pi}_1(U\vert_{\pi_1=x_1})$ is nonempty, let us pick $(y_2,\ldots,y_d)$ in it. This implies that $(x_1,y_2\ldots,y_d)\in U$. Then,
\[(x_1,\hat{y}_2,y_3,\ldots,y_d)\in \widehat{\pi}_2(U\vert_{\pi_2=y_2})=\widehat{\pi}_2(U\vert_{\pi_2=x_2}),\]
which implies that $(x_1,x_2,y_3,\ldots,y_d)\in U$. Repeating this argument for the remaining coordinates, we conclude that $(x_1,x_2,\ldots, x_d)\in U$.
\end{proof}

Recall Remark~\ref{rem:properties of finitely generated modules}~\ref{item:partially contained block}. Let \(I\in \Int(\R^d)\setminus \Ical_M\) satisfy $I\subset q_\ast^\uparrow$, and let \(B=\Pi_{1\leq j\leq d}[a_j,b_j)\in \Bcal\) such that $\emptyset \subsetneq B\cap I\subsetneq B$. By the preceding proposition, there exists $1\leq i\leq d$ such that
\begin{equation}\label{eq:intersection slices are differnt}
\widehat{\pi}_i((B\cap I)\vert_{\pi_i=t_1}) \neq \widehat{\pi}_i((B\cap I)\vert_{\pi_i=t_2})
\end{equation}
for some $t_1<t_2\in \pi_i(B)$. We claim that \(\widehat{\pi}_i(I\vert_{\pi_i=t_1}) \neq \widehat{\pi}_i( I\vert_{\pi_i=t_2})\): Note that for $\ell=1,2$, \[(B\cap I)\vert_{\pi_i=t_\ell}=(B\vert_{\pi_i=t_\ell})\cap( I\vert_{\pi_i=t_\ell}),\]
and that $\widehat{\pi}_i$ is injective on $(B\cap I)\vert_{\pi_i=t_\ell}$. Hence, Equation~\eqref{eq:intersection slices are differnt} can be read as
\[\widehat{\pi}_i(I\vert_{\pi_i=t_1})\cap  \prod_{\substack{1\leq j\leq d\\j\neq i}}[a_j,b_j) \neq  \widehat{\pi}_i(I\vert_{\pi_i=t_2})\cap  \prod_{\substack{1\leq j\leq d\\j\neq i}}[a_j,b_j),\]
which implies \(\widehat{\pi}_i(I\vert_{\pi_i=t_1}) \neq \widehat{\pi}_i( I\vert_{\pi_i=t_2})\), as desired. 

In sum, the preceding proposition guarantees the existence of $1\leq i\leq d$ such that for some $t_1<t_2\in \pi_i(B)$,  the set \(\widehat{\pi}_i(I\vert_{\pi_i=t_1})\) is different from \( \widehat{\pi}_i( I\vert_{\pi_i=t_2})\). Without loss of generality, we assume throughout the remainder of this section that $i=1$.
Accordingly, for each \(t\in \mathbb R\), we write
\(
S_t:=\widehat{\pi}_1(I\vert_{\pi_1=t})
\) so that $S_{t_1}\neq S_{t_2}$.
{The following proposition further guarantees the existence of a certain
``critical value'' $t^\ast\in\pi_1(B)$ such that the slice of $I$ on the
hyperplane $\pi_1=t^\ast$ is distinguished from all nearby parallel slices on
one side.}
\begin{proposition}\label{prop:slice unequal to one side}
 Let \(I\in \Int(\R^d)\setminus \Ical_M\) such that $I\subset q_\ast^\uparrow$, and let \(B\in \Bcal\) with $\emptyset \subsetneq B\cap I\subsetneq B$. 
  Let \([a_1,b_1):=\pi_1(B)\). %
  Then at least one of the following holds: 
    \begin{itemize}[wide]
        \item[\textbf{Case 1)}] There exists \(t^\ast\in [a_1,b_1)\) such that, for all sufficiently small \(\epsilon>0\), \(
       S_{t^\ast}
        \neq
        S_{t^\ast+\epsilon}
        \)
        or
        \item[\textbf{Case 2)}] There exists \(t^\ast\in (a_1,b_1)\) such that, for all sufficiently small \(\epsilon>0\),
        \(
        S_{t^\ast}
        \neq
        S_{t^\ast-\epsilon}
        \)
    \end{itemize}
\end{proposition}

\begin{proof}
 Let
    \(s_1<s_2\in \pi_1(B)\) be  such that
      \(
   S_{s_1}
    \neq
    S_{s_2}
    \).  Let $p\in \R^{d-1}$ be contained in only one of $S_{s_1}$ and $S_{s_2}$, and let
    \[
    A_p:=\{s\in \mathbb R : p\in S_s\}.
    \]
    This set includes either $s_1$ or $s_2$ and thus nonempty. We claim that $A_p$ is an interval in $\R$.
    Let $t_1<t_2<t_3$ in $\R$, with $t_1,t_3\in A_p$.
    By definition, this means $(t_1,p), (t_3,p)\in I$, and so by convexity of $I$, $(t_2,p)\in I$, which implies $t_2\in A_p$. 
     Since $A_p$ is a real interval that contains one of $s_1$ and $s_2$, the intersection $A_p\cap [s_1,s_2]$ is a nonempty interval in $\R$.

    Assume \(p\in S_{s_1}\setminus S_{s_2}\). Then \(s_1\in A_p\)
    and \(s_2\notin A_p\), which implies the existence of $t^\ast\in \R$ such that the interval $A_p\cap [s_1,s_2]$ is equal to either $[s_1,t^\ast]$ or $[s_1,t^\ast)$.
    In the former case, we have
    \(p\in S_{t^\ast}\) whereas, 
    for all \(\epsilon\in (0,s_2-t^\ast]\), \(p\notin S_{t^\ast+\epsilon}\). This gives Case 1. 
    In the latter case, we have \(p\notin S_{t^\ast}\) but \(p\in S_{t^\ast-\epsilon}\) for all $\epsilon\in (0,t^\ast-s_1]$, which gives Case~2.

    Assume \(p\in S_{s_2}\setminus S_{s_1}\). Then,  $A_p\cap [s_1,s_2]$ is equal to either $[t^\ast,s_2]$ or $(t^\ast,s_2]$ for some $t^\ast\in \R$.
    The former and latter cases give Case 2 and Case 1, 
    respectively. 
\end{proof}

\paragraph{Deforming intervals.}    Let $I\in \Int(\R^d)$. Let $\CC(I)$ be the set of connected components of the intersections $B\cap I$, where $B$ ranges over $\Bcal$. 
    If $I\subset q_{\ast}^\uparrow$, then $\CC(I)$ is a partition of $I$, which refines or equals the partition $\{B\cap I: B\in \Bcal\}$ of $I$.
    
    For $X,Y\subset \R^d$, we write $X\leq Y$ if there exists $(x,y)\in X\times Y$ such that $x\leq y$. 
    Define $$\CC^2_{\leq}(I) := \{(X,Y)\in \CC(I)\times \CC(I) \, : \, X\leq Y\}.$$

\medskip

The following lemma underlies the proof that a certain deformation of an interval leaves the rank of $M$ over the interval unchanged.

    \begin{proposition}\label{prop:connected components to rank general}
    Let $I,J\in \Int(\R^d)$ with $I\subset q_\ast^\uparrow$, and suppose there exists a map $\Psi:\CC(I)\to \CC(J)$ satisfying the following:%
    \begin{enumerate}[label=(\roman*)]
        \item For each $X\in \CC(I)$, $X$ and $\Psi(X)$ belong to the same unit block. \label{item:Psi condition 1}
        \item For any $(X_1,X_2)\in \CC^2_{\leq}(I)$, the pair $(\Psi(X_1),\Psi(X_2))$ belongs to the set $\CC^2_{\leq}(J)$ (in other words, $\Psi$ induces a map $\CC^2_{\leq}(I)\to \CC^2_{\leq}(J)$.).  \label{item:Psi condition 3}
    \end{enumerate}
    Then $\rk_M(I)\geq \rk_M(J)$.
\end{proposition}
When such a map $\Psi$ exists, we denote the induced map $\CC^2_{\leq}(I)\to \CC^2_{\leq}(J)$ by $\Psi^2$. 

{We remark that, the preceding proposition subsumes monotonicity of the GRI in the following sense: Assume that $I,J\in \Int(\R^d)$ with  $I\subset J$. If $I\not\subset q_\ast^\uparrow$, then we trivially obtain $0=\rk_M(I)=\rk_M(J)$. Assume that $I\subset q_\ast^\uparrow$. Then the map $\Psi:\CC(I)\to \CC(J)$ defined by sending each $X\in \CC(I)$ to the unique $Y\in\CC(J)$ that contains $X$ satisfies the two conditions in the proposition, and thus the proposition gives $\rk_M(I)\geq \rk_M(J)$.
}

\begin{proof}
 It suffices to construct linear maps $\alpha$ and $\beta$ that make the diagram below commute, where $\psi_J$ and $\psi_I$ the canonical limit-to-colimit maps so that $\rk(\psi_J)=\rk_M(J)$ and $\rk(\psi_I)=\rk_M(I)$:
    \begin{equation}\label{eq:alpha beta}\begin{tikzcd}
    \varprojlim M\vert_J \arrow[r, "\psi_J"]\arrow[d,"\alpha"] & \varinjlim M\vert_J \\
    \varprojlim M\vert_I \arrow[r, "\psi_I"] & \varinjlim M\vert_I \arrow[u,"\beta"].
    \end{tikzcd}
    \end{equation}

    We first construct $\alpha$.    Note that, as $I\subset q_\ast^\uparrow$,  $\CC(I)$ is a partition of $I$.
    For any $\mathbf{\ell}:= (\ell_p)_{p\in J}\in \varprojlim M\vert_J$ (cf. Proposition~\ref{prop:limit and colimit}~\ref{item:limit}), let 
    $\alpha(\mathbf{\ell}):= (\alpha(\ell)_p)_{p\in I}$ where $p\in X\in \CC(I)$ and $\alpha(\ell)_p := \ell_{p'}$ for any $p'\in \Psi(X)\in \CC(J)$. %
    We claim that $\alpha(\ell)_p$ is well-defined, i.e., $\ell_p = \ell_{p'}$ for all $p,p'\in \Psi(X)$. This is clear since $\Psi(X)$ is a connected component of $J$ inside a unit block. {We also note that, by the assumption \ref{item:Psi condition 1}, $\alpha(\ell)_p=\ell_{p'}\in M_{p'}=M_p$.}
    
      Next, we claim that $\alpha(\ell)\in \varprojlim M\vert_I$, i.e., for any $p\leq q\in I$, $M(p\leq q)(\alpha(\ell)_p) = \alpha(\ell)_q$. 
     This follows from assumptions \ref{item:Psi condition 1} and \ref{item:Psi condition 3}, together with Remark~\ref{rem:properties of finitely generated modules}~\ref{item:properties of finitely generated modules-maps between blocks}. 
    
    It is clear by construction that $\alpha$ is a linear map.

    Now we now prepare to define the map $\beta$.
    Recall the formulation of the colimit $\varinjlim M\vert_J$ from Proposition \ref{prop:limit and colimit} \ref{item:colimit}.
    We define maps $\delta_p:M_p\to \varinjlim M\vert_J$ for all $p\in I$ as follows.
    Since $\CC(I)$ is a partition of $I$, there exists a unique $X\in \CC(I)$ that contains $p$. 
    Let $p'\in \Psi(X)\subset J$.
    Since $X$ and $\Psi(X)$ are contained in the same unit block, $M_p = M_{p'}$.
    For any $v_p\in M_p$, let $v_{p'}:=v_p\in M_{p'}$.
    Now define $\delta_p(v_p):= [v_{p'}]\in \varinjlim M\vert_J$. 
    {The choice of $p'$ does not affect the well-definedness of $\delta_p$, since $\Psi(X)$ is a connected subposet contained in a unit block on which $M$ is constant (cf.~Proposition~\ref{prop:path characterization of zero in colim}).}
    It is also immediate that $\delta_p$ is  linear. Consequently, the family $(\delta_p)_{p\in P}$ canonically induces a map 
    \[\beta':\bigoplus_{p\in P} M_p\to \varinjlim M\vert_{J}.\]
    
    For $p\in I$, let $\bar{\iota}_p:M_p\to \bigoplus_{p\in I} M_p$ be the canonical injection.
    Recall that $\varinjlim M\vert_I$ is the quotient space $(\bigoplus_{p\in I} M_p)/V$, where $V$ is generated by vectors $\bar{\iota}_p(v_p) - \bar{\iota}_q(v_q)$ over all $p\leq q$ in $I$, with $v_p\in M_p$, $v_q\in M_q$, and $\varphi_M(p,q)(v_p)=v_q$. 
    {By showing that $V\subset \ker(\beta')$, we prove that $\beta'$ induces the linear map
     \[\left(\bigoplus_{p\in I} M_p\right)/V = \varinjlim M\vert_I\to \varinjlim M\vert_J.\] 
    The map $\beta$ is defined to be this map.}
    
    To show $V\subset \ker(\beta')$, it suffices to show that each generator of $V$, $\bar{\iota}_p(v_p) - \bar{\iota}_q(v_q)$ for $p\leq q$ in $I$, is in the kernel.
    Fix $p\leq q$ in $I$, and let $X_p,X_q\in \CC(I)$ such that $p\in X_p$ and $q\in X_q$.
    Since $p\leq q$, by assumption \ref{item:Psi condition 3} we have that $\Psi(X_p)\leq \Psi(X_q)$, so take $p'\in \Psi(X_p)$ and $q'\in \Psi(X_q)$ such that $p'\leq q'$.
    Then, \[\beta'(\bar{\iota}_p(v_p) - \bar{\iota}_q(v_q)) = \beta'(\bar{\iota}_p(v_p)) - \beta'(\bar{\iota}_q(v_q)) = [v_{p'}] - [v_{q'}].\]
    {
    By Remark~\ref{rem:properties of finitely generated modules}~\ref{item:properties of finitely generated modules-maps between blocks}, $\varphi_M(p,q):M_p\to M_q$ and $\varphi_M(p',q'):M_{p'}\to M_{q'}$ are identical maps.}
    Further, since $v_p=v_{p'}\in M_p=M_{p'}$ and $v_q=v_{q'}\in M_q=M_{q'}$, the assumption that $\varphi_M(p,q)(v_p)=v_q$ implies that $\varphi_M(p',q')(v_{p'})=v_{q'}$.
    Hence, $[v_{p'}] - [v_{q'}] = 0$, and thus $V\subseteq \ker(\beta')$, as desired.

    It remains to show that Diagram~\eqref{eq:alpha beta} commutes.
    Let $\mathbf{\ell}\in \varprojlim M\vert_J$, and let $p\in X\in \CC(I)$ and $p'\in Y=\Psi(X)$.
    By definitions, we have that $\psi_J(\mathbf{\ell})$ = $[\ell_{p'}]$, and that
    \[\beta(\psi_I(\alpha(\mathbf{\ell}))) = \beta([\alpha(\ell)_{p}]) = [\ell_{p'}],\]
   as desired.
\end{proof}

\begin{corollary}\label{cor:connected components to rank}
    Let $I, J\in \Int(\R^d)$ with $I,J\subset q_\ast^\uparrow$, and suppose there exists a map $\Psi:\CC(I)\to \CC(J)$ satisfying all the conditions in Proposition~\ref{prop:connected components to rank general}.
    If both $\Psi$ and $\Psi^2$ are bijections, then $\rk_M(I) = \rk_M(J)$.
\end{corollary}

Next, we will construct certain parametrized families of maps on $\R$ and $\R^d$. 
We use them later to construct families of intervals. 
To this end, fix some $t^\ast\in \R$ and $\epsilon^\ast>0$.
For each $\epsilon\in [0,\epsilon^\ast)$, we consider the order-preserving piecewise linear functions
$g_\epsilon,f_\epsilon:\mathbb{R}\to\mathbb{R}$ whose graphs are depicted as in Figure~\ref{fig:g epsilon and f epsilon}.
Note that $g_0$ and $f_0$ are both defined to be the identity map on $\R$.
Now we define the maps $G_\epsilon, F_\epsilon:\R^d\to \R^d$ via:
\begin{align*}
    G_\epsilon((x_1,x_2,\ldots,x_d)):= (g_\epsilon(x_1),x_2,\ldots,x_d) \ \mbox{and } 
    F_\epsilon((x_1,x_2,\ldots,x_d)):= (f_\epsilon(x_1),x_2,\ldots,x_d).
\end{align*}

\begin{figure}
\begin{subfigure}[c]{0.4\textwidth}
\begin{center}
\begin{tikzpicture}

\begin{axis}[
    axis lines=middle,
    axis equal image,
    xmin=-0.5, xmax=5.5,
    ymin=-0.5, ymax=5.5,
    xlabel={$t$},
    ylabel={$g_\epsilon(t)$},
    ylabel style={at={(axis description cs:0.5,0.9)}},
    ytick={1, 4},
    yticklabels={$t^\ast$,$t^\ast+\epsilon^\ast$},
    xtick={1,2,4},
    xticklabels={$t^\ast$,$t^\ast+\epsilon$,$t^\ast+\epsilon^\ast$},
    tick label style={font=\small},
]

\addplot[thick] coordinates {(-0.5,-0.5) (1,1)};

\addplot[thick] coordinates {(1,1) (2,1)};

\addplot[thick,domain=2:4,samples=2]
    {1 + 3/2*(x-2)};

\addplot[thick] coordinates {(4,4) (5.5,5.5)};

\addplot[only marks,mark=*]
coordinates {(1,1) (2,1) (4,4)};

\addplot[dashed] coordinates {(1,0) (1,6)};
\addplot[dashed] coordinates {(2,0) (2,5)};
\addplot[dashed] coordinates {(4,0) (4,5)};
\addplot[dashed] coordinates {(0,1) (5,1)};
\addplot[dashed] coordinates {(0,4) (5,4)};
\end{axis}

\end{tikzpicture}
\end{center}
\end{subfigure}\hspace{15mm}%
\begin{subfigure}[c]{0.4\textwidth}
\begin{center}
\begin{tikzpicture}

\begin{axis}[
    axis lines=middle,
    axis equal image,
    xmin=-0.5, xmax=5.5,
    ymin=-0.5, ymax=5.5,
    xlabel={$t$},
    ylabel={$f_\epsilon(t)$},
    ylabel style={at={(axis description cs:0.4,0.9)}},
    ytick={1,2,4},
    yticklabels={$t^\ast$,$t^\ast+\epsilon$,$t^\ast+\epsilon^\ast$},
    xtick={1,4},
    xticklabels={$t^\ast$,$t^\ast+\epsilon^\ast$},
    tick label style={font=\small},
]

\addplot[thick] coordinates {(-0.5,-0.5) (1,1)};

\addplot[thick,domain=1:4,samples=2]
    {2 + (2/3)*(x-1)};

\addplot[thick] coordinates {(4,4) (5.5,5.5)};

\addplot[only marks,mark=*]
coordinates {(1,1) (4,4)};

\addplot[only marks,mark=o]
coordinates {(1,2) (4,4)};

\addplot[dashed] coordinates {(1,0) (1,6)};
\addplot[dashed] coordinates {(4,0) (4,5)};
\addplot[dashed] coordinates {(0,1) (5,1)};
\addplot[dashed] coordinates {(0,2) (5,2)};
\addplot[dashed] coordinates {(0,4) (5,4)};
\end{axis}

\end{tikzpicture}

\end{center}
\end{subfigure}
\caption{
The graphs of the functions $g_\epsilon,f_{\epsilon}:\R\to\R$. %
We note that the composition
$g_\epsilon \circ f_\epsilon$ is the identity on $\R$, while $f_\epsilon\circ g_\epsilon$ is the identity on $\R\setminus (t^\ast, t^\ast+\epsilon]$.
}
\label{fig:g epsilon and f epsilon}
\end{figure}

The following remark and propositions contain a sequence of observations that will be used later.

As $g_\eps\circ f_\eps$ is the identity map on $\R$, we immediately have:

\begin{remark}\label{rem:GF is identity}
\begin{enumerate}[label=(\roman*)]
\item $G_\eps\circ F_\eps$ is the identity map on $\R^d$. \label{item:GF is identity1}
\item For any $I\in \Int(\R^d)$, we have $I=(G_\eps\circ F_{\eps})^{-1}(I)=F_\eps^{-1}\left(G_\eps^{-1}(I)\right)$. \label{item:GF is identity2}
\end{enumerate}
\end{remark}

\begin{proposition}\label{prop:Feps and Geps maps B to B} Let $B:=\prod_{i=1}^d[a_i,b_i)\subset \R^d$ with $[t^\ast,t^\ast+\eps^\ast]\subset [a_1,b_1)$. 
Then, the images of $B$ under the maps $G_\eps$ and $F_\eps$ are contained in $B$.
\end{proposition} 

\begin{proof} 
Let $x\in B$ and let $x_1$ be the first coordinate of $x$. It suffices to show that the images of $x_1$ under the maps $g_\eps$ and $f_\eps$ are in $[a_1,b_1)$. 
This follows from the fact that $x_1\in [a_1,b_1)$ and $[t^\ast, t^\ast+\epsilon^\ast]\subset [a_1,b_1)$, together with the observation that both $g_\epsilon$ and $f_\epsilon$ map $[t^\ast, t^\ast+\epsilon^\ast]$ into itself and fix every point outside $[t^\ast, t^\ast+\epsilon^\ast]$.
\end{proof}

\begin{proposition}\label{prop:G, F preserve qast}
Let $p,q_\ast\in \R^d$ and assume that the first coordinate of $q_\ast$ is at most $t^\ast$.
Then the following are equivalent:
\begin{enumerate*}[label=(\roman*)]
\item $p\geq q_\ast$ \label{item:p}
\item $G_\epsilon(p)\geq q_\ast$  \label{item:Geps p}
\item $F_\epsilon(p)\geq q_\ast$\label{item:Feps p}
\end{enumerate*} 
\end{proposition}

\begin{proof}Observe that the images and preimages of $[t^\ast,\infty)$ under $g_\eps$ and $f_{\eps}$ are contained in $[t^\ast,\infty)$ and that, on the interval $(-\infty,t^\ast)$, $g_\eps$ and $f_\eps$ are the identity map. These observations imply that, for $t\in \R$ and the first coordinate $q_1$ of $q_\ast$, the following are equivalent: $t\geq q_1$, $g_\eps(t)\geq q_1$, and $f_\eps(t)\geq q_1$. This immediately implies the desired statement.
\end{proof}

For $J\subset \R^d$, for ease of notation, we write
\[
J_\epsilon^G:=G_\epsilon^{-1}(J)\qquad\mbox{and}
\qquad
J_\epsilon^F:=F_{\epsilon}^{-1}(J).
\]

\begin{proposition}\label{prop:Ie is an interval for all e}
    For $J\in \Int(\R^d)$, 
    $J_\epsilon^G$ is an interval. 
\end{proposition}
\begin{proof}
Since $G_\eps$ is surjective, $J_\epsilon^G$ is nonempty. By Proposition~\ref{prop:each component is an interval}, $J_\epsilon^G$ is convex.
    It remains to show that $J_\epsilon^G$ is connected.
    Let $p,q\in J_\epsilon^G$. Then $p':=G_\epsilon(p)$ and $q':=G_\epsilon(q)$ belong to $J$.
    As $J$ is connected, there exists a path %
    in $J$ connecting $p'$ and $q'$.
      As $F_\epsilon$ is order-preserving, the image of the path under $F_\epsilon$
    is a path in $\R^d$ connecting $F_\eps(p')$ and $F_\eps(q')$. Furthermore, observing that if $r\in J$, then $F_\epsilon(r)\in J_\epsilon^G$ as $G_\epsilon(F_\epsilon(r)) = r$ (cf. Remark~\ref{rem:GF is identity}~\ref{item:GF is identity1}), the image path is actually in $J_\epsilon^G$. 
    Now we claim that prepending $p$ and appending $q$ to the image path yields a path in $J_\epsilon^G$ connecting $p$ and $q$, as desired.  Indeed, as $p$ and $F_\epsilon(p')$ 
    can only differ in the first coordinate,
    they are comparable. Likewise,
    $q$ and $F_\epsilon(q')$ are comparable.
\end{proof}

Since $J_\eps^G$ is an interval, $\rk_M(J_\eps^G)$ is well-defined. Now we claim:

\begin{proposition}\label{prop:ranks of I and Ie are equal}
    For 
    $J\in \Int(\R^d)$ with $J\subset q_\ast^\uparrow$, we have $\rk_M(J) = \rk_M(J_\epsilon^G)$. 
\end{proposition}

\begin{proof}
We prove the statement by utilizing Corollary~\ref{cor:connected components to rank}. 
To this end, we need to show that $J_\epsilon^G$ is contained in $q_\ast^\uparrow$. This directly follows from Proposition~\ref{prop:G, F preserve qast} and the assumption $J\subset q_\ast^\uparrow$.   
Next, we define the map
\begin{align*}
\Psi:\CC(J_\epsilon^G)&\to \CC(J)\\
X&\mapsto Y,
\end{align*}
where $Y$ contains $G_\epsilon(x)$ for some $x\in X$. We will show that $\Psi$ satisfies the assumptions of Corollary~\ref{cor:connected components to rank}.
  First we verify that $\Psi$ is well-defined, i.e., not affected by the choice of $x\in X$.
    Let $x'\in X$ that is distinct from $x$. Since $X$ is connected, there exists a path %
    in $X$ connecting $x$ and $x'$.  
    As $G_\epsilon$ is order-preserving, the image of the path under $G_\epsilon$ is a path    %
    connecting $G_\epsilon(x)$ and $G_\epsilon(x')$ in $J$.
    Let $B\in \Bcal$ be the unit block containing $X$.
    Then Proposition~\ref{prop:Feps and Geps maps B to B} guarantees that the image path is also in $B$.
    Thus, $G_\epsilon(x)$ and $G_\epsilon(x')$ are in the same element of $\CC(J)$, as desired.

    Proposition~\ref{prop:Feps and Geps maps B to B} also guarantees that both $X$ and $Y$ belong to the same unit block $B$, and thus $\Psi$ satisfies assumption~\ref{item:Psi condition 1} in Proposition~\ref{prop:connected components to rank general}. 
    Next, we claim that $\Psi$ satisfies assumption~\ref{item:Psi condition 3} in the same proposition.
    Suppose $(X_1,X_2)\in \CC^2_{\leq}(J_\epsilon^G)$, i.e., there exists
    a pair $x_1\leq x_2$ with $x_1\in X_1$ and $x_2\in X_2$.
    Then $G_\epsilon(x_1)\leq G_\epsilon(x_2)$, $G_\epsilon(x_1)\in \Psi(X_1)$, and $G_\epsilon(x_2)\in \Psi(X_2)$, so $(\Psi(X_1),\Psi(X_2))\in \CC^2_{\leq}(J)$, as desired.

    Next, we claim $\Psi$ is injective.
   Suppose $X_1,X_2\in \CC(J_\epsilon^G)$ such that $\Psi(X_1) = \Psi(X_2)=:Y$. This implies the following two facts: (1) As $\Psi$ satisfies assumption~\ref{item:Psi condition 1}, $X_1$ and $X_2$ are contained in the same unit block, say $B$. (2) For any $x_1\in X_1$ and $x_2\in X_2$,  $G_\epsilon(x_1),G_\epsilon(x_2)\in Y\in \CC(J)$.
    Thus, there exists a path in $J$ %
    connecting $G_\epsilon(x_1)$ to $G_\epsilon(x_2)$.
    As we saw in the proof of Proposition~\ref{prop:Ie is an interval for all e}, applying $F_\eps$ to this path induces a path connecting $x_1$ and $x_2$ in $J_\epsilon^G$. In particular, by Proposition~\ref{prop:Feps and Geps maps B to B}, the induced path is in $B$, implying that $X_1=X_2$.

    Next, we claim $\Psi$ is surjective.
    Suppose $Y\in \CC(J)$, and let $y\in Y$.
    As $G_\epsilon$ is surjective, there exists $x\in \R^d$ with $G_\epsilon(x) = y$.
    By definition, this means $x\in J_\epsilon^G$, so $x\in X$ for some $X\in \CC(J_\epsilon^G)$, and then by definition of $\Psi$, this implies $\Psi(X) = Y$. 

    Lastly, we claim that the induced map $\Psi^2:\CC^2_{\leq}(J_\epsilon^G)\to \CC^2_{\leq}(J)$ is bijective. Since $\Psi$ is injective, the injectivity of $\Psi^2$ is readily clear. To prove surjectivity, let $(Y_1,Y_2)\in \CC^2_{\leq}(J)$, and let $y_1\in Y_1, y_2\in Y_2$ such that $y_1\leq y_2$.
    Define $x_1:=F_\epsilon(y_1)$ and $x_2:=F_\epsilon(y_2)$.
    Since $G_\epsilon\circ F_\epsilon$ is the identity on $\R^d$ (cf. Remark~\ref{rem:GF is identity}~\ref{item:GF is identity1}), we have $G_\epsilon(x_1)=y_1$ and $G_\epsilon(x_2)=y_2$, implying that $x_1,x_2\in J_\epsilon^G$.
    Then there exist $X_1, X_2\in \CC(J_\epsilon^G)$ that contain $x_1$ and $x_2$, respectively.
    Further, as $F_\epsilon$ is order-preserving, we have $x_1\leq x_2$, and thus $(X_1,X_2)\in \CC^2_{\leq}(J_\epsilon^G)$. Then, we have $\Psi^2((X_1,X_2)) = (\Psi(X_1),\Psi(X_2)) = (Y_1,Y_2)$, as desired.
\end{proof}

We now proceed to the proof of Theorem \ref{thm:unique candidate}. 

\paragraph{Proof of Theorem \ref{thm:unique candidate}.}
Let $\Ical:= \supp(\dgm^\Int_M)$.  
We first show that proving $\Ical\subset \Ical_M$ suffices to establish Theorem~\ref{thm:unique candidate}:

\resetclaimcounter

    \begin{claim}\label{claim:Ical subset implies equal dgm}
        If $\Ical\subset \Ical_M$, then $\dgm^\Int_M = \dgm_M$.
    \end{claim}
    \begin{claimproof}
        By definition, the restrictions $\dgm_M\vert_{\Ical_M}$ and $\dgm^\Int_M\vert_{\Ical_M}$ are both the M\"obius inversion of $\rk^{\Ical_M}_M$ over the finite poset $(\Ical_M,\supset)$.
        Hence, by the uniqueness of M\"obius inversion, $\dgm^\Int_M\vert_{\Ical_M} = \dgm_M\vert_{\Ical_M}$.
        The assumption $\Ical\subset \Ical_M$ implies that $\dgm^\Int_M = \dgm_M$ on $\Int(\R^d)$.
    \end{claimproof}
Thanks to this claim, our next goal is to show that every element of $\Ical$ belongs to $\Ical_M$, or equivalently, is subordinate to the minimal grid $Q_M$. We begin by making a basic observation about the intervals in $\Ical$. Recall that $q_\ast$ is the minimal point of $Q_M$.

        \begin{claim}\label{claim:positive I is in the upset generated by the minimal point}
If $I\in \Ical$, then $I\subset q_\ast^\uparrow$.
    \end{claim}

    \begin{claimproof}
        If an interval $I$ contains a point $p\in \R^d$ outside $q_\ast^\uparrow$, then $M_p = 0$. Hence, monotonicity of the GRI implies $\rk_M(I) = 0$, which in turn implies by Proposition \ref{prop:properties of GGPD} \ref{item:GGPD is more parsimonious} that $\dgm^\Int_M(I) = 0$, i.e., $I\notin \Ical$.
    \end{claimproof}
\textbf{For the remainder of the proof, we assume that $\Ical\not \subset \Ical_M$, and show that this yields a contradiction with the assumption of Theorem~\ref{thm:unique candidate} that $\rk_M^\Int$ is M\"obius invertible.} 
See Claim~\ref{claim:rkm is not mobius invertible} and its subsequent paragraph.

In the next two claims, we show that any maximal interval in $\Ical\setminus \Ical_M$ has positive $\dgm^\Int_M$ and such an interval is \emph{rank-critical} (see Definition~\ref{def:rank-critical} below).

    \begin{claim}\label{claim:negative dgm interval implies larger positive}
     Any maximal interval in $\Ical\setminus \Ical_M$ has positive $\dgm^\Int_M$.
    \end{claim}
    \begin{claimproof}
    It suffices to show that, assuming that $\dgm^\Int_M(I)<0$ for some $I\in \Ical\setminus \Ical_M$, there exists an element of $\Ical\setminus \Ical_M$ that properly contains $I$ and its $\dgm^\Int_M$ value is positive.
    
    If $I$ is maximal in $\Ical$, then $\rk_M(I) = \dgm^\Int_M(I) < 0$, a contradiction. 
    Thus, there exists $J\in \Ical$ that properly contains $I$ and has positive $\dgm^\Int_M$. 
    Let $J^\ast$ be the connected component of $\bigcap \Ical^{+}_{\supset I}$ that contains $I$. Since  $\bigcap \Ical^{+}_{\supset I}$ is convex, $J^\ast$ is an interval. 

We claim that  $J^\ast = I$.      We have:
        \begin{align*}
        \rk_M(J^\ast) =  \sum_{\substack{K\supset J^\ast\\K\in \Ical}} \dgm^\Int_M(K) = & \ \sum_{\substack{K\supset I\\K\in \Ical}} \dgm^\Int_M(K) - \sum_{\substack{K\supset I, K\not\supset J^\ast \\K\in \Ical}} \dgm^\Int_M(K)\\ 
        = & \ \rk_M(I) - \sum_{\substack{K\supset I, K\not\supset J^\ast \\K\in \Ical}} \dgm^\Int_M(K) . 
        \end{align*}
        As $J^\ast$ contains $I$, monotonicity of the GRI implies that 
        \[\sum_{\substack{K\supset I, K\not\supset J^\ast \\K\in \Ical}} \dgm^\Int_M(K) \geq 0. \]
        However, by construction of $J^\ast$, any $K\in \Ical$ that contains $I$ and does not contain $J^\ast$ must have negative $\dgm^\Int_M$.
       Thus, the sum above must be empty. However, if $I\subsetneq J^\ast$, then $\dgm_M^\Int(I)$ appears as a summand in the sum. This implies that $J^\ast=I$.

    Now, by Claim~\ref{claim:positive I is in the upset generated by the minimal point}, $I\subset q_\ast^\uparrow$.
    Further, as $I\notin \Ical_M$, by Remark~\ref{rem:properties of finitely generated modules} \ref{item:partially contained block}, there exists a unit block $B$ such that $\emptyset \subsetneq B\cap I\subsetneq B$. We fix such a block $B$.
 
    We claim that there  exists $J\in \Ical^{+}_{\supset I}$ such that  $\emptyset\subsetneq B\cap J\subsetneq B$. Since every $J\in \Ical^{+}_{\supset I}$ contains $I$, we have that $\emptyset\subsetneq B\cap J$. Suppose by way of contradiction that each $J\in \Ical^{+}_{\supset I}$ contains $B$, i.e., $B\subset \bigcap \Ical^{+}_{\supset I}$.
    As $B$ is connected, $B$ must be inside a single connected component of $ \bigcap \Ical^{+}_{\supset I}$, and by the assumption $B\cap I\neq \emptyset$, this implies $B\subset J^\ast$ or equivalently
     $B\cap J^\ast = B$. However, since $J^\ast=I$, this is a contradiction to the assumption $B\cap I \subsetneq B$.
    
    Now fix any $J\in \Ical^{+}_{\supset I}$  such that $\emptyset \subsetneq B\cap J\subsetneq B$.
    By Claim~\ref{claim:positive I is in the upset generated by the minimal point}, $J\subset q_\ast^\uparrow$ and thus   Remark~\ref{rem:properties of finitely generated modules} \ref{item:partially contained block} implies that  $J\notin \Ical_M$.
    By construction, $\dgm^\Int_M(J) > 0$ with $J\supsetneq I$, as desired.
    \end{claimproof}

We introduce a useful terminology. 
    \begin{definition}\label{def:rank-critical} We call $I\in \Int(\R^d)$ \textbf{rank-critical} (with respect to $M$) if for all $J\in \Int(\R^d)$ properly containing $I$, $\rk_M(J)$ is strictly less than $\rk_M(I)$.
    \end{definition}
    
    \begin{claim}\label{claim:rank-critical}
        Any maximal interval in $\Ical\setminus \Ical_M$ is rank-critical.
    \end{claim}
    \begin{claimproof}
       Let $\Iast\in \Ical\setminus \Ical_M$ be a maximal interval. 
        Since $\dgm^\Int_M(\Iast)\neq 0$, Proposition \ref{prop:properties of GPD} \ref{item:GPD is more parsimonious} implies $\rk_M(\Iast) > 0$. Assume by way of contradiction that there exists some $K\in \Int(\R^d)$ with $\Iast \subsetneq K$ and $\rk_M(K) = \rk_M(\Iast)$. Then, we have $\rk_M(K)>0$ and
        therefore $\Iast\subsetneq K\subset q_\ast^\uparrow$.
        
        As $\Iast\subsetneq K$ and $L$ is maximal in $\Ical\setminus \Ical_M$, any $J\in \Ical$ containing $K$ belongs to $\Ical_M$. Thus, by Proposition~\ref{prop:on convex hulls and IcalM}, 
      $J$ contains $\mathbb{K}:=\lfloor - \rfloor\inv\left(\overline{\lfloor K \rfloor}\right)$. 
        This implies the equality $\{J\in \Ical \, : \, J\supset K\} = \{J\in \Ical \, : \, J\supset \mathbb{K}\}$ and thus
        \[
        \sum_{\substack{J\supset K\\ J\in \Ical}} \dgm^\Int_M(J) = \sum_{\substack{J\supset\mathbb{K}\\ J\in \Ical}} \dgm^\Int_M(J).\] 
        Since the LHS equals $\rk_M(K)$  and the RHS equals  $\rk_M\left(\mathbb{K}\right)$, we have that
        \[\rk_M(\Iast) = \rk_M(K) = \rk_M\left(\mathbb{K}\right).\]
        Then, since $\Iast\subset \mathbb{L}:=\lfloor - \rfloor\inv\left(\overline{\lfloor \Iast \rfloor}\right)\subset \mathbb{K}$,
        monotonicity of the GRI implies \[\rk_M(\Iast) = \rk_M(\mathbb{\Iast}).\]
        Then, by the definition of $\dgm_M^\Int$, the LHS and RHS are respectively equal to \[\sum_{\substack{J\supset \Iast\\J\in \Ical}} \dgm^\Int_M(J) \qquad\mbox{ and } \sum_{\substack{J\supset \mathbb{\Iast}\\ J\in \Ical}}\dgm_M^\Int(J),\]
        and thus \[\sum_{\substack{J\supset \Iast,\ J\not\supset \mathbb{\Iast}\\J\in \Ical}} \dgm^\Int_M(J)=0.\]
        We claim that the preceding sum  consists solely of the single summand $\dgm_M^\Int(L)$, which leads to the equality $\dgm_M^\Int(L)=0$, a contradiction. To see the claim, note that, by maximality of $\Iast\in \Ical\setminus \Ical_M$ and Proposition~\ref{prop:on convex hulls and IcalM},
        \begin{align*}
        \{J\in \Ical \, : \, J\supset \Iast\} &= \{J\in \Ical\setminus  \Ical_M \, : \, J\supset \Iast\}&\cup&\qquad\{J\in \Ical\cap\Ical_M \, : \, J\supset \Iast\}& \\&=\{\Iast\}&\cup&\qquad\{J\in \Ical\cap\Ical_M \, : \, J\supset \mathbb{\Iast}\}.&
        \end{align*}
        {Note also that $\Iast\neq \mathbb{\Iast}$ because  $L\notin\Ical_M$ whereas $\mathbb{\Iast}\in \Ical_M$.
         Consequently, the only $J$ satisfying the condition to be a summand in the preceding sum is $L$, completing the proof.}  
    \end{claimproof}

Let $\Iast\in \Ical\setminus \Ical_M$ be a maximal and thus rank-critical interval. Then $L\subset q_\ast^\uparrow$ by Claim~\ref{claim:positive I is in the upset generated by the minimal point}.
Also, as $L\notin \Ical_M$, by Remark~\ref{rem:properties of finitely generated modules}~\ref{item:partially contained block}, there exists a unit block $B\in \Bcal$ with $\emptyset\subsetneq B\cap L\subsetneq B$. Hence, $L$ satisfies the assumption of Proposition \ref{prop:slice unequal to one side}.
Assume that Case~1 in the conclusion of Proposition~\ref{prop:slice unequal to one side} holds. This means that,
letting   \[
    S_t:=\widehat{\pi}_1(L|_{\pi_1=t}),
    \]
there exist $t^\ast\in [a_1,b_1):=\pi_1(B)$ and $\eps^\ast>0$ so that 
\begin{equation}\label{eq:closed interval around tast}
[t^\ast,t^{\ast}+\eps^\ast]\subset [a_1,b_1)
\end{equation}
and
\begin{equation}\label{eq:slices are different}
\mbox{for all $\eps \in (0, \epsilon^\ast]$,}\qquad S_{t^\ast} \neq  S_{t^\ast+\epsilon}.
\end{equation}

\begin{remark}\label{rem:q1 is at most tast}
By construction, $B\subset q_\ast^\uparrow$, and thus \[\mbox{(The first coordinate of $q_\ast$)}\leq a_1 \leq t^\ast.\]
\end{remark}
\textbf{From Claim~\ref{claim:pairwise difference} through Claim~\ref{claim:rkm is not mobius invertible} below, we continue to use the same $t^\ast$ and $\epsilon^\ast$ as above.}

Our next goal is to show that $\{\Iast_\epsilon^G\}_{\eps\in [0,\eps^\ast)}$
is an infinite family of intervals without pairwise containment. 
Later, this family will be utilized to show the existence of an infinite subcollection $\Jcal\subset\Ical$ with nontrivial intersection, which contradicts M\"obius invertibility of $\rk_M^\Int$.

First, by Proposition~\ref{prop:Ie is an interval for all e}, we have that $\Iast_\epsilon^G$ is an interval for each $\eps$.

\begin{claim}\label{claim:pairwise difference}
For $\epsilon<\epsilon'\in [0,\epsilon^\ast)$, we have that    $\Iast_{\epsilon}^G\neq \Iast_{\epsilon'}^G$.
\end{claim}
We remark that the maximality or rank-criticality of $L$ does not play any role in the following proof.
\begin{claimproof} 
   From Figure~\ref{fig:g epsilon and f epsilon}, observe that %
   \begin{equation}\label{eq:strictly increasing}
   t^\ast =g_{\eps'}(t^\ast+\eps') < g_{\epsilon}(t^\ast+\epsilon')\leq t^\ast + \epsilon^\ast.\end{equation} %
   Also, observe that, by the definitions of $L_\epsilon^G$ and $S_{g_\eps(t)}$,
    \begin{equation}\label{eq:some equation}
\mbox{ for all $t\in \R$,  } \qquad\widehat{\pi}_1({L_\epsilon^G}\vert_{\pi_1=t}) = S_{g_\eps(t)}.
\end{equation}
    By Equations~\eqref{eq:slices are different}, \eqref{eq:strictly increasing}, and \eqref{eq:some equation},
    \[\widehat{\pi}_1(\Iast_{\epsilon'}^G\vert_{\pi_1=t^\ast+\epsilon'}) = S_{g_{\eps'}(t^\ast+\eps')}=S_{t^\ast} \neq S_{g_{\epsilon}(t^\ast+\epsilon')} = \widehat{\pi}_1(\Iast_{\epsilon}^G\vert_{\pi_1=t^\ast+\epsilon'}) .\]
    {This implies $\Iast_{\epsilon'}^G\vert_{\pi_1=t^\ast+\epsilon'}\neq \Iast_{\epsilon}^G\vert_{\pi_1=t^\ast+\epsilon'}$, and thus} 
    we have $\Iast_{\epsilon}^G\neq \Iast_{\epsilon'}^G$.\end{claimproof} %

\begin{claim}\label{claim:slice condition implications IFe}
    Let $I\in \Int(\R^d)$, and $\epsilon \geq 0$ such that for all $t\in [t^\ast, t^\ast + \epsilon]$ it holds that $\widehat{\pi}_1(I\vert_{\pi_1=t}) = \widehat{\pi}_1(I\vert_{\pi_1=t^\ast})$.
    Then we have the following:
    \begin{enumerate}[label=(\roman*)]
            \item $I = (I_{\epsilon}^F)_{\epsilon}^G$ (note that, by Remark~\ref{rem:GF is identity}~\ref{item:GF is identity2}, we readily know $I = (I_{\epsilon}^G)_{\epsilon}^F$). \label{item:F then G equals original}  
        \item $I_{\epsilon}^F$ is an interval. 
        \label{item:IFe is an interval} 
    \end{enumerate}
\end{claim}
Before proving the claim, it will be useful to note that Item \ref{item:F then G equals original} implies the following. 
\begin{equation}\label{eq:x in I then Gex in}
   \mbox{If $x\in I$, then $G_\epsilon(x)\in I_{\epsilon}^F$}. 
\end{equation}

\begin{claimproof}
\noindent    \ref{item:F then G equals original}:
    First, we show that $I\subset (I_{\epsilon}^F)_{\epsilon}^G$.
    It suffices to show that, if $x\in I$, then $F_\eps ( G_{\eps}(x))\in I$. Let $x =: (x_1,x')\in I$ with $x_1\in \R$ and $x'\in \R^{d-1}$. 
    Observe that if $x_1\in [t^\ast, t^\ast+\epsilon]$, then $f_{\epsilon}(g_{\epsilon}(x_1)) = t^\ast$, and if not, then $f_{\epsilon}(g_{\epsilon}(x_1))=x_1$ (cf. the caption of Figure~\ref{fig:g epsilon and f epsilon}).
    Hence, we have
    \begin{equation}\label{eq:FG cases}
    F_{\epsilon}(G_{\epsilon}(x)) = 
    \begin{cases}
    (t^\ast,x') & \mbox{if } x_1\in [t^\ast, t^\ast+\epsilon],\\
        x & \mbox{otherwise.}
    \end{cases}
    \end{equation}
    Thus, if $x_1\notin [t^\ast, t^\ast+\epsilon]$, then $F_{\epsilon}(G_{\epsilon}(x)) = x\in I$ as desired.
    Now assume that $x_1\in [t^\ast, t^\ast+\epsilon]$. Since $x\in I$, the point $x'$ belongs to the set $ \widehat{\pi}_1(I\vert_{\pi_1=x_1})$ which is, by %
    assumption, %
    identical to $ \widehat{\pi}_1(I\vert_{\pi_1=t^\ast})$.
    Thus, 
    $(t^\ast,x')\in I$.
    By Equation \eqref{eq:FG cases}, $F_\epsilon(G_\epsilon(x)) = (t^\ast,x')$, and thus $F_\epsilon(G_\epsilon(x))\in I$. 
    
     Next, we show that $ (I_{\epsilon}^F)_{\epsilon}^G\subset I$. 
     Suppose $x\in (I_{\epsilon}^F)_{\epsilon}^G$, i.e., there exists $y\in I$ such that $F_\eps(G_\eps(x))=y$. 
     By Equation \eqref{eq:FG cases}, if $x_1\notin [t^\ast, t^\ast + \epsilon]$, then $F_\epsilon(G_\epsilon(x)) = x$, which implies $y=x$, and thus $x\in I$ since $y$ is.    If $x_1\in [t^\ast, t^\ast+\epsilon]$, then $y = F_\epsilon(G_\epsilon(x)) = (t^\ast, x')$.    Since $y\in I$, the point $x'$ belongs to the set $ \widehat{\pi}_1(I\vert_{\pi_1=t^\ast})$, which is by our assumption identical to $\widehat{\pi}_1(I\vert_{\pi_1=x_1})$. This shows that $x=(x_1,x')\in I$. %
   
\noindent   \ref{item:IFe is an interval}:
    First, we show that $I_{\epsilon}^F\neq \emptyset$.
    As $I\in \Int(\R^d)$, $I\neq \emptyset$, so pick any $x=(x_1,x')\in I$ with $x_1\in \R$ and $x'\in \R^{d-1}$.
    Suppose
    $x_1\in \R\setminus [t^\ast, t^\ast+\epsilon]$.
    Then $f_\epsilon(x) = x$, so $x\in I_{\epsilon}^F$, implying that $I_{\epsilon}^F\neq \emptyset$.
       Suppose 
    $x_1\in [t^\ast, t^\ast + \epsilon]$.
    Then we have $x'\in \widehat{\pi}_1(I\vert_{\pi_1=x_1})=\widehat{\pi}_1(I\vert_{\pi_1=t^\ast})$, by our assumption. 
    Thus, $(t^\ast, x')\in I$.
    As $f_{\epsilon}(t^\ast) = t^\ast$, we have that $F_\epsilon((t^\ast,x')) = (t^\ast, x')\in I$, and thus $(t^\ast, x')\in I_{\epsilon}^F$. This proves that $I_{\epsilon}^F\neq \emptyset$.

    By Proposition \ref{prop:each component is an interval}, $I_{\epsilon}^F$ is convex, and it remains to show connectivity.
    Let $p,q\in I_{\epsilon}^F$, which implies that  $F_\epsilon(p), F_\epsilon(q)\in I$.
    As $I$ is an interval, there is a path in $I$ connecting $F_\epsilon(p)$ and $F_\epsilon(q)$.
    By the fact that $G_\epsilon$ is order-preserving, alongside Equation~\eqref{eq:x in I then Gex in}, the image of this path under $G_\epsilon$ 
    is a path in $I_{\epsilon}^F$ 
    connecting $G_\epsilon(F_\epsilon(p)) = p$ to $G_\epsilon(F_\epsilon(q)) = q$ (cf. Remark~\ref{rem:GF is identity}~\ref{item:GF is identity1}), as desired.
\end{claimproof}

\begin{claim}\label{claim:Ie distinct for different e}
    For $\epsilon<\epsilon'\in [0,\epsilon^\ast)$, we have that  $\Iast_{\epsilon}^G\not\subset \Iast_{\epsilon'}^G$.
\end{claim}
\begin{claimproof}
    Assume 
    $\Iast_{\epsilon}^G\subset \Iast_{\epsilon'}^G$ for some fixed $\eps<\eps'\in [0,\eps^\ast)$, which will lead to a contradiction to the rank-criticality of $\Iast$.   
    {Namely, we show that $(\Iast_{\epsilon'}^G)_\epsilon^F$ is an interval that properly contains $L$, but has the same $\rk_M$ value as $L$.} We first note the following:
   
   \begin{itemize} 
   \item By Claim \ref{claim:pairwise difference}, we have $\Iast_{\epsilon}^G\neq \Iast_{\epsilon'}^G$, and thus 
    \begin{equation}\label{eq:IastGeps is smaller}
    \Iast_{\epsilon}^G\subsetneq \Iast_{\epsilon'}^G.
    \end{equation}
  \item   As $L\subset q_\ast^\uparrow$, Proposition~\ref{prop:G, F preserve qast} and Remark~\ref{rem:q1 is at most tast}
     imply that $\Iast_{\epsilon'}^G\subset q_\ast^\uparrow$, and in turn 
    \begin{equation}\label{eq:IastGF is in the quadrant} (\Iast_{\epsilon'}^G)_\epsilon^F\subset q_\ast^\uparrow.
    \end{equation}
    
\end{itemize}

    Now we claim that $(\Iast_{\epsilon'}^G)_\epsilon^F\supsetneq \Iast$. Let $x\in L$.
    By Remark~\ref{rem:GF is identity}~\ref{item:GF is identity2}, $(\Iast_{\epsilon}^G)_\epsilon^F=\Iast$,  and thus $x\in (\Iast_{\epsilon}^G)_\epsilon^F$.
    Then $F_{\epsilon}(x)\in \Iast_{\epsilon}^G\subset \Iast_{\epsilon'}^G$, which by definition implies $x\in (\Iast_{\epsilon'}^G)_\epsilon^F$, as desired.
    Next, we show $(\Iast_{\epsilon'}^G)_\epsilon^F\neq \Iast$.
    Note that for all $t\in [t^\ast, t^\ast+\epsilon]$,  $g_\epsilon(t) = g_{\epsilon'}(t) = t^\ast$ {(here, we use the assumption $\eps<\eps'$)}, and hence Equation \eqref{eq:some equation} implies 
    \begin{equation}\label{eq:both are Stast}\mbox{for all $t\in [t^\ast,t^\ast+\eps]$,}\qquad\widehat{\pi}_1(\Iast_{\epsilon}^G\vert_{\pi_1=t}) = \widehat{\pi}_1(\Iast_{\epsilon'}^G\vert_{\pi_1=t}) = S_{t^\ast}.\end{equation}    
    Hence for any $y=(y_1,\ldots,y_d)\in \Iast_{\epsilon'}^G\setminus \Iast_{\epsilon}^G$ (which exists by Equation~\eqref{eq:IastGeps is smaller}), we have $y_1\notin [t^\ast, t^\ast + \epsilon]$; otherwise $(y_2,\ldots,y_d)\in \widehat{\pi}_1(\Iast_{\epsilon'}^G\vert_{\pi_1=y_1}) = \widehat{\pi}_1(\Iast_{\epsilon}^G\vert_{\pi_1=y_1})$, implying that $y\in L_\eps^G$, a contradiction.
    As a result, $y_1\in \im(f_{\epsilon})$, and in turn $y\in \im(F_{\epsilon})$.
    Thus, there exists $y'\in \R^d$ with $F_{\epsilon}(y')=y$, which by definition implies $y'\in (\Iast_{\epsilon'}^G)_\epsilon^F$. To show  $(\Iast_{\epsilon'}^G)_\epsilon^F\neq \Iast$, it suffices to show that $y'\not\in L$.
    Since $G_\eps\circ F_\eps$ is the identity map on $\R^d$ (cf. Remark~\ref{rem:GF is identity}~\ref{item:GF is identity1}), we have  $y' = G_{\epsilon}(F_{\epsilon}(y')) = G_{\epsilon}(y)$. But, as $y\notin \Iast_{\epsilon}^G$, we have $y'=G_{\epsilon}(y)\notin \Iast$, as desired. 

    Next we claim that $(\Iast_{\epsilon'}^G)_\epsilon^F$ is an interval. Note, by Proposition~\ref{prop:Ie is an interval for all e}, that $L_{\eps'}^G$ is an interval. Hence,  Equation~\eqref{eq:both are Stast} %
and Claim~\ref{claim:slice condition implications IFe} \ref{item:IFe is an interval} imply the claim.

Lastly, we claim that $L$ and $(\Iast_{\epsilon'}^G)_\epsilon^F$ have the same $\rk_M$ value. Indeed,
\begin{align*}
\rk_M((\Iast_{\epsilon'}^G)_\epsilon^F)&=\rk_M(((\Iast_{\epsilon'}^{G})_{\epsilon}^{F})_{\epsilon}^{G})&\mbox{by Equation~\eqref{eq:IastGF is in the quadrant} and  Proposition~\ref{prop:ranks of I and Ie are equal}}
\\&=\rk_M(\Iast_{\epsilon'}^G)&\mbox{Claim \ref{claim:slice condition implications IFe} \ref{item:F then G equals original} implies $((\Iast_{\epsilon'}^{G})_{\epsilon}^{F})_{\epsilon}^{G} = \Iast_{\epsilon'}^G$}
\\&=\rk_M(\Iast)&\mbox{by the fact $L\subset q_\ast^\uparrow$ and Proposition~\ref{prop:ranks of I and Ie are equal},}
\end{align*}
completing the proof.
\end{claimproof}

\smallskip

\begin{claim}\label{claim:the other containment}
 For $\epsilon<\epsilon'\in [0,\epsilon^\ast)$, we have that  $\Iast_{\epsilon'}^G\not\subset \Iast_{\epsilon}^G$.
\end{claim}
\begin{claimproof}
    We assume 
    that $\Iast_{\epsilon'}^G\subset \Iast_{\epsilon}^G$ for some fixed $\eps<\eps'\in [0,\eps^\ast)$ and deduce a contradiction to the rank criticality of $L$.
   Let $t'\in (t^\ast + \epsilon, t^\ast + \epsilon']$.
    By the definitions of the functions $g_{\epsilon}$ and $g_{\epsilon'}$, we have that 
    \begin{equation}\label{eq:y and tast}        
    t^\ast = g_{\epsilon'}(t') < g_{\epsilon}(t') =: t \leq t^\ast+\epsilon^\ast.
    \end{equation}
    Then, we obtain the following sequence of inclusions: %
    \begin{align*}
    &\widehat{\pi}_1(\Iast_{\epsilon'}^G\vert_{\pi_1=t'})
    \subset \widehat{\pi}_1(\Iast_{\epsilon}^G\vert_{\pi_1=t'})  &\mbox{since $\Iast_{\epsilon'}^G\subset \Iast_{\epsilon}^G$}
\\ &\implies
    S_{g_{\epsilon'}(t')}
    \subset S_{g_{\epsilon}(t')}&\mbox{by Equation~\eqref{eq:some equation}}
\\&\implies S_{t^\ast}
 \subset S_t &\mbox{by Equation~\eqref{eq:y and tast}.}
\end{align*}
  Further,
  by Equations~\eqref{eq:slices are different} and \eqref{eq:y and tast}, 
  we have $S_{t^\ast}\neq S_{t}$, 
  and thus $S_{t^\ast}\subsetneq S_{t}$. This implies that the set
     \[I:=\Iast \cup \left([t^\ast,t]\times S_{t}\right)\] 
    properly contains $L$; namely, $\emptyset\subsetneq\{t^\ast\}\times (S_t\setminus S_{t^\ast})\subset I\setminus L$. Hence, it suffices to show that $I$ is an interval and $\rk_M(I)=\rk_M(L)$. 

    We define the function $h:\R\to \R$ as 
    \[h(s):=\begin{cases}
      t, & \mbox{if } s\in [t^\ast, t).\\
        s, & \mbox{otherwise.}
    \end{cases}\]
    We also define $H:\R^d\to \R^d$ by $H((x_1,x_2,\ldots,x_d)):=(h(x_1),x_2,\ldots,x_d)$.
    It is clear that for $x\in \R^d$, $x$ and $H(x)$ can differ in only the first coordinate and that $x\leq H(x)$.

    We claim that the image of $I$ under $H$ is contained in $L$.
    Let $x\in I$.

    ($\ast$) Suppose $x\notin \Iast$. Then by construction of $I$, $x = (x_1, x')$ for some  $x_1\in [t^\ast, t]$ and  $x'\in S_{t}$.
    Then $H(x) = (t, x')$.
    As $x'\in S_t$, we have that $H(x)\in \Iast$, as desired.
    
    Now suppose $x\in \Iast$.
    If $H(x)=x$, then trivially, $H(x)\in \Iast$. Let us assume that $H(x)\neq x$.
    This implies that the first coordinate of $x$ is not fixed by $h$, i.e., $x = (x_1,x')$ for some $x_1\in [t^\ast, t)$ and $x'\in S_{x_1}$.
    If $x_1=t^\ast$, then we have $H(x) = (t,x')\in \Iast$ from the fact $S_{t^\ast}\subset S_{t}$.
    Now, assume $x_1\in (t^\ast, t)$.
    Recall that $t = g_{\epsilon}(t')$ for $t'\in (t^\ast + \epsilon, t^\ast+\epsilon']$.
    Then consider $g_{\epsilon}(f_{\epsilon'}(x_1))$.
    As $t^\ast<x_1$, 
    \[t' \leq t^\ast+\epsilon' \leq f_{\epsilon'}(x_1),\]
    and as $g_{\epsilon}$ is order-preserving, this gives 
    \[t = g_{\epsilon}(t') \leq s := g_{\epsilon}(f_{\epsilon'}(x_1)).\]
    By the assumption $\Iast_{\epsilon'}^G\subset \Iast_{\epsilon}^G$, 
    we have that $\widehat{\pi}_1(\Iast_{\epsilon'}^G\vert_{\pi_1=f_{\epsilon'}(x_1)}) \subset \widehat{\pi}_1(\Iast_{\epsilon}^G\vert_{\pi_1=f_{\epsilon'}(x_1)})$, 
    which %
    by Equation \eqref{eq:some equation} and the fact $g_{\epsilon'}(f_{\epsilon'}(x_1)) = x_1$, implies $S_{ x_1}\subset S_s$.
    Hence, as $x = (x_1, x')\in \Iast$,     %
    we have that $(s, x')\in \Iast\vert_{\pi_1=s}\subset L$.
    Also, since $x_1<t\leq s$, we have \[(x_1,x')\leq (t, x')=H(x) \leq (s,x').\] Then, convexity of $\Iast$ and the fact that $(x_1,x'),(s,x')\in \Iast$ imply  $H(x) \in \Iast$. 
    We have proved that   
    the image of $I$ under $H$ is contained in $L$, as desired.

    Now, we will show that $I$ is an interval. 
    To show connectivity of $I$, suppose $p,q\in I$. 
    Then $H(p),H(q)\in \Iast$.
    As $\Iast$ is connected, there exists a path 
    in $\Iast$ (hence in $I$) connecting 
    $H(p)$ and $H(q)$.
    As $p \leq H(p)$ and $q\leq H(q)$, 
    by prepending $p$ and appending $q$ to the path,
    we obtain a path in $I$ connecting $p$ to $q$.
    Next we show convexity of $I$. 
    Let $p<q<r\in \R^d$ with $p,r\in I$, and we show that $q\in I$. %
    First, observe that since $H$ is order-preserving, we have 
    \[(\ast\ast)\qquad \mbox{$H(p)\leq H(q)\leq H(r)$ in $\R^d$, and thus convexity of $L$ implies that $H(q)\in L$.}\]
    The claim that $q\in I$ follows by considering the following cases.
    \begin{itemize}
    \item Assume $p,r\in \Iast$. Then $q\in \Iast$ by convexity of $\Iast$, and thus $q\in I$.
    \item Assume $p\in \Iast$ and $r\in I\setminus \Iast$. Since $r\leq H(r)$ and in $H(r)\in \Iast$, we have that $p<q<H(r)$ with $p,H(r)\in \Iast$, and thus convexity of $L$ implies $q\in \Iast\subset I$. 
    \item Assume $p,r\in I\setminus \Iast$.  We saw from ($\ast$) that $H(p),H(r)\in \Iast\vert_{\pi_1=t}$. Recall also ($\ast\ast$). 
    As the first coordinates of $H(p)$ and $H(r)$ are both equal to $t$, this also holds for $H(q)$, i.e.,
    $H(q)\in \Iast\vert_{\pi_1=t}$.
    {Then, since $q$ and $H(q)$ can differ only in the first coordinate, $\widehat{\pi}_1(q)=\widehat{\pi}_1(H(q))\in S_t$. 
    In addition, since $H(q)$ is obtained by applying $h$ to the first coordinate of $q$ while leaving the other coordinates unchanged, the fact that $h^{-1}(t)=[t^\ast,t]$ 
    implies that the first coordinate of $q$ is in $[t^\ast,t]$. In sum,} %
\[
q\in [t^\ast,t]\times S_t\subset I.
\]%
    \item Assume $p\in I\setminus \Iast$ and $r\in L$.
    Let us write $p=(p_1,p'), q=(q_1,q')$ and $r=(r_1,r')$ with $p_1,q_1,r_1\in \R$ and $p',q',r'\in \R^{d-1}$.
    As $p<q<r$, we have $p_1\leq q_1\leq r_1$ and $p'\leq q'\leq r'$.
    Observe that the assumption $p\in I\setminus \Iast$ implies $p\in [t^\ast,t)\times S_t$, so in particular, we have $t^\ast \leq p_1<t$.
    Since $p_1\leq q_1$, this implies $t^\ast \leq q_1$.
    
    Recall $(\ast\ast)$.
    If $q_1\in [t,\infty)$,
    then $q=H(q)\in \Iast$.
    Assume $q_1\in (-\infty,t)$. 
    Then we have $t^\ast \leq q_1<t$, and so $H(q) = (t,q')$.
    As $H(q)\in \Iast$, we have $H(q)\in \Iast\vert_{\pi_1=t}$, which implies $q'\in S_t$. In sum, we have
    $q\in [t^\ast,t]\times S_t\subset I$, as desired.

\end{itemize}

    Finally, we prove that $\rk_M(I)=\rk_M(\Iast)$. 
    By monotonicity of the GRI, it suffices to show that $\rk_M(I)\geq \rk_M(\Iast)$. 
    We do this by showing $I\subset q^\ast$ and constructing a map $\Psi:\CC(I) \to \CC(\Iast)$  satisfying all assumptions of Proposition~\ref{prop:connected components to rank general}.
        First, we claim that $I\subset q_\ast^\uparrow$.
    Since $\Iast\subset q_\ast^\uparrow$, it suffices to show that $[t^\ast, t]\times S_t\subset q_\ast^\uparrow$. 
    {Then, by Remark~\ref{rem:q1 is at most tast}, it suffices to show that $S_t\in \widehat{\pi}_1(q_\ast^\uparrow)$. This follows from
    \begin{align*}
        S_t&=\widehat{\pi}_1(L|_{\pi=t})&\mbox{by definition}
        \\&\subset \widehat{\pi}_1(L)
        \\&\subset \widehat{\pi}_1(q_\ast^\uparrow)&\mbox{since $L\subset q_\ast^\uparrow$.}
    \end{align*}
    }

    Now we define $\Psi:\CC(I) \to \CC(\Iast)$ by $\Psi(X)= Y$ if $X\cap Y\neq \emptyset$. 
    {We claim that for all $X\in \CC(I)$,  $\Psi(X)$ exists in $\CC(\Iast)$. 
    Fix $X\in \CC(I)$. 
    Verifying that $H(X)\subset X$ suffices: Indeed, since $H(I)\subset \Iast\subset q_\ast^\uparrow$, the containment $H(X)\subset X$ implies that for any $x\in X$, there exists $Y\in \CC(\Iast)$ such that $H(x)\in X\cap Y$, and thus $X\cap Y\neq \emptyset$. Now observe that 
    \[(\ast\ast\ast)\qquad \mbox{$X$ and $H(X)$ are contained in the same unit block}\] because (1) $H$ only affects the first coordinate of its input via $h$, (2)  $h([t^\ast,t])=\{t\}\subset [t^\ast,t]\subset \pi_1(B)$ (cf. Equations~\eqref{eq:closed interval around tast} and \eqref{eq:y and tast}), and (3) $h$ fixes points outside $[t^\ast,t]$. Also, since $H$ is order-preserving and $X$ is connected, $H(X)$ is also connected. Furthermore, since  $x\leq H(x)$ for each $x\in X$, and $X$ is a connected component (within the unit block containing $X$), we have $H(X)\subset X$.}  
    
    Now we claim that $\Psi$ is well-defined, i.e., $\Psi(X)$ is unique for each $X\in \CC(I)$.
    Assume $x,x'\in X\in \CC(I)$ with $x\in Y$ and $x'\in Y'$ for some $Y,Y'\in \CC(\Iast)$.
    We show $Y=Y'$.
    Since 
    both $x$ and $H(x)$ belong to $L$, and  $x\leq H(x)$, the fact in $(\ast\ast\ast)$ and the assumption $x\in Y$ imply that 
    $H(x)\in Y$. Similarly, $H(x')\in Y'$.
   As $x,x'\in X$, 
    there is a path 
    in $X$
    connecting $x$ and $x'$.
    Then the image of the path connects $H(x)$ and $H(x')$ in $\Iast$.
    As $H$ is order-preserving, the fact $(\ast\ast\ast)$ implies that $H(x)$ and $H(x')$ belong to the same element of $\CC(L)$, i.e.,
    $Y=Y'$ as desired.

    {It only remains to show that $\Psi$ satisfies assumptions \ref{item:Psi condition 1} and \ref{item:Psi condition 3} in Proposition~\ref{prop:connected components to rank general}. %
    Assumption \ref{item:Psi condition 1}: Since any $X\in \CC(I)$ and its image $\Psi(X)\in \CC(\Iast)$ have nonempty intersection, by the definitions of $\CC(I)$ and $\CC(\Iast)$, $X$ and $\Psi(X)$ are in the same unit block.
    \newline \noindent Assumption \ref{item:Psi condition 3}: Let $(X,X')\in \CC^2_{\leq}(I)$, i.e., there exists $x\in X$ and $x'\in X'$ such that $x\leq x'$. Recall that $\Psi(X)$ and $\Psi(X')$ contain $H(x)$ and $H(x')$, respectively. Since $H$ is order-preserving, we have $H(x)\leq H(x')$ as desired.}
\end{claimproof}

Let $\epsilon<\epsilon'\in [0,\eps^\ast)$, and let $J\subset \R^d$ be the convex hull of $\Iast_{\epsilon}^G\cup \Iast_{\epsilon'}^G$ (cf. Remark~\ref{rem:properties of finitely generated modules}~\ref{item:convex hull is smallest interval}).
The following claim regarding slices of $J$ will be useful in the sequel:
\begin{claim} For all  $t\in (-\infty, t^\ast+\epsilon]$,
     {\begin{equation}\label{eq:slices of J equal slices of Ie simpler}
        J\vert_{\pi_1=t} = \Iast_{\epsilon}^G\vert_{\pi_1=t} = \Iast_{\epsilon'}^G\vert_{\pi_1=t}.
    \end{equation}  } 
\end{claim}
{We note that, by Remark~\ref{rem:slice projection vs. slice}, this claim is equivalent to:}
    \begin{equation}\label{eq:slices of J equal slices of Ie}
        \mbox{for all } t\in (-\infty, t^\ast+\epsilon],\hspace{5mm} 
        \widehat{\pi}_1(J\vert_{\pi_1=t}) = \widehat{\pi}_1(\Iast_{\epsilon}^G\vert_{\pi_1=t}) = \widehat{\pi}_1(\Iast_{\epsilon'}^G\vert_{\pi_1=t}).        
    \end{equation}
     
\begin{claimproof}
 First, from Figure~\ref{fig:g epsilon and f epsilon}, observe that 
    \begin{equation}\label{eq:identity map within a range}\mbox{for all $t\in \R\setminus (t^\ast,t^\ast+\eps^\ast)$, }\ \  g_\epsilon(t) = g_{\epsilon'}(t) = t.\end{equation}

  We prove the right equality in Equation \eqref{eq:slices of J equal slices of Ie} {(and thus the one in Equation~\eqref{eq:slices of J equal slices of Ie simpler})}.
      If $t \in (-\infty,t^\ast)$, by Equation~\eqref{eq:identity map within a range}, we have $g_\epsilon(t)=g_{\epsilon'}(t) = t$, and thus Equation \eqref{eq:some equation} implies that both $\widehat{\pi}_1(\Iast_{\epsilon}^G\vert_{\pi_1=t})$ and $\widehat{\pi}_1(\Iast_{\epsilon'}^G\vert_{\pi_1=t})$ equal $S_t$, establishing the right equality.
     If $t\in [t^\ast, t^\ast+\epsilon]$, we have $g_\epsilon(t) = g_{\epsilon'}(t) = t^\ast$, and hence Equation \eqref{eq:some equation} again implies the right equality.
    
    We show the left equality in Equation \eqref{eq:slices of J equal slices of Ie simpler} {(and thus the one in Equation~\eqref{eq:slices of J equal slices of Ie})}.
    {From the definition $J=\overline{\Iast_{\epsilon}^G\cup \Iast_{\epsilon'}^G}$, it is immediate that $J\vert_{\pi_1=t}\supset \Iast_{\epsilon}^G\vert_{\pi_1=t}$. We show  $J\vert_{\pi_1=t}\subset\Iast_{\epsilon}^G\vert_{\pi_1=t}(=\Iast_{\epsilon'}^G\vert_{\pi_1=t})$.}
    Let $p\in J\vert_{\pi_1=t}$. From the definition of convex hull, $p\in J$ implies that there exist $q,r\in \Iast_{\epsilon}^G\cup \Iast_{\epsilon'}^G$ such that $q\leq p\leq r$.
    \begin{itemize}
    \item 
    If $q, r\in \Iast_{\epsilon}^G$, then $p\in \Iast_{\epsilon}^G$ by convexity of $\Iast_{\epsilon}^G$ (cf. Proposition~\ref{prop:Ie is an interval for all e}), and as $\pi_1(p)=t$, we have $p\in \Iast_{\epsilon}^G\vert_{\pi_1=t}$. 
    \item 
    Similarly, if $q,r\in \Iast_{\epsilon'}^G$, then $p\in \Iast_{\epsilon'}^G|_{\pi_1=t}$.
    \item If $q\in \Iast_{\epsilon}^G$ and $r\in \Iast_{\epsilon'}^G$, then let $q=:(q_1, q')$ with $q_1\in \R$, $q'\in \R^{d-1}$. Also, let $p=:(t, p')$ for some $p'\in \R^{d-1}$.
    Since $q\leq p$, we have $q_1\leq t\leq t^\ast + \epsilon$. 
    {Thus, the right equality in Equation  \eqref{eq:slices of J equal slices of Ie} gives
    $\widehat{\pi}_1(\Iast_\eps^G\vert_{\pi_1=q_1}) = \widehat{\pi}_1(\Iast_{\eps'}^G\vert_{\pi_1=q_1})$. Hence, as $q'\in \widehat{\pi}_1(\Iast^G_\eps\vert_{\pi_1=q_1})$, we have $q'\in \widehat{\pi}_1(L^G_{\eps'}\vert_{\pi_1=q_1})$, which implies $q\in \Iast^G_{\eps'}$.
    Thus, $q,r\in \Iast^G_{\eps'}$, the case considered in the previous item, and thus $p\in \Iast^G_{\eps}|_{\pi_1=t}$.
    }
   \item {If $q\in \Iast_{\epsilon'}^G$ and $r\in \Iast_{\epsilon}^G$, then
$q\in \Iast_{\epsilon}^G$ by interchanging the roles of $\epsilon$ and
$\epsilon'$ in the preceding argument. Hence, we are again in the situation
described in the first bullet point.} %
\end{itemize}
{We have proved the left equality in Equation \eqref{eq:slices of J equal slices of Ie simpler}, completing the proof.}    
\end{claimproof}   

We proceed by demonstrating that if an interval in $\R^d$ contains two distinct intervals from $\{\Iast_{\epsilon}^G\}_{\eps \in [0,\eps^\ast)}$, then it must have strictly smaller rank than $L$:

\begin{claim}\label{claim:rank of interval containing two decreases}
    Let $\epsilon<\epsilon'\in [0,\epsilon^\ast)$.
    For any $I\in \Int(\R^d)$ containing both $ \Iast_{\epsilon}^G$ and $\Iast_{\epsilon'}^G$,
    we have that $\rk_M(I) < \rk_M(\Iast)$.
\end{claim}
\begin{claimproof}
   Let $J\subset \R^d$ be the convex hull of $\Iast_{\epsilon}^G\cup \Iast_{\epsilon'}^G$ (cf. Remark~\ref{rem:properties of finitely generated modules}~\ref{item:convex hull is smallest interval}). It suffices to show that $J$ is the smallest interval containing $\Iast_{\epsilon}^G\cup \Iast_{\epsilon'}^G$ and that $\rk_M(J)<\rk_M(L)$. 
   
   By Remark~\ref{rem:properties of finitely generated modules}~\ref{item:convex hull is smallest interval}, proving the former reduces to showing that $\Iast_{\epsilon}^G\cup \Iast_{\epsilon'}^G$ is connected.
    Taking $t=t^\ast$ and $t=t^\ast+\epsilon^\ast$ in Equation~\eqref{eq:identity map within a range}, we obtain from Equation~\eqref{eq:some equation} that
\[\widehat{\pi}_1(\Iast_{\epsilon}^G\vert_{\pi_1=t^\ast}) = \widehat{\pi}_1(\Iast_{\epsilon'}^G\vert_{\pi_1=t^\ast}) = S_{t^\ast} \mbox{ and } \widehat{\pi}_1(\Iast_{\epsilon}^G\vert_{\pi_1=t^\ast+\eps^\ast}) = \widehat{\pi}_1(\Iast_{\epsilon'}^G\vert_{\pi_1=t^\ast+\eps^\ast}) = S_{t^\ast+\eps^\ast}.\]
  Then, by Remark~\ref{rem:slice projection vs. slice}, we have
    \begin{equation}\label{eq:two equations}
    \Iast_{\epsilon}^G\vert_{\pi_1=t^\ast} = \Iast_{\epsilon'}^G\vert_{\pi_1=t^\ast} = \Iast\vert_{\pi_1=t^\ast}
    \mbox{ and } 
    \Iast_{\epsilon}^G\vert_{\pi_1=t^\ast+\epsilon^\ast} = \Iast_{\epsilon'}^G\vert_{\pi_1=t^\ast + \epsilon^\ast} = \Iast\vert_{t^\ast + \epsilon^\ast}.
    \end{equation}
    By Equation~\eqref{eq:slices are different}, we have $S_{t^\ast}\neq S_{t^\ast + \epsilon^\ast}$  
    , and thus at least one of $S_{t^\ast}$ and $S_{t^\ast + \epsilon^\ast}$ is nonempty, which in turn implies that at least one of $\Iast\vert_{\pi_1=t^\ast}$ and $\Iast\vert_{\pi_1=t^\ast+\eps^\ast}$ is nonempty. Let us pick a point $p\in \left(\Iast\vert_{\pi_1=t^\ast}\right)\cup \left( \Iast\vert_{\pi_1=t^\ast+\epsilon^\ast}\right)$. 
    By Equation~\eqref{eq:two equations}, $p$ belongs to both $\Iast_{\epsilon}^G$ and $\Iast_{\epsilon'}^G$. %
    As $\Iast_{\epsilon}^G$ and $\Iast_{\epsilon'}^G$ are both intervals (cf. Proposition~\ref{prop:Ie is an interval for all e}), they are connected.
    Thus, by Remark~\ref{rem:properties of finitely generated modules}~\ref{item:union of connected sets with overlap is connected}, $\Iast_{\epsilon}^G\cup \Iast_{\epsilon'}^G$ is connected, as desired.

    It remains to show that $\rk_M(J)<\rk_M(L)$. 
    {By the rank-criticality of $L$, it suffices to establish:}
    
    {\emph{Subclaim 1}. $J_\eps^F\in \Int(\R^d)$,}

    {\emph{Subclaim 2}. $\rk_M(J)=\rk_M(J_\eps^F)$, and}

    {\emph{Subclaim 3}. $J_\eps^F$ properly contains $L$.}

  \noindent  {In order to show \emph{Subclaim 2}, by Proposition~\ref{prop:ranks of I and Ie are equal}, it suffices to show:}

   {\emph{Subclaim 4}. $J_\eps^F\subset q_\ast^\uparrow$, and}
  
 {\emph{Subclaim 5}. $J=(J_\eps^F)_\eps^G$.}

\noindent{ \emph{Proof of Subclaims 1 and 5.}}      Recall that, for $t\in [t^\ast, t^\ast + \epsilon]$, $g_{\epsilon}(t) = t^\ast$, so Equation \eqref{eq:slices of J equal slices of Ie} alongside Equation \eqref{eq:some equation} imply that, 
    for all $t\in [t^\ast, t^\ast + \epsilon]$,  
    $\widehat{\pi}_1(J\vert_{\pi_1=t}) = \widehat{\pi}_1(J\vert_{\pi_1=t^\ast})$, as both equal $S_{t^\ast}$.
    Thus, by Claim \ref{claim:slice condition implications IFe} \ref{item:IFe is an interval}, we have that $J_{\epsilon}^F$ is an interval, and Claim \ref{claim:slice condition implications IFe} \ref{item:F then G equals original} implies that $(J_{\epsilon}^F)_{\epsilon}^G = J$.
    
    \noindent{{\emph{Proof of Subclaim 3.}}} 
    Let $x\in \Iast$. As $G_{\epsilon}(F_{\epsilon}(x)) = x$ (cf. Remark~\ref{rem:GF is identity}~\ref{item:GF is identity1}), we have $F_{\epsilon}(x)\in \Iast_{\epsilon}^G\subset J$. 
    From the definition of $J_{\epsilon}^F$, this implies $x\in J_{\epsilon}^F$, so $J_{\epsilon}^F\supset \Iast$.
    By Claim \ref{claim:the other containment}%
    , we have that there exists $x\in J$ such that $x\notin \Iast_{\epsilon}^G$.
    By definition, this means that $G_{\epsilon}(x)\notin \Iast$.
    As established above, $x\in J$ implies $G_{\epsilon}(x)\in J_{\epsilon}^F$, so $G_{\epsilon}(x)\in J_{\epsilon}^F\setminus \Iast$, implying $J_{\epsilon}^F\neq \Iast$, which alongside $J_{\epsilon}^F\supset \Iast$ demonstrates that $J_{\epsilon}^F\supsetneq \Iast$.

    \noindent\emph{{Proof of Subclaim 4.}} Lastly, note that $\Iast_{\epsilon}^G, \Iast_{\epsilon'}^G\subset q_\ast^\uparrow$ implies that $J\subset q_\ast^\uparrow$ as every point in $J$ is equal to or greater than some points in $\Iast^G_\epsilon\cup \Iast^G_{\epsilon'}\subset q_\ast^\uparrow$(c.f. Remark \ref{rem:properties of finitely generated modules} \ref{item:convex hull is smallest interval}) %
    Proposition~\ref{prop:G, F preserve qast} then implies $J_{\epsilon}^F\subset q_\ast^\uparrow$.
\end{claimproof}

    Fix a pair $\epsilon< \epsilon'$ in $[0, \epsilon^\ast)$. %
    By Claim~\ref{claim:rank of interval containing two decreases} {and Proposition~\ref{prop:ranks of I and Ie are equal}}, we have
\[
\rk_M(\Iast_{\epsilon}^G) >
\rk_M\left(\overline{\Iast_{\epsilon}^G\cup \Iast_{\epsilon'}^G}\right),
\]
or equivalently,
\[\sum_{\substack{J\supset \Iast_{\epsilon}^G\\ J\in \Ical}}
    \dgm^\Int_M(J)
    \;>\;
    \sum_{\substack{J\supset \overline{\Iast_{\epsilon}^G\cup \Iast_{\epsilon'}^G}\\ J\in \Ical}}
    \dgm^\Int_M(J).
    \]
    Therefore, we have
\begin{align}
    \sum_{\substack{J\supset \Iast_{\epsilon}^G,\,
    J\not\supset \overline{\Iast_{\epsilon}^G\cup \Iast_{\epsilon'}^G}\\
    J\in \Ical}}
    \dgm^\Int_M(J)
    \;>&\;
    0.
\end{align}
This implies that the set $\Jcal_{\epsilon,\epsilon'}:= \{J\in \Ical \, : \, J\supset \Iast_{\epsilon}^G \mbox{ and } J\not\supset \Iast_{\epsilon'}^G\}$ is nonempty.
For {$\epsilon\in [0,\epsilon^\ast)$,}
let $\Jcal_\epsilon:=\bigcup_{\epsilon'\in (\epsilon, \epsilon^\ast)} \Jcal_{\epsilon,\epsilon'}$.
Then, we define $\Jcal:=\bigcup_{\epsilon\in [0,\epsilon^\ast)} \Jcal_\epsilon$, which is a subcollection of $\Ical$ by construction. 

\begin{claim}\label{claim:Jcal is infinite pixelized}
    $|\Jcal| = \infty$.
\end{claim}
\begin{claimproof}
    To avoid triviality,
    assume that $|\Jcal_\epsilon|<\infty$ for every $\epsilon\in [0, \epsilon^\ast)$.
    Assume by way of contradiction that $|\Jcal|<\infty$.
    As $\Jcal$ is a union of infinitely many finite sets, there must be two sets in the union, denote them $\Jcal_{\epsilon}$ and $\Jcal_{\epsilon'}$ with $\epsilon< \epsilon'$, such that $\Jcal_{\epsilon}=\Jcal_{\epsilon'}$.
    Then $\emptyset\subsetneq \Jcal_{\epsilon,\epsilon'}\subsetneq\Jcal_{\epsilon}$, and thus we can
    take an element $J\in \Jcal_{\epsilon,\epsilon'}$.  %
    By definition, $J\in \Int(\R^d)$ with $J\supset \Iast_{\epsilon}^G$, but $J\not\supset \Iast_{\epsilon'}^G$. %
    But by definition of $\Jcal_{\epsilon'}$, this means $J\notin \Jcal_{\epsilon'}$, and so $\Jcal_{\epsilon}\neq \Jcal_{\epsilon'}$, a contradiction.
\end{claimproof}

The next claim shows that assuming $\Ical\not\subset\Ical_M$ and Case 1 of Proposition \ref{prop:slice unequal to one side} yields a contradiction.
\begin{claim}\label{claim:rkm is not mobius invertible} $\rk_M^\Int$ is not M\"obius invertible.
\end{claim}
\begin{claimproof}We begin by noting that, by Equation~\eqref{eq:slices are different}, $S_{t^\ast}\neq S_{t^\ast+\epsilon^\ast}$ and hence at least one of $S_{t^\ast}$ and $S_{t^\ast+\epsilon^\ast}$ is nonempty. This implies, by Remark \ref{rem:slice projection vs. slice}, that at least one of $\Iast\vert_{\pi_1=t^\ast}$ and $\Iast\vert_{\pi_1=t^\ast+\epsilon^\ast}$ is nonempty. Let us fix a point $p$ in the union of these two sets. 
Next, we claim that, for all $\epsilon
\in [0,\epsilon^\ast)$,  \[\Iast_{\epsilon}^G\vert_{\pi_1=t^\ast} = \Iast\vert_{\pi_1=t^\ast} \mbox{ and } \Iast_\epsilon^G\vert_{\pi_1=t^\ast + \epsilon^\ast} = \Iast\vert_{\pi_1=t^\ast+\epsilon^\ast}.\] 
Indeed, since $g_\epsilon(t^\ast) = t^\ast$ and $g_\epsilon(t^\ast + \epsilon^\ast) = t^\ast + \epsilon^\ast$, Equation \eqref{eq:some equation} together with Remark \ref{rem:slice projection vs. slice} implies the two equalities. 
Hence, the point $p$ belongs to
\[\Iast_{\epsilon}^G\vert_{\pi_1=t^\ast}\cup \Iast_\epsilon^G\vert_{\pi_1=t^\ast + \epsilon^\ast},\]
which is a subset of $\Iast_\epsilon^G.$ This implies that $p$ belongs to every interval in $\Jcal$. Then, since $\Jcal$ is an infinite subcollection of $\Ical=\supp(\dgm_M^\Int)$, 
 we have shown that the set $\{J\supset \{p\} \, : \, J\in \Ical\}$ is infinite. But, since $\{p\}\in \Int(\R^d)$, this contradicts M\"obius invertibility of  $\rk^\Int_M$ (Definition~\ref{def:mobius invertible general} \ref{item:mobius invertible general 1}).
\end{claimproof}

We have demonstrated that assuming Case 1 of Proposition \ref{prop:slice unequal to one side} yields a contradiction with the initial assumption of Theorem \ref{thm:unique candidate} that $\rk^\Int_M$ is M\"obius invertible.
Similarly, assuming Case 2 of Proposition \ref{prop:slice unequal to one side} yields a contradiction with the same assumption.
This proof is highly analogous due to symmetry between Cases 1 and 2, and follows the same argument structure as Case 1. 
Accordingly, we omit it from the main text and present the full details in Appendix \ref{sec:case 2 details}.

\newcommand{\etalchar}[1]{$^{#1}$}

\appendix

\section{{Details on Remark~\ref{rem:conflict}:} On Errors Present in the Literature
}\label{sec:details on prior errors}

In this section, 
we show that 
\cite[Lemma 3.12, Corollaries 3.13, 3.15, and 3.18]{botnan2024signed} are erroneous, as demonstrated by the counterexamples described below.

Recall that a poset is a \emph{lattice} if every pair $p,q$ of its elements has a join $p\vee q$ and a meet $p \wedge q$ in it. A \emph{lattice homomorphism} (or a \emph{morphism of lattices}) is an order-preserving map between lattices that preserve joins and meets. 
\begin{lemma}[{\cite[Lemma 3.12]{botnan2024signed}}]\label{lem:BOO 3.12}
    Let $f:\catP\to \catQ$ be a morphism of lattices. 
    Then for any interval $I$ in $\catP$, and $j$ in $\overline{f(I)}$, the convex hull of $f(I)$, the posets $\{i\in I \, :\, f(i)\leq j\}$ and $\{i\in I \,:\, f(i)\geq j\}$ are connected.
\end{lemma}

\begin{proof}[Proof from \cite{botnan2024signed}] For $i_1, i_2 \in \{i \in I : f(i) \leq j\}$ also their join lies in this set. The proof for the other poset is dual.
\end{proof}
\begin{proof}[Refutation] We note that $i_1\vee i_2$ does not necessarily belong to $I$; while we have $f(i_1\vee i_2)=f(i_1)\vee f(i_2)\leq j$, it is not necessarily the case that $i_1\vee i_2$ is an element of $I$.

For such an example, consider the finite grid $\catQ = \{0,1,2,3\}\times \{0,1\}\times \{0,1,2,3\}\subset \R^3$ from Example~\ref{ex:another R3 module}.
Let $\catP:=\{(x,y,z)\in \R^3:x,y,z\geq 0\}$ and $f:\catP\to \catQ$ be the floor map $\lfloor - \rfloor_{\catQ}:\catP\to \catQ$, sending each $p\in P$ to the maximal $q\in Q$ with $q\leq p$. This $f$ is a lattice homomorphism.  
Consider the interval $I\subset \catP$ depicted in Figure~\ref{fig:interval not mapped to interval}.
Let $i_1:=(0, 0.5, 1),i_2:=(1, 0.5, 0)\in I$, and $j:=(1, 0, 1)\in \catQ$. 
From  Remark~\ref{rem:difference for fp modules in dims 2 and 3+}, we readily observed that $j$ belongs to the convex hull of $f(I)$. One can also see that 
\begin{itemize}
\item $f(i_1) = (0, 0, 1)$ and $f(i_2)= (1, 0, 0)$, and thus $i_1,i_2\in \{i\in I \, : \, f(i) \leq j\}$.
\item  $f( i_1\vee i_2) = f ((1, 0.5, 1))  =(1, 0, 1) \leq j$.
\item However, $i_1\vee i_2 = (1, 0.5, 1)\notin I$.\qedhere
\end{itemize}

\end{proof}
The error in Lemma~\ref{lem:BOO 3.12} propagates through the following subsequent results in \cite{botnan2024signed}, which rely on Lemma~\ref{lem:BOO 3.12}.

\begin{corollary}[{\cite[Corollary 3.13]{botnan2024signed}}]\label{cor: Boo's Cor 3.13}
    Let $f:\catP\to \catQ$ be a morphism of lattices. 
    Then for any interval $I$ in $\catP$, and any $\catQ$-module $N$, we have
    \[\rk_{f^\ast N}(I) = \rk_N\left(\overline{f(I)}\right).\]
\end{corollary}
\begin{proof}[Refutation] For a counterexample to this claim, let $\catP$, $\catQ$, $f$, and $I$ be the same as before.
For the $\R^3$-module $M$ in Example~\ref{ex:another R3 module} 
(which is also depicted in Figure~\ref{fig:ia for fp not imply mi}), we restrict its indexing poset to $\catQ$; we denote the resulting $\catQ$-module by $N$.
Then, the $\catP$-module $f^\ast N$ coincides with the restriction of $M$ to $\catP$.  Consider intervals $I_a\subset \R^3$ for $a\in (0,1)$ from Example~\ref{ex:another R3 module}
and Figure~\ref{fig:ia for fp not imply mi}. Recall that $\rk_M(I_{a})=1$ for each $a$.
As $I\subset I_{0.5}$, monotonicity of the GRI implies $\rk_M(I) \geq 1$.  
Also, since $(2,0,0)\in I$, $1=\dim M_{(2,0,0)}\geq \rk_M(I)$, and thus $\rk_M(I) =\rk_{f^\ast N}(I) =1$.
However, we claim that $ \rk_{N}\left(\overline{f(I)}\right)\neq 1$, contradicting Corollary~\ref{cor: Boo's Cor 3.13}: Note that $\overline{f(I)}=f(I)\cup \{(1, 0, 1)\}\subset \catQ$ and thus the diagram $M|_{\overline{f(I)}}$ contains Diagram~\eqref{eq:sub-diagram} as part of it.
Therefore, we have
\[0=\rk\left(M\vert_{\overline{f(I)}}\right)=\rk_{M}\left(\overline{f(I)}\right)= \rk_{N}\left(\overline{f(I)}\right). 
\]
\end{proof}

    Given an order-preserving map $f:\catP\to \catQ$, any two subsets $\Ical\subset \Int(\catP)$ and $\Jcal\subset \Int(\catQ)$ are said to be \emph{compatible} \cite[Definition 3.14]{botnan2024signed} if:
    \begin{itemize}
        \item for any $I\in \Ical$, the convex hull of the image $f(I)$ is in $\Jcal$, and
        \item for any $J\in \Jcal$, the preimage $f\inv(J)$ is in $\Ical$, or is empty.
    \end{itemize}

\begin{corollary}[{\cite[Corollary 3.15]{botnan2024signed}}]\label{cor:BOO Cor 3.15}
Let $f:\catP\to \catQ$ be a lattice homomorphism, and let $\Ical_\catP$ and $\Ical_\catQ$ be a compatible collection of intervals  in $\catP$ and $\catQ$, respectively.
Let $N$ be a $\catQ$-module with a rank decomposition $(\Rcal, \Scal)$ over $\Ical_\catQ$.
Then $(f\inv(\Rcal),f\inv(\Scal))$ is a rank decomposition of $f^\ast N$ over $\Ical_\catP$, where $f\inv(\Rcal):=\{f\inv(R) \, :\, R\in \Rcal, R\cap f(\catP)\neq \emptyset\}$, and similarly for $\Scal$.
\end{corollary}

\begin{proof}[Refutation]Let $\catP$, $\catQ$, $f$, $I$,
and $N$ be the same as before. The map $f$ is surjective. Also, let 
$\Ical_\catP:= \Int(\catP)$ 
and $\Jcal_\catQ := \Int(\catQ)$, 
which form a compatible collection of intervals; for the connectedness of the preimages under $f$, see the proof of Proposition~\ref{prop:preimages under the floor map are connected}. Since $\Ical_\catQ$ is finite, the $\Int$-GRI of $N$ is M\"obius invertible.
Invoking Theorem~\ref{thm:rk decomposition is GPD}, the claim of Corollary \ref{cor:BOO Cor 3.15} implies that the $\Int$-GRI of $f^\ast N$ is also M\"obius invertible. However, this is a contradiction; in Example~\ref{ex:another R3 module}, we saw that
the $\Int$-GRI of the $\R^3$-module $M$ is not M\"obius invertible. That proof also shows that the $\Int$-GRI of  $M|_P$, which is isomorphic to $f^\ast N$, is also not M\"obius invertible. 
\end{proof}

{An interval $I$ in $\R^d$ is called \emph{half open} if
\begin{itemize}
    \item for every $x\in I$, there exists $\varepsilon>0$ such that
    $
    [ x,x+(\varepsilon,\ldots,\varepsilon)]\subseteq I,
    $
    and
    \item for every $x\notin I$, there exists $\varepsilon>0$ such that
    $
    [ x,x+(\varepsilon,\ldots,\varepsilon)]\cap I=\varnothing.
    $
\end{itemize}
Intervals in $\R^d$ of the form $\Pi_{i=1}^d[a_i,b_i)$ are called \emph{half-open rectangles}.
}

\begin{corollary}[{\cite[Corollary 3.18]{botnan2024signed}}]\label{cor:BOO's Cor 3.18}
Let $N$ be a finitely presented $\R^d$-module.
Let $\Ical$ be either the collection of half-open rectangles or the collection of all half-open intervals.
Then there are unique disjoint multi-sets $\Rcal$ and $\Scal$ of elements of $\Ical$, such that 
\[\rk^\Ical_N = \rk^\Ical_{\kf_\Rcal} - \rk^\Ical_{\kf_\Scal}. \]
\end{corollary}

\begin{proof}[Refutation] 
{
From Example~\ref{ex:another R3 module}, 
recall the $\R^3$-module $M$ whose Int-GRI
is not M\"obius invertible.
Observe that the support-maximal intervals $I_a$ for $a\in (0,1)$, depicted in Figure~\ref{fig:ia for fp not imply mi} (right), are all half-open and that $B_{1,0}=[(1,0,0),(2,1,1))$ is a half-open rectangle such that $B_{1,0}\subset \bigcap_{a\in (0,1)} I_a$.}
Thus, if we let $\Ical$ be the collection of half-open intervals in $\R^3$, Corollary \ref{cor:intersection} shows that $\rk^\Ical_M$ is not M\"obius invertible.
As $\rk^\Ical_M$ is not M\"obius invertible, Theorem \ref{thm:rk decomposition is GPD} implies that $\rk^\Ical_M$ has no rank decomposition over $\Ical$.
\end{proof}

However, {when $d\leq 2$, the statement for half-open intervals in the previous corollary does hold: since half-open intervals are intervals in the sense of Definition~\ref{def:intervals}, Remark~\ref{rem:monotonicity of mobius invertibility new}~\ref{item:monotonicity of mobius invertibility} and Theorem~\ref{thm:sufficient conditions for invertibility}~\ref{item:finitely presentable invertible} directly imply the statement.} 

{In addition,} the statement for half-open rectangles in the previous corollary does hold in every dimension $d\geq 1$.
The statements and proofs in \cite[Lemma 3.12--Corollary 3.18]{botnan2024signed} remain valid if each occurrence of ``interval'' is replaced by ``half-open rectangle''.
Alternatively, the proof of %
Theorem~\ref{thm:sufficient conditions for invertibility}~\ref{item:finitely presentable invertible}
can be easily adapted to prove the statement for half-open rectangles.

Furthermore, this same proof outline trivially extends to prove a stronger statement below. 
Let $Q$ be the minimal grid for $N$.
Let $\Ical$ be the collection of all segments \emph{and} half-open rectangles in $\R^d$.
Let \[\Jcal := \left\{\lfloor - \rfloor_{\catQ}\inv(I) \, : \, I\in \seg(\catQ)\right\},\]
a set of half-open rectangles in $\R^d$. Then, we have:

\begin{proposition}\label{prop:invertible over half-open}
    $\rk^\Ical_N$ is M\"obius invertible over $\Jcal$. 
\end{proposition}
Note that by Remark \ref{rem:monotonicity of mobius invertibility new} \ref{item:monotonicity of mobius invertibility}, 
this implies the statement that the GRI of $N$ over half-open rectangles is M\"obius invertible over half-open rectangles. 

\section{Additional details for Section \ref{sec:sufficient conditions}}

{In this section, we provide an example of an $\Int$-covering that is not a covering, as promised immediately before Theorem~\ref{thm:pull-back GPD} (Example~\ref{ex:counterexample}).}

\begin{example}\label{ex:counterexample} { 
Let $\catP$ be the union of five vertical or horizontal line segments in $\R^2$, depicted below. 
Let $\pi$ be the projection of $\catP$ to the $x$-axis, and let $\catQ$ be its image. 
Let both $\catP$ and $\catQ$ inherit the order of $\R^2$. Then it is not difficult to check that $\pi$ is an $\Int$-covering. However, $\pi$ is not a covering: 
Let $J:=\{(0,0),(1,0)\}\in \Con(\catQ)$. Its preimage under $\pi$ consists of two connected components, but neither of them map onto $J$.}

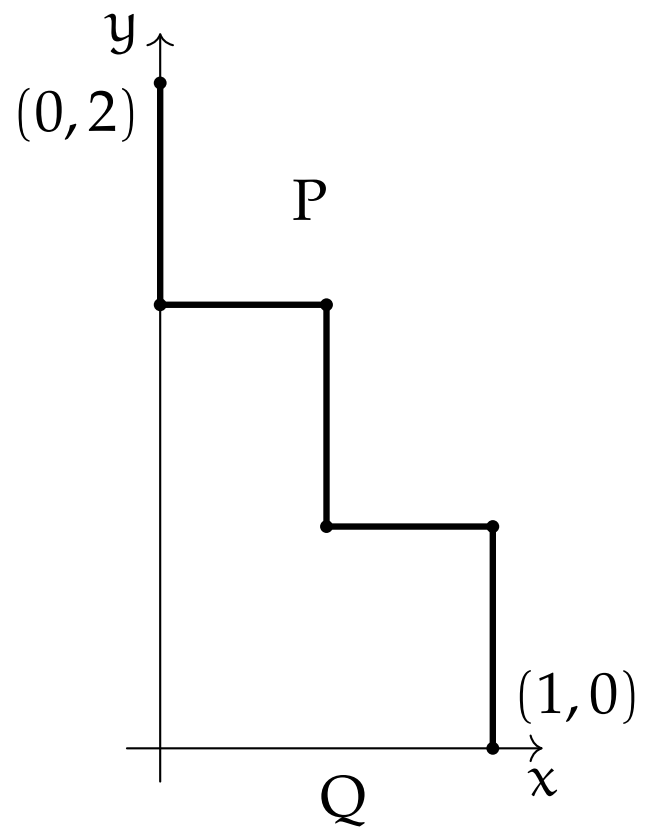
\begin{center}
\begin{tikzpicture}[scale=2.2, line cap=round, line join=round]
\draw[->] (-0.1,0) -- (1.15,0) node[below] {$x$};
\draw[->] (0,-0.1) -- (0,2.15) node[left] {$y$};

\def\ya{2/3}
\def\yb{4/3}

\draw[very thick] (0,\yb) -- (0,2);          %
\draw[very thick] (0,\yb) -- (1/2,\yb);      %
\draw[very thick] (1/2,\ya) -- (1/2,\yb);    %
\draw[very thick] (1/2,\ya) -- (1,\ya);      %
\draw[very thick] (1,0) -- (1,\ya);          %

\foreach \P in {(0,\yb),(0,2),(1/2,\yb),(1/2,\ya),(1,\ya),(1,0)}{
  \fill \P circle (0.55pt);
}

\node[left]  at (0,1.9)   {$(0,2)$};
\node[right] at (1,0.15)   {$(1,0)$};

\node at (0.45,1.65) {$\catP$};

\node[below] at (0.55,-0.02) {$\catQ$};

\end{tikzpicture}
\end{center}
 \end{example}

\section{Continuation of the Proof of Theorem \ref{thm:unique candidate} (Case 2 of Proposition \ref{prop:slice unequal to one side})
}\label{sec:case 2 details}

Recall that{, for the rank-critical interval $L$ chosen after the proof of Claim~\ref{claim:rank-critical},} we showed the existence of a unit block $B\in \Bcal$ with $\emptyset\subsetneq B\cap \Iast\subsetneq B$, thus the assumption of Proposition \ref{prop:slice unequal to one side} holds.
We have already dealt with Case 1 of the conclusion of the Proposition, so now assume that
Case 2 
holds, i.e., there exist $t^\ast\in [a_1,b_1):=\pi_1(B)$ and $\eps^\ast>0$ so that 
\begin{equation}\label{eq:closed interval around tast 2}
[t^\ast -\eps^\ast,t^{\ast}]\subset [a_1,b_1)
\end{equation}
and
\begin{equation}\label{eq:slices are different 2}
\mbox{for all $\eps \in (0, \epsilon^\ast]$,}\qquad S_{t^\ast} \neq  S_{t^\ast-\epsilon}.
\end{equation}

The following {remarks, propositions, and claims}
are analogues of those used in handling Case 1. {To facilitate comparison, we number them by suffixing the numbers of their counterparts  with the section letter C. For example, Remark~\ref{rem:q1 is at most tast 2} below is analogous to Remark~\ref{rem:q1 is at most tast}.}

\setcounter{proposition}{44}
\renewcommand{\theremark}{C.3.\arabic{remark}}

\begin{remark}\label{rem:q1 is at most tast 2}
By construction, $B\subset q_\ast^\uparrow$, and thus \[\mbox{(The first coordinate of $q_\ast$)}\leq a_1 \leq t^\ast-\eps^\ast.\]
\end{remark}
\textbf{For the remainder of this section,
, we continue to use the same $t^\ast$ and $\epsilon^\ast$ as above.}

Next, as before we will construct certain parametrized families of maps on $\R$ and $\R^d$ which we will use to construct families of intervals. 
For each $\epsilon\in [0,\epsilon^\ast)$, we consider the %
piecewise linear functions $g_\epsilon,f_\epsilon:\mathbb{R}\to\mathbb{R}$ whose graphs are depicted as in Figure~\ref{fig:g epsilon and f epsilon2}.
Note that both $g_0$ and $f_0$ are defined to be the the identity map on $\R$.%

As before, define the maps $G_\epsilon, F_\epsilon:\R^d\to \R^d$ via:
\begin{align*}
    G_\epsilon((x_1,x_2,\ldots,x_d)):= (g_\epsilon(x_1),x_2,\ldots,x_d) \ \mbox{and } 
    F_\epsilon((x_1,x_2,\ldots,x_d)):= (f_\epsilon(x_1),x_2,\ldots,x_d).
\end{align*}

\begin{figure}
\begin{subfigure}[c]{0.4\textwidth}
\begin{center}
\begin{tikzpicture}

\begin{axis}[
    axis lines=middle,
    axis equal image,
    xmin=-0.5, xmax=5.5,
    ymin=-0.5, ymax=5.5,
    xlabel={$t$},
    ylabel={$g_\epsilon(t)$},
    ylabel style={at={(axis description cs:0.5,0.9)}},
    ytick={1, 4},
    yticklabels={$t^\ast-\epsilon^\ast$,$t^\ast$},
    xtick={1,3,4},
    xticklabels={$t^\ast-\epsilon^\ast$,$t^\ast-\epsilon$,$t^\ast$},
    tick label style={font=\small},
]

\addplot[thick] coordinates {(-0.5,-0.5) (1,1)};

\addplot[thick] coordinates {(3,4) (4,4)};

\addplot[thick,domain=1:3,samples=2]
    {3/2 * x - 1/2};

\addplot[thick] coordinates {(4,4) (5.5,5.5)};

\addplot[only marks,mark=*]
coordinates {(1,1) (3,4) (4,4)};

\addplot[dashed] coordinates {(1,0) (1,6)};
\addplot[dashed] coordinates {(3,0) (3,5)};
\addplot[dashed] coordinates {(4,0) (4,5)};
\addplot[dashed] coordinates {(0,1) (5,1)};
\addplot[dashed] coordinates {(0,4) (5,4)};
\end{axis}

\end{tikzpicture}
\end{center}
\end{subfigure}\hspace{15mm}%
\begin{subfigure}[c]{0.4\textwidth}
\begin{center}
\begin{tikzpicture}

\begin{axis}[
    axis lines=middle,
    axis equal image,
    xmin=-0.5, xmax=5.5,
    ymin=-0.5, ymax=5.5,
    xlabel={$t$},
    ylabel={$f_\epsilon(t)$},
    ylabel style={at={(axis description cs:0.4,0.9)}},
    ytick={1,3,4},
    yticklabels={$t^\ast-\epsilon^\ast$,$t^\ast-\epsilon$,$t^\ast$},
    xtick={1,4},
    xticklabels={$t^\ast-\epsilon^\ast$,$t^\ast$},
    tick label style={font=\small},
]

\addplot[thick] coordinates {(-0.5,-0.5) (1,1)};

\addplot[thick,domain=1:4,samples=2]
    {1 + (2/3)*(x-1)};

\addplot[thick] coordinates {(4,4) (5.5,5.5)};

\addplot[only marks,mark=*]
coordinates {(1,1) (4,4)};

\addplot[only marks,mark=o]
coordinates {(4,3) (4,4)};

\addplot[dashed] coordinates {(1,0) (1,6)};
\addplot[dashed] coordinates {(4,0) (4,5)};
\addplot[dashed] coordinates {(0,1) (5,1)};
\addplot[dashed] coordinates {(0,3) (5,3)};
\addplot[dashed] coordinates {(0,4) (5,4)};
\end{axis}

\end{tikzpicture}

\end{center}
\end{subfigure}
\caption{
The graphs of the functions $g_\epsilon,f_{\epsilon}:\R\to\R$.
The composition
$g_\epsilon \circ f_\epsilon$ is the identity on $\R$, while $f_\epsilon\circ g_\epsilon$ is the identity on $\R\setminus (t^\ast-\epsilon^\ast, t^\ast-\epsilon]$.
}
\label{fig:g epsilon and f epsilon2}
\end{figure}

{The key difference from the functions $g_\epsilon$ and $f_\epsilon$ used in Case~1 of Proposition~\ref{prop:slice unequal to one side} is that, whereas previously we had
$g_\epsilon(x)\leq x\leq f_\epsilon(x)$
for all $x\in\R$, the new functions satisfy
$
g_\epsilon(x)\geq x\geq f_\epsilon(x).
$
The key similarities are that $g_\epsilon$ and $f_\epsilon$ remain order-preserving, and that $G_\epsilon$ and $F_\epsilon$ still map points in $B$ into $B$ (Proposition~\ref{prop:Feps and Geps maps B to B 2} below).}

As $g_\eps\circ f_\eps$ is the identity map on $\R$, we immediately have:

\setcounter{proposition}{38}
\renewcommand{\theremark}{C.3.\arabic{remark}}
\renewcommand{\theproposition}{C.3.\arabic{proposition}}

\begin{remark}\label{rem:GF is identity 2}
\begin{enumerate}[label=(\roman*)]
\item $G_\eps\circ F_\eps$ is the identity map on $\R^d$. \label{item:GF is identity1 2}
\item For any $I\in \Int(\R^d)$, we have $I=(G_\eps\circ F_{\eps})^{-1}(I)=F_\eps^{-1}\left(G_\eps^{-1}(I)\right)$. \label{item:GF is identity2 2}
\end{enumerate}
\end{remark}

\begin{proposition}\label{prop:Feps and Geps maps B to B 2} Let $B:=\prod_{i=1}^d[a_i,b_i)\subset \R^d$ with $[t^\ast-\eps^\ast,t^\ast]\subset [a_1,b_1)$. 
Then, the images of $B$ under the maps $G_\eps$ and $F_\eps$ are contained in $B$.
\end{proposition} 

\begin{proof} 
Let $x\in B$ and let $x_1$ be the first coordinate of $x$. 
It suffices to show that the images of $x_1$ under the maps $g_\eps$ and $f_\eps$ are in $[a_1,b_1)$. 
This follows from the fact that $x_1\in [a_1,b_1)$ and $[t^\ast - \epsilon^\ast, t^\ast]\subset [a_1,b_1)$, together with the observation that both $g_\epsilon$ and $f_\epsilon$ map $[t^\ast-\epsilon^\ast,t^\ast]$ into itself and fix every point outside $[t^\ast-\epsilon^\ast, t^\ast]$.
\end{proof}

\begin{proposition}\label{prop:G, F preserve qast 2}
Let $p,q_\ast\in \R^d$ and assume that the first coordinate of $q_\ast$ is at most $t^\ast - \epsilon^\ast$.
Then the following are equivalent:
\begin{enumerate*}[label=(\roman*)]
\item $p\geq q_\ast$ \label{item:p 2}
\item $G_\epsilon(p)\geq q_\ast$  \label{item:Geps p 2}
\item $F_\epsilon(p)\geq q_\ast$\label{item:Feps p 2}
\end{enumerate*} 
\end{proposition}

\begin{proof}
Observe that the images and preimages of $[t^\ast-\epsilon^\ast,\infty)$ under $g_\eps$ and $f_{\eps}$ are contained in $[t^\ast-\epsilon^\ast,\infty)$, and that, on the interval $(-\infty,t^\ast-\eps^\ast)$, $g_\eps$ and $f_{\eps}$ are the identity map. These observations imply that, for $t\in \R$ and the first coordinate $q_1$ of $q_\ast$, the following are equivalent: $t\geq q_1$, $g_\eps(t)\geq q_1$, and $f_\eps(t)\geq q_1$. This immediately implies the desired statement.
\end{proof}

\begin{proposition}\label{prop:Ie is an interval for all e2}
    For $J\in \Int(\R^d)$, 
    $J_\epsilon^G$ is an interval. 
\end{proposition}
The proof exactly as written for Proposition~\ref{prop:Ie is an interval for all e} holds again in the context of the new $G_\epsilon$, so we omit it here.

\begin{proposition}\label{prop:ranks of I and Ie are equal 2}
    For $J\in \Int(\R^d)$ with $J\subset q_\ast^\uparrow$, we have $\rk_M(J) = \rk_M(J_\epsilon^G)$. 
\end{proposition}
{The proof of Proposition~\ref{prop:ranks of I and Ie are equal} carries over verbatim to the new map $G_\epsilon$ and the corresponding auxiliary results, so we omit the details.} As in the setting of Case 1 of Proposition \ref{prop:slice unequal to one side}, our next goal is to show that $\{\Iast_\epsilon^G\}_{\eps\in [0,\eps^\ast)}$ is an  infinite family of intervals without pairwise containment. 

\setcounter{claim}{4}   %
\renewcommand{\theclaim}{C.\arabic{claim}}

\begin{claim}\label{claim:pairwise difference 2}
For $\epsilon<\epsilon'\in [0,\epsilon^\ast)$, we have that $\Iast_{\epsilon}^G\neq \Iast_{\epsilon'}^G$. 
\end{claim}
We remark that the rank-criticality of $L$ does not play any role in the following proof.
\begin{claimproof} 
   From Figure \ref{fig:g epsilon and f epsilon2}, observe that \begin{equation}\label{eq:strictly decreasing}
   t^\ast =g_{\eps'}(t^\ast-\eps') > g_{\epsilon}(t^\ast-\epsilon')\geq t^\ast - \epsilon^\ast.\end{equation} 
   Also, recall that, by the definitions of $L^G_\eps$ and $S_{g_\eps(t)}$,
    \begin{equation}\label{eq:some equation2}
\mbox{ for all $t\in \R$,  } \qquad\widehat{\pi}_1({L_\epsilon^G}\vert_{\pi_1=t}) = S_{g_\eps(t)}.
\end{equation}
    By Equations~\eqref{eq:slices are different 2}, \eqref{eq:strictly decreasing}, and \eqref{eq:some equation2},
    \[\widehat{\pi}_1(\Iast_{\epsilon'}^G\vert_{\pi_1=t^\ast-\epsilon'}) = S_{g_{\eps'}(t^\ast-\eps')}=S_{t^\ast} \neq S_{g_{\epsilon}(t^\ast-\epsilon')} = \widehat{\pi}_1(\Iast_{\epsilon}^G\vert_{\pi_1=t^\ast-\epsilon'}) .\]
    This implies $\Iast_{\epsilon'}^G\vert_{\pi_1=t^\ast-\epsilon'}\neq \Iast_{\epsilon}^G\vert_{\pi_1=t^\ast-\epsilon'}$, and thus we have $\Iast_{\epsilon}^G\neq \Iast_{\epsilon'}^G$.
\end{claimproof}

\begin{claim}\label{claim:slice condition implications IFe 2}
    Let $I\in \Int(\R^d)$, and $\epsilon \geq 0$ such that for all $t\in [t^\ast-\epsilon, t^\ast]$ %
    it holds that $\widehat{\pi}_1(I\vert_{\pi_1=t^\ast}) = \widehat{\pi}_1(I\vert_{\pi_1=t})$.
    Then we have the following:
    \begin{enumerate}[label=(\roman*)]
            \item $I = (I_{\epsilon}^F)_{\epsilon}^G$ (note that, by Remark~\ref{rem:GF is identity 2}~\ref{item:GF is identity2 2}, we readily know $I = (I_{\epsilon}^G)_{\epsilon}^F$). \label{item:F then G equals original 2}  
        \item $I_{\epsilon}^F$ is an interval. 
        \label{item:IFe is an interval 2} 
    \end{enumerate}
\end{claim}
Before proving the claim, it will be useful to note that Item \ref{item:F then G equals original 2} implies the following. 
\begin{equation}\label{eq:x in I then Gex in 2}
   \mbox{If $x\in I$, then $G_\epsilon(x)\in I_{\epsilon}^F$}. 
\end{equation}

\begin{claimproof}
\noindent    \ref{item:F then G equals original 2}:
    First, we show that $I\subset (I_{\epsilon}^F)_{\epsilon}^G$.
    It suffices to show that, if $x\in I$, then $F_\eps ( G_{\eps}(x))\in I$. 
    Let $x =: (x_1,x')\in I$ with $x_1\in \R$ and $x'\in \R^{d-1}$. 
    Observe that if $x_1\in [t^\ast-\epsilon, t^\ast]$, then $f_{\epsilon}(g_{\epsilon}(x_1)) = t^\ast$, and if not, then $f_{\epsilon}(g_{\epsilon}(x_1))=x_1$ (cf. the caption of Figure~\ref{fig:g epsilon and f epsilon2}).
    Hence, we have
    \begin{equation}\label{eq:FG cases 2}
    F_{\epsilon}(G_{\epsilon}(x)) = 
    \begin{cases}
    (t^\ast,x') & \mbox{if } x_1\in [t^\ast-\epsilon, t^\ast],\\
        x & \mbox{otherwise.}       
    \end{cases}
    \end{equation}
    Thus, if $x_1\notin [t^\ast-\eps, t^\ast]$, then $F_{\epsilon}(G_{\epsilon}(x)) = x\in I$ as desired.
    Now assume that $x_1\in [t^\ast-\eps, t^\ast]$. Since $x\in I$, the point $x'$ belongs to the set $ \widehat{\pi}_1(I\vert_{\pi_1=x_1})$ which is, by
    assumption, 
    identical to $ \widehat{\pi}_1(I\vert_{\pi_1=t^\ast})$.
    Thus, $(t^\ast,x')\in I$.
    By Equation \eqref{eq:FG cases 2}, $F_\epsilon(G_\epsilon(x)) = (t^\ast,x')$, and thus $F_\epsilon(G_\epsilon(x))\in I$. 
    
     Next, we show that $ (I_{\epsilon}^F)_{\epsilon}^G\subset I$. 
     Suppose $x\in (I_{\epsilon}^F)_{\epsilon}^G$, i.e., there exists $y\in I$ such that $F_\eps(G_\eps(x))=y$. 
     By Equation \eqref{eq:FG cases 2}, if $x_1\notin [t^\ast-\eps, t^\ast]$, then $F_\epsilon(G_\epsilon(x)) = x$, which implies $y=x$, and thus $x\in I$ since $y$ is.    
     If $x_1\in [t^\ast-\eps, t^\ast]$, then $y = F_\epsilon(G_\epsilon(x)) = (t^\ast, x')$.    
     Since $y\in I$, the point $x'$ belongs to the set $ \widehat{\pi}_1(I\vert_{\pi_1=t^\ast})$, which is by our assumption identical to $\widehat{\pi}_1(I\vert_{\pi_1=x_1})$. This shows that $x=(x_1,x')\in I$. 
   
\noindent   \ref{item:IFe is an interval 2}:
    First, we show that $I_{\epsilon}^F\neq \emptyset$.
    As $I\in \Int(\R^d)$, $I\neq \emptyset$, so pick any $x=(x_1,x')\in I$ with $x_1\in \R$ and $x'\in \R^{d-1}$.
    Suppose $x_1\in \R\setminus [t^\ast-\eps, t^\ast]$.
    Then $f_\epsilon(x) = x$, so $x\in I_{\epsilon}^F$, implying that $I_{\epsilon}^F\neq \emptyset$.
    Suppose $x_1\in [t^\ast-\eps, t^\ast]$.
    Then we have $x'\in \widehat{\pi}_1(I\vert_{\pi_1=x_1})=\widehat{\pi}_1(I\vert_{\pi_1=t^\ast})$, by our assumption. 
    Thus, $(t^\ast, x')\in I$.
    As $f_{\epsilon}(t^\ast) = t^\ast$, we have that $F_\epsilon((t^\ast,x')) = (t^\ast, x')\in I$, and thus $(t^\ast, x')\in I_{\epsilon}^F$. This proves that $I_{\epsilon}^F\neq \emptyset$.

    By Proposition \ref{prop:each component is an interval}, $I_{\epsilon}^F$ is convex, and it remains to show connectivity.
    Let $p,q\in I_{\epsilon}^F$, which implies that  $F_\epsilon(p), F_\epsilon(q)\in I$.
    As $I$ is an interval, there is a path in $I$ connecting $F_\epsilon(p)$ and $F_\epsilon(q)$.
    By the fact that $G_\epsilon$ is order-preserving, alongside Equation~\eqref{eq:x in I then Gex in 2}, the image of this path under $G_\epsilon$ is a path in $I_{\epsilon}^F$ 
    connecting $G_\epsilon(F_\epsilon(p)) = p$ to $G_\epsilon(F_\epsilon(q)) = q$ (cf. Remark~\ref{rem:GF is identity 2}~\ref{item:GF is identity1 2}), as desired.
\end{claimproof}

\begin{claim}\label{claim:Ie distinct for different e 2}
    For $\epsilon<\epsilon'\in [0,\epsilon^\ast)$, we have that  $\Iast_{\epsilon}^G\not\subset \Iast_{\epsilon'}^G$.
\end{claim}
\begin{claimproof}
    Assume 
    $\Iast_{\epsilon}^G\subset \Iast_{\epsilon'}^G$ for some fixed $\eps<\eps'\in [0,\eps^\ast)$, which will lead to a contradiction to the rank-criticality of $\Iast$.   
    {Namely, we show that $(\Iast_{\epsilon'}^G)_\epsilon^F$ is an interval that properly contains $L$, but has the same $\rk_M$ value as $L$.} We first note the following:
   
   \begin{itemize} 
   \item By Claim \ref{claim:pairwise difference 2}, we have $\Iast_{\epsilon}^G\neq \Iast_{\epsilon'}^G$, and thus 
    \begin{equation}\label{eq:IastGeps is smaller 2}
    \Iast_{\epsilon}^G\subsetneq \Iast_{\epsilon'}^G.
    \end{equation}
  \item   As $L\subset q_\ast^\uparrow$, Remark~\ref{rem:q1 is at most tast 2} and
    Proposition~\ref{prop:G, F preserve qast 2} imply that $\Iast_{\epsilon'}^G\subset q_\ast^\uparrow$, and in turn 
    \begin{equation}\label{eq:IastGF is in the quadrant 2} (\Iast_{\epsilon'}^G)_\epsilon^F\subset q_\ast^\uparrow.
    \end{equation}
\end{itemize}

    Now we claim that $(\Iast_{\epsilon'}^G)_\epsilon^F\supsetneq \Iast$.
    Let $x\in L$.
    By Remark~\ref{rem:GF is identity 2}~\ref{item:GF is identity2 2}, $(\Iast_{\epsilon}^G)_\epsilon^F=\Iast$,  and thus $x\in (\Iast_{\epsilon}^G)_\epsilon^F$.
    Then $F_{\epsilon}(x)\in \Iast_{\epsilon}^G\subset \Iast_{\epsilon'}^G$, which by definition implies $x\in (\Iast_{\epsilon'}^G)_\epsilon^F$, as desired.
    Next, we show $(\Iast_{\epsilon'}^G)_\epsilon^F\neq \Iast$.
    Note that for all $t\in [t^\ast-\eps, t^\ast]$,  $g_\epsilon(t) = g_{\epsilon'}(t) = t^\ast$ {(here, we use the assumption $\eps<\eps'$)}, and hence Equation \eqref{eq:some equation2} implies 
    
    \begin{equation}\label{eq:both are Stast 2}\mbox{for all $t\in [t^\ast-\eps,t^\ast]$,}\qquad\widehat{\pi}_1(\Iast_{\epsilon}^G\vert_{\pi_1=t}) = \widehat{\pi}_1(\Iast_{\epsilon'}^G\vert_{\pi_1=t}) = S_{t^\ast}.\end{equation}    
    Hence for any $y=(y_1,\ldots,y_d)\in \Iast_{\epsilon'}^G\setminus \Iast_{\epsilon}^G$ (which exists by Equation~\eqref{eq:IastGeps is smaller 2}), we have $y_1\notin [t^\ast-\eps, t^\ast]$; otherwise $(y_2,\ldots,y_d)\in \widehat{\pi}_1(\Iast_{\epsilon'}^G\vert_{\pi_1=y_1}) = \widehat{\pi}_1(\Iast_{\epsilon}^G\vert_{\pi_1=y_1})$, implying that $y\in L_\eps^G$, a contradiction.
    As a result, $y_1\in \im(f_{\epsilon})$, and in turn $y\in \im(F_{\epsilon})$.
    Thus, there exists $y'\in \R^d$ with $F_{\epsilon}(y')=y$, which by definition implies $y'\in (\Iast_{\epsilon'}^G)_\epsilon^F$. To show  $(\Iast_{\epsilon'}^G)_\epsilon^F\neq \Iast$, it suffices to show that $y'\not\in L$.
    Since $G_\eps\circ F_\eps$ is the identity map on $\R^d$ (cf. Remark~\ref{rem:GF is identity 2}~\ref{item:GF is identity1 2}), we have  $y' = G_{\epsilon}(F_{\epsilon}(y')) = G_{\epsilon}(y)$. But, as $y\notin \Iast_{\epsilon}^G$, we have $y'=G_{\epsilon}(y)\notin \Iast$, as desired. 

    Next we claim that $(\Iast_{\epsilon'}^G)_\epsilon^F$ is an interval. Note, by Proposition~\ref{prop:Ie is an interval for all e2}, that $L_{\eps'}^G$ is an interval. Hence,  Equation~\eqref{eq:both are Stast 2} 
and Claim~\ref{claim:slice condition implications IFe 2} \ref{item:IFe is an interval 2} imply the claim.

Lastly, we claim that $L$ and $(\Iast_{\epsilon'}^G)_\epsilon^F$ have the same $\rk_M$ value. Indeed,
\begin{align*}
\rk_M((\Iast_{\epsilon'}^G)_\epsilon^F)&=\rk_M(((\Iast_{\epsilon'}^{G})_{\epsilon}^{F})_{\epsilon}^{G})&\mbox{by Equation~\eqref{eq:IastGF is in the quadrant 2} and  Proposition~\ref{prop:ranks of I and Ie are equal 2}}
\\&=\rk_M(\Iast_{\epsilon'}^G)&\mbox{Claim \ref{claim:slice condition implications IFe 2} \ref{item:F then G equals original 2} implies $((\Iast_{\epsilon'}^{G})_{\epsilon}^{F})_{\epsilon}^{G} = \Iast_{\epsilon'}^G$}
\\&=\rk_M(\Iast)&\mbox{by the fact $L\subset q_\ast^\uparrow$ and Proposition~\ref{prop:ranks of I and Ie are equal 2},}
\end{align*}
completing the proof.\qedhere
\end{claimproof}

\begin{claim}\label{claim:the other containment 2}
For $\epsilon<\epsilon'\in [0,\epsilon^\ast)$, we have that  $\Iast_{\epsilon'}^G\not\subset \Iast_{\epsilon}^G$.
\end{claim} 
\begin{claimproof}
    We assume that $\Iast_{\epsilon'}^G\subset \Iast_{\epsilon}^G$ for some fixed $\eps<\eps'\in [0,\eps^\ast)$ and deduce a contradiction to the rank criticality of $L$.
    Let $t'\in [t^\ast - \epsilon', t^\ast - \epsilon)$.
    By the definitions of the functions $g_{\epsilon}$ and $g_{\epsilon'}$, we have that 
    \begin{equation}\label{eq:y and tast 2}        
    t^\ast = g_{\epsilon'}(t') > g_{\epsilon}(t') =: t \geq t^\ast-\epsilon^\ast.
    \end{equation}
    Then, we obtain the following sequence of inclusions: 
    \begin{align*}
    &\widehat{\pi}_1(\Iast_{\epsilon'}^G\vert_{\pi_1=t'})
    \subset \widehat{\pi}_1(\Iast_{\epsilon}^G\vert_{\pi_1=t'})  &\mbox{since $\Iast_{\epsilon'}^G\subset \Iast_{\epsilon}^G$}
\\ &\implies
    S_{g_{\epsilon'}(t')}
    \subset S_{g_{\epsilon}(t')}&\mbox{by Equation~\eqref{eq:some equation2}}
\\&\implies S_{t^\ast}
 \subset S_t &\mbox{by Equation~\eqref{eq:y and tast 2}.}
\end{align*}
  Further, 
  by Equations~\eqref{eq:slices are different 2} and \eqref{eq:y and tast 2}, 
  we have $S_{t^\ast}\neq S_{t}$, 
  and thus $S_{t^\ast}\subsetneq S_{t}$. 
  This implies that the set \[I:=\Iast \cup \left([t,t^\ast]\times S_{t}\right)\] 
    properly contains $L$; namely, $\emptyset\subsetneq\{t^\ast\}\times (S_t\setminus S_{t^\ast})\subset I\setminus L$. Hence, it suffices to show that $I$ is an interval and $\rk_M(I)=\rk_M(L)$. 
     
    We define the function $h:\R\to \R$ as 
    \[h(s):=\begin{cases}
      t, & \mbox{if } s\in (t,t^\ast].\\
        s, & \mbox{otherwise.}
    \end{cases}\]
    We also define $H:\R^d\to \R^d$ by $H((x_1,x_2,\ldots,x_d)):=(h(x_1),x_2,\ldots,x_d)$.
    It is clear that for $x\in \R^d$, $x$ and $H(x)$ can differ in only the first coordinate and that $x\geq H(x)$.

    We claim that the image of $I$ under $H$ is contained in $L$.
    Let $x\in I$.

    ($\ast$) Suppose $x\notin \Iast$. Then by construction of $I$, $x = (x_1, x')$ for some  $x_1\in [t, t^\ast]$ and  $x'\in S_{t}$.
    Then $H(x) = (t, x')$.
    As $x'\in S_t$, we have that $H(x)\in \Iast$, as desired.
    
    Now suppose $x\in \Iast$.
    If $H(x)=x$, then trivially, $H(x)\in \Iast$. Let us assume that $H(x)\neq x$.
    This implies that the first coordinate of $x$ is not fixed by $h$, i.e., $x = (x_1,x')$ for some $x_1\in (t, t^\ast]$ and $x'\in S_{x_1}$.
    If $x_1=t^\ast$, then we have $H(x) = (t,x')\in \Iast$ from the fact $S_{t^\ast}\subset S_{t}$.
    Now, assume $x_1\in (t, t^\ast)$.
    Recall that $t = g_{\epsilon}(t')$ for $t'\in [t^\ast-\eps', t^\ast-\eps)$.
    Then consider $g_{\epsilon}(f_{\epsilon'}(x_1))$.
    As $t^\ast>x_1$, 
    \[t' \geq t^\ast-\epsilon' \geq f_{\epsilon'}(x_1),\]
    and as $g_{\epsilon}$ is order-preserving, this gives 

    \[t = g_{\epsilon}(t') \geq s := g_{\epsilon}(f_{\epsilon'}(x_1)).\]
    By the assumption $\Iast_{\epsilon'}^G\subset \Iast_{\epsilon}^G$, 
    we have that $\widehat{\pi}_1(\Iast_{\epsilon'}^G\vert_{\pi_1=f_{\epsilon'}(x_1)}) \subset \widehat{\pi}_1(\Iast_{\epsilon}^G\vert_{\pi_1=f_{\epsilon'}(x_1)})$, 
    which by Equation \eqref{eq:some equation} and the fact $g_{\epsilon'}(f_{\epsilon'}(x_1)) = x_1$, implies $S_{ x_1}\subset S_s$.
    Hence, as $x = (x_1, x')\in \Iast$,     %
    we have that $(s, x')\in \Iast\vert_{\pi_1=s}\subset L$.
    Also, since $x_1>t\geq s$, we have \[(x_1,x')\geq (t, x')=H(x) \geq (s,x').\] 
    Then, convexity of $\Iast$ and the fact that $(x_1,x'),(s,x')\in \Iast$ imply  $H(x) \in \Iast$. 
    We have proved that the image of $I$ under $H$ is contained in $L$, as desired.
  
    Now, we will show that $I$ is an interval. 
    To show connectivity of $I$, suppose $p,q\in I$. 
    Then $H(p),H(q)\in \Iast$.
    As $\Iast$ is connected, there exists a path in $\Iast$ (hence in $I$) connecting     $H(p)$ and $H(q)$.
    As $p \geq H(p)$ and $q\geq H(q)$, by prepending $p$ and appending $q$ to the path,
    we obtain a path in $I$ connecting $p$ to $q$. 
    Next we show convexity of $I$. 
    Let $p<q<r\in \R^d$ with $p,r\in I$, and we show that $q\in I$. %
    First, observe that since $H$ is order-preserving, we have 
    \[(\ast\ast)\qquad \mbox{$H(p)\leq H(q)\leq H(r)$ in $\R^d$, and thus convexity of $L$ implies that $H(q)\in L$.}\]
    The claim that $q\in I$ follows by considering the following cases.

    \begin{itemize}
    \item Assume $p,r\in \Iast$. Then $q\in \Iast$ by convexity of $\Iast$, and thus $q\in I$.

    \item Assume $p\in I\setminus \Iast$ and $r\in \Iast$.
    Since $H(p)\leq p$ and in $H(p)\in \Iast$, we have that $H(p)<q<r$ with $H(p),r\in \Iast$, and thus convexity of $L$ implies $q\in \Iast\subset I$. 
    
    \item Assume $p,r\in I\setminus \Iast$.  
    We saw from ($\ast$) that $H(p),H(r)\in \Iast\vert_{\pi_1=t}$. Recall also ($\ast\ast$). 
    As the first coordinates of $H(p)$ and $H(r)$ are both equal to $t$, this also holds for $H(q)$, i.e.,
    $H(q)\in \Iast\vert_{\pi_1=t}$.
    Then, since $q$ and $H(q)$ can differ only in the first coordinate, $\widehat{\pi}_1(q)=\widehat{\pi}_1(H(q))\in S_t$. 
    In addition, since $H(q)$ is obtained by applying $h$ to the first coordinate of $q$ while leaving the other coordinates unchanged, the fact that $h^{-1}(t)=[t,t^\ast]$ 
    implies that the first coordinate of $q$ is in $[t,t^\ast]$. In sum, 
\[
q\in [t,t^\ast]\times S_t\subset I.
\]

    \item Assume $p\in \Iast$ and $r\in I\setminus \Iast$.
    Let us write $p=(p_1,p'), q=(q_1,q')$ and $r=(r_1,r')$ with $p_1,q_1,r_1\in \R$ and $p',q',r'\in \R^{d-1}$.
    As $p<q<r$, we have $p_1\leq q_1\leq r_1$ and $p'\leq q'\leq r'$.
    Observe that the assumption $r\in I\setminus \Iast$ implies $r\in (t,t^\ast]\times S_t$, so in particular, we have $t^\ast \geq r_1>t$.
    Since $r_1\geq q_1$, this implies $t^\ast \geq q_1$.
    
    Recall $(\ast\ast)$.
    If $q_1\in (-\infty, t]$,
    then $q=H(q)\in \Iast$.
    Assume $q_1\in (t,\infty)$. 
    Then we have $t^\ast \geq q_1>t$, and so $H(q) = (t,q')$.
    As $H(q)\in \Iast$, we have $H(q)\in \Iast\vert_{\pi_1=t}$, which implies $q'\in S_t$. In sum, we have 
    $q\in [t,t^\ast]\times S_t\subset I$, as desired.
\end{itemize}

    Finally, we prove that $\rk_M(I)=\rk_M(\Iast)$. 
    By monotonicity of the GRI, it suffices to show that $\rk_M(I)\geq \rk_M(\Iast)$. 
    We do this by showing $I\subset q^\ast$ and constructing a map $\Psi:\CC(I) \to \CC(\Iast)$  satisfying all assumptions of Proposition~\ref{prop:connected components to rank general}.
        First, we claim that $I\subset q_\ast^\uparrow$.
    Since $\Iast\subset q_\ast^\uparrow$, it suffices to show that $[t, t^\ast]\times S_t\subset q_\ast^\uparrow$. 
    {Then, by Remark~\ref{rem:q1 is at most tast 2}, it suffices to show that $S_t\in \widehat{\pi}_1(q_\ast^\uparrow)$. This follows from
    \begin{align*}
        S_t&=\widehat{\pi}_1(L|_{\pi=t})&\mbox{by definition}
        \\&\subset \widehat{\pi}_1(L)
        \\&\subset \widehat{\pi}_1(q_\ast^\uparrow)&\mbox{since $L\subset q_\ast^\uparrow$.}
    \end{align*}
    }

    Now we define $\Psi:\CC(I) \to \CC(\Iast)$ by $\Psi(X)= Y$ if $X\cap Y\neq \emptyset$. 
    {We claim that for all $X\in \CC(I)$,  $\Psi(X)$ exists in $\CC(\Iast)$. 
    Fix $X\in \CC(I)$. 
    Verifying that $H(X)\subset X$ suffices: Indeed, since $H(I)\subset \Iast\subset q_\ast^\uparrow$, the containment $H(X)\subset X$ implies that for any $x\in X$, there exists $Y\in \CC(\Iast)$ such that $H(x)\in X\cap Y$, and thus $X\cap Y\neq \emptyset$. 
    Now observe that 
    \[(\ast\ast\ast)\qquad \mbox{$X$ and $H(X)$ are contained in the same unit block}\] because (1) $H$ only affects the first coordinate of its input via $h$, (2)  $h([t,t^\ast])=\{t\}\subset [t,t^\ast]\subset \pi_1(B)$ (cf. Equations~\eqref{eq:closed interval around tast 2} and \eqref{eq:y and tast 2}), and (3) $h$ fixes points outside $[t,t^\ast]$. Also, since $H$ is order-preserving and $X$ is connected, $H(X)$ is also connected. Furthermore, since  $x\geq H(x)$ for each $x\in X$, and $X$ is a connected component (within the unit block containing $X$), we have $H(X)\subset X$.}  

    Now we claim that $\Psi$ is well-defined, i.e., $\Psi(X)$ is unique for each $X\in \CC(I)$.
    Assume $x,x'\in X\in \CC(I)$ with $x\in Y$ and $x'\in Y'$ for some $Y,Y'\in \CC(\Iast)$.
    We show $Y=Y'$.
    Since 
    both $x$ and $H(x)$ belong to $L$, and  $x\geq H(x)$, the fact in $(\ast\ast\ast)$ and the assumption $x\in Y$ imply that 
    $H(x)\in Y$. Similarly, $H(x')\in Y'$.
   As $x,x'\in X$, 
    there is a path 
    in $X$
    connecting $x$ and $x'$.
    Then the image of the path connects $H(x)$ and $H(x')$ in $\Iast$.
    As $H$ is order-preserving, the fact $(\ast\ast\ast)$ implies that $H(x)$ and $H(x')$ belong to the same element of $\CC(L)$, i.e.,
    $Y=Y'$ as desired.

    {It only remains to show that $\Psi$ satisfies assumptions \ref{item:Psi condition 1} and \ref{item:Psi condition 3} in Proposition~\ref{prop:connected components to rank general}. %
    Assumption \ref{item:Psi condition 1}: Since any $X\in \CC(I)$ and its image $\Psi(X)\in \CC(\Iast)$ have nonempty intersection, by the definitions of $\CC(I)$ and $\CC(\Iast)$, $X$ and $\Psi(X)$ are in the same unit block.
    \newline \noindent Assumption \ref{item:Psi condition 3}: Let $(X,X')\in \CC^2_{\leq}(I)$, i.e., there exists $x\in X$ and $x'\in X'$ such that $x\leq x'$. Recall that $\Psi(X)$ and $\Psi(X')$ contain $H(x)$ and $H(x')$, respectively. Since $H$ is order-preserving, we have $H(x)\leq H(x')$ as desired.}
\end{claimproof}

Let $\epsilon<\epsilon'\in [0,\eps^\ast]$, and let $J\subset \R^d$ be the convex hull of $\Iast_{\epsilon}^G\cup \Iast_{\epsilon'}^G$ (cf. Remark~\ref{rem:properties of finitely generated modules}~\ref{item:convex hull is smallest interval}).
The following claim regarding slices of $J$ will be useful in the sequel:

\begin{claim} For all  $t\in [t^\ast - \epsilon, \infty)$,
     {\begin{equation}\label{eq:slices of J equal slices of Ie simpler 2}
        J\vert_{\pi_1=t} = \Iast_{\epsilon}^G\vert_{\pi_1=t} = \Iast_{\epsilon'}^G\vert_{\pi_1=t}.
    \end{equation}  } 
\end{claim}
{We note that, by Remark~\ref{rem:slice projection vs. slice}, this claim is equivalent to:}
    \begin{equation}\label{eq:slices of J equal slices of Ie 2}
        \mbox{for all } t\in [t^\ast-\eps,\infty),\hspace{5mm} 
        \widehat{\pi}_1(J\vert_{\pi_1=t}) = \widehat{\pi}_1(\Iast_{\epsilon}^G\vert_{\pi_1=t}) = \widehat{\pi}_1(\Iast_{\epsilon'}^G\vert_{\pi_1=t}).        
    \end{equation}
     
\begin{claimproof}
 First, from Figure~\ref{fig:g epsilon and f epsilon2}, observe that 
    \begin{equation}\label{eq:identity map within a range 2}\mbox{for all $t\in \R\setminus (t^\ast-\eps^\ast,t^\ast)$, }\ \  g_\epsilon(t) = g_{\epsilon'}(t) = t.\end{equation}

  We prove the right equality in Equation \eqref{eq:slices of J equal slices of Ie 2} {(and thus the one in Equation~\eqref{eq:slices of J equal slices of Ie simpler 2})}.
      If $t \in (t^\ast,\infty)$, by Equation~\eqref{eq:identity map within a range 2}, we have $g_\epsilon(t)=g_{\epsilon'}(t) = t$, and thus Equation \eqref{eq:some equation2} implies that both $\widehat{\pi}_1(\Iast_{\epsilon}^G\vert_{\pi_1=t})$ and $\widehat{\pi}_1(\Iast_{\epsilon'}^G\vert_{\pi_1=t})$ equal $S_t$, establishing the right equality.
     If $t\in [t^\ast-\eps, t^\ast]$, we have $g_\epsilon(t) = g_{\epsilon'}(t) = t^\ast$, and hence Equation \eqref{eq:some equation2} again implies the right equality.
    
    We show the left equality in Equation \eqref{eq:slices of J equal slices of Ie simpler 2} {(and thus the one in Equation~\eqref{eq:slices of J equal slices of Ie 2})}.
    {From the definition $J=\overline{\Iast_{\epsilon}^G\cup \Iast_{\epsilon'}^G}$, it is immediate that $J\vert_{\pi_1=t}\supset \Iast_{\epsilon}^G\vert_{\pi_1=t}$. We show  $J\vert_{\pi_1=t}\subset\Iast_{\epsilon}^G\vert_{\pi_1=t}(=\Iast_{\epsilon'}^G\vert_{\pi_1=t})$.}
    Let $p\in J\vert_{\pi_1=t}$. From the definition of convex hull, $p\in J$ implies that there exist $q,r\in \Iast_{\epsilon}^G\cup \Iast_{\epsilon'}^G$ such that $q\leq p\leq r$.
    \begin{itemize}
    \item 
    If $q, r\in \Iast_{\epsilon}^G$, then $p\in \Iast_{\epsilon}^G$ by convexity of $\Iast_{\epsilon}^G$ (cf. Proposition~\ref{prop:Ie is an interval for all e2}), and as $\pi_1(p)=t$, we have $p\in \Iast_{\epsilon}^G\vert_{\pi_1=t}$. 
    \item 
    Similarly, if $q,r\in \Iast_{\epsilon'}^G$, then $p\in \Iast_{\epsilon'}^G|_{\pi_1=t}$.
    \item If $q\in \Iast_{\epsilon}^G$ and $r\in \Iast_{\epsilon'}^G$, then let $r=:(r_1, r')$ with $r_1\in \R$, $r'\in \R^{d-1}$. 
    Also, let $p=:(t, p')$ for some $p'\in \R^{d-1}$.
    Since $r\geq p$, we have $r_1\geq t\geq t^\ast - \epsilon$. 
    {Thus, the right equality in Equation  \eqref{eq:slices of J equal slices of Ie 2} gives
    $\widehat{\pi}_1(\Iast_\eps^G\vert_{\pi_1=r_1}) = \widehat{\pi}_1(\Iast_{\eps'}^G\vert_{\pi_1=r_1})$. Hence, as $r'\in \widehat{\pi}_1(\Iast^G_{\eps'}\vert_{\pi_1=r_1})$, we have $r'\in \widehat{\pi}_1(L^G_{\eps}\vert_{\pi_1=r_1})$, which implies $r\in \Iast^G_{\eps}$.
    Thus, $q,r\in \Iast^G_{\eps}$, the case considered in the first bullet, and thus $p\in \Iast^G_{\eps}|_{\pi_1=t}$.
    }
   \item {If $q\in \Iast_{\epsilon'}^G$ and $r\in \Iast_{\epsilon}^G$, then
$r\in \Iast_{\epsilon'}^G$ by interchanging the roles of $\epsilon$ and
$\epsilon'$ in the preceding argument. Hence, we are again in the situation
described in the second bullet point.} %
\end{itemize}
We have proved the left equality in Equation \eqref{eq:slices of J equal slices of Ie simpler 2}, completing the proof.    
\end{claimproof}   

We proceed by demonstrating that if an interval in $\R^d$ contains two distinct intervals from $\{\Iast_{\epsilon}^G\}_{\eps \in [0,\eps^\ast)}$, then it must have strictly smaller rank than $L$:

\begin{claim}\label{claim:rank of interval containing two decreases 2}
    Let $\epsilon<\epsilon'\in [0,\epsilon^\ast)$.
    For any $I\in \Int(\R^d)$ containing both $ \Iast_{\epsilon}^G$ and $\Iast_{\epsilon'}^G$,
    we have that $\rk_M(I) < \rk_M(\Iast)$.
\end{claim}
\begin{claimproof}
   Let $J\subset \R^d$ be the convex hull of $\Iast_{\epsilon}^G\cup \Iast_{\epsilon'}^G$ (cf. Remark~\ref{rem:properties of finitely generated modules}~\ref{item:convex hull is smallest interval}). It suffices to show that $J$ is the smallest interval containing $\Iast_{\epsilon}^G\cup \Iast_{\epsilon'}^G$ and that $\rk_M(J)<\rk_M(L)$. 
   
   By Remark~\ref{rem:properties of finitely generated modules}~\ref{item:convex hull is smallest interval}, proving the former reduces to showing that $\Iast_{\epsilon}^G\cup \Iast_{\epsilon'}^G$ is connected.
    Taking $t=t^\ast-\eps$ and $t=t^\ast$ in Equation~\eqref{eq:identity map within a range 2}, we obtain from Equation~\eqref{eq:some equation2} that
\[\widehat{\pi}_1(\Iast_{\epsilon}^G\vert_{\pi_1=t^\ast-\eps^\ast}) = \widehat{\pi}_1(\Iast_{\epsilon'}^G\vert_{\pi_1=t^\ast-\eps^\ast}) = S_{t^\ast-\eps^\ast} \mbox{ and } \widehat{\pi}_1(\Iast_{\epsilon}^G\vert_{\pi_1=t^\ast}) = \widehat{\pi}_1(\Iast_{\epsilon'}^G\vert_{\pi_1=t^\ast}) = S_{t^\ast}.\]
  Then, by Remark~\ref{rem:slice projection vs. slice}, we have
    \begin{equation}\label{eq:two equations 2}
    \Iast_{\epsilon}^G\vert_{\pi_1=t^\ast-\eps^\ast} = \Iast_{\epsilon'}^G\vert_{\pi_1=t^\ast-\eps^\ast} = \Iast\vert_{\pi_1=t^\ast-\eps^\ast}
    \mbox{ and } 
    \Iast_{\epsilon}^G\vert_{\pi_1=t^\ast} = \Iast_{\epsilon'}^G\vert_{\pi_1=t^\ast} = \Iast\vert_{t^\ast}.
    \end{equation}
    By Equation~\eqref{eq:slices are different 2}, we have $S_{t^\ast-\eps^\ast}\neq S_{t^\ast}$, and thus at least one of $S_{t^\ast-\eps^\ast}$ and $S_{t^\ast}$ is nonempty, which in turn implies that at least one of $\Iast\vert_{\pi_1=t^\ast-\eps^\ast}$ and $\Iast\vert_{\pi_1=t^\ast}$ is nonempty.
    Let us pick a point $p\in \left(\Iast\vert_{\pi_1=t^\ast-\eps^\ast}\right)\cup \left( \Iast\vert_{\pi_1=t^\ast}\right)$. 
    By Equation~\eqref{eq:two equations 2}, $p$ belongs to both $\Iast_{\epsilon}^G$ and $\Iast_{\epsilon'}^G$. %
    As $\Iast_{\epsilon}^G$ and $\Iast_{\epsilon'}^G$ are both intervals (cf. Proposition~\ref{prop:Ie is an interval for all e2}), they are connected.
    Thus, by Remark~\ref{rem:properties of finitely generated modules}~\ref{item:union of connected sets with overlap is connected}, $\Iast_{\epsilon}^G\cup \Iast_{\epsilon'}^G$ is connected, as desired.

    It remains to show that $\rk_M(J)<\rk_M(L)$. 
    {By the rank-criticality of $L$, it suffices to establish:}
    
    {\emph{Subclaim 1}. $J_\eps^F\in \Int(\R^d)$,}

    {\emph{Subclaim 2}. $\rk_M(J)=\rk_M(J_\eps^F)$, and}

    {\emph{Subclaim 3}. $J_\eps^F$ properly contains $L$.}

  \noindent  {In order to show \emph{Subclaim 2}, by Proposition~\ref{prop:ranks of I and Ie are equal 2}, it suffices to show:}

   {\emph{Subclaim 4}. $J_\eps^F\subset q_\ast^\uparrow$, and}
  
 {\emph{Subclaim 5}. $J=(J_\eps^F)_\eps^G$.}

\noindent{ \emph{Proof of Subclaims 1 and 5.}}      Recall that, for $t\in [t^\ast-\eps^\ast, t^\ast]$, $g_{\epsilon}(t) = t^\ast$, so Equation \eqref{eq:slices of J equal slices of Ie 2} alongside Equation \eqref{eq:some equation2} imply that, 
    for all $t\in [t^\ast-\eps, t^\ast]$,  
    $\widehat{\pi}_1(J\vert_{\pi_1=t}) = \widehat{\pi}_1(J\vert_{\pi_1=t^\ast})$, as both equal $S_{t^\ast}$.
    Thus, by Claim \ref{claim:slice condition implications IFe 2} \ref{item:IFe is an interval 2}, we have that $J_{\epsilon}^F$ is an interval, and Claim \ref{claim:slice condition implications IFe 2} \ref{item:F then G equals original 2} implies that $(J_{\epsilon}^F)_{\epsilon}^G = J$.
    
    \noindent{{\emph{Proof of Subclaim 3.}}} 
    Let $x\in \Iast$. As $G_{\epsilon}(F_{\epsilon}(x)) = x$ (cf. Remark~\ref{rem:GF is identity 2}~\ref{item:GF is identity1 2}), we have $F_{\epsilon}(x)\in \Iast_{\epsilon}^G\subset J$. 
    From the definition of $J_{\epsilon}^F$, this implies $x\in J_{\epsilon}^F$, so $J_{\epsilon}^F\supset \Iast$.
    By Claim \ref{claim:the other containment 2}%
    , we have that there exists $x\in J$ such that $x\notin \Iast_{\epsilon}^G$.
    By definition, this means that $G_{\epsilon}(x)\notin \Iast$.
    As established above, $x\in J$ implies $G_{\epsilon}(x)\in J_{\epsilon}^F$, so $G_{\epsilon}(x)\in J_{\epsilon}^F\setminus \Iast$, implying $J_{\epsilon}^F\neq \Iast$, which alongside $J_{\epsilon}^F\supset \Iast$ demonstrates that $J_{\epsilon}^F\supsetneq \Iast$.

    \noindent\emph{{Proof of Subclaim 4.}} Lastly, note that $\Iast_{\epsilon}^G, \Iast_{\epsilon'}^G\subset q_\ast^\uparrow$ implies that $J\subset q_\ast^\uparrow$ as every point in $J$ is equal to or greater than some points in $\Iast^G_\epsilon\cup \Iast^G_{\epsilon'}\subset q_\ast^\uparrow$(c.f. Remark \ref{rem:properties of finitely generated modules} \ref{item:convex hull is smallest interval}) %
    Proposition~\ref{prop:G, F preserve qast 2} then implies $J_{\epsilon}^F\subset q_\ast^\uparrow$.
\end{claimproof}

    Fix a pair $\epsilon< \epsilon'$ in $[0, \epsilon^\ast)$.
    By Claim~\ref{claim:rank of interval containing two decreases 2} {and Proposition~\ref{prop:ranks of I and Ie are equal 2}}, we have
\[
\rk_M(\Iast_{\epsilon}^G) >
\rk_M\left(\overline{\Iast_{\epsilon}^G\cup \Iast_{\epsilon'}^G}\right),
\]
or equivalently,
\[\sum_{\substack{J\supset \Iast_{\epsilon}^G\\ J\in \Ical}}
    \dgm^\Int_M(J)
    \;>\;
    \sum_{\substack{J\supset \overline{\Iast_{\epsilon}^G\cup \Iast_{\epsilon'}^G}\\ J\in \Ical}}
    \dgm^\Int_M(J).
    \]
    Therefore, we have
\begin{align*}
    \sum_{\substack{J\supset \Iast_{\epsilon}^G,\,
    J\not\supset \overline{\Iast_{\epsilon}^G\cup \Iast_{\epsilon'}^G}\\
    J\in \Ical}}
    \dgm^\Int_M(J)
    \;>&\;
    0.
\end{align*}
This implies that the set $\Jcal_{\epsilon,\epsilon'}:= \{J\in \Ical \, : \, J\supset \Iast_{\epsilon}^G \mbox{ and } J\not\supset \Iast_{\epsilon'}^G\}$ is nonempty.
For {$\epsilon\in [0,\epsilon^\ast)$,}
let $\Jcal_\epsilon:=\bigcup_{\epsilon'\in (\epsilon, \epsilon^\ast)} \Jcal_{\epsilon,\epsilon'}$.
Then, we define $\Jcal:=\bigcup_{\epsilon\in [0,\epsilon^\ast)} \Jcal_\epsilon$, which is a subcollection of $\Ical$ by construction. 

\begin{claim}\label{claim:Jcal is infinite pixelized 2}
    $|\Jcal| = \infty$.
\end{claim}
The proof follows exactly as the proof of Claim \ref{claim:Jcal is infinite pixelized}, so we omit the details here.

The next claim shows that assuming $\Ical\not\subset\Ical_M$ and Case 2 of Proposition \ref{prop:slice unequal to one side} yields a contradiction.
\begin{claim}\label{claim:rkm is not mobius invertible 2} $\rk_M^\Int$ is not M\"obius invertible.
\end{claim}
\begin{claimproof}
We begin by noting that, by Equation~\eqref{eq:slices are different 2}, $S_{t^\ast-\eps^\ast}\neq S_{t^\ast}$ and hence at least one of $S_{t^\ast-\eps^\ast}$ and $S_{t^\ast}$ is nonempty. 
This implies, by Remark \ref{rem:slice projection vs. slice}, that at least one of $\Iast\vert_{\pi_1=t^\ast-\eps^\ast}$ and $\Iast\vert_{\pi_1=t^\ast}$ is nonempty. 
Let us fix a point $p$ in the union of these two sets. 
Next, we claim that, for all $\epsilon
\in [0,\epsilon^\ast)$,  
\[\Iast_{\epsilon}^G\vert_{\pi_1=t^\ast-\eps^\ast} = \Iast\vert_{\pi_1=t^\ast-\eps^\ast} \mbox{ and } \Iast_\epsilon^G\vert_{\pi_1=t^\ast} = \Iast\vert_{\pi_1=t^\ast}.\] 
Indeed, since $g_\epsilon(t^\ast-\eps^\ast) = t^\ast-\eps^\ast$ and $g_\epsilon(t^\ast) = t^\ast$, Equation \eqref{eq:some equation2} together with Remark \ref{rem:slice projection vs. slice} implies the two equalities. 
Hence, the point $p$ belongs to
\[\Iast_{\epsilon}^G\vert_{\pi_1=t^\ast-\eps^\ast}\cup \Iast_\epsilon^G\vert_{\pi_1=t^\ast},\]
which is a subset of $\Iast_\epsilon^G.$ This implies that $p$ belongs to every interval in $\Jcal$. Then, since $\Jcal$ is an infinite subcollection of $\Ical=\supp\left(\dgm_M^\Int\right)$, 
 we have shown that the set $\{J\supset \{p\} \, : \, J\in \Ical\}$ is infinite. 
 But, since $\{p\}\in \Int(\R^d)$, this contradicts M\"obius invertibility of  $\rk^\Int_M$ (Definition~\ref{def:mobius invertible general} \ref{item:mobius invertible general 1}).
\end{claimproof}

This completes the proof of Theorem~\ref{thm:unique candidate}.

\end{document}